\documentclass[12pt,oneside,openany,leqno]{report}

\ifdefined\pdfminorversion
\fi

\newif\ifPDF
\ifx\pdfoutput\undefined\PDFfalse
\else \ifnum \pdfoutput > 0 \PDFtrue
        \else \PDFfalse
        \fi
\fi

\usepackage[centertags]{amsmath}
\usepackage{amsfonts}
\usepackage{mathrsfs}
\usepackage{textcomp}
\usepackage{amssymb}
\usepackage{amsthm}
\usepackage{microtype}
\usepackage{needspace}
\usepackage[nottoc]{tocbibind}
\usepackage{newlfont}
\usepackage[all]{xy}

\ifPDF
  \usepackage[pdftex]{color, graphicx}
  \usepackage[pdftex, bookmarks, colorlinks]{hyperref}
  \hypersetup{colorlinks=false}

\else
  \usepackage{color}
  \usepackage[dvips]{graphicx}
  \usepackage[dvips]{hyperref}
\fi

\usepackage[scale=0.8]{geometry}

\let\reportnormalsize\normalsize
\renewcommand{\normalsize}{%
  \reportnormalsize\fontsize{12}{14}\selectfont}
\normalsize

\usepackage[pagewise, mathlines]{lineno}

\newtheorem{thm}{Theorem}[chapter]
\newtheorem{cor}[thm]{Corollary}
\newtheorem{lem}[thm]{Lemma}
\newtheorem{prop}[thm]{Proposition}
\newtheorem{prp}[thm]{Proposition}

\newtheorem{introthm}{Theorem}

\newtheorem*{overviewthm}{Theorem}
\theoremstyle{definition}
\newtheorem{defn}[thm]{Definition}
\newtheorem{dfn}[thm]{Definition}
\newtheorem{introdefn}[introthm]{Definition}
\newtheorem{cnv}[thm]{Convention}
\newtheorem{ntn}[thm]{Notation}
\theoremstyle{remark}
\newtheorem{rem}[thm]{Remark}

\newtheorem{example}[thm]{Example}

\numberwithin{equation}{chapter}
\newcommand{\norm}[1]{\left\Vert#1\right\Vert}
\newcommand{\abs}[1]{\left\vert#1\right\vert}

\newcommand{\Int}{\mathbb Z}
\newcommand{\Comp}{\mathbb C}

\newcommand{\eps}{\varepsilon}

\newcommand{\F}{\mathcal{F}}
\DeclareRobustCommand{\Cstar}{\ensuremath{\mathrm C^*}}
\DeclareRobustCommand{\Wstar}{\ensuremath{\mathrm W^*}}
\DeclareRobustCommand{\Ksym}{\ensuremath{\mathrm K}}
\DeclareRobustCommand{\KKsym}{\ensuremath{\mathrm{KK}}}
\usepackage{mathtools}

\newcommand{\beq}{\begin{equation}}
\newcommand{\eeq}{\end{equation}}
\newcommand{\beqr}{\begin{eqnarray*}}
\newcommand{\eeqr}{\end{eqnarray*}}
\newcommand{\bal}{\begin{align*}} 
\newcommand{\eal}{\end{align*}}
\newcommand{\bei}{\begin{itemize}}
\newcommand{\eei}{\end{itemize}}

\newcommand{\Gm}{\Gamma}

\newcommand{\QED}{\rule{0.4em}{2ex}}

\begin{document}


\title{Weak tracial approximation and weak Rokhlin property for \Cstar{}-dynamical systems}

\author{George A.~Elliott, Chun Guang Li, Zhuang Niu, and Qingyun Wang}


\date{}



\hypersetup{pageanchor=false}
\maketitle

\begin{abstract}
We extend results on tracial approximation and regularity from topological dynamical systems to \Cstar{}-dynamical systems. For actions of countable discrete amenable groups on unital separable \Cstar{}-algebras, we introduce a weak Rokhlin property which yields weak tracial approximation of crossed products by finite direct sums of matrix algebras over hereditary subalgebras of the coefficient algebra. Examples arise from central topological towers and from \Wstar{}-central freeness.

\par
\par
\medskip
\noindent\textbf{Keywords:} \Cstar{}-dynamical system, Rokhlin property,
weak tracial approximation.

\medskip
\noindent\textbf{Acknowledgements.}
The first draft dates from 2020--2021.  The recent development of ChatGPT
has assisted the subsequent revision of the arguments and exposition
and the completion of the present version.

G.A.E.~is supported by an NSERC Discovery Grant. C.G.L.~is supported by an NNSF Grant of China (No. 12671150). Z.N.~is supported by a Simons Foundation Grant (MP-TSM-00002606).

\end{abstract}

\hypersetup{pageanchor=true}
\pagenumbering{roman}
\tableofcontents

\clearpage
\pagenumbering{arabic}

\chapter{Introduction and preliminaries}\label{chapter-introduction-preliminaries}

\section{Introduction}\label{section-introduction}

The main motivation of this paper is to extend the correspondence
between Rokhlin properties and weak tracial approximation from
topological dynamical systems to \Cstar{}-dynamical systems.  In
\cite{Niu-MD-Z-absorbing,Niu-MD-Z},
Rokhlin towers provide the local approximation structure through
which the geometry of an action enters the study of its crossed
product.  To describe this starting point, let a countably infinite
discrete amenable group $\Gamma$ act freely and minimally on a compact
metrizable space $X$.  When the action has the uniform Rokhlin property
(URP), one can choose finite families of disjoint open towers with
sufficiently invariant shapes whose union has measure arbitrarily
close to one, uniformly over all invariant probability measures.
The tower levels give matrix coordinates in $\mathrm C(X)\rtimes\Gamma$.
Controlling the boundaries of the tower shapes gives almost central
positive contractions whose compressions of prescribed finite subsets
are approximated by finite direct sums of possibly nonunital
homogeneous \Cstar{}-algebras
\cite[Theorem~3.9]{Niu-MD-Z}.  Uniform coverage by the towers gives
a compression defect whose dimension functions are uniformly small
on all tracial states of the crossed product.  Together with Cuntz
comparison of open sets (COS), this approximation connects the mean
dimension of the action to comparison and Jiang--Su stability of
its crossed product.

For an action on a general unital \Cstar{}-algebra $A$, we seek an intrinsic
formulation of the same mechanism.  An open tower base $U\subseteq X$
is represented algebraically by the hereditary subalgebra
$\mathrm C_0(U)\subseteq\mathrm C(X)$; this suggests using hereditary
\Cstar{}-subalgebras of $A$ as the coefficients of the local matrix models.
Positive contractions describe tower bases, orthogonality of their
translates replaces disjointness of level sets, and invariant tracial
states replace invariant probability measures in measuring the
complement.  In the noncommutative setting, the tower bases must also
commute asymptotically with the coefficient algebra.  This condition
is automatic for commutative coefficients and is essential for the
resulting matrix models to approximate compressions of general
elements of the crossed product.  Our formulation of the weak Rokhlin
property (WRP) incorporates this centrality and uniform tracial
control.  Its principal structural consequence is weak tracial
approximation of $A\rtimes_r\Gamma$ by finite direct sums of matrix
algebras over hereditary \Cstar{}-subalgebras of $A$, with the compressing
positive elements still belonging to $A$
(Theorem~\ref{WTA-C-dyn-1}).  Allowing nonunital models and positive
compressions is essential already in the topological setting, since
open tower bases need not be compact.  The resulting approximation
structure provides a common setting for studying how comparison and
regularity of the coefficient algebra influence the crossed product.

\subsection{Weak tracial approximation}

The preceding tower construction suggests comparing a finite subset
of an ambient algebra with a model after compression by a positive
contraction.  The model need not contain the unit of the ambient
algebra.  We therefore use a positive contraction $h$ in the model,
require it to commute approximately with the prescribed elements,
and approximate their compressions by elements of the model.
The defect $1-h$ is controlled either uniformly over tracial states
or by comparison with a prescribed nonzero positive element.
These two forms of control lead to the following approximation
properties.

For a unital \Cstar{}-algebra $A$, write $\mathrm T(A)$ for its tracial
state space and put
\[
 \mathrm d_\tau(a)=\lim_{n\to\infty}\tau(a^{1/n}),\qquad a\in A_+.
\]
We write $a\precsim b$ for Cuntz subequivalence.

\begin{introdefn}[Definition~\ref{defn_WAS}]\label{intro-defn_WAS}
Let $\mathcal S$ be a class of \Cstar{}-algebras, and let $A$ be unital.
We say that $A$ is \emph{weakly tracially approximated by $\mathcal S$},
and write $A\in\mathrm{TA}_0(\mathcal S)$, if
$\mathrm T(A)\neq\emptyset$ and, for every finite set
$\mathcal F=\{x_1,\ldots,x_n\}\subseteq A$ and every $\eps>0$, there
are a \Cstar{}-subalgebra $C\subseteq A$ with $C\in\mathcal S$, elements
$x_1',\ldots,x_n'\in A$, and a positive contraction $h\in C$ with
$\|h\|=1$ such that
\begin{enumerate}\setlength{\itemsep}{0pt}\setlength{\parsep}{0pt}
\item[(1)] $\|x_i-x_i'\|<\eps$, $1\leq i\leq n$;
\item[(2)] $\|[h,x_i']\|<\eps$, $1\leq i\leq n$;
\item[(3)] $hx_i'h\in C$, $1\leq i\leq n$;
\item[(4)] $\mathrm d_\tau(1-h)<\eps$ for every $\tau\in\mathrm T(A)$.
\end{enumerate}

We say that $A$ is \emph{tracially approximated by $\mathcal S$},
and write $A\in\mathrm{TA}_1(\mathcal S)$, if, for every such
$\mathcal F$ and $\eps$ and every nonzero positive contraction
$a_0\in A$, the same data can be chosen to satisfy (1)--(3) and
\[
 \text{(4')}\qquad 1-h\precsim a_0.
\]
The definition of $\mathrm{TA}_1(\mathcal S)$ does not require
$\mathrm T(A)\neq\emptyset$.
\end{introdefn}

Condition (4) measures the support of the defect.  For example, if
$B=\mathrm C(X)\rtimes\Gamma$, $h\in\mathrm C(X)$ is a positive
contraction, and $\tau_\mu=\mu\circ\mathbb E$ is the trace induced
by an invariant probability measure $\mu$, then
\[
 \mathrm d_{\tau_\mu}(1-h)=\mu\bigl(\{x\in X:h(x)<1\}\bigr).
\]
Thus the uniform measure estimate furnished by Rokhlin towers becomes
exactly the defect estimate in weak tracial approximation.
For a fixed $h$, a bound on $\mathrm d_\tau(1-h)$ is stronger than a bound
on $\tau(1-h)$.  At the level of the approximation property,
however, the two formulations are equivalent after a
functional-calculus cutoff and a change of tolerance
(Remark~\ref{trace-rank-ta0}).

The perturbed elements $x_i'$ in Definition~\ref{intro-defn_WAS}
ensure that the compressed elements belong to the model exactly.
Equivalently, one may work directly with the prescribed finite set
and require
\[
 \|[h,x]\|<\eps,\qquad
 \operatorname{dist}(hxh,C)<\eps,\qquad x\in\mathcal F,
\]
together with the corresponding defect condition; see
Lemma~\ref{alt-defn-was}.  This is the form naturally supplied by
the tower constructions.

The $\mathrm{TA}_1$ condition expresses the smallness of the defect
by Cuntz comparison and therefore does not require tracial states.
When $A$ is simple and has tracial states,
$\mathrm{TA}_1(\mathcal S)$ implies
$\mathrm{TA}_0(\mathcal S)$; strict comparison by tracial states
on $A_+$ gives the converse
(Remarks~\ref{ta1-to-ta0} and~\ref{comparison-gta}).
For such $A$, Remark~\ref{comparison-gta} also identifies the Cuntz-controlled
version with nonunital generalized tracial approximation
\cite{Elliott-Fan-Fang-nonunital-GTA} and explains its relation
to the formulation with unital models in
\cite{Elliott-Fan-Fang-GTA}.
For dynamical applications, the next step is to construct these
positive compressions and local models from Rokhlin towers in the
coefficient algebra.

\subsection{Weak Rokhlin property}

The classical Rokhlin lemma gives arbitrarily tall measurable towers
for an aperiodic invertible measure-preserving transformation of a standard
probability space, with complement of arbitrarily small measure.
For topological actions, the uniform Rokhlin property replaces
measurable bases by open sets and requires the complement to be small
uniformly over invariant probability measures.  For free minimal
actions of countably infinite discrete amenable groups on compact
metrizable spaces, URP together
with COS bounds the radius of comparison of the crossed product by
one half of the mean dimension \cite{Niu-MD-Z}, gives stable rank one
\cite{LN-sr1}, and gives Jiang--Su stability when the mean dimension
is zero \cite{Niu-MD-Z-absorbing}.  These results motivate a tower
condition for \Cstar{}-dynamical systems that retains uniform tracial
coverage and produces the positive compressions used in weak
tracial approximation.

Let $A$ be a separable unital \Cstar{}-algebra and let
$\alpha\colon\Gamma\to\operatorname{Aut}(A)$ be an action of a
countable discrete amenable group.
Fix a free ultrafilter
$\omega$ on $\mathbb N$, and write
\[
 A_\omega=\ell^\infty(\mathbb N,A)/
 \{(a_n):\lim_{n\to\omega}\|a_n\|=0\}.
\]
We view $A$ as the constant sequences in $A_\omega$ and extend the
action coordinatewise.  Denote by $\mathrm{T}^\alpha(A)$ the space of invariant
tracial states, and put
\[
 J_{\mathrm{T}^\alpha(A)}=\left\{[(a_n)]\in A_\omega:
   \lim_{n\to\omega}\sup_{\rho\in \mathrm{T}^\alpha(A)}\rho(a_n^*a_n)=0\right\}.
\]
The supremum is taken before the ultralimit, so the estimate is
uniform over the invariant traces.  Recall that an exact tiling of
$\Gamma$ is a partition
\[
 \Gamma=\bigsqcup_{s=1}^S\bigsqcup_{t\in C_s}D_st
\]
with finitely many finite nonempty shapes $D_1,\ldots,D_S$ and
center sets $C_s\subseteq\Gamma$.  We make the following definition.
\begin{introdefn}[Definition~\ref{D:weak-Rokhlin-property}]\label{intro-D:weak-Rokhlin-property}
Let $A$ be a unital separable \Cstar{}-algebra with
$\mathrm{T}(A)\ne\emptyset$, and let $\alpha$ be an action of a countable
discrete amenable group $\Gamma$.  The action has the
\emph{weak Rokhlin property} (WRP) if, for every exact tiling of
$\Gamma$ with shapes $D_1,\ldots,D_S$, there are positive
contractions $c_s\in A_\omega\cap A'$, $1 \leq s \leq S$, such that
\begin{enumerate}
\item the levels $\alpha_g(c_s)$, $g\in D_s$, $1\leq s\leq S$, are mutually orthogonal, and
\item
\[
 1-\sum_{s=1}^S\sum_{g\in D_s}\alpha_g(c_s)
       \in J_{\mathrm{T}^\alpha(A)}.
\]
\end{enumerate}
\end{introdefn}

Amenability ensures that $\mathrm T^\alpha(A)$ is nonempty when $\mathrm T(A)$ is
nonempty.  It also provides exact tilings whose shapes have small
boundary relative to any prescribed finite subset of $\Gamma$
\cite{DHZ-tiling}.  The quantifier over prescribed tilings lets us
choose these shapes to fit the Fourier supports of the elements to
be approximated.

For $A=\mathrm C(X)$, WRP agrees with the topological
uniform Rokhlin property (URP): disjoint open towers cover all but
a uniformly small amount of every invariant measure
(Proposition~\ref{top-urp-wrp}).  For general coefficient algebras,
the finite-stage uniform Rokhlin condition in
Definition~\ref{URP-defn-C-dyn} implies WRP
(Proposition~\ref{urp-implies-wrp}).

For comparison, the literature contains several Rokhlin properties
with different tower and remainder requirements.  For finite group
actions, the strict Rokhlin property uses approximately central
projection towers whose sum is the unit
\cite[Definition~1.1]{Phillips-tracial-Rokhlin}.
The tracial Rokhlin property retains projection towers but allows a
complement which is Murray--von Neumann equivalent to a projection
in the hereditary subalgebra generated by a prescribed nonzero
positive element of $A$
\cite[Definition~1.2]{Phillips-tracial-Rokhlin}.
Weak tracial Rokhlin properties replace these projections by
positive contractions while retaining Cuntz control of the
complement in $A$.  Such formulations are considered for amenable
groups in \cite{Wang-WRP,OPW-Strict-comparison-crossed-products}
and, including nonunital coefficient algebras, for finite groups in
\cite{Forough-Golestani-weak-tracial}.

There are also positive-tower formulations based directly on
tracial estimates, as in \cite{Sato-Rohlin-Z}.  The weak Rohlin
property of \cite[Definition~2.5(3)]{Matui-Sato-strongly-outer-II}
uses a central sequence and one sufficiently invariant tower shape,
with the traces of the base converging uniformly over all of
$\mathrm T(A)$ to the reciprocal of the tower size.  The amenable-group
formulation of \cite[Definition~2.9]{Gardella-Hirshberg-Vaccaro} allows a finite
family of shapes, requires the complement to vanish uniformly over
all traces, and requires the levels within each tower to have one
common positive tracial value, independent of the trace.  Thus even
among positive-tower conditions,
the trace space tested and the normalization of the levels matter.

A different relaxation is finite Rokhlin dimension
\cite{HWZ-RD}, extended to residually finite groups in
\cite{SWZ-Z}.  Here finitely many colours of positive towers cover
the unit in norm, with orthogonality required within each colour.
The version with commuting towers also requires approximate
commutation between different colours.  Our WRP requires all tower
levels to be mutually orthogonal, but permits a complement invisible
to invariant traces.  These are different ways to relax projection
Rokhlin towers.  In particular, the permanence results for finite
nuclear dimension and, with commuting towers, Jiang--Su stability
in \cite{HWZ-RD} arise from this coloured norm-covering structure.

Our formulation dispenses with projections, Cuntz control in $A$,
and trace estimates outside $\mathrm T^\alpha(A)$.  It is closest to the
weak Rokhlin property with respect to invariant traces in
\cite[Definition~5.2]{Gardella-Geffen-Naryshkin-Vaccaro-Z}.
Our WRP implies that condition, since sufficiently invariant exact
tilings may be chosen.  The additional requirement here is that
towers exist for every prescribed exact tiling; in the cited
formulation the sufficiently invariant shapes are chosen as part
of the construction.  This requirement on the shapes is useful
for fitting the towers to the finite sets occurring in the
approximation theorem below.

The link with weak tracial approximation is the following structure
theorem.  It approximates the crossed product
by matrix algebras over hereditary subalgebras of the coefficient
algebra, with the compression itself still in $A$.

\Needspace{12\baselineskip}
\begin{introthm}[Theorem~\ref{WTA-C-dyn-1}]\label{intro-WTA-C-dyn-1}
Let $A$ be a unital separable \Cstar{}-algebra with
$\mathrm T(A)\ne\emptyset$, and let $\alpha$ be an action of a countable
discrete amenable group $\Gamma$ with WRP.  Put
$B=A\rtimes_{\alpha,r}\Gamma$.  For every finite
$\mathcal F\subseteq B$ and every $\eps>0$, there are a subalgebra
\[
 C\cong\bigoplus_{s=1}^S \mathrm M_{n_s}(H_s)\subseteq B,
 \qquad H_s\subseteq A\ \text{hereditary},
\]
and a norm-one positive contraction $h\in A\cap C$ such that
\[
 \|[h,x]\|<\eps,\qquad \operatorname{dist}(hxh,C)<\eps
       \quad(x\in\mathcal F),
 \qquad \mathrm d_\tau(1-h)<\eps\quad(\tau\in \mathrm T(B)).
\]
When $x\in\mathcal F\cap A$, its compressed approximation may be
chosen in $A\cap C$.  In particular,
$B\in\mathrm{TA}_0(\mathcal S)$, where $\mathcal S$ consists of
finite direct sums of matrix algebras over hereditary subalgebras of $A$.
\end{introthm}

The matrix coordinates in these models record the levels of a tower;
the hereditary algebra records its base.  In the commutative case,
$H_s=\mathrm C_0(U_s)$ for an open subset $U_s\subseteq X$, recovering the
models supplied by topological towers.  In the general case, the
proof constructs matrix cones from central sequence towers and lifts
them while keeping their diagonals in $A$.  Enlarging the cones gives
the hereditary matrix models, and a spectral cutoff changes the
uniform tracial $2$-norm estimate into the rank estimate on $1-h$.
The conclusion requires neither simplicity nor nuclearity of $A$.
WRP also identifies the traces of $B$ with the invariant traces of
$A$; when the action is minimal, it makes $B$ simple
(Proposition~\ref{trace-factor} and Corollary~\ref{wrp-simple}).

The location of the defect is important here: $1-h\in A$, but WRP
alone gives its smallness through traces.  When $A$ is simple and
infinite-dimensional, $B$ is simple, and the action also has COS,
the same approximation can be chosen so that, for any prescribed
$0\ne z\in B_+$,
\[
                  1-h\precsim_B z
\]
(Proposition~\ref{WRP-COS-TA1}).  This upgrades the conclusion to
$B\in\mathrm{TA}_1(\mathcal S)$.  The Cuntz comparison takes place
in the crossed product $B$, even though the compressing element and
its defect belong to $A$.  This is the form of remainder control
used below to transfer comparison and rank properties through the
tower models.

\subsection[Sources of weak tracial approximation and the weak Rokhlin property]{Sources of weak tracial approximation\\and the weak Rokhlin property}

Chapter~\ref{section-examples} collects the examples of weak tracial
approximation and WRP discussed here, together with their proofs
and constructions.  We begin with interval approximation for
classifiable \Cstar{}-algebras, which does not require an action, and then
describe two sources of WRP: central topological towers and towers
obtained from invariant traces.  By
Theorem~\ref{intro-WTA-C-dyn-1}, each of these WRP examples also
gives weak tracial approximation of the corresponding crossed product.

The first result shows that the Cuntz-controlled approximation in
Definition~\ref{intro-defn_WAS} holds throughout the unital simple
infinite-dimensional stably finite classifiable setting.  The local
models can be chosen to be matrix algebras over half-open or open
intervals.

\begin{introthm}[Theorem~\ref{classifiable-interval-approximation}]\label{intro-classifiable-interval-approximation}
Let $A$ be a unital, separable, simple, infinite-dimensional, nuclear,
stably finite, $\mathcal Z$-stable \Cstar{}-algebra satisfying the UCT.
Let $\mathcal I_{\mathrm h}$ and $\mathcal I_{\mathrm o}$ be the classes
of finite direct sums
\[
 \bigoplus_{j=1}^{r}\mathrm M_{n_j}(\mathrm C_0([0,1)))\quad\hbox{and}\quad
 \bigoplus_{j=1}^{r}\mathrm M_{n_j}(\mathrm C_0((0,1))),
 \qquad r,n_1,\ldots,n_r\in\mathbb N,
\]
respectively.  Then
\[
 A\in\mathrm{TA}_1(\mathcal I_{\mathrm h})
       \cap\mathrm{TA}_1(\mathcal I_{\mathrm o}).
\]
In particular, $A$ can be weakly tracially approximated by either class.
\end{introthm}

These interval models are projectionless.  They approximate finite
subsets after compression by a positive contraction, as in
Definition~\ref{intro-defn_WAS}, so the ambient algebra may still
have nontrivial projections and \Ksym{}-theory.
The proof starts from the classification models of
\cite{GLN-TAS-1,GLN-TAS-2} and constructs half-open interval
subalgebras with uniformly small tracial defect.  Strict comparison
gives the required Cuntz control, and $\mathcal Z$-absorption gives
the open interval version.

For actions, a first source of WRP is a central topological
subsystem.  Suppose that $\mathrm C(X)\subseteq Z(A)$ is invariant under
the action.  Tower functions in $\mathrm C(X)$ then commute with all of $A$,
and invariant traces on $A$ restrict to invariant probability
measures on $X$.  Thus open towers on $X$ supply algebraic towers
with the centrality and uniform coverage required in the preceding
subsection.  This gives the following criterion.

\begin{introthm}[Theorem~\ref{central-topological-urp}]\label{intro-central-topological-urp}
Let $A$ be a unital separable \Cstar{}-algebra with
$\mathrm T(A)\ne\emptyset$, and let $\alpha$ be an action of a
countable discrete amenable group $\Gamma$ on $A$.
Suppose that $A$ contains an $\alpha$-invariant unital subalgebra
$\mathrm C(X)\subseteq Z(A)$, where $X$ is compact Hausdorff.
If the induced action on $X$ has topological URP, then
$(A,\Gamma)$ has the URP of Definition~\ref{URP-defn-C-dyn}
and hence WRP.

In particular, this applies to a unital separable \Cstar{}-algebra that is
a continuous field over $X$ and has a nonempty tracial state space, provided that
the action preserves its canonical central copy of $\mathrm C(X)$ and the
induced action on $X$ has topological URP.
\end{introthm}

A second source comes from the von Neumann algebras associated with
invariant traces.  Each $\rho\in \mathrm T^\alpha(A)$ gives a GNS von
Neumann algebra $M_\rho=\pi_\rho(A)''$ with an induced action.
The action on $A$ is called \emph{\Wstar{}-centrally free} when all these
induced actions are centrally free
(Definition~\ref{wstar-freeness-cstar}).  The construction begins
with Rokhlin towers in these von Neumann algebras and lifts them to
positive towers in the norm central sequence algebra of $A$.

When $\mathrm T^\alpha(A)$ has only finitely many extreme points,
\Wstar{}-central freeness suffices for WRP without a nuclearity assumption
on $A$ (Theorem~\ref{finite-invariant-traces-wrp}).  This includes
the case of a unique invariant trace
(Theorem~\ref{unique-invariant-trace-wrp}).  For a larger family
of invariant traces, nuclearity allows the construction to be made
uniform over a compact finite-dimensional extreme boundary.
The next theorem states this criterion; recall that a Bauer simplex
is a Choquet simplex whose extreme boundary is compact.

\begin{introthm}[Theorem~\ref{finite-dimensional-trace-wrp}]\label{intro-finite-dimensional-trace-wrp}
Let $A$ be a unital separable nuclear \Cstar{}-algebra with
$\mathrm T(A)\ne\emptyset$, and let
$\alpha\colon\Gamma\to\operatorname{Aut}(A)$ be an action
of a countable discrete amenable group.
Suppose that $\alpha$ is \Wstar{}-centrally free and that
$S=\mathrm T^\alpha(A)$ is a Bauer simplex with
$\dim(\partial_eS)<\infty$.
If $\Gamma$ is infinite, then $\alpha$ has WRP.
The same conclusion holds for finite $\Gamma$ when $A$ has
no finite-dimensional quotients.
\end{introthm}

The dimension assumption concerns the extreme boundary of the
\emph{invariant} trace simplex, which need not be a face of $\mathrm T(A)$.
For infinite groups, finite-dimensional quotients of $A$ are allowed.
This is important for commutative coefficients: a free action of a
countably infinite discrete amenable group on a compact metrizable
space has WRP, and hence topological URP, whenever its
invariant-measure simplex is Bauer with finite-dimensional extreme
boundary (Corollary~\ref{commutative-bauer-urp}).  The theorem
therefore also supplies topological towers from information about
invariant measures.

\subsection{Comparison and Jiang--Su stability}

The preceding constructions give local models on large positive
compressions.  To use these models to study the whole algebra, one
must also control the part left outside the compression.  Comparison
provides this control by absorbing the defect into a prescribed
positive element.  Chapter~\ref{section-comparison} develops the
comparison results, and Chapter~\ref{section-Z-stability} uses them
to obtain Jiang--Su stability.

For dynamical applications, the first comparison condition concerns
positive elements of the coefficient algebra.  It extends the Cuntz
comparison of open sets used for topological actions and retains
the name COS.  Let $\alpha$ be an action of a discrete group
$\Gamma$ on a unital \Cstar{}-algebra $A$, suppose that
$\mathrm T^\alpha(A)\ne\emptyset$, and put $B=A\rtimes_{\alpha,r}\Gamma$.
For $\lambda\in(0,1]$ and $m\in\mathbb N$, the condition
$(\lambda,m)$-COS requires
\[
 \mathrm d_\rho(a)<\lambda \mathrm d_\rho(b)\quad(\rho\in \mathrm T^\alpha(A))
 \quad\Longrightarrow\quad
 a\precsim_B b^{\oplus m},
 \qquad a,b\in A_+.
\]
The inequalities are tested on invariant traces of $A$, and the
comparison takes place in $B$.  The constants allow a fixed loss
in the trace inequality and a fixed number of copies of the target.
For a simple infinite-dimensional coefficient algebra with simple
crossed product, WRP and COS turn the tracially small defect of
Theorem~\ref{intro-WTA-C-dyn-1} into a Cuntz-small defect in $B$.
Thus $B$ satisfies the $\mathrm{TA}_1$ version of the tower
approximation described in the preceding subsections
(Proposition~\ref{WRP-COS-TA1}).

COS can also follow from regularity of the coefficient algebra.
If $A$ is simple and unital, has strict comparison by tracial states,
and has almost divisible Cuntz semigroup, then every discrete group
action admitting an invariant tracial state has $(1,1)$-COS
(Proposition~\ref{pre-COM}).  In particular, this applies when
$A$ is simple, unital, exact, stably finite, and $\mathcal Z$-stable
(Corollary~\ref{Z=>COS}).  The amenable-group case of this criterion
was established in
\cite[Proposition~5.4]{OPW-Strict-comparison-crossed-products}.
The conclusion compares coefficient elements, including their finite
matrix amplifications; comparison of arbitrary elements of $B$
uses the approximation argument as well.

The abstract form of that argument is
Theorem~\ref{WAD=>SCP}.  Let $A$ be simple, separable, exact, unital,
and infinite-dimensional, with $\mathrm T(A)\ne\emptyset$ and
$A\in\mathrm{TA}_1(\mathcal S)$.  If $\mathrm{Cu}(C)$ is almost
unperforated for every $C\in\mathcal S$, then $A$ has strict
comparison.  The theorem also transfers quantitative comparison
bounds.  For nonunital models, the admissible gap is expressed by
a finite-valued continuous affine function on the cone of densely
finite traces, including zero, where the function vanishes.
Such bounds pass to $A$ when they are compatible with the chosen
compressions.  This compatibility is part of the hypothesis:
traces on a model need not extend to traces on $A$.

For Jiang--Su stability, comparison is accompanied by a second
requirement: sufficiently many positive elements must be available
to realize prescribed rank functions.  A positive element
$x\in(A\otimes\mathcal K)_+$ determines the function
$\tau\mapsto \mathrm d_\tau(x)$ on $\mathrm T(A)$.  Write
$\mathrm{LAff}_{++}(\mathrm T(A))$ for the strictly positive lower
semicontinuous affine functions on $\mathrm T(A)$, allowing the value
$\infty$, and $\mathrm V(A)$ for the semigroup of projection classes.
The next theorem shows that, under the stated local approximation
hypothesis, strict comparison and realization of all these functions
are exactly what is needed for $\mathcal Z$-absorption.

\begin{introthm}[Theorem~\ref{TW-conj}]\label{intro-TW-conj}
Let $\mathcal S$ be a class such that every $C\in\mathcal S$ can be
locally approximated by \Cstar{}-subalgebras of finite nuclear dimension.  Let
$A$ be a separable, simple, unital, infinite-dimensional, nuclear
\Cstar{}-algebra with $A\in\mathrm{TA}_0(\mathcal S)$.  Then
$A\cong A\otimes\mathcal Z$ if and only if
\begin{enumerate}
\item $A$ has strict comparison of positive elements; and
\item for every $f\in\mathrm{LAff}_{++}(\mathrm T(A))$ there is
      $x\in(A\otimes\mathcal K)_+$ such that
      $\mathrm d_\tau(x)=f(\tau)$ for every $\tau\in \mathrm T(A)$.
\end{enumerate}
Consequently, within this class, $\mathcal Z$-stability is also equivalent
to the canonical ordered-semigroup description
\[
 \mathrm{Cu}(A)\cong
 \mathrm V(A)\sqcup\mathrm{LAff}_{++}(\mathrm T(A)).
\]
\end{introthm}

In the displayed description of $\mathrm{Cu}(A)$, $\mathrm V(A)$ records
projection classes, and the other classes are described by their
tracial rank functions.  The hypothesis on the models requires no
common bound on the nuclear dimensions of their local approximants.
For basic homogeneous models, the approximation and strict comparison
already supply the required rank realization.  This gives the
following criterion.

\begin{introthm}[Theorem~\ref{homogeneous-TA-Z}]\label{intro-homogeneous-TA-Z}
Let
\[
 \mathcal H=\big\{\mathrm M_n(\mathrm C_0(X)):n\in\mathbb N,
       \ X\text{ locally compact and metrizable}\big\}
\]
be the class of basic homogeneous \Cstar{}-algebras.  Let $A$ be a separable,
simple, unital, nuclear, non-type-I \Cstar{}-algebra.  If
$A\in\mathrm{TA}_0(\mathcal H)$, then
\[
 A\cong A\otimes\mathcal Z
 \quad\Longleftrightarrow\quad
 A\text{ has strict comparison of positive elements}.
\]
\end{introthm}

For crossed products, the hereditary models in
Theorem~\ref{intro-WTA-C-dyn-1} allow regularity of the coefficient
algebra to enter the argument.  Corollary~\ref{WRP-WTZ} shows that,
for a separable, simple, unital, infinite-dimensional $A$ with simple
crossed product $B$, WRP transfers weak tracial
$\mathcal Z$-absorption from $A$ to $B$.  This supplies approximately
central order zero maps from matrix algebras with tracially small
complement.  Strict comparison in $B$ strengthens this to Cuntz
control of the complement, and nuclearity then gives
$\mathcal Z$-stability.

If $A$ is nuclear and $\mathcal Z$-stable, the required comparison
in $B$ follows from the results above: the action has COS, the hereditary
matrix models have almost unperforated Cuntz semigroups, and the
$\mathrm{TA}_1$ comparison theorem gives strict comparison in $B$.
This yields the following permanence theorem.

\begin{introthm}[Theorem~\ref{WRP-finite-nucdim-Z}]\label{intro-WRP-finite-nucdim-Z}
Let $A$ be a separable, simple, unital, nuclear, $\mathcal Z$-stable
\Cstar{}-algebra, and let $\Gamma$ be a countable discrete amenable group
acting on $A$.  If $(A,\Gamma)$ has WRP, then
\[
 B=A\rtimes_{\alpha,r}\Gamma
\]
is simple and $\mathcal Z$-stable.
In particular, the conclusion holds for a separable, simple, unital,
infinite-dimensional coefficient algebra of finite nuclear dimension
with WRP.
\end{introthm}

The theorem requires no separate strict comparison assumption on
the crossed product.  Preservation of $\mathcal Z$-stability under
the weak tracial Rokhlin property was obtained in
\cite[Corollary~8.13]{OPW-Strict-comparison-crossed-products}.
That formulation controls the defect in the coefficient algebra in
the Cuntz order; our WRP uses uniform control over invariant traces.
For the present WRP, the conclusion also follows from
\cite[Corollary~5.6]{Gardella-Geffen-Naryshkin-Vaccaro-Z}, using
the implication between the two Rokhlin conditions described above.
The sources of WRP in the preceding subsection can therefore be
combined with $\mathcal Z$-stable coefficients to obtain simple
$\mathcal Z$-stable crossed products.

\subsection{Real rank and stable rank}

Comparison also allows rank properties to pass from the local
models to the ambient algebra.  Recall that real rank zero means
that self-adjoint elements can be approximated by finite-spectrum
self-adjoint elements, and stable rank one means that invertible
elements are dense in the unitization.  Chapter~\ref{section-ranks}
first proves permanence results for these two properties, and then
uses the tower structure to obtain stable rank one for crossed
products.

For the permanence results, weak tracial approximation supplies
the desired approximations on a large positive compression.
Strict comparison controls the remaining defect by comparison
with positive elements in the model.  This allows the
approximation to be completed in the whole algebra.

\begin{introthm}[Theorem~\ref{rr0-twa}]\label{intro-rr0-twa}
Let $\mathcal S$ be a class of \Cstar{}-algebras of real rank zero, and let $A$
be a simple unital \Cstar{}-algebra in $\mathrm{TA}_0(\mathcal S)$.  If $A$
has strict comparison of positive elements, then $A$
has real rank zero.
\end{introthm}

The real-rank conclusion requires no stable-rank-one assumption on
the models.  The corresponding stable-rank permanence theorem is
as follows.

\begin{introthm}[Theorem~\ref{sr1-twa}]\label{intro-sr1-twa}
Let $\mathcal S$ be a class of \Cstar{}-algebras of stable rank one.  Let $A$ be
a simple unital \Cstar{}-algebra in $\mathrm{TA}_0(\mathcal S)$.  If $A$ has
strict comparison of positive elements, then $A$ has stable rank one.
\end{introthm}

The models in these permanence theorems need not be unital or
simple.  For crossed products, the dynamics give a further
stable-rank result without requiring stable rank one of the
coefficient algebra or of the hereditary tower models.  Minimality
and WRP provide orthogonal positive coefficient elements that are
equivalent in the crossed product.  Choosing enough such elements
allows the fixed comparison multiplicity in COS to be absorbed.
Together with the matrix structure of the towers, this gives the
approximations by invertibles needed for the following theorem.
The argument extends the stable-rank methods of \cite{LN-sr1}
to noncommutative coefficient algebras.

\begin{introthm}[Theorem~\ref{sr1-wrp-cos-D}]\label{intro-sr1-wrp-cos-D}
Let $(A,\Gamma)$ be a minimal \Cstar{}-dynamical system with WRP, where
$A$ is unital and separable, $\mathrm T(A)\ne\varnothing$, and
$\Gamma$ is countably infinite, discrete, and amenable.  If
$(A,\Gamma)$ has $(\lambda,m)$-COS for some $\lambda\in(0,1]$
and $m\in\mathbb N$, then
\[
              \operatorname{tsr}(A\rtimes_{\alpha,r}\Gamma)=1.
\]
\end{introthm}

The coefficient algebra in Theorem~\ref{intro-sr1-wrp-cos-D} may be
nonsimple; minimality excludes nonzero proper invariant ideals.  The theorem
requires neither nuclearity nor strict comparison of the whole
crossed product.  Its comparison hypothesis is COS, tested on
coefficient elements and implemented in the crossed product.

Rank properties also combine with the Jiang--Su stability results
to give classification structure.  Corollary~\ref{WRP-TAF-AH}
shows that a WRP crossed product is tracially AF, and is isomorphic
to a simple unital AH algebra with slow dimension growth and real
rank zero, under its additional hypotheses.  These include real
rank zero and weak tracial $\mathcal Z$-absorption for the simple
coefficient algebra, and strict comparison, nuclearity, and the UCT
for the crossed product.  The stable-rank theorem will also be used
with the quantitative comparison estimates in the next subsection.

\subsection{Orbit radius of comparison}

To obtain quantitative comparison estimates, we introduce the
\emph{orbit radius of comparison}, $\operatorname{orc}(A,\alpha)$,
in Section~\ref{orbit-radius-section}.  It can be regarded as a
\emph{\Cstar{}-algebraic analogue of mean dimension}: comparison of
positive matrices over $A$ takes the place of dimension estimates
for open covers along finite orbit sets.  The rank functions are
normalized by the size of the orbit set, just as covering dimensions
are normalized in the definition of mean dimension.  The comparison
tests are made over $A$ before the Rokhlin towers are chosen.

The tower approximation transfers these tests to the crossed
product.  Let a countably infinite discrete amenable group act
minimally on a unital separable \Cstar{}-algebra $A$ with
$\mathrm T(A)\ne\emptyset$.  If the action has WRP and COS and
$B=A\rtimes_{\alpha,r}\Gamma$ is exact, then
\[
 \operatorname{rc}(A\rtimes_{\alpha,r}\Gamma)
       \leq\operatorname{orc}(A,\alpha)
 \qquad\text{(Theorem~\ref{orbit-radius-crossed-product}).}
\]
For $A=\mathrm C(X)$, with $X$ compact and metrizable, the relation to mean
dimension is made explicit by
\[
 \operatorname{orc}(\mathrm C(X),\alpha)
       \leq\tfrac12\operatorname{mdim}(X,\Gamma)
 \qquad\text{(Proposition~\ref{orbit-radius-mdim}).}
\]
This bound itself requires neither freeness nor minimality nor a
Rokhlin hypothesis.  In the noncommutative setting, a fixed finite
coefficient comparison radius becomes negligible after normalization
by growing orbit size: if $A$ is stably finite and exact and
$\operatorname{rc}(A)<\infty$, then
$\operatorname{orc}(A,\alpha)=0$ for an infinite amenable group
(Proposition~\ref{orbit-finite-radius-zero}).  This parallels the
vanishing of mean dimension for actions on finite-dimensional
compact spaces.  The averaging statement does not require WRP
or COS; those hypotheses control the passage to the crossed product.

There are also tower conditions that supply the required complement
comparison directly.  WRP with Cuntz control in the crossed product
supplements WRP by requiring sufficiently invariant central sequence towers whose
complement is Cuntz below any prescribed nonzero positive element
of $B$, with the comparison made in $B_\omega$
(Definition~\ref{cuntz-controlled-wrp-definition}).  Strong URP
uses central projection towers covering the unit, so their
complement is zero
(Definition~\ref{strong-urp-definition} and
Lemma~\ref{strong-urp-cuntz-controlled-wrp}).  For the minimal
actions above, either condition together with zero orbit radius
gives $(1,1)$-COS and, when $B$ is exact, strict comparison
(Proposition~\ref{cuntz-controlled-wrp-orbit-zero-cos}).  Thus these
stronger tower conditions can replace a separate COS assumption.

Combining these comparison results with
Theorem~\ref{intro-sr1-wrp-cos-D} gives strict comparison and stable
rank one under the hypotheses below.  Stable rank one supplies
positive elements realizing the strictly positive lower
semicontinuous affine rank functions.  Locally finite nuclear dimension then allows the
Jiang--Su stability criterion of the preceding subsections to be
applied.  The coefficient algebra need not already be
$\mathcal Z$-stable.

\begin{introthm}[Corollary~\ref{wrp-finite-radius-Z}]
\label{intro-wrp-finite-radius-Z}
Let $A$ be a separable, unital, stably finite \Cstar{}-algebra with
$\mathrm T(A)\neq\varnothing$, locally finite nuclear dimension, and
$\operatorname{rc}(A)<\infty$.  Let a countably infinite discrete
amenable group $\Gamma$ act minimally on $A$.  
Suppose that one
of the following holds:
\begin{enumerate}
\item the action has WRP and $(\lambda,m)$-COS for some
      $\lambda\in(0,1]$ and $m\in\mathbb N$;
\item the action has strong URP;
\item the action has WRP with Cuntz control in the crossed product.
\end{enumerate}
Then
\[
 B=A\rtimes_{\alpha,r}\Gamma
\]
is simple and $\mathcal Z$-stable.
\end{introthm}

In particular, Theorem~\ref{intro-wrp-finite-radius-Z} gives a
\emph{noncommutative extension of the finite-dimensional-space
case} of \cite[Theorem~4.9]{Niu-MD-Z-absorbing}.  If $A=\mathrm C(X)$
with $X$ finite-dimensional, compact, and metrizable, then $A$ has
finite nuclear dimension and finite radius of comparison, by
\cite[Proposition~2.4]{WZ-ndim} and \cite{EN-RC}.
Topological URP agrees with WRP, so the theorem recovers
$\mathcal Z$-stability for free minimal actions of countably
infinite discrete amenable groups with URP and COS.  The general
statement replaces finite covering dimension by finite coefficient
comparison radius together with locally finite nuclear dimension.

The coefficient algebra may therefore be nonsimple, have positive
finite radius of comparison, and have infinite nuclear dimension;
no common bound is imposed on the nuclear dimensions of its local
approximating algebras.  For infinite groups, the theorem also
extends Theorem~\ref{intro-WRP-finite-nucdim-Z}: its
$\mathcal Z$-stable coefficient algebras have finite nuclear
dimension and comparison radius zero, simplicity makes the action
minimal, and Corollary~\ref{Z=>COS} supplies COS.
Theorem~\ref{intro-WRP-finite-nucdim-Z} additionally covers finite
groups.

\section{Notation and preliminaries}\label{section-preliminaries}

We fix the notation for the approximation and comparison arguments
that follow.  We begin with positive elements, traces, and rank
functions, then specify the conventions for actions and central
sequences.  The final subsection recalls the rank and regularity
properties used in the main results.  Hypotheses such as simplicity, nuclearity, and
the existence of tracial states are specified in each statement.  The coefficient algebras in our dynamical results are
unital; the local models and hereditary subalgebras may be nonunital.

\subsection{\Cstar{}-algebras, cutoffs, and order zero maps}

For a \Cstar{}-algebra $A$, let $A_+$ and $A_{\mathrm{sa}}$ denote its
positive cone and self-adjoint part.  A positive contraction is an
element of $A_+$ with norm at most one.  We write $\widetilde A=A$
when $A$ is unital and use the minimal unitization otherwise.
For a unital algebra, $\mathrm U(A)$ and $\mathrm{GL}(A)$ are the unitary
group and the group of invertible elements.  The hereditary
subalgebra generated by $a\in A_+$ is
\[
 \operatorname{Her}(a)=\overline{aAa}.
\]
We write $\mathcal K$ for the compact operators on a separable
infinite-dimensional Hilbert space, $\mathrm M_n(A)$ for the algebra of
$n$ by $n$ matrices over $A$, and
$\mathrm M_\infty(A)=\bigcup_n\mathrm M_n(A)$ for the algebraic union under corner
embeddings $x\mapsto\operatorname{diag}(x,0)$.  Tensor products are
minimal unless specified otherwise.  The notation $x\in_\eps C$
means $\operatorname{dist}(x,C)<\eps$, and $[x,y]=xy-yx$.
The multiplier algebra of $D$ is denoted by $\operatorname{Mult}(D)$.

For $a\in A_+$ and $\eps>0$, the positive part $(a-\eps)_+$ is
defined by continuous functional calculus.  The approximation
arguments repeatedly use the following continuous cutoff, for $s>0$:
\begin{equation}\label{defn-g}
 g_s(t)=
 \begin{cases}
  0,&0\leq t\leq s/2,\\
  2t/s-1,&s/2\leq t\leq s,\\
  1,&t\geq s.
 \end{cases}
\end{equation}
In particular, $g_s(a)\in\operatorname{Her}(a)$ is a positive
contraction.  When $A$ is nonunital, functional calculus by a
function vanishing at zero takes values in $A$.

The matrix structure of positive compressions will often be
expressed through order zero maps.  A map is \emph{c.p.c.} if it is
completely positive and contractive.  A completely positive map
$\phi:D\to A$ has \emph{order zero} if
$\phi(x)\phi(y)=0$ whenever $x,y\in D_+$ and $xy=0$.
For a c.p.c. order zero map with unital domain, the structure theorem
of \cite[Theorem~3.3]{WZ-order-zero} gives
\[
 \phi(x)=h\pi(x)=\pi(x)h,\qquad h=\phi(1_D),
\]
where $\pi:D\to\operatorname{Mult}(\mathrm C^*(\phi(D)))$ is a supporting
*-homomorphism.  Equivalently, such a map determines a
*-homomorphism from $\mathrm C_0((0,1])\otimes D$, taking $t\otimes x$ to
$\phi(x)$, where $t$ is the coordinate function on $(0,1]$
\cite[Corollary~4.1]{WZ-order-zero}.  A unital order zero map is a
*-homomorphism.  We use the order zero and cone descriptions
interchangeably.

\subsection{Traces and Cuntz comparison}

The comparison estimates use ranks measured by traces.  We write
$\mathrm T(A)$ for the tracial state space of $A$.  For unital $A$
this space is weak* compact, and it is metrizable when $A$ is
separable.  We extend a tracial state to $\mathrm M_n(A)$ using the
\emph{unnormalized} matrix trace, and to $A\otimes\mathcal K$ as
the associated lower semicontinuous trace.  The resulting rank
function is
\[
 \mathrm d_\tau(a)=\lim_{m\to\infty}\tau(a^{1/m})\in[0,\infty],
 \qquad a\in(A\otimes\mathcal K)_+.
\]
Thus $\mathrm d_\tau(1_{\mathrm M_n(A)})=n$ when $A$ is unital.  For
$a\in \mathrm M_n(A)_+$, the rank is at most $n$; for a projection $p$,
$\mathrm d_\tau(p)=\tau(p)$.  Every tracial state on a simple unital algebra
is faithful.

We compare positive elements using Cuntz subequivalence.  For
$a,b\in(A\otimes\mathcal K)_+$, write $a\precsim b$ if there are
$v_j\in A\otimes\mathcal K$ such that $v_j^*bv_j\to a$ in norm.
Write $a\sim b$ when both comparisons hold, and use a subscript to
specify the ambient algebra when needed.  The Cuntz semigroup
$\mathrm{Cu}(A)$ consists of these equivalence classes, with order
induced by $\precsim$ and addition by orthogonal direct sum.  Its
finite-matrix version is denoted by $\mathrm W(A)$;
$\mathrm V(A)$ denotes the Murray--von Neumann semigroup of
projections in $\mathrm M_\infty(A)$.  We use $[a]$ for a Cuntz class and,
when the context specifies $\mathrm V(A)$, for a projection class.

The following cutoff facts allow norm approximations to be used
in comparison:
\[
 \|a-b\|<\eps\ \Longrightarrow\ (a-\eps)_+\precsim b,
 \qquad
 a\precsim b\ \Longleftrightarrow\
 (a-\eps)_+\precsim b\ \text{for every }\eps>0;
\]
see \cite[Proposition~2.1(ii)--(iii)]{BPT-Cuntz}.
Increasing sequences in $\mathrm{Cu}(A)$ have suprema.  The relation
$x\ll y$ means that, whenever $y\leq\sup_n y_n$ for an increasing
sequence, one has $x\leq y_n$ for some $n$.  In particular,
$[(a-\eps)_+]\ll[a]$.  We refer to \cite{CEI-CuntzSG} for the
ordered-semigroup conventions.

Two comparison properties will be used repeatedly.  An ordered
semigroup is \emph{almost unperforated} if $(n+1)x\leq ny$ implies
$x\leq y$ for every positive integer $n$.  A Cuntz semigroup $S$ is
\emph{almost divisible} if, for every $x'\ll x$ and $k\geq1$,
there is $y\in S$ such that
\[
 ky\leq x,\qquad x'\leq(k+1)y.
\]
It is \emph{pure} if it has both properties.

\begin{defn}[Strict comparison]\label{strict-comparison-definition}
Let $A$ be a unital \Cstar{}-algebra with $\mathrm T(A)\ne\emptyset$.
We say that $A$ has \emph{strict comparison by tracial states} if
\[
 \mathrm d_\tau(a)<\mathrm d_\tau(b)\quad(\tau\in\mathrm T(A))
 \quad\Longrightarrow\quad a\precsim b,
 \qquad a,b\in \mathrm M_\infty(A)_+.
\]
The phrase \emph{strict comparison by tracial states on $A_+$}
denotes the same implication restricted to $a,b\in A_+$.
Unless expressly qualified, strict comparison in this paper is
understood to be by tracial states.
\end{defn}

The usual quasitracial formulation replaces $\mathrm T(A)$ by the
normalized bounded $2$-quasitraces on $A$ and uses the same
unnormalized matrix extensions.  For exact algebras these
$2$-quasitraces are traces \cite[Theorem~5.11]{Haagtrace}, so the
two formulations agree.  For a unital stably finite exact algebra
with nonempty tracial state space, the radius of comparison
$\operatorname{rc}(A)$ is the infimum of $r\geq0$ such that
\[
 \mathrm d_\tau(a)+r<\mathrm d_\tau(b)\quad(\tau\in\mathrm T(A))
 \quad\Longrightarrow\quad a\precsim b,
 \qquad a,b\in \mathrm M_\infty(A)_+.
\]
The infimum is $\infty$ if no such $r$ exists.  Gap functions on the
unnormalized trace cone, and the pointed comparison radius for
nonunital models, are introduced in
Chapter~\ref{section-comparison}.

\subsection{Affine functions and rank functions}

Rank functions connect comparison with affine functions on the
trace space.  For a compact convex set $S$, $\mathrm{Aff}(S)$ is the
real Banach space of continuous affine functions, with the uniform
norm.  The set $\mathrm{LAff}_{++}(S)$ consists of strictly positive
lower semicontinuous affine functions $S\to(0,\infty]$.
The function $\widehat a(\tau)=\mathrm d_\tau(a)$ is the rank function of a
positive element; for a projection $p$,
$\widehat p(\tau)=\tau(p)$.
We denote the \Ksym{}-groups by $\mathrm K_0(A)$ and $\mathrm K_1(A)$.  The trace pairing
on $\mathrm K_0(A)$ is
\[
 \rho_A([p]-[q])(\tau)=\tau(p)-\tau(q).
\]
We write $\partial_eS$ for the extreme boundary of $S$ and $\dim$
for covering dimension.  A Bauer simplex is a Choquet simplex whose
extreme boundary is compact.

For a simple unital stably finite \Cstar{}-algebra with strict comparison
by tracial states, the expression
$\mathrm V(A)\sqcup\mathrm{LAff}_{++}(\mathrm T(A))$ carries the
usual addition and order on each summand, together with
\[
 \begin{aligned}
 [p]+f&=\widehat p+f,\\
 [p]\leq f&\ \Longleftrightarrow\
       \widehat p(\tau)<f(\tau)\quad\text{for every }\tau\in\mathrm T(A),\\
 f\leq[p]&\ \Longleftrightarrow\ f\leq\widehat p.
 \end{aligned}
\]
When we assert a \emph{canonical} Cuntz-semigroup identification
with this model, projection classes retain their projection data
and nonprojection classes are represented by their rank functions.
The existence of such an identification, including its comparison
and rank-realization conclusions, is proved under the hypotheses
of the relevant theorem.
See \cite{BPT-Cuntz,CEI-CuntzSG} for the underlying conventions.

\subsection{Actions, crossed products, and tilings}

We now fix the dynamical conventions used to form the tower models.
Groups $\Gamma$ are discrete, with identity $e$.  The Rokhlin and
orbit-radius results concern countable amenable groups, whereas the
coefficient comparison criterion in Section~\ref{cos-section}
allows arbitrary discrete groups.  An action is written
$\alpha:\Gamma\to\operatorname{Aut}(A)$, with
$\alpha_g\alpha_k=\alpha_{gk}$.  We also write $g(a)=\alpha_g(a)$.
For the reduced crossed product $B=A\rtimes_{\alpha,r}\Gamma$,
the canonical unitaries satisfy
\[
 u_gu_k=u_{gk},\qquad u_gau_g^*=\alpha_g(a),
 \qquad
 \mathbb E\Big(\sum_g a_gu_g\Big)=a_e.
\]
Here $\mathbb E:B\to A$ is the canonical faithful conditional
expectation.  The algebra $\mathrm C_{\mathrm c}(\Gamma,A)$ of finite Fourier sums is
dense in $B$.  For amenable groups the full and reduced crossed
products agree, and we sometimes omit the subscript $r$.
The action is \emph{minimal} if $A$ has no nonzero proper
$\alpha$-invariant closed two-sided ideals.

The traces relevant to the action form the invariant trace space,
for which we use the notation
\[
 S=\mathrm T^\alpha(A)=\mathcal M_1(A,\Gamma)
   =\{\tau\in\mathrm T(A):\tau\circ\alpha_g=\tau
                                      \text{ for all }g\in\Gamma\}.
\]
If $A$ is unital, $\mathrm T(A)\ne\emptyset$, and $\Gamma$ is
amenable, then $S\ne\emptyset$.  For a finite nonempty
$F\subseteq\Gamma$, the averaging convention is
\[
 \tau_F=\frac1{|F|}\sum_{g\in F}\tau\circ\alpha_{g^{-1}}.
\]
For a left action on a compact space $X$, the induced action on $\mathrm C(X)$
is $\alpha_g(f)(x)=f(g^{-1}x)$, and invariant tracial states are
invariant Borel probability measures.

To specify tower shapes, for finite $K,F\subseteq\Gamma$ with
$e\in K$, put
\[
 \operatorname{Int}_K(F)=\{g\in F:Kg\subseteq F\}.
\]
A nonempty finite $F$ is $(K,\eps)$-invariant if
$|F\setminus\operatorname{Int}_K(F)|<\eps|F|$.
An exact tiling with shapes $D_1,\ldots,D_R$ is a disjoint
decomposition
\[
 \Gamma=\bigsqcup_{s=1}^R\bigsqcup_{c\in\mathcal C_s}D_sc,
\]
where the shapes are finite and nonempty and the sets
$\mathcal C_s$ are their centers.  We may normalize shapes to contain
$e$.  Thus the tiles are \emph{right} translates, whereas the
invariance and boundary estimates use \emph{left} multiplication.
The prescribed-tiling quantifier is part of WRP in
Definition~\ref{D:weak-Rokhlin-property}.

\subsection{Ultrapowers and trace-kernel ideals}

The tower constructions use central sequences and uniform tracial
estimates.  Fix a free ultrafilter $\omega$ on $\mathbb N$.
Throughout the paper, $A_\omega$ denotes the \emph{norm} ultrapower
\[
 A_\omega=\ell^\infty(\mathbb N,A)/c_\omega(A),\qquad
 c_\omega(A)=\{(a_n):\lim_{n\to\omega}\|a_n\|=0\}.
\]
Constant sequences identify $A$ with a subalgebra of $A_\omega$;
its central sequence algebra is $A_\omega\cap A'$.
Actions induce coordinatewise actions on these algebras.
For nonempty $S\subseteq\mathrm T(A)$, define
\[
 \|a\|_{2,S}=\sup_{\tau\in S}\tau(a^*a)^{1/2},\qquad
 \|[(a_n)]\|_{2,S}=\lim_{n\to\omega}\|a_n\|_{2,S},
 \qquad
 J_S=\{x\in A_\omega:\|x\|_{2,S}=0\}.
\]
Taking the supremum before the ultralimit makes this condition
uniform over $S$, including sequences of varying traces.
The set $J_S$ is a closed two-sided ideal, and
$q_S:A_\omega\to A_\omega/J_S$ is the quotient map.
Each sequence $\tau_n\in S$ defines a limit trace
$\tau_\omega([(a_n)])=\lim_{n\to\omega}\tau_n(a_n)$, which
annihilates $J_S$.

The lifting arguments also use the following ideal property.
An ideal $I\triangleleft D$ is a $\sigma$-ideal if, for every
separable subalgebra $C\subseteq D$, there is a positive contraction
$z\in I\cap C'$ with $zb=b$ for $b\in I\cap C$.
The trace-kernel ideals $J_S$ have this property
\cite[Remark~4.7]{KR-CenSeq}; it permits lifts from quotient relative
commutants that commute with a prescribed separable subalgebra
\cite[Proposition~4.5(iii)]{KR-CenSeq}.
The detailed lifting statements are collected in
Proposition~\ref{trace-kernel-sigma}.

For von Neumann algebras we use a tracial ultrapower in the finite
case.  If $M$ is finite with faithful normal tracial state $\tau$,
then
\[
 M^\omega=\ell^\infty(\mathbb N,M)/
       \{(x_n):\lim_{n\to\omega}\|x_n\|_{2,\tau}=0\},
 \qquad M_\omega=M'\cap M^\omega.
\]
For general von Neumann algebras, $M^\omega$ denotes the Ocneanu
ultrapower and $M_\omega$ the asymptotic centralizer, as defined in
Definition~\ref{ocneanu-ultrapowers}.
For $\tau\in\mathrm T(A)$, write $M_\tau=\pi_\tau(A)''$ and
$\bar\tau$ for its faithful normal trace.
If $\tau\in\mathrm T^\alpha(A)$, the action induces a normal
trace-preserving action $\alpha_\tau$ on $M_\tau$, with
$(\alpha_\tau)_g(\pi_\tau(a))=\pi_\tau(\alpha_g(a))$.
\Wstar{}-central freeness means central freeness of every such action;
see Definition~\ref{wstar-freeness-cstar}.

\subsection{Rank and regularity conventions}

Finally, we recall the regularity properties used in the rank and
absorption results.  Real rank zero, written $\mathrm{RR}(A)=0$,
means that finite-spectrum self-adjoint elements are dense in
$\widetilde A_{\mathrm{sa}}$.  Stable rank one, written
$\operatorname{tsr}(A)=1$, means that
$\mathrm{GL}(\widetilde A)$ is dense in $\widetilde A$.
The Jiang--Su algebra is denoted by $\mathcal Z$, and
$\mathcal Z$-stability means $A\cong A\otimes\mathcal Z$.
The UCT is the universal coefficient theorem for \KKsym{}-theory.

We use nuclear dimension in the sense of
\cite[Definition~2.1]{WZ-ndim}.  The condition
$\dim_{\mathrm{nuc}}(A)\leq d$ means that finite sets can be
approximated by completely positive factorizations
\[
 A\xrightarrow{\ \psi\ }F_0\oplus\cdots\oplus F_d
   \xrightarrow{\ \phi\ }A,
\]
where each $F_i$ is a finite-dimensional \Cstar{}-algebra, $\psi$ is c.p.c., and each
$\phi|_{F_i}$ is c.p.c. and order zero.  The sum map $\phi$ need
not be contractive.  An algebra has \emph{locally finite nuclear
dimension}, or is \emph{locally approximated by finite nuclear
dimension algebras}, if, for every finite set and tolerance, a
\Cstar{}-subalgebra of finite nuclear dimension approximates every
element of the set within that tolerance.  The dimension bound
may depend on the finite set and the tolerance.

\chapter{Weak tracial approximation properties}\label{section-approximation}

\section{Definitions and equivalent formulations}

For a unital \Cstar{}-algebra $A$, write $\mathrm T(A)$ for its tracial state
space.  If $\tau\in\mathrm T(A)$ and $a\in A_+$, set
\[
 \mathrm d_\tau(a)=\lim_{m\to\infty}\tau(a^{1/m}).
\]
For $a,b\in A_+$, the notation $a\precsim b$ means Cuntz
subequivalence: there are $v_j\in A$ such that $v_j^*bv_j\to a$ in
norm.  We write $x\in_\eps C$ when
$\operatorname{dist}(x,C)<\eps$.  Classes of \Cstar{}-algebras are understood
up to isomorphism; the approximating subalgebras need not be unital.

\begin{defn}\label{defn_WAS}
Let $\mathcal S$ be a class of \Cstar{}-algebras, and let $A$ be unital.
We say that $A$ is \emph{weakly tracially approximated by $\mathcal S$},
and write $A\in\mathrm{TA}_0(\mathcal S)$, if
$\mathrm T(A)\neq\emptyset$ and, for every finite set
$\mathcal F=\{x_1,\ldots,x_n\}\subseteq A$ and every $\eps>0$, there
are a \Cstar{}-subalgebra $C\subseteq A$ with $C\in\mathcal S$, elements
$x_1',\ldots,x_n'\in A$, and a positive contraction $h\in C$ with
$\|h\|=1$ such that
\begin{enumerate}
\item[(1)] $\|x_i-x_i'\|<\eps$, $1\leq i\leq n$;
\item[(2)] $\|[h,x_i']\|<\eps$, $1\leq i\leq n$;
\item[(3)] $hx_i'h\in C$, $1\leq i\leq n$;
\item[(4)] $\mathrm d_\tau(1-h)<\eps$ for every $\tau\in\mathrm T(A)$.
\end{enumerate}

We say that $A$ is \emph{tracially approximated by $\mathcal S$},
and write $A\in\mathrm{TA}_1(\mathcal S)$, if, for every such
$\mathcal F$ and $\eps$ and every nonzero positive contraction
$a_0\in A$, the same data can be chosen to satisfy (1)--(3) and
\[
 \text{(4')}\qquad 1-h\precsim a_0.
\]
The definition of $\mathrm{TA}_1(\mathcal S)$ does not require
$\mathrm T(A)\neq\emptyset$.
\end{defn}

\begin{example}[Motivation from topological dynamics]\label{Rokhlin-example}
Topological Rokhlin towers give the positive compressions in
Definition~\ref{defn_WAS}.
Let a countably infinite amenable group $\Gamma$ act on a nonempty
compact metrizable space $X$, and let $\mathcal C_X$ be the class
of finite direct sums of algebras
$\mathrm M_n(\mathrm C_0(U))$, where $U\subseteq X$ is open.
If the action has the topological uniform Rokhlin property, then
\[
 \mathrm C(X)\rtimes\Gamma\in\mathrm{TA}_0(\mathcal C_X).
\]
By Theorem~\ref{central-topological-urp} and
Proposition~\ref{urp-implies-wrp}, topological URP implies WRP.
Theorem~\ref{WTA-C-dyn-1} then gives approximation by matrix algebras
over hereditary subalgebras of $\mathrm C(X)$.
For $b\in\mathrm C(X)_+$,
\[
 \overline{b\mathrm C(X)b}
 =\mathrm C_0\bigl(\{x\in X:b(x)>0\}\bigr).
\]
Amenability ensures the existence of invariant probability measures,
which induce tracial states of the crossed product.  The earlier
tracial approximation construction appears in
\cite[Theorem~3.9]{Niu-MD-Z}.

An open tower base need not be compact, and its algebra $\mathrm C_0(U)$
may therefore be nonunital.  Positive functions supported on the
towers define the compression.  For each invariant probability
measure, the dimension function of the defect is the measure of
the set on which the compression is less than one.
Thus Definition~\ref{defn_WAS} accommodates
these nonunital models and positive compressions.

If the action is also free and minimal and has Cuntz comparison of
open sets (COS), then
$\mathrm C(X)\rtimes\Gamma\in\mathrm{TA}_1(\mathcal C_X)$.
To prove the required comparison, put
$B=\mathrm C(X)\rtimes\Gamma$, fix $0\neq a_0\in B_+$, and
let $(\lambda,m)$ be COS constants as in
\cite[Definition~4.1]{Niu-MD-Z}.
The algebra $B$ is simple and infinite dimensional.
Choose $m$ nonzero orthogonal Cuntz equivalent positive elements
$b_1,\ldots,b_m\in\overline{a_0Ba_0}$, using
\cite[Lemma~2.4]{NCP-LAlg}.
By \cite[Lemma~4.2]{Niu-MD-Z}, there is $0\neq c\in\mathrm C(X)_+$
with $c\precsim b_1$.  Consequently
\[
 c\otimes1_m\precsim\sum_{j=1}^m b_j\precsim a_0.
\]
For $V=\{c>0\}$, minimality gives
$\beta=\inf_\mu\mu(V)>0$, where $\mu$ ranges over invariant
probability measures: finitely many translates of $V$ cover $X$.
Choose the approximation data in Theorem~\ref{WTA-C-dyn-1} so that
$\mu(\{h<1\})<\lambda\beta$ for every such $\mu$.
This is possible because $\tau_\mu=\mu\circ\mathbb E$ is a tracial
state of $B$ and
$\mathrm d_{\tau_\mu}(1-h)=\mu(\{h<1\})$.
Since $h\in\mathrm C(X)$, COS then gives
\[
 1-h\precsim c\otimes1_m\precsim a_0.
\]
The $\mathrm{TA}_1$ conclusion follows from
Lemma~\ref{alt-defn-was}.
\end{example}

\begin{example}[Classifiable \Cstar{}-algebras]\label{classifiable-interval-example}
Let $\mathcal I_{\mathrm h}$ and $\mathcal I_{\mathrm o}$ be the
classes of finite direct sums of matrix algebras over $\mathrm C_0([0,1))$
and $\mathrm C_0((0,1))$, respectively.
Every unital, separable, simple, infinite-dimensional, nuclear,
stably finite, $\mathcal Z$-stable \Cstar{}-algebra $A$ satisfying the UCT
belongs to
\[
 \mathrm{TA}_0(\mathcal I_{\mathrm h})
       \cap\mathrm{TA}_0(\mathcal I_{\mathrm o}).
\]
Thus every such classifiable \Cstar{}-algebra admits weak tracial
approximation by either half-open or open interval algebras.
Theorem~\ref{classifiable-interval-approximation} gives the stronger
conclusion
\[
 A\in\mathrm{TA}_1(\mathcal I_{\mathrm h})
       \cap\mathrm{TA}_1(\mathcal I_{\mathrm o}),
\]
and Remark~\ref{ta1-to-ta0} gives the stated weak tracial approximation.
The proof uses the Gong--Lin--Niu models, whose traces concentrate
on the point quotients.  Compression to endpoint collars gives
half-open interval models, and $\mathcal Z$-absorption gives open
interval models; see Section~\ref{interval-models-section}.
\end{example}

\begin{rem}[Comparison with generalized tracial approximation]\label{comparison-gta}
Let $A$ be simple and unital with $\mathrm T(A)\neq\emptyset$.
Then
\[
 A\in\mathrm{TA}_1(\mathcal S)
 \quad\Longleftrightarrow\quad A\in\mathrm{gTA}(\mathcal S),
\]
where $\mathrm{gTA}$ is the nonunital generalized tracial approximation
property of Elliott--Fan--Fang
\cite[Definition~1.2 and Theorem~2.11]{Elliott-Fan-Fang-nonunital-GTA}.
When $\mathcal S$ consists of unital algebras, these conditions also
agree with their $\mathrm{WTA}(\mathcal S)$
\cite[Definition~1.1]{Elliott-Fan-Fang-GTA}; see also
\cite[Theorem~2.12]{Elliott-Fan-Fang-nonunital-GTA}.
In that formulation, the approximating algebra has unit $p=1_C$,
and the compression is $p-g$, with $0\leq g\leq p$.
The property called weak tracial approximation in that reference
therefore corresponds to $\mathrm{TA}_1$.

We verify the equivalence, including the one-sided approximation
and norm preservation conditions in $\mathrm{gTA}$.  Include the
prescribed nonzero positive contraction $a_0$ in the finite set.
By Lemma~\ref{almost-pres-norm}, the $\mathrm{TA}_1$ approximation
data can be chosen with $\|ha_0h\|$ arbitrarily close to
$\|a_0\|$.
For $0<s<1$, put $k=g_s(h)$ and $\ell=g_{s/2}(h)$.
Then $k\ell=k$ and
\[
 kx_i'\ell=
 \frac{g_s(h)}{h}\,(hx_i'h)\,
 \frac{g_{s/2}(h)}{h}\in C,
\]
where the quotients are continuous functions extended by zero at
the origin.  Likewise $\ell x_i'k\in C$.
It follows that the distances of $kx_i'$ and $x_i'k$ from $C$ are
bounded by $\|[\ell,x_i']\|$.  These commutators are small by
continuous functional calculus.
Moreover,
\[
 1-k\leq(1-s)^{-1}(1-h),\qquad
 \|ka_0k\|\geq\|ha_0h\|-2s\|a_0\|.
\]
The first estimate gives $1-k\precsim a_0$.  The second follows
from $\|h(1-k)\|\leq s$ and proves the required norm estimate.
For sufficiently small $s$ and initial tolerance, these estimates
give the formulation of $\mathrm{gTA}$ in
\cite[Theorem~2.11]{Elliott-Fan-Fang-nonunital-GTA}.
They also imply the conditions of Definition~1.2, since, for every
$y\in A_+$,
\[
 (y^2-yky-\eps)_+\precsim y(1-k)y\precsim1-k\precsim a_0.
\]
Rescaling gives the same conclusion for arbitrary nonzero positive
comparison elements.

Conversely, suppose that $A$ satisfies Definition~1.2.
Fix $0\neq a_0\in A_+$, put
$a=a_0/\|a_0\|$, and use $y=1_A$ in Definition~1.2 with a
tolerance $0<\delta<1/4$.  There are an approximating algebra
$C\in\mathcal S$ and a positive contraction $k_0\in C$ such that
\[
 (1-k_0-\delta)_+\precsim a,\qquad
 \|k_0ak_0\|\geq1-\delta,
\]
and the one-sided approximation and commutator estimates hold
with tolerance $\delta$.  Set
\[
 h=\min\{k_0/(1-\delta),1_A\}.
\]
Continuous functional calculus defines this minimum, and the
scalar function vanishes at zero.  Thus $h\in C$, including when
$C$ is nonunital.  Since
$\|k_0\|\geq\sqrt{1-\delta}>1-\delta$, we have $\|h\|=1$.
Moreover, $\|h-k_0\|\leq\delta$, and
\[
 1-h=(1-\delta)^{-1}(1-k_0-\delta)_+\precsim a\sim a_0.
\]
Multiplication by $k_0$ of an approximation to $k_0x$ in $C$ gives
$\operatorname{dist}(k_0xk_0,C)<\delta$.  Hence, for every
$x$ in the prescribed finite set,
\[
 \|[h,x]\|<\delta(1+2\|x\|),\qquad
 \operatorname{dist}(hxh,C)<\delta(1+2\|x\|).
\]
Taking $\delta$ sufficiently small and applying
Lemma~\ref{alt-defn-was} proves $A\in\mathrm{TA}_1(\mathcal S)$.
When $C$ is unital, setting $g=p-h$ identifies these conditions
with those of $\mathrm{WTA}$.

The defect in $\mathrm{TA}_0$ is controlled by dimension functions.
Remark~\ref{ta1-to-ta0} gives
$\mathrm{TA}_1\Rightarrow\mathrm{TA}_0$.
If $A$ has strict comparison by traces, the reverse implication
also holds. For a nonzero positive contraction $a_0$, faithfulness
and compactness give
\[
 m=\min_{\tau\in\mathrm T(A)}\tau(a_0)>0,\qquad
 \mathrm d_\tau(a_0)\geq\tau(a_0)\geq m.
\]
Choosing $h$ so that $\mathrm d_\tau(1-h)<m$ for every $\tau$ gives
$1-h\precsim a_0$.  Thus strict comparison by tracial states on
$A_+$ suffices for the reverse implication.
\end{rem}

\begin{rem}\label{ta1-to-ta0}
If $A$ is simple and $\mathrm T(A)\neq\emptyset$, then
\[
 A\in\mathrm{TA}_1(\mathcal S)
 \quad\Longrightarrow\quad A\in\mathrm{TA}_0(\mathcal S).
\]
Moreover, the $\mathrm{TA}_1$ approximation data can be chosen
with uniformly small dimension-function defect and
$1-h\precsim a_0$ for a prescribed $a_0$.

Suppose first that $A$ is infinite dimensional.  Fix $a_0\neq0$ and
$\eta>0$, and choose $N$ with $1/N<\eta$.
Such an $A$ is not of type I, since a simple unital type-I algebra
is a matrix algebra.
By \cite[Lemma~2.4]{NCP-LAlg}, there are nonzero, pairwise orthogonal,
Cuntz equivalent positive elements $b_1,\ldots,b_N$ in
$\overline{a_0Aa_0}$.  Thus $b_1\precsim a_0$ and
\[
 N \mathrm d_\tau(b_1)
 =\mathrm d_\tau\left(\sum_{j=1}^N b_j\right)\leq1,
 \qquad \tau\in\mathrm T(A).
\]
After rescaling $b_1$ to be a contraction, apply the
$\mathrm{TA}_1$ property with $b_1$ in place of $a_0$.
The resulting defect satisfies both $1-h\precsim a_0$ and
$\mathrm d_\tau(1-h)<\eta$ for every $\tau$.

Now suppose $A\cong\mathrm M_r$.
Choose a finite family $\mathcal U$ of unitaries whose conjugation
average is $x\mapsto\operatorname{tr}_r(x)1_A$, where
$\operatorname{tr}_r$ is the normalized matrix trace.
Apply the approximation property to $\mathcal F\cup\mathcal U$ with
tolerance $\delta<\min\{\eps,1/12\}$.  Conditions (1) and (2) give
$\|[h,u]\|<3\delta$ for $u\in\mathcal U$, so
\[
 \|h-\operatorname{tr}_r(h)1_A\|<3\delta.
\]
Since $\|h\|=1$, we have $\operatorname{tr}_r(h)>1-3\delta$,
and therefore $h>(1-6\delta)1_A>\frac12 1_A$.
Consequently $g_{1/2}(h)=1_A\in C$, and $h^{-1}\in C$.
Condition (3) then gives $x_i'\in C$ for the original finite set.
Replacing $h$ by $1_A$ makes the defect zero and preserves the
required approximation.
\end{rem}

\begin{rem}\label{trace-rank-ta0}
In Definition~\ref{defn_WAS}, the $\mathrm{TA}_0$ condition (4) can
equivalently be replaced by
\[
 \tau(1-h)<\eps,\qquad \tau\in\mathrm T(A).
\]
The equivalence allows a decrease in the initial tolerance; the
two estimates need not be equivalent for a fixed $h$.

The dimension-function estimate implies the trace estimate because
$0\leq1-h\leq1_A$.  Conversely, choose the approximation data
with a smaller trace tolerance and put $k=g_{1/2}(h)$.  Then
$k\in C$ and
$\|k\|=1$.  In the tracial von Neumann algebra associated to $\tau$,
\[
 \mathrm d_\tau(1-k)
 =\tau\bigl(1_{[0,1/2)}(h)\bigr)\leq2\tau(1-h).
\]
Define $r(t)=g_{1/2}(t)/t$ for $t>0$, with $r(0)=0$.
This is continuous on $[0,1]$, and
\[
 kx_i'k=r(h)(hx_i'h)r(h)\in C.
\]
The commutators $[k,x_i']$ are small when $[h,x_i']$ are sufficiently
small.  This follows by polynomial approximation to $g_{1/2}$,
using the bound $\|x_i'\|\leq\|x_i\|+1$.
A sufficiently small initial tolerance therefore makes $k$ satisfy
the required conditions.
\end{rem}

\begin{example}
Requiring $C$ to be unital and taking $h=1_C$ gives compression
by the projection $p=1_C\in A$.  With the Cuntz comparison
condition of $\mathrm{TA}_1$, this recovers the tracial approximation
property of \cite{EN-Tapprox}.  The $\mathrm{TA}_0$ condition gives
the corresponding property with tracial control of the defect.
Allowing general positive contractions also permits nonunital
approximating subalgebras, as in Example~\ref{Rokhlin-example}.
\end{example}

The following equivalent formulation uses the original finite set
rather than the auxiliary perturbations $x_i'$.

\begin{lem}\label{alt-defn-was}
Let $A$ be a unital \Cstar{}-algebra.  Then
$A\in\mathrm{TA}_0(\mathcal S)$ if and only if
$\mathrm T(A)\neq\emptyset$ and, for every finite set
$\{x_1,\ldots,x_n\}\subseteq A$ and every $\eps>0$, there are a
\Cstar{}-subalgebra $C\subseteq A$ with $C\in\mathcal S$ and a positive
contraction $h\in C$ with $\|h\|=1$ such that
\[
 \|[h,x_i]\|<\eps,\qquad
 \operatorname{dist}(hx_ih,C)<\eps,\qquad 1\leq i\leq n,
\]
and $\mathrm d_\tau(1-h)<\eps$ for every $\tau\in\mathrm T(A)$.

Likewise, $A\in\mathrm{TA}_1(\mathcal S)$ if and only if, for every
such finite set and $\eps$ and every nonzero positive contraction
$a_0\in A$, the same data can be chosen with $1-h\precsim a_0$ in
place of the dimension-function estimate.
\end{lem}
\begin{proof}
The approximation data in Definition~\ref{defn_WAS} with
tolerance~$\delta$ satisfy
\[
 \|[h,x_i]\|<3\delta,\qquad
 \operatorname{dist}(hx_ih,C)<\delta.
\]
The defect condition is unchanged.

Conversely, fix the finite set and $\eps>0$, and put
$\rho=\eps/16$.  Choose $C$ and $k$ satisfying the stated
conditions with tolerance $\delta$, where $\delta$ is sufficiently
small.
Set
\[
 u=g_{1/4}(k),\qquad v=g_{1/2}(k)^{1/2},\qquad
 w=u^{1/2},\qquad z=(1-u)^{1/2}.
\]
Polynomial approximation in continuous functional calculus allows
$\delta$ to be chosen so that
\[
 \|[w,x_i]\|<\rho,\qquad
 \|[z,x_i]\|<\rho,\qquad
 \|[v,x_i]\|<\rho
\]
for all $i$.  Define the continuous function
\[
 r(t)=
 \begin{cases}
 g_{1/4}(t)^{1/2}/t,&t>0,\\
 0,&t=0.
 \end{cases}
\]
Choose $y_i\in C$ with $\|kx_ik-y_i\|<\delta$ and put
$c_i=r(k)y_i r(k)\in C$.  Since $r(k)k=w$,
\[
 \|c_i-wx_iw\|\leq\|r\|^2\delta<\rho
\]
after decreasing $\delta$ further.  Set
\[
 x_i'=zx_i z+c_i.
\]
Using $w^2+z^2=1_A$ gives
\[
 \|wx_iw+zx_i z-x_i\|
 \leq\|[w,x_i]\|+\|[z,x_i]\|<2\rho.
\]
Hence $\|x_i'-x_i\|<3\rho$ and
$\|[v,x_i']\|<7\rho<\eps$.
Moreover, $vz=0$, so $vx_i'v=vc_iv\in C$, and
$v\in C$ is a positive contraction of norm one.

For $\mathrm{TA}_0$, the zero set of $1-v$ contains that of $1-k$,
so
\[
 \mathrm d_\tau(1-v)\leq \mathrm d_\tau(1-k)<\delta.
\]
For $\mathrm{TA}_1$, the scalar inequality
$1-g_{1/2}(t)^{1/2}\leq2(1-t)$ on $[0,1]$ gives
\[
 1-v\precsim1-k\precsim a_0.
\]
Taking also $\delta<\eps$ proves both converse implications.
\end{proof}

\section{Norm preservation under tracial cut-downs}

A small tracial defect gives a small compression error in the
tracial $2$-norm.  In a simple algebra with tracial states, this
estimate yields norm preservation for the compressions in
Definition~\ref{defn_WAS}.  The corresponding results for
asymptotically order zero maps and relatively central compressions
apply to individual Rokhlin towers.

For a trace $\tau$, write $\|x\|_{2,\tau}=\tau(x^*x)^{1/2}$.

\begin{lem}\label{small-tr-cut}
Let $A$ be unital, let $h\in A$ be a positive contraction, and let
$\tau\in\mathrm T(A)$.  For every $x\in A$,
\[
 \|hxh-x\|_{2,\tau}
 \leq2\|x\|\tau(1-h)^{1/2}
\]
and
\[
 |\tau(hxh)-\tau(x)|
 \leq\|x\|\tau(1-h^2)\leq2\|x\|\tau(1-h).
\]
\end{lem}
\begin{proof}
Write $hxh-x=hx(h-1)+(h-1)x$.
The inequalities
$\|ab\|_{2,\tau}\leq\|a\|\|b\|_{2,\tau}$ and
$\|ab\|_{2,\tau}\leq\|a\|_{2,\tau}\|b\|$ give
\[
 \|hxh-x\|_{2,\tau}
 \leq2\|x\|\|1-h\|_{2,\tau}
 \leq2\|x\|\tau(1-h)^{1/2}.
\]
By traciality,
$\tau(hxh)-\tau(x)=-\tau(x(1-h^2))$.
For $b\geq0$, the positive functional
$x\mapsto\tau(b^{1/2}xb^{1/2})$ has norm $\tau(b)$.
Taking $b=1-h^2$ proves the second estimate, since
$1-h^2\leq2(1-h)$.
\end{proof}

\begin{lem}\label{almost-pres-norm}
Let $A$ be simple and unital with
$\mathrm T(A)\neq\emptyset$, and suppose
$A\in\mathrm{TA}_0(\mathcal S)$ or
$A\in\mathrm{TA}_1(\mathcal S)$.
The approximation data in Definition~\ref{defn_WAS} can also be
chosen so that
\[
 \|hx_i'h\|>\|x_i\|-\eps
 \quad\hbox{and}\quad
 \|hx_i'h\|>\|x_i'\|-\eps,\qquad 1\leq i\leq n.
\]
In the $\mathrm{TA}_1$ case, comparison with the prescribed
$a_0$ is retained.
\end{lem}
\begin{proof}
Fix $\tau_0\in\mathrm T(A)$, which is faithful by simplicity.
Faithfulness implies, for every $x\in A$,
\[
 \lim_{p\to\infty}\tau_0((x^*x)^p)^{1/(2p)}=\|x\|.
\]
Thus, for each $i$, choose an integer $p_i\geq1$ so that
\[
 \tau_0((x_i^*x_i)^{p_i})^{1/(2p_i)}
 >\|x_i\|-\eps/4.
\]

Choose the approximation data with tolerance $\delta$.
In the $\mathrm{TA}_0$ case, $\tau_0(1-h)<\delta$ follows from
condition (4).  In the $\mathrm{TA}_1$ case,
Remark~\ref{ta1-to-ta0} gives the same estimate while preserving
$1-h\precsim a_0$.  Put $y_i=hx_i'h$.
Lemma~\ref{small-tr-cut} gives
\[
 \|y_i-x_i\|_{2,\tau_0}
 \leq\delta+2\|x_i\|\delta^{1/2}.
\]
The elements $x_i$ and $y_i$ lie in a fixed norm-bounded set,
on which the function $x\mapsto\tau_0((x^*x)^p)$ is continuous
in the tracial $2$-norm:
if $\|x\|,\|y\|\leq M$ and $p\geq1$ is an integer, then
\[
 \big|\tau_0((x^*x)^p)-\tau_0((y^*y)^p)\big|
       \leq2pM^{2p-1}\|x-y\|_{2,\tau_0}.
\]
Indeed, $\|x^*x-y^*y\|_{2,\tau_0}\leq2M\|x-y\|_{2,\tau_0}$;
expanding the difference of the $p$th powers and applying
Cauchy--Schwarz to each term proves the estimate.
Therefore, choosing $\delta$ sufficiently small gives
\[
 \|y_i\|\geq
 \tau_0((y_i^*y_i)^{p_i})^{1/(2p_i)}
 >\|x_i\|-\eps/2
\]
for all $i$.  Choosing also $\delta<\eps/2$ and using
$\|x_i'\|<\|x_i\|+\delta$ proves both required inequalities.
\end{proof}

A completely positive map has \emph{order zero} if it preserves
orthogonality of positive elements.  A sequence $(\phi_n)$ of
completely positive contractions is asymptotically order zero if
$\|\phi_n(a)\phi_n(b)\|\to0$ whenever $a,b\geq0$ and $ab=0$.

\begin{lem}\label{almost-ord0-lem}
Let $A$ be unital, let $\omega$ be a free ultrafilter on $\mathbb N$,
and let $\phi_n:A\to D_n$ be completely positive contractions.
Suppose $\tau_n\in\mathrm T(D_n)$ and
\begin{equation}\label{almost-ord0}
 \lim_{n\to\infty}\|\phi_n(a)\phi_n(b)\|=0
 \quad\hbox{whenever }a,b\in A_+\hbox{ and }ab=0,
\end{equation}
\begin{equation}\label{state-asumption}
 \lim_{n\to\infty}\tau_n(\phi_n(1_A))=1.
\end{equation}
If the induced state
\[
 \tau(a)=\lim_{n\to\omega}\tau_n(\phi_n(a)),\qquad a\in A,
\]
is faithful, then
\[
 \lim_{n\to\omega}\|\phi_n(a)\|=\|a\|,\qquad a\in A.
\]
The state $\tau$ is automatically tracial.
If every tracial state of $A$ is faithful (in particular, if $A$
is simple), then the ordinary limit satisfies
\[
 \lim_{n\to\infty}\|\phi_n(a)\|=\|a\|,\qquad a\in A.
\]
\end{lem}
\begin{proof}
Let $\widetilde D_n=D_n$ if $D_n$ is unital, and let
$\widetilde D_n$ be its unitization otherwise.
Extend $\tau_n$ to a tracial state $\widetilde\tau_n$ on
$\widetilde D_n$, and form the norm ultraproduct
\[
 Q_\omega=
 \prod_n\widetilde D_n\big/
 \{(d_n):\lim_{n\to\omega}\|d_n\|=0\}.
\]
It carries the tracial state
\[
 \tau_\omega([(d_n)])
 =\lim_{n\to\omega}\widetilde\tau_n(d_n).
\]
Condition~\eqref{almost-ord0} makes the induced completely positive
contraction $\Phi:A\to Q_\omega$ an order zero map.
Put $h=\Phi(1_A)$.  By \eqref{state-asumption},
$\tau_\omega(h)=1$, and
\[
 \tau_\omega((1-h)^2)\leq\tau_\omega(1-h)=0.
\]
Since $\tau_\omega$ is tracial, its GNS representation
$\pi_{\tau_\omega}$ annihilates $1-h$.
Thus $\Pi=\pi_{\tau_\omega}\circ\Phi$ is a unital completely
positive order zero map, hence a *-homomorphism by
\cite[Theorem~3.3]{WZ-order-zero}.
In particular, $\tau=\tau_\omega\circ\Phi$ is tracial.
Faithfulness of $\tau$ makes $\Pi$ injective, so
\[
 \|a\|=\|\Pi(a)\|\leq\|\Phi(a)\|
 =\lim_{n\to\omega}\|\phi_n(a)\|\leq\|a\|.
\]
This proves the first conclusion.
If all tracial states of $A$ are faithful, the argument applies to
every free ultrafilter.  Since a bounded scalar sequence whose
limit along every free ultrafilter is the same must converge to
that value, the ordinary norm limit follows.

The state normalization also controls all powers of
$h_n=\phi_n(1_A)$: for $k\geq1$,
\begin{equation}\label{small-trace-bd}
 0\leq\tau_n(h_n-h_n^k)
 \leq(k-1)\bigl(1-\tau_n(h_n)\bigr)\longrightarrow0.\qedhere
\end{equation}
\end{proof}

\begin{lem}\label{central-cut}
Let $A$ and $B_n$ be unital \Cstar{}-algebras, and suppose every tracial
state of $A$ is faithful.
Let $\iota_n:A\to B_n$ be unital *-homomorphisms and let
$h_n\in B_n$ be positive contractions such that
\[
 \|[h_n,\iota_n(a)]\|\longrightarrow0,\qquad a\in A.
\]
Suppose there are $\tau_n\in\mathrm T(B_n)$ with
$\mathrm d_{\tau_n}(h_n)>0$ and
\begin{equation}\label{small-bdr-cond}
 \frac{\tau_n(h_n^2)}{\mathrm d_{\tau_n}(h_n)}
 \longrightarrow1.
\end{equation}
Then
\[
 \lim_{n\to\infty}\|h_n\iota_n(a)h_n\|=\|a\|,
 \qquad a\in A.
\]
\end{lem}
\begin{proof}
Set
\[
 D_n=\overline{h_nB_nh_n},\qquad
 \phi_n(a)=h_n\iota_n(a)h_n,\qquad
 \widetilde\tau_n=
 \frac{\tau_n|_{D_n}}{\mathrm d_{\tau_n}(h_n)}.
\]
The norm of $\tau_n|_{D_n}$ is $\mathrm d_{\tau_n}(h_n)$, as is seen by
evaluating it on the approximate unit $h_n^{1/j}$.
Thus $\widetilde\tau_n$ is a tracial state of $D_n$.
Each $\phi_n$ is completely positive and contractive.
If $a,b\in A_+$ and $ab=0$, then
\[
 \phi_n(a)\phi_n(b)
 =h_n\iota_n(a)[h_n^2,\iota_n(b)]h_n,
\]
whose norm tends to zero by relative centrality.
Finally, $\phi_n(1_A)=h_n^2$, and~\eqref{small-bdr-cond} gives
\[
 \widetilde\tau_n(\phi_n(1_A))
 =\frac{\tau_n(h_n^2)}{\mathrm d_{\tau_n}(h_n)}
 \longrightarrow1.
\]
Lemma~\ref{almost-ord0-lem} now gives the asserted norm limit.
\end{proof}

\begin{rem}
Condition~\eqref{small-bdr-cond} also implies
\[
 \frac{\tau_n(h_n^k)}{\mathrm d_{\tau_n}(h_n)}\longrightarrow1
 \quad\hbox{for every }k\geq1.
\]
Indeed, integrate
$0\leq1-t^k\leq k(1-t^2)$ over the nonzero spectrum of $h_n$
and divide by $\mathrm d_{\tau_n}(h_n)$.
\end{rem}

\begin{rem}
In Lemma~\ref{central-cut}, the support traces
$\mathrm d_{\tau_n}(h_n)$ may tend to zero.  The normalized traces on the
hereditary subalgebras must concentrate at $1$, as required for
individual Rokhlin towers.  For compressions whose defects have
uniformly small trace, Lemma~\ref{almost-pres-norm} gives the
corresponding norm estimate on a finite set.
\end{rem}
\chapter{Rokhlin properties and weak tracial approximation}\label{section-rokhlin}

Rokhlin properties for actions of discrete amenable groups on
unital \Cstar{}-algebras give matrix models over hereditary coefficient
subalgebras.  The weak Rokhlin property thereby gives weak tracial
approximation of the crossed product.  We also establish the
freeness conditions and lifting results underlying these tower
constructions.

Throughout this chapter, we write the action as
\[
 \alpha_g(a)=g(a)=u_gau_g^*,\qquad a\in A,\ g\in\Gamma,
\]
where $(u_g)_{g\in\Gamma}$ are the canonical unitaries and
$u_gu_k=u_{gk}$.  For finite subsets $K,F\subseteq\Gamma$, put
\[
 \operatorname{Int}_K(F)=\{g\in F:Kg\subseteq F\}.
\]
All invariance and collar estimates below therefore use left multiplication.

\section{Rokhlin properties for \Cstar{}-dynamical systems}


\subsection{A finite-stage sufficient condition: uniform Rokhlin property}
For a topological action on $X$, we use the left action convention
$\alpha_g(f)(x)=f(g^{-1}x)$ on $\mathrm C(X)$.  Invariance is understood
with respect to left multiplication, as above.
\begin{defn}[Topological URP]\label{topological-urp-definition}
Let $(X, \Gamma)$ be a topological dynamical system, where $X$ is a
compact Hausdorff space and $\Gamma$ is a discrete amenable group.
The system has the \emph{uniform Rokhlin property} (URP) if, for every
$\eps>0$ and every finite set $K\subseteq\Gamma$, there are
$(K, \eps)$-invariant finite sets
$\Gamma_1, \Gamma_2, ..., \Gamma_S\subseteq \Gamma$ and open sets
$$U_s\subseteq X,\quad s=1, ..., S $$
such that
\begin{enumerate}
\item $\gamma U_s$, $\gamma\in \Gamma_s$, $s=1, 2, ..., S$, are mutually disjoint, and
\item $\mu(X\setminus\bigcup_{s=1}^S\bigcup_{\gamma\in\Gamma_s} \gamma U_s)<\eps$ for every invariant probability measure $\mu$.
\end{enumerate}
\end{defn}
The corresponding property for \Cstar{}-dynamical systems is defined as follows.
\begin{defn}\label{URP-defn-C-dyn}
Let $(A, \Gamma)$ be a \Cstar{}-dynamical system, where $A$ is unital,
$\Gamma$ is countable and amenable, and
$\mathrm{T}(A) \neq \varnothing$.  A finite set $F\subseteq\Gamma$
is $(K,\eps)$-invariant if
$\abs{F\setminus\{g\in F:Kg\subseteq F\}}<\eps\abs F$;
we may assume $e\in K$.  The action has the \emph{uniform Rokhlin
property} (URP) if, for every $\eps>0$ and all finite sets
$\F\subseteq A$ and $\mathcal K \subseteq \Gamma$, there are
$(\mathcal K, \eps)$-invariant finite sets
$\Gamma_1, \Gamma_2, ..., \Gamma_S\subseteq \Gamma$ and positive
elements
$$c_s\in A,\quad s=1, ..., S $$
such that
\begin{enumerate}
\item $\norm{c_s} = 1$, $s=1, ..., S$,
\item $\gamma(c_s)$, $\gamma\in \Gamma_s$, $s=1, 2, ..., S$, are mutually orthogonal,
\item\label{almost-com-urp} $\displaystyle\sum_{s=1}^S\sum_{\gamma\in\Gamma_s}\norm{[\gamma(c_s), f]} < \eps$ for every $f\in\F$,
\item\label{small-leftover-urp} $$\mathrm{d}_\tau(1-\sum_{s=1}^S\sum_{\gamma\in\Gamma_s}\gamma(c_s)) < \eps, \quad \tau\in\mathrm{T}(A \rtimes \Gamma).$$
\end{enumerate}
 \end{defn}
 
\begin{rem}\label{rem-urp}
The aggregate commutator estimate in \eqref{almost-com-urp} is
stronger than a levelwise bound, since the number of tower levels
is not bounded in advance.  For $A=\mathrm C(X)$, all commutators
vanish.  Cutting finite URP towers along an exact tiling gives
central sequence towers with the prescribed shapes
(Proposition \ref{urp-implies-wrp}).
\end{rem}

\subsection{Weak Rokhlin property}
Weak Rokhlin-type properties for positive towers were studied by
Sato for automorphisms of the Jiang--Su algebra
\cite{Sato-Rohlin-Z} and by Wang for actions of countable discrete
amenable groups \cite{Wang-WRP}. We use the following
prescribed-tiling formulation, with the remainder controlled
uniformly over invariant tracial states.


\begin{defn}\label{D:weak-Rokhlin-property}
Let $(A, \Gamma)$ be a \Cstar{}-dynamical system, where $A$ is a unital
separable \Cstar{}-algebra with $\mathrm{T}(A) \neq \varnothing$ and
$\Gamma$ is a countable discrete amenable group. 
The action has the \emph{weak Rokhlin property} (WRP) if, 
for every exact tiling $\mathcal{T}$ of $\Gm$ with shapes
$K_1, \dots, K_S$, there are positive contractions
$c_s\in A' \cap A_\omega, s=1, ..., S$ such that
\begin{enumerate}
    \item $\gamma(c_s)$, $\gamma\in K_s$, $s=1, 2, ..., S$, are mutually orthogonal, and 
    \item $$1-\sum_{s=1}^S\sum_{\gamma\in K_s}\gamma(c_s) \in J_{\mathcal M_1(A, \Gamma)}, $$ where ${\mathcal M_1(A, \Gamma)}$ is the simplex of invariant
    tracial states, and $J_{\mathcal M_1(A, \Gamma)}$ is the
    corresponding trace-kernel ideal (Definition \ref{D:trace-kernel}).
\end{enumerate}
\end{defn}

The weak Rokhlin property admits the following local formulation
for each exact tiling.

\begin{lem}\label{local-wrp}
WRP is equivalent to the following condition: for every exact tiling with shapes
$D_1,\ldots,D_R$, every finite $\F\subseteq A$, and every $\eps>0$,
there are positive contractions $c_s\in A$ such that
\begin{enumerate}
\item $\norm{\gamma(c_s)\gamma'(c_{s'})}<\eps$ for distinct levels
$(s,\gamma)\ne(s',\gamma')$, with $\gamma\in D_s$ and $\gamma'\in D_{s'}$;
\item $\norm{[c_s,f]}<\eps$ for $s\leq R$ and $f\in\F$;
\item $\displaystyle\sup_{\rho\in\mathcal M_1(A,\Gamma)}
\rho\left(\left(1-\sum_{s=1}^R\sum_{\gamma\in D_s}
\gamma(c_s)\right)^2\right)<\eps$.
\end{enumerate}
\end{lem}
\begin{proof}
Suppose first that WRP holds.  Lift the positive contractions in
$A_\omega$ to sequences of positive contractions.  The norm conditions
and the uniform tracial $2$-norm condition hold simultaneously at
some coordinate.

Conversely, apply the local condition to an increasing sequence of
finite sets with dense union in $A$ and to tolerances tending to zero.
The resulting sequences define central positive contractions in
$A_\omega$ with orthogonal translates and complement in the
trace-kernel ideal.  At a finite stage, their sum need not be at most
one; condition (3) therefore uses the square of the complement.
\end{proof}

\begin{prop}\label{urp-implies-wrp}
For separable $A$, the URP of Definition
\ref{URP-defn-C-dyn} implies the prescribed-tiling WRP of Definition
\ref{D:weak-Rokhlin-property}.
\end{prop}
\begin{proof}
Fix an exact tiling
$\Gamma=\bigsqcup_{s=1}^R\bigsqcup_{t\in C_s}D_st$.
After translating each shape and its centers, assume $e\in D_s$.
Set $L=\bigcup_s D_sD_s^{-1}$, and choose $\eta_n\downarrow0$ and
increasing finite sets $\F_n\subset A$ with dense union. Apply URP with
$K=L$, $\F=\F_n$, and tolerance $\eta_n$. Write the resulting
shapes as $F_{r,n}$, their bases as $b_{r,n}$, and their levels as
$p_{r,g,n}=g(b_{r,n})$. Define
\[
 c_{s,n}=\sum_r\sum_{\substack{t\in C_s\\D_st\subset F_{r,n}}}
                         t(b_{r,n})\in A_+.
\]
These finite sums are contractions because $e\in D_s$ and the
summands are distinct levels of the URP towers.  For distinct $(s,k)$,
$k\in D_s$, the levels $k(c_{s,n})$ are orthogonal: within each
$F_{r,n}$ they comprise disjoint tiles of the prescribed tiling.
The aggregate URP estimate gives
$\norm{[c_{s,n},f]}<\eta_n$ for $f\in\F_n$.

Put $q_n=\sum_{r,g\in F_{r,n}}p_{r,g,n}$ and
$q'_n=\sum_{s,k\in D_s}k(c_{s,n})$. If
$g\in F_{r,n}$ satisfies $Lg\subset F_{r,n}$, its unique tiling
piece $D_st$ lies in $F_{r,n}$: writing $g=kt$, $k\in D_s$,
one has $D_st\subset Lg$. Hence at most $\eta_n|F_{r,n}|$
levels of each URP tower are omitted from $q'_n$.
For any invariant trace $\rho$, $\rho\circ E$ is a trace on the
crossed product, and subadditivity of dimension functions yields
\begin{align*}
\mathrm d_\rho(1-q'_n)&\leq \mathrm d_\rho(1-q_n)
 +\sum_r\sum_{g\text{ omitted from }F_{r,n}}\mathrm d_\rho(p_{r,g,n})\\
&<\eta_n+\eta_n\sum_r|F_{r,n}|\mathrm d_\rho(b_{r,n})\leq2\eta_n.
\end{align*}
Thus $\sup_\rho\rho((1-q'_n)^2)\leq2\eta_n$.
The classes $[(c_{s,n})]\in A_\omega\cap A'$ are the required
WRP bases for the prescribed tiling.
\end{proof}

\begin{prop}\label{top-urp-wrp}
Let a countable discrete amenable group $\Gamma$ act on a nonempty
compact metrizable space $X$. Then topological URP is equivalent to
WRP for the induced action on $\mathrm C(X)$.
\end{prop}
\begin{proof}
Suppose first that the topological system has URP, and discard
empty bases.  For the finite family of open towers, choose open
sets $V_s$ with $\overline{V_s}\subset U_s$ such that the union
of their translates has measure greater than $1-\eps$ for every
invariant measure.  To make this choice uniformly, regularity
gives such sets for each measure, and lower semicontinuity of the
measure of an open set gives neighborhoods covering the compact
invariant measure space.  A finite subcover and unions within
each $U_s$ give the required $V_s$.  Choose norm-one positive contractions
$c_s\in\mathrm C(X)$ supported in $U_s$ and equal to one on
$\overline{V_s}$.  The translated functions are orthogonal,
and the dimension function of $1-\sum_{s,g}g(c_s)$ is below
$\eps$ for every invariant measure. The aggregate commutator
condition is automatic.  Thus URP holds, and Proposition
\ref{urp-implies-wrp} applies.

Conversely, given $K$ and $\eps$, choose an exact tiling with
$(K,\eps)$-invariant shapes $D_s$ and use WRP for this tiling.
For infinite $\Gamma$, this follows from
\cite[Theorem~4.3 and Lemma~2.6]{DHZ-tiling}; for finite
$\Gamma$, use the single shape $D_1=\Gamma$.
Let $c_{s,n}\in\mathrm C(X)_+$ be representatives by positive contractions,
put $M=\sum_s|D_s|$, and choose $t<1/(2M)$. At a suitable coordinate,
products of distinct translated functions have norm below $t^2$ and
\[
\sup_\mu\mu\left(\left(1-\sum_{s,g\in D_s}g(c_{s,n})\right)^2\right)
<\eps/4.
\]
The open sets $\{c_{s,n}>t\}$ form disjoint towers. Outside their
union the sum of all translated functions is at most $Mt<1/2$,
so the complement has measure below $\eps$ for every invariant measure.
\end{proof}

\subsection{Invariant traces}

For a \Cstar{}-dynamical system $(A, \Gamma)$ with WRP, restriction to
$A$ identifies the tracial state simplex of $A\rtimes\Gamma$
with $\mathcal M_1(A, \Gamma)$.

\begin{prop}\label{trace-factor}
Let $(A, \Gamma)$ be a \Cstar{}-dynamical system, where $A$ is a unital
separable \Cstar{}-algebra with $\mathrm{T}(A) \neq \varnothing$ and
$\Gamma$ is a countable discrete amenable group. 
Let $\mathbb E \colon A \rtimes \Gm \rightarrow A$ be the
canonical conditional expectation.  If $(A, \Gamma)$ has WRP,
then the map
$$
\theta: \mathcal M_1(A, \Gamma) \ni \tau \mapsto \tau\circ\mathbb E \in \mathrm{T}(A \rtimes \Gm), 
$$
is an affine homeomorphism.
\end{prop}
\begin{proof}
The restriction of a trace $\sigma\in\mathrm T(A\rtimes\Gamma)$ to $A$
is invariant; conversely each invariant trace $\rho$ extends to
$\rho\circ\mathbb E$. Fix $k\ne e$, $a\in A$, and $\eta>0$.
Choose an exact tiling with shapes $D_s$ satisfying
$\abs{D_s\setminus\{g\in D_s:kg\in D_s\}}<\eta\abs{D_s}$.
Let $c_s\in A_\omega\cap A'$ be WRP bases,
$p_{s,g}=g(c_s)$, and $q=\sum_{s,g\in D_s}p_{s,g}$.
The trace-kernel condition gives $\sigma_\omega(1-q)=0$ for
$\sigma\in\mathrm T(A\rtimes\Gamma)$, where $\sigma_\omega$ is
its ultralimit trace.  Hence $q$ acts as the unit in the GNS
representation associated with $\sigma_\omega$.  Off-diagonal
terms vanish because
$\sigma_\omega(p_{s,g}a u_kp_{t,j})
=\sigma_\omega(p_{t,j}p_{s,g}a u_k)=0$ for distinct levels.
Thus orthogonality and traciality give
\[
\sigma(a u_k)=\sigma_\omega(q a u_k q)
 =\sum_{s,g\in D_s}\sigma_\omega(p_{s,g}a u_kp_{s,g}).
\]
For $kg\in D_s$, centrality with $a$ and orthogonality of
$p_{s,g}$ and $k(p_{s,g})$ make the summand zero. On the
remaining levels, Cauchy--Schwarz and invariance of
$\sigma|_A$ give
\[
\abs{\sigma(a u_k)}\leq\norm a\sum_s
 \abs{D_s\setminus\{g:kg\in D_s\}}\sigma_\omega(c_s)
 \leq\eta\norm a\sigma_\omega(q)\leq\eta\norm a.
\]
Let $\eta$ decrease to zero. Fourier polynomials are dense, so
$\sigma=\sigma|_A\circ\mathbb E$. The affine bijection is a
homeomorphism by compactness.
\end{proof}

\begin{cor}\label{separate-traces}
Let $(A, \Gamma)$ be a \Cstar{}-dynamical system with WRP. If projections of
$A$ separate its tracial states, then projections of $A\rtimes\Gamma$
separate its tracial states.
\end{cor}
\begin{proof}
Tracial states of $A\rtimes\Gamma$ are determined by their restrictions
to $A$, by Proposition \ref{trace-factor}. Distinct restrictions are
therefore distinguished by a projection in $A\subseteq A\rtimes\Gamma$.
%
\end{proof}

\begin{cor}\label{min-faithful-trace}
Let $(A, \Gamma)$ be a minimal \Cstar{}-dynamical system with WRP.
Then every tracial state of $A\rtimes\Gamma$ is faithful.
\end{cor}
\begin{proof}
Let $\tau\in\mathrm{T}(A\rtimes\Gamma)$, and let $a\in A\rtimes\Gamma$
be a positive element with $\tau(a)=0$.

By Proposition \ref{trace-factor},
$0=\tau(a) = \tau(\mathbb E(a))$.  The set
$I=\{x\in A: \tau(xx^*) = 0\}$ is a closed two-sided ideal of
$A$, since $\tau$ is tracial.  Invariance of $\tau|_A$ makes $I$
$\Gamma$-invariant.  Since $\tau$ is nonzero, minimality gives
$I = \{0\}$, so $\mathbb E(a) = 0$.  Faithfulness of
$\mathbb E$ now implies $a = 0$.
\end{proof}

\begin{rem}
Under the same hypotheses, $A\rtimes\Gamma$ is simple
(Corollary \ref{wrp-simple}).
\end{rem}

\section{Weak Rokhlin property and weak tracial approximation}
\label{wrp-approximation-section}

For separable $A$, Proposition \ref{urp-implies-wrp} shows that URP
implies WRP.  The approximation theorem in this section therefore
also applies to systems with URP.

Let $(A, \Gamma)$ be a \Cstar{}-dynamical system with WRP, and let
$\Gamma_1, ..., \Gamma_S$ be the shapes of an exact tiling of
$\Gamma$, with $e\in\Gamma_s$, $s=1, ..., S$.  There are
$c_s \in A' \cap A_\omega$, $s=1, ..., S$, such that
$$\gamma(c_s),\quad \gamma\in\Gamma_s,\ s=1, ..., S, $$
are mutually orthogonal.
For each $s$, index the matrix units of
$\mathrm M_{\abs{\Gamma_s}}$ by $\Gamma_s$.  The elements
\[
 v_{s,\gamma}:=c_su_\gamma^*,\qquad
 \gamma\in\Gamma_s,\quad 1\leq s\leq S,
\]
determine a *-homomorphism
\[
 \phi\colon
 \bigoplus_{s=1}^S
 \mathrm M_{\abs{\Gamma_s}}(\mathrm C_0((0,1]))
 \longrightarrow (A\rtimes\Gamma)_\omega,
 \qquad
 \phi(t\otimes e_{e,\gamma}^{(s)})=v_{s,\gamma}.
\]
Its image is the \Cstar{}-subalgebra generated by these elements.

Projectivity of the cone gives a lift of $\phi$.  The following
relative form of Loring's construction gives a lift whose diagonal
belongs to the coefficient algebra.

\begin{lem}[Relative cone lifting]\label{lift-d-drop}
Let $\Gamma_0=\{\gamma_1=e,\gamma_2,\ldots,\gamma_n\}$, let
$c\in A_\omega$ be a positive contraction, and suppose that
\[
  \gamma_i(c):=u_{\gamma_i}cu_{\gamma_i}^*,\qquad 1\leq i\leq n,
\]
are mutually orthogonal.  Let
\[
 \phi\colon \mathrm M_n(\mathrm C_0((0,1]))\longrightarrow
 (A\rtimes\Gamma)_\omega
 \quad\text{be given by}\quad
 \phi(t\otimes e_{1i})=cu_{\gamma_i}^*.
\]
Then $\phi$ has a lift
\[
 \widetilde\phi\colon \mathrm M_n(\mathrm C_0((0,1]))\longrightarrow
 \ell^\infty(A\rtimes\Gamma)
\]
such that
\[
 \widetilde\phi(t\otimes e_{ii})\in\ell^\infty(A),\qquad
 1\leq i\leq n.
\]
Moreover, the elements
$\widetilde\phi(t\otimes e_{i1})$ are generalized normalizers of
$\ell^\infty(A)$ in $\ell^\infty(A\rtimes\Gamma)$.
\end{lem}

\begin{proof}
Put $B=A\rtimes\Gamma$, $E=\ell^\infty(B)$, and
$D=\ell^\infty(A)\subseteq E$, and denote the quotient map by
$q\colon E\to B_\omega$.  We use the same notation for the constant
sequences $u_{\gamma_i}\in E$.

The standard lifting lemma for a finite family of orthogonal positive
contractions gives
mutually orthogonal positive contractions $k_1,\ldots,k_n\in D$ such
that
\[
 q(k_i)=\gamma_i(c),\qquad 1\leq i\leq n.
\]
This lifting lemma follows by applying the usual
finite orthogonalization perturbation to arbitrary positive lifts at
each coordinate; the perturbation tends to zero because the products
of distinct lifts tend to zero along $\omega$.  This is also the
commutative instance of the cone-lifting argument in
\cite[Section~10]{Lor-Mon}.

Write $k_i=(k_{i,m})_m$ and choose numbers $\eta_m\searrow0$.  Set
\[
 h_{i,m}=(k_{i,m}-\eta_m)_+,
 \]
and let $g_{i,m}=f_m(k_{i,m})$, where $0\leq f_m\leq1$, $f_m=0$ on
$[0,\eta_m/2]$, and $f_m=1$ on $[\eta_m,1]$.  With
$h_i=(h_{i,m})_m$ and $g_i=(g_{i,m})_m$, we have
\begin{equation}\label{relative-cone-collars}
 \begin{split}
 &q(h_i)=\gamma_i(c),\qquad g_ih_i=h_i,\\
 &h_ih_j=g_ig_j=0\quad (i\ne j).
 \end{split}
\end{equation}

Put $p_{n+1}=h_1^2$ and,
recursively for $i=n,n-1,\ldots,2$, define
\begin{equation}\label{relative-cone-recursion}
 b_i=u_{\gamma_i}p_{i+1}^{1/2},\qquad
 a_i=g_i b_i,\qquad p_i=a_i^*a_i.
\end{equation}
Since conjugation by $u_{\gamma_i}$ preserves $D$, we have $p_i\in D$
and
\begin{equation}\label{relative-cone-moduli}
 0\leq p_2\leq p_3\leq\cdots\leq p_{n+1}=h_1^2.
\end{equation}
Moreover, an induction using \eqref{relative-cone-collars} gives
\begin{equation}\label{relative-cone-quotients}
 q(a_i)=u_{\gamma_i}c,\qquad q(p_i)=c^2,\qquad 2\leq i\leq n.
\end{equation}
Indeed, if $q(p_{i+1})=c^2$, then $q(b_i)=u_{\gamma_i}c$; and
$u_{\gamma_i}c=\gamma_i(c)u_{\gamma_i}$, while
$q(g_i)\gamma_i(c)=\gamma_i(c)$.

All the $p_i$ belong to $\overline{h_1Dh_1}$.  Since $g_jh_1=0$ for
$j\geq2$, it follows that
\begin{equation}\label{relative-cone-zero}
 b_i g_j=0\quad (2\leq i,j\leq n),
 \qquad a_i^*a_j=0\quad (i\ne j).
\end{equation}
Each $a_i=g_i u_{\gamma_i}p_{i+1}^{1/2}$ is also a
generalized normalizer of $D$ in $E$:
\begin{equation}\label{relative-cone-normalizer}
 a_iDa_i^*\subseteq D,\qquad a_i^*Da_i\subseteq D.
\end{equation}

For $\lambda>0$, define
\[
 z_i(\lambda)=a_i(\lambda+p_i)^{-1/2}p_2^{1/2},\qquad 2\leq i\leq n.
\]
Because $p_2\leq p_i$, the elements $z_i(\lambda)$ converge in norm as
$\lambda\searrow0$; denote the limit by $z_i$.  Indeed, the elements
\[
 e_\lambda=p_i(\lambda+p_i)^{-1}
\]
form an approximate identity for $\overline{p_iDp_i}$, which contains
$p_2$.  Thus
\[
 z_i(\lambda)^*z_i(\lambda)
   =p_2^{1/2}e_\lambda p_2^{1/2}\longrightarrow p_2.
\]
The same calculation with two parameters shows that
$z_i(\lambda)$ is Cauchy: its squared difference is obtained by
sandwiching
\[
 e_\lambda+e_\mu-2(e_\lambda e_\mu)^{1/2}
   =(e_\lambda^{1/2}-e_\mu^{1/2})^2
\]
between two copies of $p_2^{1/2}$.  Consequently,
\begin{equation}\label{relative-cone-common-modulus}
 z_i^*z_i=p_2,\qquad 2\leq i\leq n.
\end{equation}
By \eqref{relative-cone-normalizer},
$z_i(\lambda)z_i(\lambda)^*\in D$, and hence
\begin{equation}\label{relative-cone-diagonal-range}
 z_iz_i^*\in D,\qquad 2\leq i\leq n.
\end{equation}

The right-hand factor
$(\lambda+p_i)^{-1/2}p_2^{1/2}$ belongs to
$\overline{h_1Dh_1}$.  Equations \eqref{relative-cone-zero} therefore
give, after taking limits,
\begin{equation}\label{relative-cone-relations}
 z_i z_j=0,\qquad z_i^*z_j=\delta_{ij}p_2,\qquad 2\leq i,j\leq n.
\end{equation}
Finally, \eqref{relative-cone-quotients} yields
\[
 \begin{split}
 q(z_i)
 &=\lim_{\lambda\searrow0}
   u_{\gamma_i}c(\lambda+c^2)^{-1/2}c
  =u_{\gamma_i}c.
 \end{split}
\]
The last convergence is in norm, since
\[
 \sup_{0\leq s\leq1}
 \left|\frac{s^2}{(\lambda+s^2)^{1/2}}-s\right|
 \longrightarrow0.
\]

Equation \eqref{relative-cone-common-modulus} also gives
$\|z_i\|^2=\|p_2\|\leq1$.  Thus the relations in
\eqref{relative-cone-relations} are Loring's cone relations: the
universal \Cstar{}-algebra generated by contractions
$x_2,\ldots,x_n$ satisfying
\[
 x_ix_j=0,\qquad x_i^*x_j=\delta_{ij}x_2^*x_2
\]
is $\mathrm M_n(\mathrm C_0((0,1]))$, with
$x_i=t\otimes e_{i1}$; see
\cite[Sections~3 and~10]{Lor-Mon}. Hence there is a
*-homomorphism $\widetilde\phi$ with
\[
 \widetilde\phi(t\otimes e_{i1})=z_i,\qquad 2\leq i\leq n.
\]
It lifts $\phi$, because these generators generate the cone and
$q(z_i)=\phi(t\otimes e_{i1})$.  Furthermore,
\[
 \widetilde\phi(t\otimes e_{11})=(z_2^*z_2)^{1/2}=p_2^{1/2}\in D
\]
and, for $2\leq i\leq n$,
\[
 \widetilde\phi(t\otimes e_{ii})=(z_iz_i^*)^{1/2}\in D.
\]
The generalized normalizer assertion follows from
\eqref{relative-cone-normalizer} and the norm-limit construction of the
$z_i$.  When $n=1$, any lift of $c$ by a positive contraction
gives the assertion.
\end{proof}

Applying the construction simultaneously to all pairs $(s,\gamma)$
gives the corresponding lift for the finite direct sum above.  Indeed,
the collar elements can be chosen mutually orthogonal across all
summands, so the images of distinct summands are orthogonal.

\begin{lem}[Relative linking algebra]\label{app-tower-alg}
Let $A\subseteq B$ be \Cstar{}-algebras and let
\[
 \psi:\mathrm M_n(\mathrm C_0((0,1]))\longrightarrow B
\]
be a *-homomorphism.  Suppose that
\[
 d_i=\psi(t\otimes e_{ii})\in A,\qquad 1\leq i\leq n,
\]
and that $x_i=\psi(t\otimes e_{i1})$ is a generalized normalizer of
$A$ in $B$ for every $i$.  Put $H=\overline{d_1Ad_1}$.  Then there is
a \Cstar{}-subalgebra $C\subseteq B$ such that
\[
 C\cong \mathrm M_n(H),\qquad d_i\in C\cap A,\qquad
 C_{ii}=\overline{d_iAd_i}.
\]
\end{lem}

\begin{proof}
Write $x_i=w_id_1$ for the polar decomposition in $B^{**}$, taking
$w_1$ to be the support projection of $d_1$.  The cone relations give
\[
 w_i^*w_j=\delta_{ij}s(d_1)
\]
on the support of $d_1$, and $w_id_1w_i^*=d_i$.  Define
\[
 \Theta:\mathrm M_n(H)\longrightarrow B^{**},\qquad
 \Theta((a_{ij}))=\sum_{i,j=1}^nw_ia_{ij}w_j^*.
\]
These relations show that $\Theta$ is an injective *-homomorphism.

To locate the image and its diagonal, note that $d_1Ad_1$ is dense
in $H$ and that, for
$b\in A$,
\[
 w_i(d_1bd_1)w_j^*=x_i b x_j^*\in B.
\]
For $i=j$ this element belongs to $A$, since $x_i$ is a generalized
normalizer.  Passing to norm limits proves
$\Theta(\mathrm M_n(H))\subseteq B$ and
$\Theta(e_{ii}\otimes H)\subseteq A$. Every element of this diagonal
corner is supported on $s(d_i)$, so the corner lies in
$\overline{d_iAd_i}$. Conversely, for $b\in A$, the inverse-normalizer
assumption gives $x_i^*bx_i\in A$. Its support lies in $s(d_1)$,
hence $x_i^*bx_i\in H$, and
\[
 \Theta(e_{ii}\otimes x_i^*bx_i)=d_i b d_i.
\]
Density proves the reverse inclusion. Taking
$C=\Theta(\mathrm M_n(H))$ proves the lemma.
\end{proof}

\begin{rem}
The relative lift places the coordinate diagonals in \(A\), and
the linking lemma enlarges each lifted matrix cone to a matrix
algebra over a hereditary \Cstar{}-subalgebra of \(A\).  The
trace-kernel condition gives a uniform tracial estimate on the
complement.  A functional-calculus cutoff then gives the
dimension-function estimate, as in Remark \ref{trace-rank-ta0}.
\end{rem}

\begin{thm}[WRP approximation]\label{WTA-C-dyn-1}
Let \((A,\Gamma)\) be a \Cstar{}-dynamical system with WRP, and put
\(B=A\rtimes_{\alpha,r}\Gamma\).  Let \(\mathcal S\) be the class of
finite direct sums of matrix algebras over hereditary \Cstar{}-subalgebras of
\(A\).  For every finite set \(\mathcal F\subseteq B\) and every
\(\eps>0\), there are a \Cstar{}-subalgebra
\[
 C\cong\bigoplus_{s=1}^S \mathrm M_{n_s}(H_s),\qquad
 H_s\subseteq A\ \text{hereditary},
\]
and a norm-one positive contraction \(h\in A\cap C\) such that
\begin{enumerate}
\item \(\|[h,x]\|<\eps\) for \(x\in\mathcal F\);
\item \(\operatorname{dist}(hxh,C)<\eps\) for \(x\in\mathcal F\), and
      the approximating element may be chosen in \(A\cap C\) when
      \(x\in A\);
\item \(\mathrm d_\tau(1-h)<\eps\) for every \(\tau\in\mathrm T(B)\).
\end{enumerate}
Consequently,
\[
 B\in\mathrm{TA}_0(\mathcal S).
\]
\end{thm}

\begin{proof}
By Fourier approximation and a preliminary reduction of the tolerance, it
is enough to treat elements \(au_k\), where \(a\) belongs to a fixed
finite subset \(\mathcal E\) of the unit ball of \(A\) and \(k\) belongs
to a finite symmetric set \(N\subseteq\Gamma\) containing \(e\).

Choose a large integer \(L\) and a small number \(\beta>0\).  If
\(\Gamma\) is infinite, the exact tiling theorem
\cite[Theorem~4.3 and Lemma~2.6]{DHZ-tiling} gives an
exact tiling whose shapes \(D_1,\ldots,D_S\) satisfy
\begin{equation}\label{wrp-small-boundary}
 |D_s\setminus\operatorname{Int}_{N^L}(D_s)|
       <\beta|D_s|,\qquad 1\leq s\leq S.
\end{equation}
Moreover, the shapes may be required to have arbitrarily large cardinality.
Indeed, enlarge the invariance set by a finite set \(K_0\) of the desired
cardinality and take the invariance tolerance below one.  Then
\(\operatorname{Int}_{K_0}(D_s)\ne\varnothing\), so some translate of
\(K_0\) is contained in \(D_s\).
After right translating a shape and its centers, we may assume
\(e\in D_s\).  If \(\Gamma\) is finite, take the one shape
\(D_1=\Gamma\), for which the boundary in
\eqref{wrp-small-boundary} is empty.

Apply WRP to this tiling.  We obtain positive contractions
\(c_s\in A_\omega\cap A'\) such that the elements
\[
 \alpha_g(c_s),\qquad g\in D_s,\quad 1\leq s\leq S,
\]
are mutually orthogonal and
\begin{equation}\label{wrp-trace-kernel-defect}
 1-\sum_{s=1}^S\sum_{g\in D_s}\alpha_g(c_s)
       \in J_{\mathcal M_1(A,\Gamma)}.
\end{equation}
They determine the cone homomorphism
\[
 \phi:\bigoplus_{s=1}^S \mathrm M_{|D_s|}(\mathrm C_0((0,1]))
       \longrightarrow B_\omega,\qquad
 \phi(t\otimes e_{e,g}^{(s)})=c_su_g^*.
\]
Let \(\widetilde\phi=(\widetilde\phi_m)_m\) be the lift supplied by
Lemma \ref{lift-d-drop}, where
\[
 \widetilde\phi_m:
 \bigoplus_{s=1}^S \mathrm M_{|D_s|}(\mathrm C_0((0,1]))\longrightarrow B
\]
is the \(m\)th coordinate homomorphism.  Put
\[
 d_{s,g,m}=\widetilde\phi_m(t\otimes e_{g,g}^{(s)})\in A,
 \qquad q_m=\sum_{s,g}d_{s,g,m}.
\]
The elements \(d_{s,g,m}\) are exactly orthogonal and
\(0\leq q_m\leq1\).  Equation
\eqref{wrp-trace-kernel-defect} says that
\begin{equation}\label{wrp-coordinate-two-norm}
 \lim_{m\to\omega}\sup_{\rho\in\mathcal M_1(A,\Gamma)}
       \rho((1-q_m)^2)=0.
\end{equation}

For each \(s,m\), apply Lemma \ref{app-tower-alg} to the \(s\)th cone
summand.  We obtain mutually orthogonal \Cstar{}-subalgebras
\[
 C_{s,m}\cong \mathrm M_{|D_s|}(H_{s,m}),\qquad
 H_{s,m}=\overline{d_{s,e,m}Ad_{s,e,m}},
\]
whose \(g\)th diagonal corner is
\(\overline{d_{s,g,m}Ad_{s,g,m}}\).  Put
\(C_m=\bigoplus_sC_{s,m}\), omitting zero summands.  Also set
\[
 y_{s,g,m}=\widetilde\phi_m(t\otimes e_{g,e}^{(s)})\in C_{s,m}.
\]
In the quotient one has
\[
 [(d_{s,g,m})_m]=\alpha_g(c_s),\qquad
 [(y_{s,g,m})_m]=u_gc_s.
\]

Fix an auxiliary tolerance \(\eta>0\).  By
\eqref{wrp-coordinate-two-norm}, centrality, and the quotient
identities, the following finite collection of conditions holds
at some coordinate \(m\):
\begin{align}
 &\sup_{\rho\in\mathcal M_1(A,\Gamma)}
       \rho((1-q_m)^2)<\eta^2,                                     
       \label{wrp-coordinate-trace}\\
 &\sum_{s,g}\|[d_{s,g,m},a]\|<\eta,
       &&a\in\mathcal E,                                          
       \label{wrp-coordinate-centrality}\\
 &\sum_{\substack{s,g\\kg\in D_s}}
   \|u_kd_{s,g,m}u_k^*-d_{s,kg,m}\|<\eta,
       &&k\in N,                                                   
       \label{wrp-coordinate-covariance}
\end{align}
and, whenever \(g,j\in D_s\), \(k\in N\), and \(g=kj\),
\begin{equation}\label{wrp-coordinate-linking}
 \sum_{s,g,j,a,k}
 \left\|d_{s,g,m}a d_{s,g,m}u_k
 -y_{s,g,m}\alpha_{g^{-1}}(a)y_{s,j,m}^*\right\|<\eta.
\end{equation}
The last sum ranges only over the indicated finite collection.  Its
summands tend to zero along \(\omega\), since their images in \(B_\omega\)
both equal
\[
 u_gc_s\alpha_{g^{-1}}(a)c_su_j^*.
\]

For \(g\in D_s\), define
\[
 r_s(g)=\max\{0\leq r\leq L:N^rg\subseteq D_s\},\qquad
 \lambda_s(g)=r_s(g)/L,
\]
and extend \(\lambda_s\) by zero off \(D_s\).  Symmetry of \(N\) gives
\begin{equation}\label{wrp-collar-lipschitz}
 |\lambda_s(kg)-\lambda_s(g)|\leq L^{-1},
 \qquad k\in N.
\end{equation}
Also, \(\lambda_s(g)>0\) implies \(Ng\subseteq D_s\), and
\(\lambda_s=1\) on \(\operatorname{Int}_{N^L}(D_s)\).  Define
\[
 h_0=\sum_{s,g}\lambda_s(g)d_{s,g,m}\in A\cap C_m.
\]
Equations \eqref{wrp-coordinate-centrality},
\eqref{wrp-coordinate-covariance}, and
\eqref{wrp-collar-lipschitz} give
\begin{equation}\label{wrp-h0-commutators}
 \|[h_0,a]\|<\eta,\qquad
 \|[h_0,u_k]\|<L^{-1}+\eta
\end{equation}
for \(a\in\mathcal E\) and \(k\in N\).

For the compression estimate, fix \(k\in N\) and put
\(\mu_{s,g}=\lambda_s(k^{-1}g)\), again using the zero extension.
Equation \eqref{wrp-coordinate-covariance} gives
\[
 \alpha_k(h_0)\approx_\eta\sum_{s,g}\mu_{s,g}d_{s,g,m}.
\]
Using exact orthogonality and \eqref{wrp-coordinate-centrality}, we get
\begin{align*}
 h_0au_kh_0
 &=h_0a\alpha_k(h_0)u_k\\
 &\approx_{2\eta}
   \sum_{s,g}\lambda_s(g)\mu_{s,g}
          d_{s,g,m}a d_{s,g,m}u_k.
\end{align*}
Writing $D_\mu=\sum_{s,g}\mu_{s,g}d_{s,g,m}$, the centrality
error is
\[
 \begin{aligned}
 &h_0aD_\mu-\sum_{s,g}\lambda_s(g)\mu_{s,g}d_{s,g,m}a d_{s,g,m}\\
 &\qquad=\sum_{s,g}\mu_{s,g}(h_0-\lambda_s(g)d_{s,g,m})[a,d_{s,g,m}],
 \end{aligned}
\]
whose norm is below $\eta$.  The covariance error is also below
$\eta$.
When a coefficient in this sum is nonzero,
\(j=k^{-1}g\in D_s\), and \eqref{wrp-coordinate-linking} replaces the
corresponding summand, within total error \(\eta\), by
\[
 y_{s,g,m}\alpha_{g^{-1}}(a)y_{s,j,m}^*\in C_{s,m}.
\]
Thus
\begin{equation}\label{wrp-h0-compression}
 \operatorname{dist}(h_0au_kh_0,C_m)<3\eta.
\end{equation}
For \(k=e\), the displayed approximating element belongs to \(A\cap C_m\).

It remains to control the defect.  If \(\tau\in\mathrm T(B)\), then
\(\rho=\tau|_A\) belongs to \(\mathcal M_1(A,\Gamma)\), and
\eqref{wrp-coordinate-trace} gives
\[
 \tau(1-q_m)\leq\tau((1-q_m)^2)^{1/2}<\eta.
\]
Under the identification
\(C_{s,m}\cong \mathrm M_{|D_s|}(H_{s,m})\), the elements
\(d_{s,g,m}\) are identical diagonal copies of the same positive element.
Consequently, they have the same value under every trace on \(B\).  Therefore
\begin{align}
 \tau(1-h_0)
 &=\tau(1-q_m)+
   \sum_{s,g}(1-\lambda_s(g))\tau(d_{s,g,m})\nonumber\\
 &<\eta+\beta\sum_{s,g}\tau(d_{s,g,m})
 \leq\eta+\beta.                                     
 \label{wrp-h0-trace-defect}
\end{align}

Finally, set
\[
 h=g_{1/2}(h_0)\in A\cap C_m.
\]
Taking \(\eta+\beta<1/2\) in \eqref{wrp-h0-trace-defect} shows that
\(\|h_0\|>1/2\), and hence \(\|h\|=1\).  Moreover,
\[
 \mathrm d_\tau(1-h)
 =\tau(1_{[0,1/2)}(h_0))
 \leq2\tau(1-h_0)<2(\eta+\beta).
\]
If \(r(t)=g_{1/2}(t)/t\), extended continuously by zero at the origin,
then
\[
 hxh=r(h_0)(h_0xh_0)r(h_0).
\]
Here \(r(h_0)\in A\cap C_m\) and \(\|r(h_0)\|\leq2\).  Consequently,
\eqref{wrp-h0-compression} gives
\(\operatorname{dist}(hxh,C_m)<12\eta\) for the prescribed Fourier
terms.  Continuous functional calculus, together with
\eqref{wrp-h0-commutators}, gives the required commutator estimate for
\(h\) after decreasing the preliminary tolerances.  When \(x\in A\),
the approximating element remains in \(A\cap C_m\).
Choose \(L\) sufficiently large, then \(\beta\) sufficiently small,
and finally \(\eta\) and the coordinate tolerances sufficiently
small.  The initial Fourier approximation gives the required
estimates for \(\mathcal F\).  The trace space of
\(B\) is nonempty because amenability supplies an invariant trace on
\(A\), whose composition with the canonical conditional expectation is
a trace on \(B\).
\end{proof}

\begin{cor}\label{wrp-rr0-affine}
Let \((A,\Gamma)\) be a \Cstar{}-dynamical system with WRP, where \(A\) has
real rank zero and \(\Gamma\) is a countable discrete amenable group.
Put \(B=A\rtimes_{\alpha,r}\Gamma\), and let
\[
 \rho_B:\mathrm K_0(B)\longrightarrow\operatorname{Aff}(\mathrm T(B))
\]
be the usual trace-pairing map.  Then
\[
 \overline{\operatorname{span}_{\mathbb R}\rho_B(\mathrm K_0(B))}
 =\operatorname{Aff}(\mathrm T(B)).
\]
If \(\Gamma\) is infinite, then, in addition, for every
projection \(p\in A\), every \(r\in[0,1]\), and every \(\eps>0\),
there is a projection \(q\in A\) such that
\[
 \sup_{\tau\in\mathrm T(B)}|\tau(q)-r\tau(p)|<\eps.
\]
Consequently, in this case,
\[
 \overline{\rho_B(\mathrm K_0(B))}=\operatorname{Aff}(\mathrm T(B)).
\]
\end{cor}
\begin{proof}
Real rank zero implies that projections of \(A\) separate its tracial
states.  Corollary \ref{separate-traces} therefore shows that the trace
functions of projections in \(A\), viewed as projections in \(B\),
separate the points of \(\mathrm T(B)\).  The real linear subspace
\[
 V=\operatorname{span}_{\mathbb R}\rho_B(\mathrm K_0(B))
 \subseteq\operatorname{Aff}(\mathrm T(B))
\]
contains the constant functions and separates points.  The affine
Stone--Weierstrass theorem, in the form
\cite[Corollary~7.4]{Goodearl}, gives
\(\overline V=\operatorname{Aff}(\mathrm T(B))\).
Alternatively, if the closure were proper, the Hahn--Banach and
Riesz representation theorems
would give a nonzero functional on $\operatorname{Aff}(\mathrm T(B))$
represented by a signed measure annihilating $V$. Since $1\in V$, its
positive and negative parts have equal positive mass. Their normalized
barycenters agree because $V$ separates points. The signed measure
therefore annihilates every continuous affine function, a contradiction.

Assume now that \(\Gamma\) is infinite.  Fix a projection \(p\in A\),
\(r\in[0,1]\), and \(\eps>0\).  Choose \(N\in\mathbb N\) and
\(\eta,\delta>0\) such that
\[
 N^{-1}+\eta+2\delta<\eps.
\]
Use the proof of Theorem \ref{WTA-C-dyn-1}, with its preliminary
tolerances sufficiently small, and require all the tiling shapes to have
cardinality
\(n_s=|D_s|\geq N\).  At the chosen coordinate that proof supplies
mutually orthogonal algebras
\[
 C_s\cong \mathrm M_{n_s}(H_s),\qquad
 H_s=\overline{d_{s,e}Ad_{s,e}},
\]
and identical diagonal copies \(d_{s,g}\), \(g\in D_s\), of
\(d_{s,e}\), such that, with
\[
 z=\sum_{s=1}^S\sum_{g\in D_s}d_{s,g},
\]
one has
\begin{equation}\label{wrp-large-coordinate-for-K0}
 \tau(1-z)<\delta,
 \qquad \tau\in\mathrm T(B).
\end{equation}

Every \(H_s\) has real rank zero and hence has an approximate identity
of projections.  Put \(M=\sum_sn_s\).  For each nonzero summand choose a
projection \(e_s\in H_s\) such that
\[
 \|d_{s,e}-e_sd_{s,e}e_s\|<\delta/M,
\]
and denote its identical diagonal copies in \(C_s\) by \(e_{s,g}\).
Set
\[
 P=\sum_{s,g}e_{s,g}.
\]
For every \(\tau\in\mathrm T(B)\),
\[
 \tau(d_{s,g})
 \leq \tau(e_{s,g}d_{s,g}e_{s,g})+\delta/M
 \leq \tau(e_{s,g})+\delta/M.
\]
After summing and using \eqref{wrp-large-coordinate-for-K0}, we obtain
\begin{equation}\label{wrp-large-projection-for-K0}
 \tau(1-P)<2\delta,
 \qquad \tau\in\mathrm T(B).
\end{equation}

Fix identifications \(C_s\cong \mathrm M_{n_s}(H_s)\), and denote the induced
isomorphism from \(H_s\) onto the \(g\)th diagonal corner by
\(\theta_{s,g}\).  Thus \(e_{s,g}=\theta_{s,g}(e_s)\).  Define
\[
 a_{s,g}=\theta_{s,g}^{-1}(e_{s,g}pe_{s,g})\in (H_s)_+,
 \qquad
 a_s=\frac{1}{n_s}\sum_{g\in D_s}a_{s,g}.
\]
These are positive contractions.  Since \(H_s\) has real rank zero,
there is a positive finite-spectrum contraction
\[
 b_s=\sum_{j=1}^{m_s}\lambda_{s,j}f_{s,j}\in H_s,
 \qquad 0<\lambda_{s,j}\leq1,
\]
such that \(\lVert b_s-a_s\rVert<\eta\), where the
\(f_{s,j}\) are mutually orthogonal projections.  For each \(s,j\),
choose an integer \(k_{s,j}\) such that
\[
 0\leq k_{s,j}\leq n_s,
 \qquad |k_{s,j}-r\lambda_{s,j}n_s|\leq1,
\]
and choose a subset \(E_{s,j}\subseteq D_s\) with
\(|E_{s,j}|=k_{s,j}\).  Then
\[
 Q=\sum_{s=1}^S\sum_{j=1}^{m_s}
       \sum_{g\in E_{s,j}}\theta_{s,g}(f_{s,j})
\]
is a projection in \(A\cap C\): each summand belongs to the diagonal
corner \(\overline{d_{s,g}Ad_{s,g}}\).

Let \(\tau\in\mathrm T(B)\).  Its restrictions to the diagonal corners
of \(C_s\cong \mathrm M_{n_s}(H_s)\) induce the same bounded trace
\(\sigma_s\) on \(H_s\).  Since \(P=\sum_{s,g}e_{s,g}\), traciality
gives
\begin{equation}\label{wrp-relative-average-for-K0}
 \tau(pP)
 =\sum_{s,g}\tau(e_{s,g}pe_{s,g})
 =\sum_s n_s\sigma_s(a_s).
\end{equation}
Moreover,
\begin{align*}
 |\tau(Q)-r\tau(pP)|
 &\leq
 \sum_{s,j}|k_{s,j}-r\lambda_{s,j}n_s|\sigma_s(f_{s,j})
 +r\sum_sn_s|\sigma_s(b_s-a_s)|\\
 &\leq \sum_{s,j}\sigma_s(f_{s,j})
       +\eta\sum_sn_s\lVert\sigma_s\rVert.
\end{align*}
The sum of the norms of the restrictions of \(\tau\) to the mutually
orthogonal summands \(C_s\) is at most one.  Under the matrix
identification, the norm of \(\tau|_{C_s}\) is
\(n_s\lVert\sigma_s\rVert\).  Hence
\[
 \sum_sn_s\lVert\sigma_s\rVert\leq1.
\]
On the other hand,
\[
 F=\sum_{s,j}\sum_{g\in D_s}\theta_{s,g}(f_{s,j})
\]
is a projection in \(B\), and therefore
\[
 \sum_{s,j}\sigma_s(f_{s,j})
 \leq \frac1N\sum_{s,j}n_s\sigma_s(f_{s,j})
 =\frac1N\tau(F)\leq\frac1N.
\]
It follows from \eqref{wrp-relative-average-for-K0} that
\begin{equation}\label{wrp-relative-rounding-for-K0}
 |\tau(Q)-r\tau(pP)|\leq N^{-1}+\eta.
\end{equation}
Finally, \(\tau(pP)=\tau(pPp)\), so
\[
 0\leq\tau(p)-\tau(pP)
 =\tau(p(1-P)p)
 =\tau((1-P)p(1-P))
 \leq\tau(1-P)<2\delta.
\]
Together with \eqref{wrp-relative-rounding-for-K0}, this gives
\[
 \sup_{\tau\in\mathrm T(B)}|\tau(Q)-r\tau(p)|
 <N^{-1}+\eta+2\delta<\eps.
\]

Set \(G=\overline{\rho_B(\mathrm K_0(B))}\).  The preceding construction shows
that \(r\widehat p\in G\) for every projection \(p\in A\) and every
\(r\in[0,1]\), where \(\widehat p(\tau)=\tau(p)\).  Since \(G\) is a
closed additive subgroup, it contains every real scalar multiple of
\(\widehat p\), and hence the real linear span of the trace functions of
projections in \(A\).  This span contains the constants and separates
the points of \(\mathrm T(B)\).  The affine Stone--Weierstrass
theorem therefore gives
\(G=\operatorname{Aff}(\mathrm T(B))\).
\end{proof}

\begin{rem}
For infinite \(\Gamma\), the relative averaging in the proof gives
approximate divisibility of projection trace functions and hence
density of the \(\mathrm K_0\)-range.  Approximation of constant functions
alone would not give this conclusion.  No non-elementarity
assumption on \(A\) is required, but infinitude of \(\Gamma\) is
essential: for \(\Gamma=\{e\}\) and \(A=\mathrm M_n\), WRP and real rank
zero hold, whereas \(B=\mathrm M_n\) and
\(\rho_B(\mathrm K_0(B))=n^{-1}\mathbb Z\), which is not dense in
\(\operatorname{Aff}(\mathrm T(B))\cong\mathbb R\).

If the hypotheses of Corollary \ref{Z=>SCP} hold and \(A\) has real rank zero,
then Theorem \ref{WTA-C-dyn-1} gives
\(B\in\mathrm{TA}_0(\mathcal S)\) with real-rank-zero models,
while Corollary \ref{Z=>SCP} gives strict
comparison for \(B\).  Theorem \ref{rr0-twa} therefore gives real rank
zero for \(B\).  This conclusion uses both approximation by
real-rank-zero models and strict comparison.
\end{rem}

\section{Towers and freeness conditions}

This section develops norm estimates for tower cut-downs and relates
WRP to Kishimoto's and Elliott's freeness conditions.  Equivariant
lifting of central sequence towers from finite von Neumann algebras
gives a further criterion for Elliott's condition.

\subsection{Towers in a \Cstar{}-dynamical system}

A tower in a topological dynamical system consists of a base and finitely
many pairwise disjoint translates.  The following central-sequence
formulation is its \Cstar{}-algebraic analogue.

\begin{defn}\label{central-sequence-tower}
Let $(A,\Gamma)$ be a unital \Cstar{}-dynamical system.  A \emph{Rokhlin
tower} is a pair $(c,F)$, where $F\subseteq\Gamma$ is finite and nonempty and
$c\in (A'\cap A_\omega)_+$ is a contraction of norm one, such that
the elements
\[
 g(c),\qquad g\in F,
\]
are mutually orthogonal.  Its support is the positive contraction
\[
 h_{(c,F)}:=\sum_{g\in F}g(c)\in A'\cap A_\omega.
\]

Let $D$ be a unital \Cstar{}-algebra containing $A_\omega$, let
$\tau\in\mathrm T(D)$ satisfy $\mathrm d_\tau(c)>0$, and let $\delta>0$.
The tower is \emph{$(\delta,\tau)$-supported in $D$} if
\begin{equation}\label{supported-tower-condition}
 \frac{\tau(c^l)}{\mathrm d_\tau(c)}>1-\delta,\qquad l\in\mathbb N.
\end{equation}
We omit the ambient algebra when it is clear from the context.
\end{defn}

A functional-calculus cut-down of the base gives
\eqref{supported-tower-condition} whenever the base has nonzero
dimension function for the chosen trace.

\begin{lem}\label{supported-tower-cut}
Let $(c,F)$ be a tower, let $D\supseteq A_\omega$ be unital, and let
$\tau\in\mathrm T(D)$ satisfy $\mathrm d_\tau(c)>0$.  For every $\delta>0$
there is a continuous function $f\colon[0,1]\to[0,1]$, with
$f(0)=0$ and $f(1)=1$, such that $(f(c),F)$ is a
$(\delta,\tau)$-supported tower.
\end{lem}

\begin{proof}
Let $\mu$ be the finite measure induced by $\tau$ on
$\mathrm C^*(1,c)$.  Then $\mathrm d_\tau(c)=\mu((0,1])$.  Choose $\eta>0$
such that
\[
 \mu([\eta,1])>(1-\delta)\mu((0,1]),
\]
and choose $f$ to be strictly positive on $(0,\eta)$ and equal to one
on $[\eta,1]$.  Thus $\mathrm d_\tau(f(c))=\mathrm d_\tau(c)$ and, for every
$l\in\mathbb N$,
\[
 \tau(f(c)^l)\geq\mu([\eta,1])
 >(1-\delta)\mathrm d_\tau(f(c)).
\]
Functional calculus preserves centrality and the orthogonality of the
translates, so $(f(c),F)$ is again a tower.
\end{proof}

\begin{rem}\label{tower-orbit-average}
A tower base also defines partial-orbit averages.  If
$\tau\in\mathrm T(A_\omega)$ and $\tau(c^2)>0$, then
\[
 \rho_{c,F,\tau}(a)
 =\frac{1}{|F|\tau(c^2)}
   \sum_{g\in F}\tau\big(c\,g^{-1}(a)c\big),\qquad a\in A,
\]
is a tracial state.  Moreover, for $k\in\Gamma$,
\begin{equation}\label{tower-average-invariance}
 \left|\rho_{c,F,\tau}(k(a))-\rho_{c,F,\tau}(a)\right|
 \leq \|a\|\frac{|F\mathbin\triangle k^{-1}F|}{|F|}.
\end{equation}
Indeed, each summand is tracial because $c$ commutes with $A$, and
\eqref{tower-average-invariance} follows by reindexing the two finite
sums.  Consequently, weak* limit points of such states along a left
F\o lner sequence are invariant.
\end{rem}

\subsubsection{Cut-downs by towers}

The following finite-set cutoff gives almost invariant tower
compressions.  Recall that
$\operatorname{Int}_K(F)=\{g\in F:Kg\subseteq F\}$.

\begin{lem}\label{folner-taper}
Let $K\subseteq\Gamma$ be finite and symmetric with $e\in K$, let
$m\in\mathbb N$, and let $F\subseteq\Gamma$ be finite.  If
$\operatorname{Int}_{K^m}(F)$ is nonempty, there is a function
$w\colon\Gamma\to[0,1]$ such that
\begin{enumerate}
\item $\operatorname{supp}(w)\subseteq\operatorname{Int}_K(F)$;
\item $w=1$ on $\operatorname{Int}_{K^m}(F)$; and
\item $|w(kg)-w(g)|\leq 1/m$ for $k\in K$ and $g\in\Gamma$.
\end{enumerate}
If $(c,F)$ is a tower and
\[
 q_w:=\sum_{g\in F}w(g)g(c),
\]
then $\|q_w\|=1$, $0\leq q_w\leq h_{(c,F)}$, and
\begin{equation}\label{taper-almost-invariant}
 \|k(q_w)-q_w\|\leq 1/m,\qquad k\in K.
\end{equation}
If $B=A\rtimes\Gamma$, $\tau\in\mathrm T(B_\omega)$, $0<\delta<1$,
the tower is $(\delta,\tau)$-supported, and $F$ is
$(K^m,\delta)$-invariant, then
\begin{equation}\label{taper-supported}
 \frac{\tau(q_w^l)}{\mathrm d_\tau(q_w)}>(1-\delta)^2,\qquad l\in\mathbb N.
\end{equation}
\end{lem}

\begin{proof}
Give $\Gamma$ the extended graph metric whose edges join $g$ to $kg$
for $k\in K$, and put
$d_F(g)=\operatorname{dist}(g,\Gamma\setminus F)$, allowing the value
$+\infty$.  The function
\[
 w(g)=\max\left\{0,\min\left\{1,\frac{d_F(g)-1}{m}\right\}\right\}
\]
is understood to be one when $d_F(g)=+\infty$.  It satisfies
(1)--(3), since distance to a fixed set is one-Lipschitz.  In particular, a nonzero coefficient of either
$q_w$ or $k(q_w)$ occurs at a level belonging to $F$.  Orthogonality
of the levels therefore gives \eqref{taper-almost-invariant}, as well
as $\|q_w\|=1$ and $q_w\leq h_{(c,F)}$.

Every trace on $B_\omega$ is invariant under the automorphisms
implemented by the constant unitaries $u_g$.  Orthogonality therefore
gives
\[
 \mathrm d_\tau(q_w)=|\operatorname{supp}(w)|\mathrm d_\tau(c)
\]
and
\[
 \tau(q_w^l)
 =\sum_{g\in F}w(g)^l\tau(c^l)
 \geq |\operatorname{Int}_{K^m}(F)|\tau(c^l).
\]
Since $|\operatorname{supp}(w)|\leq |F|$, the last two displays,
$(K^m,\delta)$-invariance, and
\eqref{supported-tower-condition} imply
\eqref{taper-supported}.
\end{proof}

\begin{prop}\label{full-tower-cut}
Let $(A,\Gamma)$ be a unital \Cstar{}-dynamical system, where $\Gamma$ is
countable, discrete, and amenable, and put $B=A\rtimes\Gamma$.  Assume
$\mathrm T(B)\ne\emptyset$ and every tracial state of $B$ is
faithful.  Let $a\in B_+$ have norm one.  For every $\eps>0$ there
are a finite set $L\subseteq\Gamma$ and $\delta>0$ such that, whenever
$(c,F)$ is a tower, $F$ is $(L,\delta)$-invariant, and $(c,F)$ is
$(\delta,\tau)$-supported for some $\tau\in\mathrm T(B_\omega)$, one
has
\[
 \|h_{(c,F)}a h_{(c,F)}\|>1-\eps.
\]
Here and below, elements of $B$ are regarded as constant sequences in
$B_\omega$.
\end{prop}
\begin{proof}
Suppose otherwise.  Then there is $\eps_0>0$ for which the conclusion
fails for every pair $(L,\delta)$.  Choose finite symmetric sets
\[
 e\in K_1\subseteq K_2\subseteq\cdots\subseteq\Gamma,
 \qquad \bigcup_nK_n=\Gamma,
\]
and numbers $\delta_n\downarrow0$.  Applying the assumed failure with
$L=K_n^n$, choose a $(K_n^n,\delta_n)$-invariant tower $(c_n,F_n)$
and $\tau_n\in\mathrm T(B_\omega)$ such that the tower is
$(\delta_n,\tau_n)$-supported and
\begin{equation}\label{tower-cut-contradiction}
 \|h_{(c_n,F_n)}a h_{(c_n,F_n)}\|\leq1-\eps_0.
\end{equation}
After discarding finitely many terms, the interiors
$\operatorname{Int}_{K_n^n}(F_n)$ are nonempty.  Apply Lemma
\ref{folner-taper}, with $K=K_n$ and $m=n$, and denote the resulting
positive contraction by $q_n$.

The elements $q_n$ commute with $A$.  By
\eqref{taper-almost-invariant}, they asymptotically commute with every
constant unitary $u_g$, and hence
\[
 \lim_{n\to\infty}\|[q_n,x]\|=0,\qquad x\in B.
\]
Equation \eqref{taper-supported} also gives, for every
$l\in\mathbb N$,
\[
 \lim_{n\to\infty}
 \frac{\tau_n(q_n^l)}{\mathrm d_{\tau_n}(q_n)}=1.
\]
Lemma \ref{central-cut}, applied to the constant embeddings
$B\to B_\omega$, now yields
\begin{equation}\label{taper-full-norm}
 \limsup_{n\to\infty}\|q_na q_n\|=1.
\end{equation}

The diagonal form of the cutoff gives a norm comparison with the
tower support.  In the bidual of the abelian \Cstar{}-algebra generated by
the orthogonal levels, there is a positive contraction $z_n$ with value
$w_n(g)$ on the $g$-summand, such that
\[
 q_n=z_nh_{(c_n,F_n)}=h_{(c_n,F_n)}z_n.
\]
Consequently,
\[
 \|q_na q_n\|
 =\|z_nh_{(c_n,F_n)}a h_{(c_n,F_n)}z_n\|
 \leq\|h_{(c_n,F_n)}a h_{(c_n,F_n)}\|,
\]
contradicting \eqref{tower-cut-contradiction} and
\eqref{taper-full-norm}.
\end{proof}

For positive elements of the coefficient algebra, the compression
by a tower support is an orthogonal sum of level compressions.
For general crossed-product elements, off-diagonal matrix
coefficients between levels may contribute to the norm.

\begin{cor}\label{full-tower-cut-cor}
Under the hypotheses of Proposition \ref{full-tower-cut}, let
$a\in A_+$ have norm one.  For every $\eps>0$ there are a finite set
$L\subseteq\Gamma$ and $\delta>0$ such that, whenever $(c,F)$ is a
$(\delta,\tau)$-supported tower for some
$\tau\in\mathrm T((A\rtimes\Gamma)_\omega)$ and $F$ is
$(L,\delta)$-invariant, there is $g\in F$ such that
\[
 \|g(c)a g(c)\|>1-\eps.
\]
\end{cor}

\begin{proof}
Put $h=h_{(c,F)}$.  Since every level $g(c)$ commutes with $a$ and
the levels are mutually orthogonal,
\[
 hah=\sum_{g\in F}g(c)a g(c)
\]
is an orthogonal sum of positive elements.  Hence
\[
 \|hah\|=\max_{g\in F}\|g(c)a g(c)\|,
\]
and the result follows from Proposition \ref{full-tower-cut}.
\end{proof}

\subsubsection{Order zero maps and towers}

\begin{lem}\label{tower-order-zero-map}
Let $(c,F)$ be a tower and enumerate
$F=\{g_1,\ldots,g_n\}$.  There is a c.p.c. order zero map
\[
 \phi\colon\mathrm M_n(\Comp)\longrightarrow
 (A\rtimes\Gamma)_\omega
\]
such that
\begin{equation}\label{tower-order-zero-formula}
 \phi(e_{ij})=u_{g_i}c u_{g_j}^*,\qquad 1\leq i,j\leq n,
\end{equation}
and $\phi(1)=h_{(c,F)}$.
\end{lem}

\begin{proof}
On elementary tensors define
\[
 \Theta(f\otimes e_{ij})=u_{g_i}f(c)u_{g_j}^*,
 \qquad f\in\mathrm C_0((0,1]).
\]
If $j\ne k$, tower orthogonality gives
\[
 f(c)u_{g_j}^*u_{g_k}\widetilde f(c)=0
 \qquad(f,\widetilde f\in\mathrm C_0((0,1])).
\]
The usual matrix multiplication therefore shows that $\Theta$ extends
to a *-homomorphism
\[
 \mathrm C_0((0,1])\otimes\mathrm M_n(\Comp)\longrightarrow
 (A\rtimes\Gamma)_\omega.
\]
If $\iota(t)=t$ is the canonical positive contraction in
$\mathrm C_0((0,1])$, then
$\phi(x)=\Theta(\iota\otimes x)$ is c.p.c. and order zero, and
\eqref{tower-order-zero-formula} follows directly.
\end{proof}

\subsection{Kishimoto's condition, Elliott's condition, and WRP}

Kishimoto-type conditions are noncommutative freeness conditions.  We
use the following simultaneous finite-set form to cut down Fourier
polynomials.

\begin{defn}\label{Kishimoto-Condition}
Let $(A,\Gamma)$ be a \Cstar{}-dynamical system, where $A\ne0$ and $\Gamma$
is discrete.  We say that $(A,\Gamma)$ satisfies \emph{Kishimoto's
condition} if, for every $a\in A_+$ with $\|a\|=1$, finite sets
$K\subseteq\Gamma$ and $\mathcal F\subseteq A$, and $\eps>0$, there
is a norm-one positive contraction $c\in A$ such that
\begin{enumerate}
\item $\|cac\|>1-\eps$; and
\item $\|cb\gamma(c)\|<\eps$ for every
$\gamma\in K\setminus\{e\}$ and $b\in\mathcal F$.
\end{enumerate}
\end{defn}

\begin{rem}
Definition \ref{Kishimoto-Condition} is a simultaneous,
norm-preserving version of the condition.  The literature also uses
``Kishimoto's condition'' for the pointwise hereditary-subalgebra
condition on each nontrivial automorphism.  No equivalence between
those formulations is needed here.
\end{rem}

Elliott's condition abstracts the conclusion of
\cite[Lemma~3.3]{Elliott-simple-crossed}.

\begin{defn}\label{Elliott-Condition}
A \Cstar{}-dynamical system $(A,\Gamma)$ satisfies \emph{Elliott's
condition} if, for every $a\in A_+$, finite sets
$\mathcal F\subseteq A$ and $K\subseteq\Gamma\setminus\{e\}$, and
$\eps>0$, there is a norm-one positive contraction $c\in A$ such that
\begin{enumerate}
\item $\|cac\|>\|a\|-\eps$;
\item $\|cb-bc\|<\eps$ for every $b\in\mathcal F$; and
\item $\|c\gamma(c)\|<\eps$ for every $\gamma\in K$.
\end{enumerate}
\end{defn}

\begin{lem}\label{Elliott-implies-Kishimoto}
Elliott's condition implies Kishimoto's condition.
\end{lem}

\begin{proof}
Fix the data in Definition \ref{Kishimoto-Condition}, and put
$M=\max(\{\|b\|:b\in\mathcal F\}\cup\{0\})$.  Apply Elliott's
condition with a tolerance
\[
 0<\eta<\frac{\eps}{1+M}.
\]
For $b\in\mathcal F$ and $\gamma\in K\setminus\{e\}$, the resulting
positive contraction $c$ satisfies
\[
 \|cb\gamma(c)\|
 \leq\|(cb-bc)\gamma(c)\|+\|bc\gamma(c)\|
 <(1+\|b\|)\eta\leq(1+M)\eta<\eps.
\]
After also requiring $\eta<\eps$, the compression estimate is the
one in Definition \ref{Kishimoto-Condition}.
\end{proof}

For separable coefficient algebras, Elliott's condition is precisely
the existence of finite central-sequence towers that preserve the norm
of a prescribed positive element.

\begin{prop}\label{Elliott-central-sequence}
Let $A$ be a nonzero, separable, unital \Cstar{}-algebra.  A \Cstar{}-dynamical
system $(A,\Gamma)$ satisfies Elliott's condition if and only if, for
every finite set $F\subseteq\Gamma$ containing $e$ and every
$a\in A_+$, there is a tower $(c,F)$ such that
\[
 \|cac\|=\|a\|.
\]
\end{prop}

\begin{proof}
Assume first that Elliott's condition holds.  Let
$\mathcal F_1\subseteq\mathcal F_2\subseteq\cdots$ be finite sets
whose union is dense in $A$.  For each $n$, apply the condition with
$\mathcal F_n$, tolerance $1/n$, and
\[
 K=F^{-1}F\setminus\{e\}.
\]
Denote the resulting norm-one positive contraction by $c_n$.  Then
$c=[(c_n)]\in A'\cap A_\omega$ is a norm-one positive contraction.
If $g,h\in F$ are distinct, then
\[
 g(c)h(c)=g\big(c\,(g^{-1}h)(c)\big)=0,
\]
so $(c,F)$ is a tower.  The compression estimates and the reverse
inequality $\|c_nac_n\|\leq\|a\|$ give
$\|cac\|=\|a\|$.

Conversely, fix $a$, $\mathcal F$, $K$, and $\eps$ as in Definition
\ref{Elliott-Condition}, and apply the stated tower property to
$F=K\cup\{e\}$.  Represent the tower base by norm-one positive
contractions $(c_n)$ in $A$.  Centrality, tower orthogonality, and
$\|cac\|=\|a\|$ imply that the three finite-stage estimates hold
on a set of indices belonging to $\omega$.  Any coordinate in
this set verifies Elliott's condition.
\end{proof}

\begin{lem}\label{WRP-Kishimoto}
Let $(A,\Gamma)$ be a minimal \Cstar{}-dynamical system with WRP, where
$A$ is unital and separable and $\Gamma$ is countable, discrete, and
amenable.  If $\mathrm T(A)\ne\emptyset$, then $(A,\Gamma)$
satisfies Elliott's condition and hence Kishimoto's condition.
\end{lem}

\begin{proof}
Put $B=A\rtimes_{\mathrm r}\Gamma$.  By amenability, the action of
$\Gamma$ on the nonempty compact convex set $\mathrm T(A)$ has a
fixed point.  Its composition with the canonical conditional
expectation is a tracial state on $B$, so $\mathrm T(B)\ne\emptyset$.
Corollary \ref{min-faithful-trace} shows that every tracial state of
$B$ is faithful.

We first prove Elliott's condition for a norm-one element $a\in A_+$.
Fix finite sets $K\subseteq\Gamma\setminus\{e\}$ and
$\mathcal F\subseteq A$, and fix $\eps>0$.  Enlarge $K\cup\{e\}$
to a finite symmetric set $K_0$.  Let $(L,\delta)$ be supplied by
Corollary \ref{full-tower-cut-cor} for $a$ and $\eps/2$, and enlarge
$L$ so that $e\in L$.  Choose $0<\eta<1$ so that
$\eta/(1-\eta)<\delta$.  By the exact tiling theorem
\cite[Theorem~4.3 and Lemma~2.6]{DHZ-tiling}, there is an exact tiling
whose shapes $D_1,\ldots,D_S$ are $(K_0L,\eta)$-invariant.  An
expansion-invariance tolerance below $\eta/|K_0L|$ gives this
interior estimate.  For finite $\Gamma$, take the tiling with
a single shape
$D_1=\Gamma$.  Set
\[
 F_s=\operatorname{Int}_{K_0}(D_s),\qquad 1\leq s\leq S.
\]
Since
\[
 \operatorname{Int}_{K_0L}(D_s)
 \subseteq\operatorname{Int}_L(F_s)
 \quad\text{and}\quad
 |F_s|\geq|\operatorname{Int}_{K_0L}(D_s)|
 >(1-\eta)|D_s|,
\]
we have
\[
 \frac{|F_s\setminus\operatorname{Int}_L(F_s)|}{|F_s|}
 <\frac{\eta}{1-\eta}<\delta.
\]
Thus every $F_s$ is $(L,\delta)$-invariant.

Apply WRP to this tiling, and let
$c_s\in A'\cap A_\omega$ be the resulting bases.  Fix
$\sigma\in\mathrm T(B)$ and let $\tau=\sigma_\omega$ be the induced
trace on $B_\omega$.  Put
\[
 h=\sum_{s=1}^S\sum_{g\in D_s}g(c_s).
\]
The WRP defect condition gives $\tau((1-h)^2)=0$.  By
Cauchy--Schwarz, $\tau(1-h)=0$, and hence
\[
 1=\tau(h)=\sum_{s=1}^S|D_s|\tau(c_s).
\]
Choose $s$ with $\tau(c_s)>0$.  In particular,
$\mathrm d_\tau(c_s)>0$ and $c_s\ne0$.  Apply Lemma
\ref{supported-tower-cut} to $c_s/\|c_s\|$.  It gives a function $f$
such that, with
\[
 q=f(c_s/\|c_s\|),
\]
the tower $(q,D_s)$ is $(\delta,\tau)$-supported.  Its restriction
$(q,F_s)$ has the same support property.  Corollary
\ref{full-tower-cut-cor} therefore gives $g_0\in F_s$ such that
\[
 \|g_0(q)a g_0(q)\|>1-\eps/2.
\]

Set $d=g_0(q)$.  Then $d$ is a norm-one positive contraction in
$A'\cap A_\omega$.  For $k\in K$, the definition of $F_s$ gives
$kg_0\in D_s$; since $k\ne e$, tower orthogonality gives
\[
 d\,k(d)=g_0(q)\,(kg_0)(q)=0.
\]
Represent $d$ by norm-one positive contractions $(d_n)$ in $A$.
Centrality, the preceding orthogonality relations, and the compression
estimate hold simultaneously on a set of indices belonging to
$\omega$.  Thus some $d_n$
satisfies the three estimates in Elliott's condition with tolerance
$\eps$.

For a nonzero $a\in A_+$ of arbitrary norm, apply the norm-one case
to $a/\|a\|$ with a smaller tolerance and rescale the compression
estimate.  When $a=0$, apply the construction to $1_A$; the
compression estimate for $a$ is automatic.  This proves
Elliott's condition;
Lemma \ref{Elliott-implies-Kishimoto} gives the final assertion.
\end{proof}

\begin{lem}\label{coefficient-cutdown}
Let $(A,\Gamma)$ be a \Cstar{}-dynamical system satisfying Kishimoto's
condition.  For every nonzero $x\in(A\rtimes_{\mathrm r}\Gamma)_+$,
there is a nonzero element $h\in A_+$ such that
\[
 h\precsim x\quad\text{in }A\rtimes_{\mathrm r}\Gamma.
\]
\end{lem}

\begin{proof}
Put $B=A\rtimes_{\mathrm r}\Gamma$, and let
$\mathbb E\colon B\to A$ be the canonical conditional expectation.
Its faithfulness gives $\mathbb E(x)\ne0$.  Multiplying $x$ by a
positive scalar does not change its Cuntz class, so we may assume
\[
 \|\mathbb E(x)\|=1.
\]
Fix $0<\eta<1/3$.  Density of the algebraic crossed product gives a
Fourier polynomial
\[
 y=\sum_{\gamma\in F}a_\gamma u_\gamma,
 \qquad e\in F,
\]
such that
\begin{equation}\label{coefficient-polynomial-approximation}
 \|x-y\|<\eta
 \quad\text{and}\quad
 a_e=\mathbb E(x).
\end{equation}
Indeed, after choosing an arbitrary sufficiently close Fourier
polynomial, replace its identity coefficient by $\mathbb E(x)$; the
contractivity of $\mathbb E$ controls the additional error.

Set $N=\max\{1,|F\setminus\{e\}|\}$.  Apply Kishimoto's condition to
$a_e$, the group set $F\setminus\{e\}$, the coefficient set
$\{a_\gamma:\gamma\in F\setminus\{e\}\}$, and tolerance $\eta/N$.
We obtain a norm-one positive contraction $c\in A$ such that
\[
 \|ca_ec\|>1-\eta
 \quad\text{and}\quad
 \|ca_\gamma\gamma(c)\|<\frac{\eta}{N},
 \qquad\gamma\in F\setminus\{e\}.
\]
Using $u_\gamma c=\gamma(c)u_\gamma$ and
\eqref{coefficient-polynomial-approximation}, we get
\begin{align*}
 \|cxc-ca_ec\|
 &\leq\|c(x-y)c\|
   +\sum_{\gamma\in F\setminus\{e\}}
      \|ca_\gamma\gamma(c)\| \\
 &<2\eta.
\end{align*}
The standard perturbation lemma for positive elements gives
\[
 h:=(ca_ec-2\eta)_+\precsim cxc.
\]
Moreover, $h\ne0$ because
$\|ca_ec\|>1-\eta>2\eta$.  Finally,
$cxc\precsim x$: indeed, $cxc=(x^{1/2}c)^*(x^{1/2}c)$ is Cuntz
equivalent to
$x^{1/2}c^2x^{1/2}\in\overline{xBx}$.  Consequently
$0\ne h\in A_+$ and $h\precsim x$, as required.
\end{proof}

\begin{cor}\label{wrp-simple}
Let $(A,\Gamma)$ be a minimal \Cstar{}-dynamical system with WRP, where
$A$ is unital and separable, $\mathrm T(A)\ne\emptyset$, and $\Gamma$
is countable, discrete, and amenable.  Then
$A\rtimes_{\mathrm r}\Gamma$ is simple.  Since $\Gamma$ is amenable,
the same statement holds without distinguishing the full and reduced
crossed products.
\end{cor}

\begin{proof}
Put $B=A\rtimes_{\mathrm r}\Gamma$, and let $I\triangleleft B$ be a
nonzero ideal.  Choose $0\ne x\in I_+$.  Lemma
\ref{WRP-Kishimoto} and Lemma \ref{coefficient-cutdown} give a
nonzero $h\in A_+$ with $h\precsim_B x$.  Since ideals are closed
under Cuntz subequivalence, $h\in I$.  Thus $I\cap A$ is a nonzero
$\Gamma$-invariant ideal of $A$.  Minimality gives $I\cap A=A$, so
$1_A\in I$ and hence $I=B$.
\end{proof}

\subsection{Lifting towers and Elliott's condition}

The following equivariant form of the central-sequence lifting
argument of Kirchberg and R{\o}rdam lifts orthogonality in a
quotient to exact orthogonality in the norm ultrapower while
preserving centrality.

\begin{lem}[Equivariant orthogonal lifting]\label{lift-towers-0}
Let $A\subseteq N$ be a unital inclusion of \Cstar{}-algebras, with $A$
separable, and let $\Gamma\curvearrowright N$ be an action which
leaves $A$ invariant.  Let $J\triangleleft\ell^\infty(N)$ be a
$\Gamma$-invariant closed ideal, and put
\[
 N^\omega=\ell^\infty(N)/J.
\]
Suppose that
\begin{enumerate}
\item $c_\omega(A)\subseteq J$ and $J\cap N=\{0\}$, where $N$ is
viewed as the algebra of constant sequences;
\item the homomorphism induced by the inclusion,
\[
 \Phi\colon A_\omega=\ell^\infty(A)/c_\omega(A)\longrightarrow
 N^\omega,
\]
is surjective; and
\item there is a norm-continuous function
$\theta\colon N\to[0,\infty)$, bounded on bounded subsets, such that
\[
 J=\{(x_n)\in\ell^\infty(N):
          \lim_{n\to\omega}\theta(x_n)=0\}.
\]
\end{enumerate}
Then the following statements hold.
\begin{enumerate}
\item[(i)] Let $p_1,\ldots,p_S$ be positive contractions in
$N'\cap N^\omega$, and let $E_1,\ldots,E_S$ be finite subsets of
$\Gamma$.  If
\[
 \{\gamma(p_s):1\leq s\leq S,\ \gamma\in E_s\}
\]
is pairwise orthogonal, then there are positive contractions
$c_1,\ldots,c_S\in A'\cap A_\omega$ such that
\[
 \Phi(c_s)=p_s\quad(1\leq s\leq S)
\]
and
\[
 \{\gamma(c_s):1\leq s\leq S,\ \gamma\in E_s\}
\]
is pairwise orthogonal.  In particular,
\[
 N'\cap N^\omega\subseteq\Phi(A'\cap A_\omega).
\]
\item[(ii)] Let $p\in N'\cap N^\omega$ be a positive contraction
and let $K\subseteq\Gamma$ be finite.  If
\[
 p\gamma(p)=0\quad(\gamma\in K),
\]
then there is a positive contraction $c\in A'\cap A_\omega$ such
that
\[
 \Phi(c)=p
 \quad\text{and}\quad
 c\gamma(c)=0\quad(\gamma\in K).
\]
\end{enumerate}
\end{lem}

\begin{proof}
Write
\[
 q_J\colon\ell^\infty(N)\longrightarrow N^\omega
 \quad\text{and}\quad
 q_A\colon\ell^\infty(A)\longrightarrow A_\omega
\]
for the quotient maps.  We prove (i).  By surjectivity and positive
lifting, for each $s$ there is a sequence
$b^{(s)}=(b_n^{(s)})$ of positive contractions in $A$ such that
$q_J(b^{(s)})=p_s$.

Let $H$ be the subgroup of $\Gamma$ generated by
$\bigcup_sE_s$.  It is countable.  The \Cstar{}-algebra
\[
 C=\mathrm C^*\big(A,\gamma(b^{(s)}):
      \gamma\in H,\ 1\leq s\leq S\big)
 \subseteq\ell^\infty(A)
\]
is separable and $H$-invariant.  Set $I=C\cap J$.  Then $I$ is an
$H$-invariant ideal of $C$.  Centrality of the $p_s$ gives
\begin{equation}
 [b^{(s)},a]\in I
 \quad(a\in A,\ 1\leq s\leq S).
 \label{lift-central-defect}
\end{equation}
and the assumed orthogonality gives
\begin{equation}
 \gamma(b^{(s)})\eta(b^{(t)})\in I
 \label{lift-orthogonal-defect}
\end{equation}
whenever $(s,\gamma)\ne(t,\eta)$, with
$\gamma\in E_s$ and $\eta\in E_t$.

Choose a quasicentral approximate identity
$(e^{(m)})_{m=1}^\infty$ for $I$ in $C$, consisting of positive
contractions.  Since $I$ and $C$ are $H$-invariant,
$(\gamma(e^{(m)}))_m$ is again a quasicentral approximate identity
for $I$ in $C$ for every $\gamma\in H$.  Consequently,
\begin{equation}
 \lim_{m\to\infty}
 \big\|[(1-e^{(m)})b^{(s)}(1-e^{(m)}),a]\big\|=0
 \quad(a\in A).
 \label{lift-central-cut}
\end{equation}
Indeed, modulo terms containing $[e^{(m)},a]$, the commutator in
\eqref{lift-central-cut} is
\[
 (1-e^{(m)})[b^{(s)},a](1-e^{(m)}),
\]
which tends to zero by \eqref{lift-central-defect}.  Similarly,
quasicentrality allows the cutoff factors to be moved past the
finitely many elements occurring in
\eqref{lift-orthogonal-defect}; since every translate of
$(e^{(m)})_m$ is an approximate identity for $I$, we obtain
\begin{equation}
 \lim_{m\to\infty}
 \big\|\gamma((1-e^{(m)})b^{(s)}(1-e^{(m)}))
       \eta((1-e^{(m)})b^{(t)}(1-e^{(m)}))\big\|=0
 \label{lift-orthogonal-cut}
\end{equation}
for all the pairs in \eqref{lift-orthogonal-defect}.

Apply Kirchberg's epsilon test \cite[Lemma~3.1]{KR-CenSeq} as
follows.  Choose a dense sequence $(a_j)_{j=1}^\infty$ in the unit
ball of $A$.  At the $n$th coordinate, the variable ranges over the
positive contractions of $A$, with test functions
\[
 x\longmapsto\theta(x),
\]
\[
 x\longmapsto
 \big\|[(1-x)b_n^{(s)}(1-x),a_j]\big\|,
\]
and, for every pair occurring in
\eqref{lift-orthogonal-defect}, the function
\[
 x\longmapsto
 \big\|\gamma((1-x)b_n^{(s)}(1-x))
       \eta((1-x)b_n^{(t)}(1-x))\big\|.
\]
For a fixed $m$, use the $n$th coordinate of $e^{(m)}$ as the test
variable.  Since $e^{(m)}\in I\subseteq J$, the ultralimit of the
first function is zero.  Equations
\eqref{lift-central-cut} and \eqref{lift-orthogonal-cut} show that,
as $m\to\infty$, the ultralimits of every prescribed finite
collection of the remaining functions tend to zero.  The epsilon
test gives a sequence $d=(d_n)$ of positive
contractions in $A$ such that
\[
 \lim_{n\to\omega}\theta(d_n)=0
\]
and all the commutator and orthogonality norms above have
ultralimit zero.

Thus $d\in J$.  Define
\[
 c_s=q_A\big((1-d)b^{(s)}(1-d)\big),\qquad 1\leq s\leq S.
\]
These are positive contractions.  The centrality tests imply that
$c_s\in A'\cap A_\omega$, while $q_J(d)=0$ gives
$\Phi(c_s)=p_s$.  The remaining tests give the required pairwise
orthogonality of the translates.  The final inclusion in (i)
follows by scaling positive elements and then writing an arbitrary
element of $N'\cap N^\omega$ as a linear combination of four
positive elements of that \Cstar{}-algebra.

For (ii), use the same construction with one lift $b$ of $p$ and
with only the orthogonality tests indexed by the pairs
$(e,\gamma)$, $\gamma\in K$.  The conclusion of those tests is
$c\gamma(c)=0$ for every $\gamma\in K$, as required.
\end{proof}

\begin{cor}\label{lift-towers}
Let $N$ be a finite von Neumann algebra with faithful normal tracial
state $\rho$, and let $A\subseteq N$ be a unital separable weakly
dense \Cstar{}-subalgebra.  Suppose that $\Gamma\curvearrowright N$ is a
$\rho$-preserving action which leaves $A$ invariant.  Form the
tracial ultrapower
\[
 N^\omega=\ell^\infty(N)/J_\rho,\qquad
 J_\rho=\{(x_n)\in\ell^\infty(N):
              \lim_{n\to\omega}\|x_n\|_{2,\rho}=0\}.
\]
If $p\in N'\cap N^\omega$ is a positive contraction and
$K\subseteq\Gamma$ is finite with
\[
 p\gamma(p)=0\quad(\gamma\in K),
\]
then there is a positive contraction $c\in A'\cap A_\omega$ such
that
\[
 \Phi(c)=p
 \quad\text{and}\quad
 c\gamma(c)=0\quad(\gamma\in K),
\]
where $\Phi\colon A_\omega\to N^\omega$ is induced by the inclusion.
\end{cor}

\begin{proof}
The ideal $J_\rho$ is $\Gamma$-invariant, contains $c_\omega(A)$,
and has trivial intersection with the constant copy of $N$.
Moreover, $\theta(x)=\|x\|_{2,\rho}$ satisfies condition (3) of
Lemma \ref{lift-towers-0}.  Finally, coordinatewise Kaplansky
density gives the surjectivity of $\Phi$: a bounded sequence in
$N$ can be approximated, coordinate by coordinate and without
increasing its norm bound, by a sequence in $A$ in the
$\|\cdot\|_{2,\rho}$-norm.  Lemma \ref{lift-towers-0}(ii) now
applies.
\end{proof}

Ocneanu's central-freeness and paving arguments give the following
local finite-set result.

\begin{lem}[Ocneanu's small-tower lemma]\label{Ocneanu-small-tower}
Let $\Gamma$ be a countable discrete group acting on a
finite von Neumann algebra $N$ with separable predual.  Suppose that
$\rho$ is a faithful normal $\Gamma$-invariant tracial state and
that the action is centrally free.  Given a finite set
$K\subseteq\Gamma\setminus\{e\}$ and $\delta>0$, there are
projections $p_0,p_1,\ldots,p_m\in N'\cap N^\omega$ such that
\[
 \sum_{j=0}^m p_j=1,\qquad
 \rho_\omega(p_0)<\delta,
\]
and
\[
 p_j\gamma(p_j)=0
 \quad(1\leq j\leq m,\ \gamma\in K).
\]
\end{lem}

\begin{proof}
If $K$ is empty, take $p_0=0$ and $p_1=1$. Otherwise replace $K$
by $K\cup K^{-1}$. The lemma in \cite[Section~5.6]{Ocneanu-thesis}
shows that the induced action on $R=N'\cap N^\omega$ is properly
outer on $S'\cap R$ for every countably generated invariant
von Neumann subalgebra $S\subseteq R$.
Apply the local small-tower construction of
\cite[Section~6.3]{Ocneanu-thesis} with remainder tolerance $\delta/2$;
see also the monograph \cite{Ocneanu}.
This gives the stated strict bound and exact orthogonality.
The local construction is valid for every countable discrete group
under the stated hypotheses.  Section~6.2 uses Popa's paving lemma
in the finite reduced von Neumann crossed product, and Section~6.3
uses maximality and countable reindexing to obtain exact
orthogonality from arbitrarily small overlaps.  These steps do not
require the amenability assumption of the full Rokhlin theorem.
\end{proof}

\begin{thm}\label{central-free-Elliott}
Let $(A,\Gamma)$ be a \Cstar{}-dynamical system, where $A$ is a nonzero
separable unital \Cstar{}-algebra and $\Gamma$ is a discrete
group. Suppose that $A$ has a faithful
$\Gamma$-invariant tracial state $\tau$.  If the induced action on
\[
 N=\pi_\tau(A)''
\]
is centrally free, then $(A,\Gamma)$ satisfies Elliott's condition.
\end{thm}

\begin{proof}
Identify $A$ with its image in $N$, and denote by $\rho$ the
faithful normal extension of $\tau$ to $N$.  Fix
$a\in A_+$, finite sets $\mathcal F\subseteq A$ and
$K\subseteq\Gamma\setminus\{e\}$, and $\eps>0$.
The subgroup $H=\langle K\rangle$ is countable, and its induced
action on $N$ retains central freeness and the invariant faithful
trace.

First suppose that $a\ne0$, and put
\[
 r=\max\{\|a\|-\eps,0\}.
\]
Choose $f\in \mathrm C([0,\|a\|])_+$ such that
\[
 0\leq f\leq1,\qquad
 f|_{[0,r]}=0,\qquad f(\|a\|)=1,
\]
and set $b=f(a)$.  Then $b\ne0$, and faithfulness gives
$\rho(b)>0$.  If $a=0$, instead put $b=1_A$.  In either case apply
Lemma \ref{Ocneanu-small-tower} with
\[
 \delta=\frac{\rho(b)}2.
\]
We obtain a partition of unity
$p_0,p_1,\ldots,p_m\in N'\cap N^\omega$ with
$\rho_\omega(p_0)<\rho(b)/2$ and
$p_j\gamma(p_j)=0$ for $j\geq1$ and $\gamma\in K$.

The projections commute with the constant copy of $N$, and hence
with $b$.  Therefore
\[
 \rho(b)=\sum_{j=0}^m\rho_\omega(p_jbp_j).
\]
Since $0\leq b\leq1$,
\[
 \rho_\omega(p_0bp_0)\leq\rho_\omega(p_0)<\frac{\rho(b)}2.
\]
It follows that $p_jbp_j\ne0$ for some $j\geq1$; fix such a $j$ and
write $p=p_j$.  If $a\ne0$, functional calculus in the corner
$pN^\omega p$ gives
\[
 p f(a)p=f(pap)\ne0.
\]
As $f$ vanishes on $[0,r]$, this implies
\begin{equation}
 \|pap\|>r\geq\|a\|-\eps.
 \label{large-norm-cut}
\end{equation}
If $a=0$, the required compression inequality holds because
$0>-\eps$.  In either case, $p$ is a nonzero projection.

Apply Corollary \ref{lift-towers} to $p$ and $K$.  It gives a
positive contraction $c\in A'\cap A_\omega$ satisfying
\[
 \Phi(c)=p,\qquad c\gamma(c)=0\quad(\gamma\in K).
\]
Since $\Phi$ is contractive and $p$ is a nonzero projection,
$\|c\|=1$.  Moreover,
\[
 \Phi(cac)=pap,
\]
and hence, by \eqref{large-norm-cut},
\[
 \|cac\|\geq\|pap\|>\|a\|-\eps.
\]

Choose a representing sequence $(c_n)$ of norm-one positive
contractions for $c$.  Such a representative exists because any
representative by positive contractions has
$\lim_{n\to\omega}\|c_n\|=1$ and may be normalized without changing
its class.  Centrality, the relations $c\gamma(c)=0$, and the last
compression estimate imply that the sets of indices on which
\[
 \|c_nx-xc_n\|<\eps\quad(x\in\mathcal F),\qquad
 \|c_n\gamma(c_n)\|<\eps\quad(\gamma\in K),
\]
and
\[
 \|c_nac_n\|>\|a\|-\eps
\]
all belong to $\omega$.  Their finite intersection is nonempty, and
any coordinate in it verifies Elliott's condition.
\end{proof}

\begin{rem}
The lifting lemmas and the local small-tower argument do not require
amenability.  Since Elliott's condition concerns only finitely many
group elements at a time, Theorem \ref{central-free-Elliott} applies
to arbitrary discrete groups.  The prescribed-tiling WRP and the
large-tower approximation arguments above use amenability.
\end{rem}

\chapter{Comparison of positive elements}\label{section-comparison}

The tower constructions of Chapter~\ref{section-rokhlin} give local
models for crossed products with uniformly small tracial complements.
Comparison allows these complements to be absorbed into prescribed
positive elements.  This chapter establishes comparison results for
tracially approximated \Cstar{}-algebras, introduces comparison of coefficient
elements through property (COS), and defines the orbit radius of
comparison.  The automatic COS criterion permits arbitrary discrete
groups; the Rokhlin and orbit-radius arguments use amenability.
Comparison of finite orbit matrices then passes through
Rokhlin towers to the crossed product.  The final section computes orbit
radii and gives examples of property (COS).

\section{Radius of comparison and strict comparison under weak tracial approximation}\label{radius-comparison-section}

\begin{lem}\label{constant-ultrapower-comparison}
Let $a,b\in A_+$.  Regarding $a$ and $b$ as constant sequences, one
has $a\precsim b$ in $A$ if and only if $a\precsim b$ in $A_\omega$.
\end{lem}

\begin{proof}
The forward implication follows from the inclusion of $A$ in
$A_\omega$.  For the converse, suppose that $a\precsim b$ in
$A_\omega$, and fix $\eps>0$.  There is $x\in A_\omega$ such that
$\|a-xbx^*\|<\eps$.  Choose a bounded lift $(x_k)_k$ of $x$.  On an
$\omega$-large set of indices,
\[
 \|a-x_kbx_k^*\|<\eps.
\]
For each such index, the norm perturbation lemma for Cuntz
comparison gives
\((a-\eps)_+\precsim x_kbx_k^*\precsim b\) in $A$.  Since $\eps$ is
arbitrary, $a\precsim b$ in $A$.
\end{proof}

\begin{lem}\label{comp-local-global}
Let $A$ be a \Cstar{}-algebra, and let $a=[(a_k)]$, $b=[(b_k)]$, and
$c=[(c_k)]$ be positive elements of $A_\omega$, represented by bounded
sequences of positive elements of $A$.  Suppose that there is an
$\omega$-large set $U\subseteq\mathbb N$ such that, for every $k\in U$,
\begin{enumerate}
\item $a_k \precsim b_k$, and
\item $b_kc_k=b_k$.
\end{enumerate}
Then $a \precsim c$ in $A_\omega$.
\end{lem}
\begin{proof}
Fix $\eps>0$.  For $k\in U$, the usual cut-down form of Cuntz
subequivalence gives $x_k\in A$ such that
\[
 (a_k-\eps)_+=x_kb_kx_k^*.
\]
Put $s_k=x_kb_k^{1/2}$ for $k\in U$ and $s_k=0$ otherwise.  The sequence
$(s_k)_k$ is bounded, since
$s_ks_k^*=(a_k-\eps)_+$ on $U$.  As $b_k$ and $c_k$ are positive,
$b_kc_k=b_k$ also gives $c_kb_k=b_k$, and hence
$b_k^{1/2}c_kb_k^{1/2}=b_k$.  Therefore, in $A_\omega$,
\[
 [(s_k)]c[(s_k)]^*=(a-\eps)_+.
\]
Thus $(a-\eps)_+\precsim c$ for every $\eps>0$, and consequently
$a\precsim c$.
\end{proof}

\begin{cnv}[Gap functions and pointed comparison]\label{pointed-comparison-convention}
Write $\mathrm T^+(C)$ for the cone of densely finite lower
semicontinuous traces on $C$, including the zero trace, with the
topology of pointwise convergence on the Pedersen ideal; see
\cite{ESR-Cuntz}.  Matrix extensions are unnormalized.  Following
\cite[Definition~5.4]{Elliott-Niu-Villadsen-II}, a \emph{gap function}
for $C$ is a finite-valued, nonnegative continuous affine function
$\gamma:\mathrm T^+(C)\to[0,\infty)$, with $\gamma(0)=0$, such that
\begin{equation}\label{eq:gap-function}
 \mathrm d_\tau(a)+\gamma(\tau)<\mathrm d_\tau(b)
 \quad(0\ne\tau\in\mathrm T^+(C))
 \quad\Longrightarrow\quad a\precsim b,
 \qquad a,b\in(C\otimes\mathcal K)_+.
\end{equation}
Denote the set of gap functions by $G_C$.  Affinity and the condition
at zero imply positive homogeneity.  The zero trace is excluded from
the strict inequality.  Densely finite traces on nonunital models may
be unbounded.  For a simple unital algebra $A$, $\mathrm T^+(A)$ consists precisely of the
nonnegative multiples of tracial states.

A \emph{pointed \Cstar{}-algebra} is a pair $(C,e)$, where
$e\in C_+$ is a norm-one full element.  Following
\cite[Definition~3.3.2]{RC-5-person}, its relative radius of
comparison $r_{C,e}$ is the infimum of the numbers $s>0$ for which
the following condition holds:
\begin{equation}\label{eq:relative-radius-R1}
 \lambda(x)+s\lambda([e])\leq\lambda(y)
 \quad(\lambda\in F(\operatorname{Cu}(C)))
 \quad\Longrightarrow\quad x\leq y.
\end{equation}
We write
$A\in\mathrm{TA}_1^\bullet(\mathcal S_\bullet)$ if the witnesses in
Definition~\ref{defn_WAS} can always be chosen so that $(C,h)$ belongs
to $\mathcal S_\bullet$.  Thus the distinguished element compresses the
finite set.  In particular, this implies $A\in\mathrm{TA}_1(|\mathcal S_\bullet|)$ for the underlying
class $|\mathcal S_\bullet|=\{C:(C,e)\in\mathcal S_\bullet\}$.
For the affine-gap assertion below, the witnesses are required to
have $e$ strictly positive in $C$, so that $C=\overline{eCe}$.
Strict positivity is stronger than fullness and permits norm errors to
be estimated using $\mathrm d_\tau(e)=\|\tau\|$ for bounded traces on $C$.
\end{cnv}

%
%
%

\begin{lem}\label{comparison-tools}
Let $A$ be a simple, separable, exact, unital, infinite-dimensional
\Cstar{}-algebra with $\mathrm T(A)\ne\emptyset$.
\begin{enumerate}
\item[(1)] If $0\ne b\in \mathrm M_l(A)_+$ and $\eta>0$, there are
orthogonal elements $b_0,b_1\in\overline{b\mathrm M_l(A)b}_+$, with
$b_1\ne0$, such that
\begin{equation}\label{eq:reserve-trace}
 \mathrm d_\tau(b_0)>\mathrm d_\tau(b)-\eta,
 \qquad \tau\in\mathrm T(A).
\end{equation}
Moreover, for every $L\in\mathbb N$ there is $0\ne z\in A_+$ such that
\begin{equation}\label{eq:reserve-cuntz}
 L[z]\leq[b_1]\quad\hbox{in }\operatorname{Cu}(A).
\end{equation}

\item[(2)] Put $e=[1_A]$.  Suppose that $u,v\in\operatorname{Cu}(A)$
are bounded by finite multiples of $e$, that $v\ne0$, and that, for
some $s>0$,
\begin{equation}\label{eq:functional-gap}
 \lambda(u)+s<\lambda(v)
\end{equation}
for every $\lambda\in F(\operatorname{Cu}(A))$ satisfying
$\lambda(e)=1$.  If $0\leq t<s$ and $u_0\ll u$, then there are
$N,M\in\mathbb N$ such that
\begin{equation}\label{eq:stable-gap}
 \frac{M}{N}>t,
 \qquad (N+1)u_0+Me\leq Nv.
\end{equation}

\item[(3)] Let $\mathcal S$ be a class of \Cstar{}-algebras and suppose that
$A\in\mathrm{TA}_1(\mathcal S)$.  Given $0\ne z\in A_+$, there are
$C_j\in\mathcal S$ and norm-one positive contractions $h_j\in C_j$
such that, with
\[
 h=[(h_j)]\in A_\omega,\qquad
 Q=\pi\left(\prod_{j=1}^\infty C_j\right)\subseteq A_\omega,
\]
one has
\begin{equation}\label{eq:central-ta-cut}
 h\in A'\cap A_\omega,\qquad
 h\mathrm M_l(A)h\subseteq \mathrm M_l(Q)\quad(l\in\mathbb N),\qquad
 1-h^2\precsim z.
\end{equation}
The c.p.c. order zero map
\[
 \Phi:A\longrightarrow Q,\qquad \Phi(a)=hah,
\]
induces an additive, order-preserving map on Cuntz classes,
denoted by $\operatorname{Cu}(\Phi)$.  Compact containment in the
argument below comes from spectral cutoffs of its images.
For $x,y\in\operatorname{Cu}(A)$
represented by positive elements in finite matrix algebras over $A$,
the following assertions hold.
If $\operatorname{Cu}(C)$ is almost unperforated for every
$C\in\mathcal S$, then
\begin{equation}\label{eq:local-almost-unperf}
 (k+1)x\leq ky
 \quad\Longrightarrow\quad
 \operatorname{Cu}(\Phi)(x)\leq\operatorname{Cu}(\Phi)(y).
\end{equation}

For the quantitative assertion, suppose that the witnesses can be
chosen from a class $\mathcal S_\bullet$ of pointed algebras $(C,e)$,
where $e$ is full, with a fixed number $r\geq0$.  Assume that every
$(C,e)\in\mathcal S_\bullet$ satisfies at least one of the following
conditions:
\begin{enumerate}
\item[(i)] $r_{C,e}\leq r$;
\item[(ii)] there is $\gamma_C\in G_C$ with
$\gamma_C(\tau)\leq r \mathrm d_\tau(e)$ for every
$\tau\in\mathrm T^+(C)$; when $r=0$, this means $\gamma_C=0$.
\end{enumerate}
If $A\in\mathrm{TA}_1^\bullet(\mathcal S_\bullet)$, put
$L=\lfloor r\rfloor+1$.  Identifying the images from
$\operatorname{Cu}(Q)$ with their images in
$\operatorname{Cu}(A_\omega)$, one then has
\begin{equation}\label{eq:local-radius}
 (k+1)x+m[1_A]\leq ky,\quad \frac mk>r
 \quad\Longrightarrow\quad
 \operatorname{Cu}(\Phi)(x)
 \leq\operatorname{Cu}(\Phi)(y)+L[1-h^2].
\end{equation}
No uniform cutoff modulus for the ranks of the compressing elements
is required.
\end{enumerate}
\end{lem}

\begin{proof}
For (1), use the unnormalized matrix extensions of traces on $A$.
The hereditary algebra
\[
 B=\overline{b\mathrm M_l(A)b}
\]
is simple and infinite-dimensional.  Every nonzero hereditary
\Cstar{}-subalgebra of a simple non-elementary \Cstar{}-algebra contains arbitrarily
long finite families of nonzero, mutually orthogonal, Cuntz equivalent
positive elements, by the cut-down argument of
\cite[Section~2]{NCP-LAlg}.
Choose an integer $q>l/\eta$ and norm-one elements
$c_1,\ldots,c_q\in B_+$ with these properties.  Then
\[
 q\mathrm d_\tau(c_1)=\sum_{i=1}^q \mathrm d_\tau(c_i)
 \leq \mathrm d_\tau(b)\leq l,
\]
and hence $\mathrm d_\tau(c_1)<\eta$ for every $\tau\in\mathrm T(A)$.

Choose $\delta>0$ such that $(c_1-\delta)_+\ne0$, and choose
$f\in \mathrm C([0,1])_+$ which is zero on $[0,\delta/2]$ and one on
$[\delta,1]$.  With the unit taken in the multiplier algebra of $B$,
put
\[
 t=1-f(c_1),\qquad b_0=t\frac{b}{\|b\|}t,
 \qquad b_1=(c_1-\delta)_+.
\]
Then $b_0,b_1\in B_+$ and $b_0b_1=0$.  In the von Neumann algebra
generated by $B$, the support projection of $b/\|b\|$ is the unit of
$B^{**}$ and the support projection of $b_0$ is the support projection
of $t$.  Since $t+f(c_1)=1$, it follows that
\[
 \mathrm d_\tau(b_0)\geq \mathrm d_\tau(b)-\mathrm d_\tau(f(c_1))
 \geq \mathrm d_\tau(b)-\mathrm d_\tau(c_1)>\mathrm d_\tau(b)-\eta.
\]
Applying the same argument in
$\overline{b_1\mathrm M_l(A)b_1}$, and then taking a nonzero diagonal
compression, gives, for any prescribed $L$, an element $0\ne z\in A_+$
with $L[z]\leq[b_1]$.  This proves (1).

For (2), choose integers $p,q$ such that
\[
 t<\frac pq<s.
\]
The separation argument in the proof of
\cite[Proposition~3.2.1]{RC-5-person}, applied to $qu_0+pe$ and $qv$
using \eqref{eq:functional-gap}, gives integers $K,K_1$ and an element
$w\in\operatorname{Cu}(A)$, with $w\leq K_3v$ for some $K_3$, such that
\[
 Kqu_0+Kpe+w+v\leq Kqv+w
\]
and then
\[
 K_1Kqu_0+K_1Kpe+v\leq K_1Kqv.
\]
Regularization as in
\cite[Lemma~2.4.1]{RC-5-person} passes from additive order-preserving
maps to functionals; the compact containment $u_0\ll u$ ensures the
required strict inequalities.  The boundedness
of $u$ and $v$ by finite multiples of $e$ ensures that all the
normalizations occurring in that separation argument are finite.
Since $A$ is
simple and $v\ne0$, $v$ is full.  Choose $K_2$ with $u_0\leq K_2v$.
Multiplying the last display by $K_2$ and using this inequality gives
\[
 (K_2K_1Kq+1)u_0+K_2K_1Kp\,e
 \leq K_2K_1Kq\,v.
\]
Thus \eqref{eq:stable-gap} holds with
$N=K_2K_1Kq$ and $M=K_2K_1Kp$.  Moreover, $M/N=p/q>t$.

For (3), choose $\theta>0$ with $(z-\theta)_+\ne0$ and increasing
finite subsets whose union is dense in the unit ball of $A$.
Apply the relevant tracial approximation property to these subsets,
with tolerances tending to zero and comparison element
$(z-\theta)_+$, rescaled to be a contraction if necessary.
Let $C_j,h_j$ be the resulting approximation data.  The compression and
commutator conditions give
\[
 h\in A'\cap A_\omega,\qquad hAh\subseteq Q.
\]
The matrix-level assertion follows entry by entry.  Moreover,
\[
 1-h_j\precsim(z-\theta)_+,\qquad
 (z-\theta)_+g_\theta(z)=(z-\theta)_+.
\]
Lemma~\ref{comp-local-global}, applied with the fixed target cutoff
$g_\theta(z)$, gives $1-h\precsim g_\theta(z)\precsim z$.
Since $1-h^2\sim1-h$, this proves
\eqref{eq:central-ta-cut}.  Centrality of $h$ makes $\Phi$ order
zero.

Suppose first that $(k+1)x\leq ky$ and the model Cuntz semigroups
are almost unperforated.  Write
\[
 X=\operatorname{Cu}(\Phi)(x),\qquad
 Y=\operatorname{Cu}(\Phi)(y),
\]
and represent these classes by bounded positive sequences
$(x_j),(y_j)$ in common finite matrix algebras over $C_j$.
Denote their images in $Q$ by $\widetilde x,\widetilde y$.
Fix $\rho>0$.  Compact containment gives a fixed $\sigma>0$ with
\[
 (k+1)[(\widetilde x-\rho/2)_+]
 \leq k[(\widetilde y-\sigma)_+].
\]
Indeed, the left side is compactly contained in $(k+1)X\leq kY$,
and the classes of the spectral cutoffs of $\widetilde y$ increase
to $Y$.  Lifting a Cuntz witness to the product and applying the norm
perturbation lemma gives, on an $\omega$-large set of indices,
\[
 (k+1)[(x_j-\rho)_+]
 \leq k[(y_j-\sigma)_+].
\]
Almost unperforation therefore gives
$(x_j-\rho)_+\precsim(y_j-\sigma)_+$ there.
The exact support identity
\[
 (y_j-\sigma)_+g_\sigma(y_j)=(y_j-\sigma)_+
\]
permits the argument of Lemma~\ref{comp-local-global} to be applied
in $\prod_j\mathrm M_l(C_j)$.  Multiplying the Cuntz witnesses by
$(y_j-\sigma)_+^{1/2}$ bounds their norms by the square root of the
source norm.  Consequently
\[
 (\widetilde x-\rho)_+\precsim g_\sigma(\widetilde y)
 \precsim\widetilde y\quad\text{in }Q.
\]
Letting $\rho$ decrease to zero proves
\eqref{eq:local-almost-unperf}.

For the quantitative assertion, choose the approximation data from
$\mathcal S_\bullet$.  Put $E=[h^2]=[h]$ and $q=m/k$.
Applying $\operatorname{Cu}(\Phi)$ to the hypothesized relation
gives
\[
 (k+1)X+mE\leq kY.
\]
Choose $s$ with $r<s<\min\{q,L\}$, and fix $0<\delta<1/2$.
For each $\rho>0$, compact containment again supplies a fixed
$\sigma>0$ such that
\[
 (k+1)[(\widetilde x-\rho/2)_+]
 +m[(h^2-\delta/2)_+]
 \leq k[(\widetilde y-\sigma)_+].
\]
Lifting a Cuntz witness and applying a perturbation smaller than
$\min\{\rho/2,\delta/2\}$ gives, on an $\omega$-large set of indices,
\begin{equation}\label{eq:pointed-coordinate-relation}
 (k+1)[(x_j-\rho)_+]+m[(h_j^2-\delta)_+]
 \leq k[(y_j-\sigma)_+].
\end{equation}

Set
\[
 e_j=h_j^2,\qquad p_j=(e_j-\delta)_+,\qquad
 v_j=e_j(1-g_{2\delta}(e_j)),\qquad a=1-2\delta>0.
\]
The supports of $p_j$ and $v_j$ cover the support of $e_j$.
Functional calculus therefore gives
\[
 [e_j]\leq[p_j]+[v_j],\qquad
 v_jg_a(1-e_j)=v_j.
\]
For every $\lambda\in F(\operatorname{Cu}(C_j))$,
\eqref{eq:pointed-coordinate-relation} implies
\begin{align*}
 \lambda([(x_j-\rho)_+])+s\lambda([h_j])
 &\leq\lambda([(x_j-\rho)_+])
       +s\lambda([p_j])+s\lambda([v_j])\\
 &\leq\lambda([(y_j-\sigma)_+])+L\lambda([v_j]).
\end{align*}
These inequalities hold in $[0,\infty]$ and require no subtraction of
infinite values.
Under condition (i), since $s>r\geq r_{C_j,h_j}$, pointed
comparison yields
\begin{equation}\label{eq:buffered-coordinate-comparison}
 (x_j-\rho)_+
 \precsim (y_j-\sigma)_+\oplus v_j^{\oplus L}.
\end{equation}

Under condition (ii), the same functional inequality holds for
$\mathrm d_\tau$.  The source rank $\mathrm d_\tau((x_j-\rho)_+)$ is finite:
the cutoff $g_\rho(x_j)$ lies in the Pedersen ideal, every densely
defined lower semicontinuous trace is finite there, and
\[
 \mathrm d_\tau((x_j-\rho)_+)\leq\tau(g_\rho(x_j)).
\]
Fullness of $h_j$ gives $\mathrm d_\tau(h_j)>0$ for each nonzero trace.
If this value is finite, $s>r$ gives the strict inequality
\[
 \mathrm d_\tau((x_j-\rho)_+)+\gamma_{C_j}(\tau)
 <\mathrm d_\tau((y_j-\sigma)_+)+L \mathrm d_\tau(v_j).
\]
If $\mathrm d_\tau(h_j)=\infty$, the target rank is infinite, while the
source rank and the gap are finite; the same strict inequality
holds.  Gap comparison again proves
\eqref{eq:buffered-coordinate-comparison}.

Finally put
\[
 b_j=(y_j-\sigma)_+\oplus v_j^{\oplus L},\qquad
 c_j=g_\sigma(y_j)\oplus g_a(1-e_j)^{\oplus L}.
\]
Here $c_j$ is a positive contraction in a matrix algebra over the
ambient algebra $A$, and $b_jc_j=b_j$ exactly.
The cutoff $g_a(1-e_j)$ need not belong to the nonunital model $C_j$,
so the resulting comparison takes place in $A_\omega$.
Applying Lemma~\ref{comp-local-global} at matrix level gives
\begin{align*}
 [(\widetilde x-\rho)_+]
 &\leq[g_\sigma(\widetilde y)]+L[g_a(1-h^2)]\\
 &\leq Y+L[1-h^2]\qquad\text{in }\operatorname{Cu}(A_\omega).
\end{align*}
The fixed target cutoff $\sigma$ and low-spectrum cutoff $a$ ensure
that the coordinate witnesses form a bounded sequence.  Letting $\rho$
decrease to zero proves
\eqref{eq:local-radius}.
\end{proof}

\Needspace{8\baselineskip}
\begin{thm}\label{WAD=>SCP}
Let $A$ be a simple, separable, exact, unital, infinite-dimensional
\Cstar{}-algebra with $\mathrm T(A)\ne\emptyset$.
\begin{enumerate}
\item Suppose that $\mathcal S$ is a class of \Cstar{}-algebras such that
$\operatorname{Cu}(C)$ is almost unperforated for every
$C\in\mathcal S$.  If $A\in\mathrm{TA}_1(\mathcal S)$, then $A$ has
strict comparison of positive elements.

\item Let $r\geq0$, and let $\mathcal S_\bullet$ be a class of
pointed \Cstar{}-algebras.  For each $(C,e)\in\mathcal S_\bullet$, suppose
that at least one of the following conditions holds:
\begin{enumerate}
\item[(a)] $r_{C,e}\leq r$;
\item[(b)] there is $\gamma_{C,e}\in G_C$ such that
\begin{equation}\label{eq:scalar-model-gap}
 \gamma_{C,e}(\tau)\leq r \mathrm d_\tau(e),
 \qquad \tau\in\mathrm T^+(C).
\end{equation}
When $r=0$, this means $\gamma_{C,e}=0$.
\end{enumerate}
If $A\in\mathrm{TA}_1^\bullet(\mathcal S_\bullet)$, then
\[
 \bigl(\tau\mapsto r\tau(1_A)\bigr)\in G_A,
 \qquad \operatorname{rc}(A)\leq r.
\]
In particular, condition (b) holds if
$\gamma_{C,e}(\tau)\leq r\tau(e^2)$.

\item Let $H:\mathrm T^+(A)\to[0,\infty)$ be continuous and
affine, with $H(0)=0$, and suppose that
$A\in\mathrm{TA}_1(\mathcal S)$ for a class $\mathcal S$ of
\Cstar{}-algebras.  Suppose that, for every $0<\eta<1$, there
is $c_\eta\in A_+$ satisfying
\begin{equation}\label{eq:ambient-gap-budget}
 \tau(c_\eta)\leq H(\tau)+\eta\tau(1_A),
 \qquad \tau\in\mathrm T^+(A),
\end{equation}
and that the $\mathrm{TA}_1(\mathcal S)$ witnesses can always be
chosen with $h$ strictly positive in $C$, a positive element
$c'\in C$, and $\gamma_C\in G_C$ satisfying
\begin{equation}\label{eq:compatible-model-gap}
 \|c'-hc_\eta h\|<\eps,
 \qquad \gamma_C(\tau)\leq\tau(c')+\eta \mathrm d_\tau(h),
 \quad \tau\in\mathrm T^+(C).
\end{equation}
These data are required for every finite set, every tolerance
$\eps>0$, and every prescribed nonzero defect-comparison element.  Then $H\in G_A$.
\end{enumerate}
\end{thm}

\begin{proof}
Every trace on $A$ is faithful, and therefore $A$ is stably finite.
Exactness implies that every bounded $2$-quasitrace is a trace.  Thus
the normalized functionals on $\operatorname{Cu}(A)$ are precisely the
dimension functions coming from $\mathrm T(A)$.

The comparison arguments reduce to finite matrix algebras.
Indeed, fix $a,b\in
(A\otimes\mathcal K)_+$ with a strict gap and $\eps>0$.
Norm approximation in finite matrix corners and spectral cutoffs give
$u_1$ with $[(a-\eps)_+]\ll u_1\ll[a]$, represented by
$a_1$ in a finite matrix algebra.  Choose $\rho>0$ with
$[(a-\eps)_+]\leq[(a_1-\rho)_+]$.  The continuous function
$\tau\mapsto\tau(g_{\rho/2}(a_1))$ is at most $\mathrm d_\tau(a)$.
Compactness of $\mathrm T(A)$, increasing cut-downs of $b$, and
finite-corner norm approximation give $v_1\ll[b]$, represented by
a finite-matrix element $b'$, with
\[
 \tau(g_{\rho/2}(a_1))+H(\tau)<\mathrm d_\tau(b')
 \quad(\tau\in\mathrm T(A)),
\]
where $H$ is zero in (1), the scalar gap in (2), and the given affine
function in (3).
The finite-matrix argument applied to
$(a_1-\rho/2)_+$ and $b'$ compares $(a-\eps)_+$ with $b$.
Letting $\eps$ decrease to zero proves the stabilized assertion.

For (1), let $a,b\in\mathrm M_l(A)_+$ satisfy
\begin{equation}\label{eq:strict-trace-gap}
 \mathrm d_\tau(a)<\mathrm d_\tau(b),\qquad \tau\in\mathrm T(A),
\end{equation}
and fix $\eps>0$.  Put
\[
 x_0=[(a-\eps)_+],\qquad x_1=[(a-\eps/2)_+].
\]
Then $x_0\ll x_1$.  Since
\[
 \mathrm d_\tau((a-\eps/2)_+)\leq\tau(g_{\eps/2}(a))
 \leq \mathrm d_\tau(a),
\]
the lower semicontinuous function
\[
 \tau\longmapsto \mathrm d_\tau(b)-\tau(g_{\eps/2}(a))
\]
is strictly positive on the compact space $\mathrm T(A)$.  Let
$3\eta>0$ be no larger than its minimum.
Lemma~\ref{comparison-tools}(1), applied to $b$ and $\eta$, gives
orthogonal
$b_0,b_1\in\overline{b\mathrm M_l(A)b}_+$ such that $b_1\ne0$ and
\[
 \mathrm d_\tau(x_1)+2\eta<\mathrm d_\tau(b_0),
 \qquad \tau\in\mathrm T(A).
\]
Lemma~\ref{comparison-tools}(2), with $s=\eta$ and $t=0$, gives an integer
$k\geq1$ such that
\begin{equation}\label{eq:stable-comparison-strict}
 (k+1)x_0\leq k[b_0].
\end{equation}
Here the positive multiple of $[1_A]$ in \eqref{eq:stable-gap} has
been omitted.

Choose $0\ne z\in A_+$ so that
\begin{equation}\label{eq:matrix-defect-reserve}
 l[z]\leq[b_1],
\end{equation}
and apply Lemma~\ref{comparison-tools}(3) with this $z$.  Write
$e=h^2$.  From \eqref{eq:stable-comparison-strict} and
\eqref{eq:local-almost-unperf},
\[
 h(a-\eps)_+h\precsim hb_0h
\]
in the appropriate matrix algebra over $A_\omega$.  Since $e$ is central,
\begin{align*}
 (a-\eps)_+
  &=h(a-\eps)_+h
    +(1-e)^{1/2}(a-\eps)_+(1-e)^{1/2},
\end{align*}
and hence, in $\operatorname{Cu}(A_\omega)$,
\begin{align*}
 [(a-\eps)_+]
 &\leq [h(a-\eps)_+h]+l[1-e]\\
 &\leq [hb_0h]+l[z]\\
 &\leq [b_0]+[b_1]\leq[b].
\end{align*}
Lemma~\ref{constant-ultrapower-comparison}, applied at matrix level,
gives the same comparison in $A$.  Since $\eps>0$ was arbitrary,
$a\precsim b$.
This proves strict comparison.

For (2), let $a,b\in\mathrm M_l(A)_+$ satisfy
\[
 \mathrm d_\tau(a)+r<\mathrm d_\tau(b),\qquad \tau\in\mathrm T(A),
\]
and fix $\eps>0$.  As in (1), the continuous upper bound
$\tau(g_{\eps/2}(a))$ and compactness give a uniform positive
margin after replacing $a$ by $(a-\eps/2)_+$.  Using
Lemma~\ref{comparison-tools}(1), choose orthogonal $b_0,b_1$ in
$\overline{b\mathrm M_l(A)b}$, with $b_1\ne0$, and numbers $s,t$
such that
\[
 r<t<s,
 \qquad \mathrm d_\tau((a-\eps/2)_+)+s<\mathrm d_\tau(b_0)
 \quad(\tau\in\mathrm T(A)).
\]
Lemma~\ref{comparison-tools}(2) gives $N,M\in\mathbb N$ with
\begin{equation}\label{eq:radius-relation-reserved}
 (N+1)[(a-\eps)_+]+M[1_A]\leq N[b_0],
 \qquad M/N>t>r.
\end{equation}
Put $L=\lfloor r\rfloor+1$, and choose $0\ne z\in A_+$ with
\[
 (l+L)[z]\leq[b_1].
\]
Lemma~\ref{comparison-tools}(3), with pointed approximation data,
gives \eqref{eq:local-radius}, and hence
\[
 [h(a-\eps)_+h]\leq[hb_0h]+L[1-h^2]
 \quad\hbox{in }\operatorname{Cu}(A_\omega).
\]
Centrality of $h$ and the positive decomposition used in (1) give
\[
 [(a-\eps)_+]\leq[hb_0h]+(l+L)[1-h^2]
 \leq[b_0]+[b_1]\leq[b].
\]
The comparison holds in $A$ by
Lemma~\ref{constant-ultrapower-comparison}.  Letting $\eps$ decrease
to zero gives $a\precsim b$.
This proves the scalar gap assertion for finite matrices and hence,
by the preceding cut-down reduction, for $A\otimes\mathcal K$.
The definition of radius of comparison then gives
$\operatorname{rc}(A)\leq r$.

For (3), suppose that $a,b\in \mathrm M_l(A)_+$ satisfy
\[
 \mathrm d_\tau(a)+H(\tau)<\mathrm d_\tau(b),\qquad \tau\in\mathrm T(A).
\]
Fix $\eps>0$, and put $a_0=(a-\eps)_+$ and
$a_1=(a-\eps/2)_+$.  Compactness of $\mathrm T(A)$, continuity of
$H$ and $\tau\mapsto\tau(g_{\eps/2}(a))$, and lower
semicontinuity of $\tau\mapsto \mathrm d_\tau(b)$ give
$0<\mu<1/8$ such that
\[
 \mathrm d_\tau(b)-\tau(g_{\eps/2}(a))-H(\tau)>8\mu
 \quad(\tau\in\mathrm T(A)).
\]
By Lemma~\ref{comparison-tools}(1), choose orthogonal
$b_0,b_1\in\operatorname{Her}_{\mathrm M_l(A)}(b)_+$, with $b_1\ne0$,
such that $\mathrm d_\tau(b_0)>\mathrm d_\tau(b)-\mu$ for all normalized traces.
Choose $0<\eta<\mu/4$, and take the element $c=c_\eta$ from the
hypothesis.  Then
\begin{equation}\label{eq:affine-ambient-budget}
 \mathrm d_\tau(a_1)+\tau(c)+6\mu<\mathrm d_\tau(b_0),
 \qquad \tau\in\mathrm T(A).
\end{equation}

Choose an integer $R\geq\max\{1,\|c\|\}$, an integer $n$ with
$2/n<\mu/4$, and $0<\beta<1/n$.  For $1\leq i\leq nR$, put
\[
 c_i=(c-i/n-\beta)_+,
 \qquad f_i(t)=\min\{1,\max\{0,n(t-i/n-\beta)\}\}.
\]
Elementary scalar inequalities give
\begin{equation}\label{eq:affine-riemann-budget}
 c\leq\frac2n1_A+\frac1n\sum_{i=1}^{nR}f_i(c),
 \qquad f_i(c)\leq n c_i,
 \qquad
 \frac1n\sum_{i=1}^{nR}\mathrm d_\tau((c-i/n)_+)\leq\tau(c).
\end{equation}
In $\operatorname{Cu}(A)$, set
\[
 u=n[a_1]+\sum_i[(c-i/n)_+],\quad
 u_0=n[a_0]+\sum_i[c_i],\quad v=n[b_0].
\]
Then $u_0\ll u$, and \eqref{eq:affine-ambient-budget} and
\eqref{eq:affine-riemann-budget} imply
$\lambda(u)+6n\mu<\lambda(v)$ for every normalized functional.
Apply Lemma~\ref{comparison-tools}(2), with $t=2n\mu$.
For the resulting integers $K,P$,
$(K+1)u_0+P[1_A]\leq Kv$ and $P/K>2n\mu$.
Set $N=Kn$ and $M=P$.  Omitting the additional positive terms gives
\begin{equation}\label{eq:affine-stable-budget}
 (N+1)[a_0]+\frac Nn\sum_i[c_i]+M[1_A]\leq N[b_0],
 \qquad q:=M/N>2\mu.
\end{equation}

Choose $0\ne z\in A_+$ with $(l+1)[z]\leq[b_1]$, and choose
$\theta>0$ with $z_\theta=(z-\theta)_+\ne0$.
Apply the compatibility hypothesis with this fixed $\eta$, comparison
element $z_\theta$, and tolerances tending to zero, to increasing
finite sets with dense union containing the matrix entries of $a_0,b_0$
and the elements $c,c_i,f_i(c)$.  This gives models
$C_j$, strictly positive contractions $h_j\in C_j$,
elements $c'_j\in(C_j)_+$, and finite-valued gap functions $\gamma_j$
satisfying the compatibility inequality.  Put $e_j=h_j^2$ and
\[
 h=[(h_j)],\quad e=h^2,\quad
 Q=\pi\Big(\prod_jC_j\Big)\subseteq A_\omega,
 \qquad \Phi(x)=hxh.
\]
Then $h$ commutes with $A$, and $\Phi:A\to Q$ is c.p.c.
order zero.  Moreover, $1-e_j\sim1-h_j\precsim z_\theta$.
Since $z_\theta g_\theta(z)=z_\theta$,
Lemma~\ref{comp-local-global} gives
\begin{equation}\label{eq:affine-defect-reserve}
 1-e\precsim g_\theta(z)\precsim z,
 \qquad (l+1)[1-e]\leq[b_1].
\end{equation}

Choose positive representatives $x_j,y_j\in \mathrm M_l(C_j)$ and
$U_{i,j}\in C_j$ for
$\Phi(a_0),\Phi(b_0),\Phi(c_i)$, respectively, and positive
contractions $F_{i,j}\in C_j$ for $\Phi(f_i(c))$.
They may be chosen to approximate their corresponding compressed
elements in norm, while $c'_j$ approximates $h_jch_j$.
Consequently \eqref{eq:affine-riemann-budget} gives numbers
$\xi_j\to0$ such that, in the unitization of $C_j$,
\[
 F_{i,j}\leq nU_{i,j}+\xi_j1,
 \qquad
 c'_j\leq\frac2n e_j+\frac1n\sum_iF_{i,j}+\xi_j1.
\]
Fix $\rho>0$ with $Rn\rho<\mu/4$.
For a bounded trace $\tau$ on $C_j$, spectral splitting of
$U_{i,j}$ at $\rho$ gives
\[
 \tau(F_{i,j})\leq \mathrm d_\tau((U_{i,j}-\rho)_+)
                    +(n\rho+\xi_j)\|\tau\|.
\]
Strict positivity of $h_j$ implies
$\|\tau\|=\mathrm d_\tau(e_j)$.  Hence, for all sufficiently large $j$,
the compatibility bound on $\gamma_j$ gives
\begin{equation}\label{eq:affine-local-gap-budget}
 \gamma_j(\tau)
 \leq\frac1n\sum_i \mathrm d_\tau((U_{i,j}-\rho)_+)
                      +\mu \mathrm d_\tau(e_j)
 \quad\text{for every bounded trace }\tau.
\end{equation}
Indeed, the coefficient of $\mathrm d_\tau(e_j)$ is at most
$\eta+2/n+Rn\rho+(R+1)\xi_j<\mu$.

Fix $0<\delta<1/4$ and $\chi>0$.
Apply $\operatorname{Cu}(\Phi)$ to
\eqref{eq:affine-stable-budget}, and compactly cut down all
terms on the left.  A fixed target cut and the usual lifting of
a Cuntz witness give $\sigma>0$ and an $\omega$-large set on
which
\begin{equation}\label{eq:affine-coordinate-budget}
 (N+1)[(x_j-\chi)_+]
 +\frac Nn\sum_i[(U_{i,j}-\rho)_+]
 +M[(e_j-\delta)_+]
 \leq N[(y_j-\sigma)_+].
\end{equation}
To obtain this relation, apply the cutoffs $\chi/2,\rho/2,\delta/2$
in $Q$ and choose a compactly contained target class.  Lifting with
norm errors smaller than the remaining spectral margins gives a
cutoff $\sigma$ independent of $j$.

Put $\ell_j=e_j(1-g_{2\delta}(e_j))$.
The supports of $(e_j-\delta)_+$ and $\ell_j$ cover the support of
$e_j$, so
\[
 \mathrm d_\tau(e_j)\leq \mathrm d_\tau((e_j-\delta)_+)+\mathrm d_\tau(\ell_j).
\]
For every nonzero bounded trace, combining
\eqref{eq:affine-local-gap-budget} and
\eqref{eq:affine-coordinate-budget} gives
\begin{equation}\label{eq:affine-local-comparison-premise}
 \mathrm d_\tau((x_j-\chi)_+)+\gamma_j(\tau)
 <\mathrm d_\tau((y_j-\sigma)_+)+\mathrm d_\tau(\ell_j).
\end{equation}
The strictness follows from $q>2\mu$, $\mu<1$, and
$\mathrm d_\tau(e_j)=\|\tau\|>0$.
The same strict inequality holds for every unbounded densely finite
lower semicontinuous trace: fixed spectral cuts have finite
dimension, while strict positivity gives $\mathrm d_\tau(e_j)=\infty$;
therefore $\mathrm d_\tau(\ell_j)=\infty$.  The left side is finite
because $\gamma_j$ is finite.  The local comparison-gap property
therefore yields
\[
 (x_j-\chi)_+\precsim (y_j-\sigma)_+\oplus\ell_j.
\]

The target in this comparison is fixed by multiplication with
\[
 g_\sigma(y_j)\oplus g_{1-2\delta}(1-e_j):
\]
the first summand is fixed by the cutoff at $\sigma$, and the second
by the support condition $e_j<2\delta$ for $\ell_j$.
Lemma~\ref{comp-local-global} thus gives
\[
 (\Phi(a_0)-\chi)_+
 \precsim g_\sigma(\Phi(b_0))\oplus
                   g_{1-2\delta}(1-e)
 \precsim\Phi(b_0)\oplus(1-e).
\]
Letting $\chi\downarrow0$ proves
$\Phi(a_0)\precsim\Phi(b_0)\oplus(1-e)$.
Finally, centrality of $h$ and \eqref{eq:affine-defect-reserve}
give
\[
 [a_0]\leq[\Phi(a_0)]+l[1-e]
 \leq[\Phi(b_0)]+(l+1)[1-e]
 \leq[b_0]+[b_1]\leq[b]
 \quad\text{in }\operatorname{Cu}(A_\omega).
\]
Comparison of constant elements descends to $A$.
Since $\eps>0$ was arbitrary, $a\precsim b$.
The stabilized reduction and homogeneity extend the comparison to all
positive elements of $A\otimes\mathcal K$ and all nonzero traces on $A$.
Thus $H\in G_A$.

\end{proof}

\begin{cor}\label{Z=>SCP}
Let \(A\) be a simple, separable, unital, infinite-dimensional
\Cstar{}-algebra with strict comparison and almost divisible Cuntz
semigroup.  Let \(\Gamma\) be a
countable discrete amenable group acting on \(A\).  Suppose that the
action has WRP and that \(B=A\rtimes_{\alpha,r}\Gamma\) is simple and
exact.  Then \(B\) has strict comparison of positive elements.
\end{cor}

\begin{proof}
Amenability gives $\mathcal M_1(A,\Gamma)\ne\varnothing$.
Proposition~\ref{pre-COM} gives (COS), and
Proposition~\ref{WRP-COS-TA1} gives \(B\in\mathrm{TA}_1(\mathcal S)\).
Exactness passes from \(B\) to its \Cstar{}-subalgebra \(A\).  Since WRP identifies
\(\mathrm T(B)\) with the nonempty compact set
\(\mathcal M_1(A,\Gamma)\), all hypotheses on the ambient algebra in
Theorem~\ref{WAD=>SCP}(1) are satisfied.  Strict comparison for the
simple exact algebra \(A\) implies that \(\operatorname{Cu}(A)\) is
almost unperforated.  Almost unperforation passes to matrix algebras,
hereditary \Cstar{}-subalgebras, and finite direct sums.  Thus every member
of \(\mathcal S\) has almost unperforated Cuntz semigroup, and
Theorem~\ref{WAD=>SCP}(1) applies to \(B\).
\end{proof}

\begin{rem}\label{gap-transfer-scope}
The strict inequality in \eqref{eq:gap-function} is essential: strict
comparison does not imply comparison from equality of rank functions.
The scalar assertion includes the endpoint $r$ because a source cutoff
produces a uniform positive margin.

The model gaps in (2) are bounded relative to the compression
scale; they need not be constant on a normalized trace space.
The affine bound $\gamma_{C,e}(\tau)\leq r\tau(e^2)$ is a sufficient
condition.  In (3), the gap functions are controlled by compressed
positive representatives of the ambient affine function.  This
compatibility does not follow from the existence or continuity of gap
functions on individual models, since model traces need not extend to
traces on $A$.

\label{gap-budget-realization}
The positive representatives required in
\eqref{eq:ambient-gap-budget} exist for every $H$ in (3).
Indeed, \cite[Theorem~4.3]{Dadarlat-Toms-Ranks} gives
$c_\eta\in A_+$ with
\[
 \tau(c_\eta)=H(\tau)+\eta\tau(1_A),
 \qquad \tau\in\mathrm T^+(A).
\]
Thus the hypothesis in (3) concerns approximation data satisfying
\eqref{eq:compatible-model-gap} for such representatives.  For example, if $c\in A_+$ and the witnesses have
$\gamma_C(\tau)\leq\tau(c')$ with $c'\approx hch$, then
$\tau\mapsto\tau(c)$ is a gap function of $A$.

Strict positivity of $h$ localizes the model to $\overline{hCh}$.
The proof also applies to these hereditary subalgebras when they belong
to the approximation class and admit the required gap functions.
A gap function on $C$ alone does not ensure these hereditary comparison
properties.
Throughout, model gap functions are finite-valued on the whole
cone of densely finite traces.  Allowing infinite gap values would
require a different comparison convention, since $\infty<\infty$
is false.
\end{rem}

\section{Property (COS) for \Cstar{}-dynamical systems}\label{cos-section}

The following condition extends comparison of open sets to
\Cstar{}-dynamical systems.
\begin{defn}\label{COS-definition}
Let $A$ be unital, let $\alpha:\Gamma\to\operatorname{Aut}(A)$ be
an action of a discrete group, and suppose that
$\mathcal M_1(A,\Gamma)\ne\varnothing$.  The system
$(A,\Gamma)$ has property (COS) if there are $\lambda\in(0,1]$ and
$m\in\mathbb N$ such that, whenever $a,b\in A_+$ satisfy
\[
 \mathrm d_\rho(a)<\lambda \mathrm d_\rho(b),\qquad
 \rho\in\mathcal M_1(A,\Gamma),
\]
one has
\[
 a\precsim_{A\rtimes_r\Gamma} b^{\oplus m}.
\]
We specify the constants by saying that the action has
$(\lambda,m)$-COS.
\end{defn}

\begin{prop}\label{WRP-COS-TA1}
Let \(A\) be a simple, separable, unital, infinite-dimensional \Cstar{}-algebra,
let \(\Gamma\) be a countable discrete amenable group acting on \(A\),
and put \(B=A\rtimes_{\alpha,r}\Gamma\).  Suppose that the action has
WRP and (COS), and that \(B\) is simple.  If \(\mathcal S\) is the class
from Theorem~\ref{WTA-C-dyn-1}, then
\[
 B\in\mathrm{TA}_1(\mathcal S).
\]
\end{prop}

\begin{proof}
Let \(\lambda\) and \(m\) be constants for (COS), and fix
\(0\ne x\in B_+\).  Since \(B\) is simple and infinite-dimensional,
\(\overline{xBx}\) contains nonzero pairwise orthogonal positive
elements \(x_1,\ldots,x_m\).  Simplicity of \(A\) makes the action
minimal.  Lemmas~\ref{WRP-Kishimoto} and \ref{coefficient-cutdown}
therefore give nonzero elements \(b_j\in A_+\) with
\(b_j\precsim_B x_j\).

A standard hereditary-subalgebra argument, applied finitely many times,
gives \(0\ne b\in A_+\) such that \(b\precsim_A b_j\) for every \(j\).
It follows that
\begin{equation}\label{common-coefficient-cutdown}
 b^{\oplus m}\precsim_B x_1\oplus\cdots\oplus x_m\precsim_B x.
\end{equation}
Every invariant trace on the simple algebra \(A\) is faithful.  Since
\(\mathcal M_1(A,\Gamma)\) is compact and
\(\rho\mapsto \mathrm d_\rho(b)\) is lower semicontinuous,
\[
 \kappa:=\min_{\rho\in\mathcal M_1(A,\Gamma)}\mathrm d_\rho(b)>0.
\]

Apply Theorem~\ref{WTA-C-dyn-1} with the prescribed finite set and with
the defect tolerance smaller than \(\lambda\kappa\).  It gives
\(C\in\mathcal S\) and a norm-one positive contraction \(h\in A\cap C\)
with the required commutator and compression estimates.  If
\(\rho\in\mathcal M_1(A,\Gamma)\), then
\(\rho\circ\mathbb E\) is a trace on \(B\); hence
\[
 \mathrm d_\rho(1-h)<\lambda\kappa\leq\lambda \mathrm d_\rho(b).
\]
Property (COS) and \eqref{common-coefficient-cutdown} give
\[
 1-h\precsim_B b^{\oplus m}\precsim_B x.
\]
The tracial-approximation version of Lemma~\ref{alt-defn-was}
therefore gives \(B\in\mathrm{TA}_1(\mathcal S)\).
\end{proof}

\begin{rem}
Every minimal free system $(\mathrm{C}(X), \Int^d)$ has URP and COS
by \cite{Niu-MD-Zd}.
\end{rem}


%
%


For a Cuntz semigroup $S$ with a $\Gamma$-action, write
\[
 S_\Gamma=\gamma(S/\Gamma).
\]
Here $S/\Gamma$ is the dynamical quotient in the category of
precompleted Cuntz semigroups, and $\gamma$ denotes Cu-completion.
By \cite[Corollary~9.8(ii)]{BPWZ-dynamical-Cu}, the canonical map
$q:S\to S_\Gamma$ is the universal invariant Cu-morphism: every
invariant Cu-morphism from $S$ to a Cuntz semigroup factors uniquely
through $q$.  Thus $S_\Gamma$ is the colimit in the category
$\mathrm{Cu}$ of the action diagram $B\Gamma\to\mathrm{Cu}$,
where $B\Gamma$ is the one-object category with morphisms $\Gamma$.
For a \Cstar{}-algebra action, this is the complete dynamical Cuntz
semigroup $\mathrm{Cu}_\Gamma(A)$ when $S=\mathrm{Cu}(A)$
\cite[Corollary~9.10]{BPWZ-dynamical-Cu}.

\begin{lem}\label{pure-dynamical-Cu}
Let $S$ be a pure Cuntz semigroup and let a discrete group $\Gamma$
act on $S$.  Then $S_\Gamma$ is pure.
\end{lem}

\begin{proof}
Put $\mathsf Z=\mathrm{Cu}(\mathcal Z)$.  By
\cite[Theorems~7.1.13 and~7.3.11]{APT-Cu-tensor}, pureness of $S$
is equivalent to the canonical map
\[
 \eta_S:S\longrightarrow\mathsf Z\otimes_{\mathrm{Cu}}S,
 \qquad x\longmapsto 1_{\mathsf Z}\otimes x,
\]
being an isomorphism.  These maps are natural in $S$, so $\eta_S$
intertwines the given action with $\mathrm{id}_{\mathsf Z}\otimes\alpha$.
The category $\mathrm{Cu}$ is closed symmetric monoidal
\cite[Theorems~5.10 and~5.11]{APT-bivariant-Cu}; hence
$\mathsf Z\otimes_{\mathrm{Cu}}-$ is a left adjoint and preserves
colimits.  Applying it to the action diagram gives
\[
 \mathsf Z\otimes_{\mathrm{Cu}}S_\Gamma
 \cong\operatorname*{colim}_{B\Gamma}
             (\mathsf Z\otimes_{\mathrm{Cu}}S)
 \cong\operatorname*{colim}_{B\Gamma}S=S_\Gamma.
\]
Under these identifications the colimit of the maps $\eta_S$ is
$\eta_{S_\Gamma}$.  This canonical map is therefore an isomorphism,
and the same characterization proves that $S_\Gamma$ is pure.
\end{proof}

The amenable-group case of the following automatic comparison
criterion was proved by Orovitz, Phillips, and Wang
\cite[Proposition~5.4]{OPW-Strict-comparison-crossed-products}.
The dynamical Cuntz semigroup permits comparison for arbitrary
discrete groups.

\begin{prop}\label{pre-COM}
Let $A$ be a simple unital \Cstar{}-algebra with strict comparison by
tracial states and with almost divisible Cuntz semigroup.
Let $\alpha:\Gamma\to\operatorname{Aut}(A)$ be an action of a
discrete group, and suppose that
$\mathcal M_1(A,\Gamma)\ne\varnothing$.
For every $n\in\mathbb N$ and $a,b\in \mathrm M_n(A)_+$,
\[
 \mathrm d_\rho(a)<\mathrm d_\rho(b),\qquad \rho\in\mathcal M_1(A,\Gamma),
 \quad\Longrightarrow\quad
 a\precsim_{A\rtimes_{\alpha,r}\Gamma}b.
\]
In particular, $(A,\Gamma)$ has $(1,1)$-COS.
\end{prop}

\begin{proof}
Strict comparison by traces makes $S=\mathrm{Cu}(A)$ almost
unperforated.  Indeed, the assertion first holds in $\mathrm W(A)$:
if $(k+1)[c]\leq k[d]$ and $d\ne0$, faithfulness of every tracial
state gives
\[
 \mathrm d_\tau(c)\leq\frac{k}{k+1}\mathrm d_\tau(d)<\mathrm d_\tau(d),
 \qquad\tau\in\mathrm T(A),
\]
and strict comparison applies.  If $d=0$, the assertion is immediate.
For $x,y\in S$ with $(k+1)x\leq ky$, fix $x'\ll x$ and choose a
finite-matrix class $z$ with $x'\leq z\ll x$.  Write
$y=\sup_j w_j$ for an increasing sequence of finite-matrix classes.
Then $(k+1)z\ll ky$, so $(k+1)z\leq kw_j$ for some $j$.
The finite-matrix case gives $x'\leq z\leq w_j\leq y$.
Thus $x\leq y$, proving almost unperforation of $S$.
Together with almost divisibility, this
makes $S$ pure.  Lemma~\ref{pure-dynamical-Cu} shows that
$T=S_\Gamma$ is pure as well.

Write $q:S\to T$ for the canonical map.  A functional on a Cuntz
semigroup is an additive, order-preserving map to $[0,\infty]$
which preserves zero and increasing suprema.  By
\cite[Proposition~11.2 and Paragraph~11.6]{BPWZ-dynamical-Cu},
pullback by $q$ identifies functionals on $T$ with invariant
functionals on $S$.  A functional normalized at $q([1_A])$
pulls back to the rank function of an invariant normalized bounded
$2$-quasitrace on $A$.  Strict comparison by tracial states makes
this quasitrace a trace
\cite[Theorem~3.6 and Remark~3.7]{Ng-Robert-pure-commutators}.
Thus functionals normalized at $q([1_A])$ correspond exactly to
$\mathcal M_1(A,\Gamma)$.

Fix $a,b\in \mathrm M_n(A)_+$ with the stated rank inequalities.
Nonemptiness of $\mathcal M_1(A,\Gamma)$ implies that $b\ne0$.
Since $A$ is simple and unital, $b$ is full.  Compact containment
in $[1_A]\ll[1_A]\leq\sup_j j[b]$ gives an integer $L$ with
$[1_A]\leq L[b]$.  Put
\[
 x=q([a]),\qquad y=q([b]),\qquad u=q([1_A]).
\]
Then $x,y\leq nu$ and $u\leq Ly$, so $x\leq nLy$.
If $\lambda$ is a functional on $T$ with $\lambda(y)=1$, then
\[
 \frac1n\leq\lambda(u)\leq L.
\]
Consequently $\lambda\circ q/\lambda(u)=\mathrm d_\rho$ for some
$\rho\in\mathcal M_1(A,\Gamma)$, and therefore
\[
 \lambda(x)=\lambda(u)\mathrm d_\rho(a)
   <\lambda(u)\mathrm d_\rho(b)=\lambda(y).
\]
Almost unperforation of $T$ now gives $x\leq y$, by
\cite[Proposition~5.2.20(3)]{APT-Cu-tensor}: $x$ belongs to the
ideal generated by $y$, and the required strict inequality holds
for every functional normalized at $y$.

Finally, the coefficient inclusion
$\iota:A\to B=A\rtimes_{\alpha,r}\Gamma$ induces an invariant
Cu-morphism, since the canonical unitaries implement $\alpha$.
The universal property gives a Cu-morphism
$\kappa:T\to\mathrm{Cu}(B)$ such that
$\mathrm{Cu}(\iota)=\kappa\circ q$.
Applying $\kappa$ to $x\leq y$ proves $a\precsim_B b$.
\end{proof}

\begin{rem}\label{COS-stabilized-divisibility}
Almost divisibility in Proposition~\ref{pre-COM} is understood in
the stabilized sense: for $x'\ll x$ and $k\geq1$, there is
$y\in\mathrm{Cu}(A)$ with $ky\leq x$ and $x'\leq(k+1)y$.
The stronger condition $ky\leq x\leq(k+1)y$ is sufficient.
The proof does not require a positive representative of $y$ inside
$A$ itself, nor orthogonal orbit translates of such representatives.
\end{rem}

\begin{rem}\label{COS-invariant-traces}
For amenable $\Gamma$, the condition $\mathrm T(A)\ne\varnothing$
ensures $\mathcal M_1(A,\Gamma)\ne\varnothing$ by averaging.
For arbitrary groups, nonemptiness of the invariant tracial state
space is an explicit hypothesis.  Neither separability of $A$ nor
countability of $\Gamma$ is needed in Proposition~\ref{pre-COM}.
\end{rem}

\begin{rem}\label{COS-relative-comparison-scope}
Proposition~\ref{pre-COM} compares coefficient elements at every
finite matrix level.  It requires neither freeness of the action
nor simplicity of the crossed product.  It does not assert strict
comparison or pureness of the entire crossed product.  The examples
in Section~\ref{cos-examples-section} also include nonsimple
coefficient algebras; their comparison follows from separate tower
or crossed-product arguments.
\end{rem}

\begin{cor}\label{Z=>COS}
Let $A$ be a simple, unital, stably finite, exact,
$\mathcal Z$-stable \Cstar{}-algebra.  Every action of a discrete group
$\Gamma$ on $A$ with $\mathcal M_1(A,\Gamma)\ne\varnothing$
has $(1,1)$-COS, including comparison at every finite matrix level.
\end{cor}

\begin{proof}
$\mathcal Z$-stability gives $\mathrm{Cu}(A)$ a
$\mathrm{Cu}(\mathcal Z)$-multiplication, so it is almost
unperforated and almost divisible
\cite[Proposition~7.1.5(2) and Theorem~7.3.11]{APT-Cu-tensor}.
Exactness ensures that $2$-quasitraces are traces, so almost
unperforation gives strict comparison by tracial states.
Proposition~\ref{pre-COM} applies.
\end{proof}

\section{Orbit radius of comparison}\label{orbit-radius-section}

The orbit radius of comparison measures comparison of finite orbit
matrices over the coefficient algebra.  It replaces the dimension--rank
ratio of topological tower models when the coefficient algebra is
noncommutative.  The comparison condition is imposed on prescribed
Fourier polynomials before the Rokhlin bases are chosen, independently
of WRP and of comparison in the crossed product.

Throughout this section, $A$ is unital, $\mathrm T(A)\ne\emptyset$,
$\alpha:\Gamma\to\operatorname{Aut}(A)$ is an action of a countable
discrete amenable group, and $B=A\rtimes_{\alpha,r}\Gamma$.
We use the left F\o lner convention of Section~\ref{section-preliminaries}.
Matrix traces are unnormalized unless a normalization is displayed.
Every finite invariance set $K$ is understood to contain $e$.

\begin{ntn}[Orbit matrices and scalar collars]\label{orbit-matrix-notation}
For a nonempty finite set $F\subseteq\Gamma$, regular compression gives
a unital completely positive map
\[
 \Pi_F:B\longrightarrow \mathrm M_F(A),\qquad
 \Pi_F(x)=\big[\mathbb E(u_g^*xu_h)\big]_{g,h\in F}.
\]
We use the same notation for its matrix amplifications.  If
$x=\sum_kx_ku_k\in \mathrm M_l(B)$ is a Fourier polynomial, then
\begin{equation}\label{orbit-matrix-formula}
 \Pi_F(x)=\big[\alpha_{g^{-1}}(x_{gh^{-1}})\big]_{g,h\in F}
 \in \mathrm M_F(\mathrm M_l(A)).
\end{equation}
In particular, positivity of $x$ implies positivity of $\Pi_F(x)$.

A \emph{scalar collar} on $F$ is a diagonal matrix
$\Lambda=\operatorname{diag}(\lambda_g)_{g\in F}$ with
$0\leq\lambda_g\leq1$.  It acts identically on the external
$\mathrm M_l$ factor.  Put
\[
 R_{F,\Lambda}(x)=\Lambda\Pi_F(x)\Lambda,
 \qquad D(\Lambda)=\{g\in F:\lambda_g\ne1\}.
\]
Thus $\Lambda$ is independent of the coefficient algebra, and, when
$A=\mathrm C(X)$, independent of the point of $X$.
For $z\in \mathrm M_F(\mathrm M_l(A))_+$ and $\tau\in\mathrm T(A)$, write
\begin{equation}\label{normalized-orbit-rank}
 \widehat{\mathrm d}_{\tau,F}(z)
 =\frac1{|F|}\mathrm d_{\operatorname{Tr}_F\otimes\operatorname{Tr}_l\otimes\tau}(z).
\end{equation}
Only the orbit factor is normalized.
\end{ntn}

\begin{lem}[Finite-section trace formula]\label{orbit-moment-lemma}
Let $a\in \mathrm M_l(B)_+$, and let $f$ be continuous on an interval
containing $[0,\|a\|]$.  Given $\zeta>0$, there are a finite
$K\subseteq\Gamma$ and $\eta>0$ such that, for every nonempty
$(K,\eta)$-invariant $F$ and every scalar collar with
$|D(\Lambda)|<\eta|F|$,
\begin{equation}\label{orbit-moment-estimate}
 \left|
 \frac{(\operatorname{Tr}_F\otimes\operatorname{Tr}_l\otimes\tau)
          (f(R_{F,\Lambda}(a)))}{|F|}
 -(\operatorname{Tr}_l\otimes\bar\tau_F)
          (\mathbb E^{(l)}(f(a)))
 \right|<\zeta
\end{equation}
for every $\tau\in\mathrm T(A)$, where
$\bar\tau_F=|F|^{-1}\sum_{g\in F}\tau\circ\alpha_{g^{-1}}$.
\end{lem}

\begin{proof}
First let $a$ be a Fourier polynomial with support in a finite
symmetric set $N$ containing the identity.  If $k\geq1$ and
$g\in\operatorname{Int}_{N^k}(F)$, all the intermediate indices in
the $g$th diagonal entry of $\Pi_F(a)^k$ stay in $F$.
Consequently that entry equals
$\alpha_{g^{-1}}(\mathbb E^{(l)}(a^k))$.
The remaining diagonal entries in both expressions have norm at most
$\|a\|^k$.  It follows that
\[
 \left|
 \frac{(\operatorname{Tr}_F\otimes\operatorname{Tr}_l\otimes\tau)
          (\Pi_F(a)^k)}{|F|}
 -(\operatorname{Tr}_l\otimes\bar\tau_F)(\mathbb E^{(l)}(a^k))
 \right|
 \leq 2l\|a\|^k
       \frac{|F\setminus\operatorname{Int}_{N^k}(F)|}{|F|}.
\]
The equality for the constant polynomial is immediate.

To insert the collar, put $P=\Pi_F(a)$ and $D=1-\Lambda$.
Then $P-\Lambda P\Lambda=DP+\Lambda PD$.
Telescoping the $k$th powers and using traciality gives
\[
 \frac1{|F|}\left|
 (\operatorname{Tr}_F\otimes\operatorname{Tr}_l\otimes\tau)
       (P^k-(\Lambda P\Lambda)^k)\right|
 \leq 2kl\|a\|^k\frac{|D(\Lambda)|}{|F|}.
\]
Here $0\leq D\leq1_{D(\Lambda)}$, and the bound
$|\sigma(Dz)|\leq\|z\|\sigma(D)$ is applied to each term.
Polynomial approximation of $f$ proves
\eqref{orbit-moment-estimate} when $a$ is a positive Fourier polynomial.
The normalized trace in the first term has norm $l$, independent of $F$.

For arbitrary positive $a$, approximate $a^{1/2}$ by Fourier
polynomials and take their squares, rescaling if necessary so their
norms do not exceed $\|a\|$.  Contractivity of regular compression and
norm continuity of continuous functional calculus on a bounded interval,
obtained by polynomial approximation, control both sides of \eqref{orbit-moment-estimate} uniformly in $F$ and
$\tau$.  This proves the assertion.
\end{proof}

\begin{defn}[Orbit radius of comparison]\label{orbit-radius-definition}
Let $r\geq0$.  The system $(A,\alpha)$ has \emph{orbit
$r$-comparison} if, for every $l\in\mathbb N$, every pair of
positive Fourier polynomials $x,y\in \mathrm M_l(B)$, every
$\eps,\delta>0$, and every $t>r$, there are a finite
$K\subseteq\Gamma$ and $\eta\in(0,1)$ with the following property.
For every nonempty $(K,\eta)$-invariant $F$ and every scalar collar
$\Lambda$ with $|D(\Lambda)|<\eta|F|$, the inequalities
\begin{equation}\label{orbit-local-rank-gap}
 \widehat{\mathrm d}_{\tau,F}((R_{F,\Lambda}(x)-\eps)_+)+t
 <\widehat{\mathrm d}_{\tau,F}((R_{F,\Lambda}(y)-\delta)_+),
 \qquad \tau\in\mathrm T(A),
\end{equation}
imply
\begin{equation}\label{orbit-local-comparison}
 (R_{F,\Lambda}(x)-2\eps)_+
 \precsim_{\mathrm M_\infty(A)}
 (R_{F,\Lambda}(y)-\delta/2)_+.
\end{equation}
The \emph{orbit radius of comparison} is
\[
 \operatorname{orc}(A,\alpha)
 =\inf\{r\geq0:(A,\alpha)\text{ has orbit }r\text{-comparison}\},
\]
with the infimum of the empty set equal to $\infty$.
\end{defn}

\begin{rem}\label{orbit-radius-conventions}
This definition uses traces and is available without exactness.
Whenever the ordinary radius of comparison is used below, exactness
ensures that bounded $2$-quasitraces are traces.
The set of admissible radii is upward closed.  Quantification over
finite invariance sets and tolerances makes the definition independent
of a selected F\o lner sequence.  An equivariant isomorphism transports
the orbit matrices, traces, and comparison witnesses, so the orbit
radius is an invariant of the \Cstar{}-dynamical system.

Only matrices arising from fixed Fourier inputs are tested; comparison
of all positive elements of the whole algebra $\mathrm M_F(A)$ is not required.
The spectral cutoffs in \eqref{orbit-local-comparison} give bounded
witnesses before hereditary localization.  The normalization therefore
does not depend on the possibly noncompact scale of a tower base.
This notion differs from the dynamical radius based on coefficient
elements and dynamical Cuntz comparison in
\cite[Definition~3.4]{Asadi-Phillips-dynamical-radius}.
\end{rem}

\begin{lem}[From invariant traces to orbit rank gaps]\label{orbit-trace-gap}
Let $a,b\in \mathrm M_l(B)_+$ and $s\geq0$ satisfy
\[
 \mathrm d_{\rho\circ\mathbb E}(a)+s<\mathrm d_{\rho\circ\mathbb E}(b),
 \qquad \rho\in\mathcal M_1(A,\Gamma).
\]
For every $\eps>0$ there are $\delta,\kappa,\nu>0$, a finite
$K\subseteq\Gamma$, and $\eta>0$ such that
\begin{equation}\label{orbit-uniform-rank-gap}
 \widehat{\mathrm d}_{\tau,F}((R_{F,\Lambda}(x)-\eps)_+)+s+\kappa
 <\widehat{\mathrm d}_{\tau,F}((R_{F,\Lambda}(y)-\delta)_+)
\end{equation}
for every $\tau\in\mathrm T(A)$, every nonempty
$(K,\eta)$-invariant $F$, every scalar collar with
$|D(\Lambda)|<\eta|F|$, and every positive $x,y\in \mathrm M_l(B)$
with $\|x-a\|,\|y-b\|<\nu$.
\end{lem}

\begin{proof}
The function $\rho\mapsto(\rho\circ\mathbb E)(g_\eps(a))$ is
continuous and no larger than $\mathrm d_{\rho\circ\mathbb E}(a)$.
The lower semicontinuous difference
\[
 \mathrm d_{\rho\circ\mathbb E}(b)
 -(\rho\circ\mathbb E)(g_\eps(a))-s
\]
has a positive minimum on the compact invariant trace space.
Choose $\kappa>0$ sufficiently small.  Since
$g_{2\delta}(b)$ increases to the support of $b$ as
$\delta\downarrow0$, compactness of the invariant trace space,
applied to the increasing open sets defined by the strict inequalities,
gives a fixed $\delta>0$ such that
\begin{equation}\label{orbit-continuous-trace-gap}
 (\rho\circ\mathbb E)(g_\eps(a))+s+4\kappa
 <(\rho\circ\mathbb E)(g_{2\delta}(b))
\end{equation}
for all invariant $\rho$.

For all sufficiently invariant $F$, the same inequality with
$3\kappa$ holds for $\bar\tau_F$, uniformly in
$\tau\in\mathrm T(A)$.  Indeed, otherwise a sequence of increasingly
invariant sets and traces violating the inequality would have
orbit-average traces converging to an invariant trace, contradicting
\eqref{orbit-continuous-trace-gap}.  Separability is not needed here:
one may take a weak* convergent subnet of these states.
The inverse convention in $\bar\tau_F$ is compatible with left
invariance.

Apply Lemma~\ref{orbit-moment-lemma} to $g_\eps$ and $g_{2\delta}$,
with error smaller than $\kappa/4$ in each application.
The inequalities
\[
 1_{(\eps,\infty)}\leq g_\eps,
 \qquad g_{2\delta}\leq1_{(\delta,\infty)}
\]
give \eqref{orbit-uniform-rank-gap} for $x=a,y=b$ with a strictly
positive margin.  Regular compression and scalar compression are
contractive.  Uniform continuity of the two fixed
continuous cutoff functions, and the trace bound $l$, give a
$\nu>0$ for which the same estimates hold for $x,y$ as stated.
The target cutoff $\delta$ is therefore fixed before these
approximations are chosen.
\end{proof}

\begin{prop}[Averaging finite radius of comparison]\label{orbit-finite-radius-zero}
Suppose that $A$ is stably finite and exact, and that
$\operatorname{rc}(A)<\infty$.
If $\Gamma$ is infinite, then every action satisfies
\[
 \operatorname{orc}(A,\alpha)=0.
\]
If $\Gamma$ is finite of cardinality $q$, then
\[
 \operatorname{orc}(A,\alpha)
 \leq\frac{\operatorname{rc}(A)}q.
\]
\end{prop}

\begin{proof}
Choose $R>\operatorname{rc}(A)$ for which $A$ has $R$-comparison.
For infinite $\Gamma$, fix $t>0$ and choose a finite set $K$ with
$|K|>R/t$.  Take $0<\eta<1$.
If $F$ is nonempty and $(K,\eta)$-invariant, then
$\operatorname{Int}_K(F)\ne\emptyset$, so $|F|\geq|K|$.
Under \eqref{orbit-local-rank-gap}, the unnormalized matrix traces give
\[
 \mathrm d_\tau((R_{F,\Lambda}(x)-\eps)_+)+R
 <\mathrm d_\tau((R_{F,\Lambda}(y)-\delta)_+),
 \qquad \tau\in\mathrm T(A).
\]
Ordinary comparison over $A$ gives the stronger conclusion
\[
 (R_{F,\Lambda}(x)-\eps)_+
 \precsim (R_{F,\Lambda}(y)-\delta)_+.
\]
This implies \eqref{orbit-local-comparison} and proves orbit
$0$-comparison.  Exactness justifies testing ordinary comparison
on traces.

For finite $\Gamma$, take $K=\Gamma$ and $\eta<1/q$.
The invariance condition forces $F=\Gamma$, and
$|D(\Lambda)|<\eta q<1$ forces $\Lambda=1$.
The same proof gives orbit $R/q$-comparison.  Let
$R\downarrow\operatorname{rc}(A)$.
\end{proof}

\begin{rem}\label{orbit-zero-analogy}
The infinite-group assertion requires no minimality, WRP, or COS.
It is the comparison analogue of the vanishing of mean dimension
for actions of infinite amenable groups on finite-dimensional compact
spaces: a fixed comparison obstruction is divided by a growing orbit
height.  Even a trivial infinite-group action can have orbit radius
zero.  The additional tower and complement-comparison hypotheses
are used in the crossed-product theorem below.
\end{rem}

\begin{prop}[The topological bound]\label{orbit-radius-mdim}
Let $X$ be a nonempty compact metrizable space, and let an infinite
countable amenable group $\Gamma$ act on $X$.  Then
\[
 \operatorname{orc}(\mathrm C(X),\alpha)
 \leq\tfrac12\operatorname{mdim}(X,\Gamma).
\]
No freeness, minimality, or Rokhlin hypothesis is required here.
\end{prop}

\begin{proof}
There is nothing to prove if the mean dimension is infinite.
Put $d=\operatorname{mdim}(X,\Gamma)$ and fix $t>d/2$,
positive Fourier polynomials $x,y\in \mathrm M_l(B)$, and $\eps,\delta>0$.
Choose $\xi>0$ with $d+\xi<2t$ and
$0<\nu<\min\{\eps/2,\delta/4\}$.

Choose a finite open cover $\mathcal U$ controlling the finitely many
coefficient functions of $x$ and $y$ so closely that membership in a
common atom of
\[
 \mathcal U_F=\bigvee_{g\in F}g^{-1}\mathcal U
\]
implies norm variation less than $\nu$ for both orbit matrices.
The estimate is independent of $F$: every row and column has at most
the fixed number of Fourier-support blocks, and each coefficient
block has uniformly small variation.  It also remains valid after
multiplication by any fixed scalar collar, since $\|\Lambda\|\leq1$.

For all sufficiently invariant $F$, choose a refinement
$\mathcal V_F=\{V_i\}_{i=1}^n$ of $\mathcal U_F$ whose nerve $P_F$
satisfies
\begin{equation}\label{orbit-nerve-dimension}
 \dim P_F<(d+\xi)|F|<2t|F|.
\end{equation}
This uniform statement follows from the F\o lner-independent limit
for the cover-dimension function
\cite[Theorem~2.1 and Lemma~3.1]{Coornaert-Krieger-mdim}.
If uniformity failed, sets violating the dimension estimate along a
countable exhaustion would form a F\o lner sequence contradicting that
limit.

Choose $z_i\in V_i$ and a partition-of-unity map
$\varphi:X\to P_F$ subordinate to $\mathcal V_F$.
Write $R_x=R_{F,\Lambda}(x)$ and $R_y=R_{F,\Lambda}(y)$, and define
positive matrix functions on $P_F$ by
\[
 A_F(v)=\sum_i v_iR_x(z_i),\qquad
 B_F(v)=\sum_i v_iR_y(z_i).
\]
These are positive because they are convex combinations of positive
matrices.  For a point $v$ in a simplex with vertices $I$, choose
$z\in\bigcap_{i\in I}V_i$.  Each $z_i$, $i\in I$, lies in a common
atom of $\mathcal U_F$ with $z$.  Hence
\[
 \|A_F(v)-R_x(z)\|<\nu,
 \qquad \|B_F(v)-R_y(z)\|<\nu.
\]
If \eqref{orbit-local-rank-gap} holds, apply it to the evaluation
trace at $z$.  Eigenvalue perturbation then gives
\[
 \operatorname{rank}((A_F(v)-3\eps/2)_+)+t|F|
 <\operatorname{rank}((B_F(v)-3\delta/4)_+).
\]
By \eqref{orbit-nerve-dimension} and homogeneous comparison
\cite[Theorem~4.6]{Toms-Comp-DS},
\[
 (A_F-3\eps/2)_+\precsim(B_F-3\delta/4)_+
 \quad\hbox{in }\mathrm M_{l|F|}(\mathrm C(P_F)).
\]

Pull this comparison back by $\varphi$.  By construction,
\[
 \|A_F\circ\varphi-R_x\|<\nu,
 \qquad \|B_F\circ\varphi-R_y\|<\nu.
\]
The estimates
\[
 \|R_x-(A_F\circ\varphi-3\eps/2)_+\|
 <\nu+3\eps/2<2\eps,
\]
and
\[
 \|B_F\circ\varphi-(R_y-\delta/2)_+\|
 <\nu+\delta/2<3\delta/4
\]
give, by the norm perturbation lemma for Cuntz comparison,
\[
 (R_x-2\eps)_+
 \precsim(A_F\circ\varphi-3\eps/2)_+
 \precsim(B_F\circ\varphi-3\delta/4)_+
 \precsim(R_y-\delta/2)_+.
\]
This proves the required orbit comparison for $r=d/2$.
This argument applies without a bound on the number of tapered levels.
The same finite-input factorization occurs in \cite[Theorem~3.9]{Niu-MD-Z}.
\end{proof}

\begin{lem}[Bounded spectral witnesses]\label{orbit-bounded-witness}
Let $p,q\in \mathrm M_n(D)_+$ for a \Cstar{}-algebra $D$, let $a\geq0$ and
$b>0$, and suppose that $(p-a)_+\precsim(q-b)_+$.
For every $\gamma>0$ there is $v\in \mathrm M_n(D)$ such that
\[
 vqv^*=(p-a-\gamma)_+,
 \qquad \|v\|\leq\sqrt{\|p\|/b}.
\]
\end{lem}

\begin{proof}
The cut-down factorization, compressed to the $n\times n$ corner when
comparison is taken in the stabilization, gives
$w\in \mathrm M_n(D)$ with $w(q-b)_+w^*=(p-a-\gamma)_+$.
Put $z=w(q-b)_+^{1/2}$; then $\|z\|^2\leq\|p\|$.
Let $f(t)=b^{-1/2}$ on $[0,b]$ and $f(t)=t^{-1/2}$ on $[b,\infty)$.
The element $v=zf(q)$ belongs to $\mathrm M_n(D)$, with functional calculus
taken in its unitization if necessary.  Its norm has the stated
bound, and
\[
 (q-b)_+^{1/2}f(q)qf(q)(q-b)_+^{1/2}=(q-b)_+
\]
proves the required identity.
\end{proof}

\begin{lem}[Hereditary localization of orbit comparison]
\label{orbit-tower-localization}
Suppose that $A$ is separable and the action has WRP.
Let $x_j,y_j$ be positive Fourier
polynomials in $\mathrm M_l(B)$ with
$x_j\to a$, $y_j\to b_0$, and $\|x_j\|\leq M$.
Fix $\eps,\delta>0$.  For each $j$, suppose that an exact tiling
has shapes $F_{s,j}$ and scalar collars $\Lambda_{s,j}$ such that
\[
 (R_{F_{s,j},\Lambda_{s,j}}(x_j)-2\eps)_+
 \precsim
 (R_{F_{s,j},\Lambda_{s,j}}(y_j)-\delta/2)_+.
\]
Assume that the collars are chosen as in
Theorem~\ref{WTA-C-dyn-1}, with boundary fractions tending to zero
and widths tending to infinity sufficiently fast for the Fourier
supports and a prescribed increasing dense family of finite subsets
of $B$, and that the shapes satisfy the invariance conditions required
for these collars.  Then there are tower algebras
$C_j=\bigoplus_s\mathrm M_{|F_{s,j}|}(H_{s,j})\subseteq B$, contractions
$h_{0,j}\in A\cap C_j$, and $T_j\in \mathrm M_l(C_j)$ such that
\begin{align}
 &\|T_j\|\leq C:=\sqrt{2M/\delta},\nonumber\\
 &\|T_jh_{0,j}y_jh_{0,j}T_j^*-h_{0,j}x_jh_{0,j}\|
       \leq3\eps+\alpha_j,
       \qquad \alpha_j\longrightarrow0.\label{orbit-tower-estimate}
\end{align}
For $h_j=g_{1/2}(h_{0,j})$ one may also arrange
\begin{equation}\label{orbit-tower-defect}
 \|[h_j,z]\|\longrightarrow0\quad(z\in B),
 \qquad
 \sup_{\rho\in\mathcal M_1(A,\Gamma)}\mathrm d_\rho(1-h_j)
       \longrightarrow0.
\end{equation}
\end{lem}

\begin{proof}
Fix $j$ and suppress it in the notation.  Let $N$ be a finite symmetric
set containing $e$, the Fourier supports of $x,y$, and the finite set
used in the collar construction.
The assumed collar construction gives
\[
 \lambda_s(g)>0\ \Longrightarrow\ Ng\subseteq F_s,
 \qquad |\lambda_s(kg)-\lambda_s(g)|\leq L^{-1}\quad(k\in N),
\]
with $L=L_j\to\infty$ and the collars equal to one off a boundary
whose relative cardinality tends to zero.  The exact tiling theorem
supplies these shapes before the orbit comparisons are applied.

Lemma~\ref{orbit-bounded-witness} gives matrices $V_s$ over $A$ with
\begin{equation}\label{orbit-shaped-witness}
 V_sR_{F_s,\Lambda_s}(y)V_s^*
       =(R_{F_s,\Lambda_s}(x)-3\eps)_+,
 \qquad \|V_s\|\leq C.
\end{equation}
With the matrices $V_s$ fixed, WRP gives bases
$c_s\in A'\cap A_\omega$.  These commute with every entry of the
matrices, and all the levels $\alpha_g(c_s)$ are mutually orthogonal.
Set
\[
 h_0=\sum_{s,g}\lambda_s(g)\alpha_g(c_s).
\]
For each $s$, the map
\[
 \psi_s(Z)=\sum_{g,h\in F_s}u_gc_sZ_{gh}c_su_h^*
\]
is c.p.c. order zero, since $c_s$ commutes with $A$ and its levels
are orthogonal.  The Fourier-interior condition gives the exact
identities
\begin{equation}\label{orbit-tower-identity}
 \sum_s\psi_s(R_{F_s,\Lambda_s}(x))=h_0xh_0,
 \qquad
 \sum_s\psi_s(R_{F_s,\Lambda_s}(y))=h_0yh_0.
\end{equation}
Indeed, in a Fourier term with index $k\in N$, a nonzero source
collar level $h\in F_s$ has $kh\in F_s$.  Orthogonality leaves only
the target level $g=kh$ in the same tower.

The relative cone lifts and linking-algebra isomorphisms of
Lemmas~\ref{lift-d-drop} and \ref{app-tower-alg} give coordinate base
diagonals $d_{s,m}$ and linking isomorphisms:
\[
 \Theta_{s,m}:\mathrm M_{F_s}(\operatorname{Her}_A(d_{s,m}))\longrightarrow C_{s,m}
\]
For $f_n(t)=\min\{1,nt\}$ define
\[
 T_{s,n,m}=\Theta_{s,m}^{(l)}
 \big((1_{F_s}\otimes1_l\otimes f_n(d_{s,m}))V_s
            (1_{F_s}\otimes1_l\otimes f_n(d_{s,m}))\big).
\]
This is an element of $\mathrm M_l(C_{s,m})$ of norm at most $C$.
Its quotient image is
\[
 \sum_{g,h\in F_s}u_gf_n(c_s)(V_s)_{gh}f_n(c_s)u_h^*.
\]
The formula follows from the linking identities because
$f_n(t)/t$ extends continuously at zero.

Let $q_s=\sum_{g\in F_s}\alpha_g(c_s)$ and denote the quotient
image of $T_{s,n,m}$ by $T_{s,n}$.  Centrality of the base and
\eqref{orbit-shaped-witness} give
\[
 T_{s,n}\psi_s(R_{F_s,\Lambda_s}(y))T_{s,n}^*
 =f_n(q_s)^4\psi_s((R_{F_s,\Lambda_s}(x)-3\eps)_+).
\]
The order zero map has support element $\psi_s(1)=q_s^2$.
Thus the truncation error on the $s$th tower is bounded by
\[
 M\sup_{0\leq t\leq1}t^2(1-f_n(t)^4)\leq M/n^2.
\]
Replacing the source by its $3\eps$ cutoff adds at most $3\eps$ to
the error.
The tower summands are orthogonal, so their combined norm error is
the maximum of their errors.  The sum of the coordinate witnesses
also has norm at most $C$.

Choose $n=n_j\to\infty$ and then one sufficiently accurate
coordinate $m=m_j$.  As in the proof of
Theorem~\ref{WTA-C-dyn-1}, put
$h_{0,j}=\sum_{s,g}\lambda_s(g)d_{s,g,m_j}$.
The linking identities and \eqref{orbit-tower-identity} give
\eqref{orbit-tower-estimate}, with $\alpha_j\to0$.
The coordinate $m_j$ can also be chosen to satisfy the commutator,
covariance, linking, and uniform $2$-norm defect estimates of that
theorem.  The shrinking boundary fractions and increasing collar widths
give \eqref{orbit-tower-defect} after applying $g_{1/2}$.  Each step
imposes only finitely many conditions, so $m_j$ can be chosen after
the witnesses and their truncations have been fixed.
\end{proof}

\begin{thm}[Comparison in the crossed product]\label{orbit-radius-crossed-product}
Let $(A,\Gamma)$ be minimal with WRP and $(\lambda,m)$-COS,
where $\lambda\in(0,1]$, $m\in\mathbb N$, $A$ is unital and
separable, $\mathrm T(A)\ne\emptyset$,
and $\Gamma$ is countably infinite, discrete, and amenable.
If $B=A\rtimes_{\alpha,r}\Gamma$ is exact, then
\[
 \operatorname{rc}(B)\leq\operatorname{orc}(A,\alpha).
\]
\end{thm}

\begin{proof}
The crossed product is simple by Corollary~\ref{wrp-simple}, and
Lemma~\ref{no-atoms} shows that every nonzero hereditary coefficient
corner is infinite-dimensional.  WRP identifies its traces with
invariant coefficient traces by Proposition~\ref{trace-factor}.
They are faithful, and $B$ is infinite-dimensional and stably finite.
Exactness identifies the normalized $2$-quasitraces with traces.

Let $r$ be an admissible orbit radius, and let $q>r$.
Suppose that $a,b\in \mathrm M_l(B)_+$ satisfy
\[
 \mathrm d_\tau(a)+q<\mathrm d_\tau(b),\qquad \tau\in\mathrm T(B).
\]
Fix $\zeta>0$ and put
$p=(a-\zeta/2)_+$.
Choose $s$ with $r<s<q$.
By Lemma~\ref{comparison-tools}(1), applied to $B$, choose orthogonal
$b_0,b_1\in\operatorname{Her}_{\mathrm M_l(B)}(b)_+$ with $b_1\ne0$ and
$\mathrm d_\tau(b_0)>\mathrm d_\tau(b)-(q-s)/2$ for all traces.
Thus
\[
 \mathrm d_\tau(p)+s<\mathrm d_\tau(b_0),\qquad \tau\in\mathrm T(B),
 \qquad b_0\oplus b_1\precsim b.
\]

A nonzero diagonal compression of $b_1$, followed by
Lemma~\ref{coefficient-cutdown}, gives $0\ne c\in A_+$ with
$c\precsim_B b_1$.
Choose $\sigma>0$ with $c'=(c-\sigma)_+\ne0$; then
$g_\sigma(c)c'=c'$.  Apply
Lemma~\ref{coefficient-multiplicity-reserves} to $c'$ and
$g_\sigma(c)$, with $lm$ columns.  Its orthogonal range
squares belong to $\operatorname{Her}_A(c)$, so the common source
$z\in A_+$ satisfies
\begin{equation}\label{orbit-complement-reserve}
 z^{\oplus lm}\precsim_B b_1.
\end{equation}
Choose $\theta>0$ with $z_\theta=(z-\theta)_+\ne0$.
Minimality makes invariant traces faithful; hence
\[
 \kappa_z=\min_{\rho\in\mathcal M_1(A,\Gamma)}\mathrm d_\rho(z_\theta)>0.
\]

Choose one $\eps>0$ with $48\eps<\zeta/2$.
Lemma~\ref{orbit-trace-gap}, applied to $p,b_0,s$, gives fixed
$\delta,\kappa,\nu>0$ and a F\o lner condition.
Choose positive Fourier polynomials $x_j\to p$ and $y_j\to b_0$
within the specified neighborhood, with $\|x_j\|\leq M$.
Fix $t$ with $s<t<s+\kappa$.
For each $j$, use orbit $r$-comparison for $x_j,y_j,\eps,\delta,t$.
Choose the shape invariance and collar boundary fractions to satisfy
this comparison condition, Lemma~\ref{orbit-trace-gap}, and the finite
centrality conditions at step $j$.
The exact tiling theorem and the collar construction of
Theorem~\ref{WTA-C-dyn-1} give compatible shapes and collars.
Their orbit rank gaps imply the internal comparisons in
Lemma~\ref{orbit-tower-localization}.

That lemma gives $h_{0,j},T_j$, with
$\|T_j\|\leq C=\sqrt{2M/\delta}$, and put
\[
 h_j=g_{1/2}(h_{0,j}),\qquad
 k_j=\frac{g_{1/2}(h_{0,j})}{h_{0,j}},\qquad S_j=k_jT_j.
\]
The quotient defining $k_j$ extends continuously by zero at the
origin; the bound $\|k_j\|\leq4$ suffices.
Consequently $\|S_j\|\leq4C$ and
\begin{align*}
 &\|S_jh_{0,j}b_0h_{0,j}S_j^*-h_jph_j\|\\
 &\quad\leq48\eps+16\alpha_j
          +16C^2\|y_j-b_0\|+\|x_j-p\|.
\end{align*}
Here $\eps$ and $\delta$ remain fixed as $j\to\infty$.
Their dependence on $\zeta$ therefore does not affect boundedness
of the witnesses as $j$ varies.

For large $j$, \eqref{orbit-tower-defect} gives
\[
 \mathrm d_\rho(1-h_j^2)=\mathrm d_\rho(1-h_j)<\lambda\kappa_z,
 \qquad \rho\in\mathcal M_1(A,\Gamma).
\]
COS yields $1-h_j^2\precsim_B z_\theta^{\oplus m}$.
Put $h=[(h_j)]\in B_\omega$, pad $1-h_j^2$ with zeroes to
an element $P_j\in \mathrm M_m(B)_+$,
and put $Z=z^{\oplus m}$.  Lemma~\ref{orbit-bounded-witness},
with the fixed target cutoff $\theta$, gives $v_j\in \mathrm M_m(B)$ with
\[
 v_jZv_j^*=(P_j-1/j)_+,\qquad \|v_j\|\leq\theta^{-1/2}.
\]
The bounded sequence $(v_j)$ therefore implements
$1-h^2\precsim_{B_\omega}z^{\oplus m}$.
Thus, in $\operatorname{Cu}(B_\omega)$,
\begin{equation}\label{orbit-defect-absorbed}
 l[1-h^2]\leq[z^{\oplus lm}]\leq[b_1].
\end{equation}

Let $h_0=[(h_{0,j})]$ and $S=[(S_j)]$.  The element $h$ commutes
with $p$.  Put
\[
 X=S h_0b_0h_0S^*,\qquad
 R=(1-h^2)^{1/2}p(1-h^2)^{1/2}.
\]
Then $X\precsim b_0$, $R\precsim(1-h^2)^{\oplus l}\precsim b_1$,
and $p=hph+R$.  The preceding norm estimate gives
\[
 \|p-(X+R)\|\leq48\eps<\zeta/2.
\]
Norm perturbation and \eqref{orbit-defect-absorbed} therefore give
\[
 (a-\zeta)_+=(p-\zeta/2)_+
 \precsim X+R\precsim X\oplus R
 \precsim b_0\oplus b_1\precsim b
 \quad\hbox{in }\mathrm M_\infty(B_\omega).
\]
Comparison between constant elements descends to $B$ by
Lemma~\ref{constant-ultrapower-comparison}, applied to matrix algebras.
Letting $\zeta\downarrow0$ proves $a\precsim b$.
Hence $\operatorname{rc}(B)\leq r$ for every admissible $r$,
and taking the infimum proves the theorem.
\end{proof}

\begin{defn}[Strong URP]\label{strong-urp-definition}
Let $A$ be unital and separable.  We say that $(A,\Gamma,\alpha)$
has \emph{strong URP} if, for every finite $K\subseteq\Gamma$
and every $\eta>0$, there are nonempty finite
$(K,\eta)$-invariant sets $F_1,\ldots,F_S$ and projections
$p_s\in A_\omega\cap A'$ such that the levels
$\alpha_g(p_s)$, $g\in F_s$, are mutually orthogonal and
\[
 \sum_{s=1}^S\sum_{g\in F_s}\alpha_g(p_s)=1
 \quad\text{in }A_\omega.
\]
The shapes are fixed for the central-sequence towers in this
condition; they may depend on $K$ and $\eta$.
\end{defn}

\begin{defn}[WRP with Cuntz control in the crossed product]
\label{cuntz-controlled-wrp-definition}
Let $A$ be unital and separable, with $\mathrm T(A)\ne\emptyset$,
and put $B=A\rtimes_{\alpha,r}\Gamma$.
We say that the action has \emph{WRP with Cuntz control in the
crossed product} if it has WRP and the following condition holds.
For every finite $K\subseteq\Gamma$, every $\eta>0$, and every
$0\ne z\in B_+$, there are nonempty finite
$(K,\eta)$-invariant shapes $F_1,\ldots,F_S$ and positive
contractions $c_s\in A_\omega\cap A'$ whose levels
$\alpha_g(c_s)$, $g\in F_s$, are mutually orthogonal, such that
\[
 q=\sum_{s=1}^S\sum_{g\in F_s}\alpha_g(c_s),\qquad
 1-q\in J_{\mathcal M_1(A,\Gamma)},\qquad
 1-q\precsim_{B_\omega}z.
\]
Here $A_\omega$ is viewed as a subalgebra of $B_\omega$, and $z$
denotes its constant image.  The shapes may depend on $K,\eta,z$,
but are fixed within the central-sequence towers.  The reserve $z$
is prescribed before choosing these towers.
\end{defn}

The additional comparison condition is imposed in the norm
ultrapower of the crossed product, in addition to the usual uniform
tracial control of the complement.  It differs from the weak tracial
Rokhlin property of \cite[Definition~2.1]{Wang-WRP}, where comparison
with a prescribed coefficient element is imposed pointwise inside $A$.

\begin{lem}[Full towers and Cuntz control]
\label{strong-urp-cuntz-controlled-wrp}
If $A$ is unital and separable with $\mathrm T(A)\ne\emptyset$,
then strong URP implies WRP with Cuntz control in the crossed product.
\end{lem}

\begin{proof}
Fix an exact tiling with prescribed shapes $D_s$.
Normalize the shapes to contain the identity and put
$L=\bigcup_sD_sD_s^{-1}$.  Apply strong URP with $K=L$ and
$\eta\downarrow0$, and cut its projection towers along this
prescribed tiling, retaining only tiles wholly contained in each
$F_r$.  The construction in the proof of Proposition~\ref{urp-implies-wrp}
applies in $A_\omega$.
The new bases are central projections and their prescribed levels
are orthogonal.  At most $\eta|F_r|$ original levels are lost
from the $r$th tower.  For every ultralimit $\rho_\omega$ of a
possibly varying sequence of invariant traces, their total trace
is therefore at most
\[
 \eta\sum_r|F_r|\rho_\omega(p_r)=\eta.
\]
Choosing representatives and diagonalizing as in
Lemma~\ref{local-wrp} gives the prescribed-tiling WRP.

The original full towers from strong URP satisfy the additional
condition in Definition~\ref{cuntz-controlled-wrp-definition}.
Their sum is $q=1$ in $A_\omega$,
so $1-q=0\precsim_{B_\omega}z$ for every prescribed nonzero
$z\in B_+$.  These towers need not be the recut towers used to
establish prescribed-tiling WRP.
\end{proof}

\begin{prop}[Comparison from a Cuntz-small tower complement]
\label{cuntz-controlled-wrp-orbit-zero-cos}
Let $(A,\Gamma,\alpha)$ be minimal, where $A$ is unital and
separable, $\mathrm T(A)\ne\emptyset$, and $\Gamma$ is countably
infinite, discrete, and amenable.  Suppose that the action has WRP
with Cuntz control in the crossed product, and assume that
$\operatorname{orc}(A,\alpha)=0$.  Then it has $(1,1)$-COS.
If $B=A\rtimes_{\alpha,r}\Gamma$ is exact, it also has strict
comparison of positive elements.
In particular, these conclusions hold when $A$ is stably finite
and exact with $\operatorname{rc}(A)<\infty$.

Strong URP is a special case, by
Lemma~\ref{strong-urp-cuntz-controlled-wrp}.  For strong URP,
the implication from zero orbit radius to $(1,1)$-COS holds even
without minimality and also for finite amenable groups.
\end{prop}

\begin{proof}
Let $a,b\in A_+$ satisfy $\mathrm d_\rho(a)<\mathrm d_\rho(b)$ for every
$\rho\in\mathcal M_1(A,\Gamma)$.  Rescale $b$ so that it is a
contraction; this does not change its dimension functions or Cuntz
class.  Fix $\zeta>0$ and put $p=(a-\zeta/2)_+$.
Since $\rho(g_{\zeta/2}(a))\leq \mathrm d_\rho(a)$, lower
semicontinuity and compactness give
\[
 \mu=\min_{\rho\in\mathcal M_1(A,\Gamma)}
       \bigl(\mathrm d_\rho(b)-\rho(g_{\zeta/2}(a))\bigr)>0.
\]
Choose $N\in\mathbb N$ with $1/N<\mu/2$.

Lemma~\ref{no-atoms} shows that every nonzero hereditary
\Cstar{}-subalgebra of $A$ is infinite-dimensional.  This permits the
choice of a reserve in the coefficient algebra.
Choose $\sigma_0>0$ with $(b-\sigma_0)_+\ne0$ and apply
Lemma~\ref{coefficient-multiplicity-reserves} to
$(b-\sigma_0)_+$ and $g_{\sigma_0}(b)$, with $N$ columns.
Their range squares $v_1,\ldots,v_N$ are nonzero, mutually
orthogonal, belong to $\operatorname{Her}_A(b)$, and are Cuntz
equivalent in $B$.  Each invariant $\rho$ extends to the trace
$\rho\circ\mathbb E$ on $B$.  Consequently
\[
 N \mathrm d_\rho(v_1)=\mathrm d_\rho\Bigl(\sum_{i=1}^Nv_i\Bigr)
 \leq \mathrm d_\rho(b)\leq1,
 \qquad \rho\in\mathcal M_1(A,\Gamma).
\]
Choose $\sigma>0$ so that $b_1=(v_1-\sigma)_+\ne0$, and put
\[
 h=g_\sigma(v_1),\qquad b_0=(1-h)b(1-h).
\]
Then $hb_1=b_1$, so $b_0b_1=0$, and
$b_0,b_1\in\operatorname{Her}_A(b)$ give
$b_0\oplus b_1\precsim_A b$.
The decomposition
\[
 b^{1/2}=(1-h)b^{1/2}+hb^{1/2}
\]
gives $b\leq2b_0+2hbh$.  Dimension-function subadditivity and
$hbh\in\operatorname{Her}_A(h)$ therefore imply
$\mathrm d_\rho(b)\leq \mathrm d_\rho(b_0)+\mathrm d_\rho(h)$.
Since $h\in\operatorname{Her}_A(v_1)$, it follows that
\[
 \mathrm d_\rho(b_0)\geq \mathrm d_\rho(b)-1/N
 >\mathrm d_\rho(p)+\mu/2.
\]
The nonzero element $b_1$ is fixed before the towers are chosen and
will absorb their complement.

Apply Lemma~\ref{orbit-trace-gap} to $p,b_0,s=0$, with
$\eps=\zeta/4$ and scalar collar $1$.
It gives $\delta,\kappa>0$ and a F\o lner condition under which
\[
 \widehat{\mathrm d}_{\tau,F}((\Pi_F(p)-\eps)_+)+\kappa
 <\widehat{\mathrm d}_{\tau,F}((\Pi_F(b_0)-\delta)_+)
 \quad(\tau\in\mathrm T(A)).
\]
Choose $0<t<\kappa$.  Since the orbit radius is zero, choose an
admissible orbit radius $r<t$ and strengthen this F\o lner condition
to satisfy orbit $r$-comparison for $p,b_0,\eps,\delta,t$.
Coefficient elements are Fourier polynomials, and their orbit
matrices are diagonal, so the resulting comparison is
\[
 \Pi_F((a-\zeta)_+)
 \precsim_{\mathrm M_\infty(A)}\Pi_F((b_0-\delta/2)_+)
\]
for every qualifying shape $F$.

Choose the towers of
Definition~\ref{cuntz-controlled-wrp-definition} for this
strengthened condition and the already fixed reserve $z=b_1$.
Write $q_s=\sum_{g\in F_s}\alpha_g(c_s)$, so
$q=\sum_s q_s$ and $1-q\precsim_{B_\omega}b_1$.
Define
\[
 \Psi_s:\mathrm M_{F_s}(A)\longrightarrow B_\omega,
 \qquad
 \Psi_s(Z)=\sum_{g,h\in F_s}
        u_gc_s^{1/2}Z_{gh}c_s^{1/2}u_h^*.
\]
The row $V_s=(u_gc_s^{1/2})_{g\in F_s}$ satisfies
$V_s^*V_s=1_{F_s}\otimes c_s$, which commutes with coefficient
matrices.  Hence $\Psi_s$ is c.p.c. and has order zero, with
$\Psi_s(1)=q_s$.  The ranges of the different $\Psi_s$ are
orthogonal.  For $x\in A$, diagonal reconstruction gives
\[
 \Psi_s(\Pi_{F_s}(x))=q_sx,
 \qquad \sum_s\Psi_s(\Pi_{F_s}(x))=qx.
\]
Order zero maps preserve Cuntz comparison
\cite[Corollary~4.5]{WZ-order-zero}.
Apply them to the orbit comparisons and take the orthogonal sum to
obtain, with $a_\zeta=(a-\zeta)_+$,
\[
 qa_\zeta\precsim q(b_0-\delta/2)_+\precsim b_0
 \quad\text{in }B_\omega.
\]
Since $q$ commutes with $A$, the remaining part is positive and
satisfies
\[
 (1-q)a_\zeta\precsim1-q\precsim b_1.
\]
Subadditivity of Cuntz classes gives
\[
 a_\zeta=qa_\zeta+(1-q)a_\zeta
 \precsim qa_\zeta\oplus(1-q)a_\zeta
 \precsim b_0\oplus b_1\precsim b.
\]
Comparison of constant elements descends from $B_\omega$ to $B$
by Lemma~\ref{constant-ultrapower-comparison}.  Letting
$\zeta\downarrow0$ proves $a\precsim_B b$, which is $(1,1)$-COS.

Under strong URP, apply the same orbit comparison and order zero
reconstruction with $b_0=b$, $b_1=0$, and full projection towers.
Fix $\zeta>0$, put
$p=(a-\zeta/2)_+$, and apply
Lemma~\ref{orbit-trace-gap} directly to $p,b$ with $s=0$,
$\eps=\zeta/4$, and scalar collar $1$.
Zero orbit radius gives the diagonal orbit comparisons above.
Choose full towers satisfying their F\o lner condition; then
$q=1$, so reconstruction gives $(a-\zeta)_+\precsim b$ in
$B_\omega$.  Lemma~\ref{constant-ultrapower-comparison} gives the
comparison in $B$.
This argument requires no coefficient reserve and proves the COS
assertion without minimality or infiniteness.

If $A$ is stably finite and exact with finite radius of comparison
and $\Gamma$ is infinite, Proposition~\ref{orbit-finite-radius-zero}
gives orbit radius zero.  Under the minimal, infinite-group
hypotheses of the proposition, exactness of $B$ permits an application of
Theorem~\ref{orbit-radius-crossed-product}, giving
$\operatorname{rc}(B)=0$.  For a strict tracial gap, fixing a
source cutoff gives a uniform positive gap by compactness and
continuous cutoff functions, as in
Lemma~\ref{orbit-trace-gap}.  Since $B$ has $t$-comparison for
every $t>0$, each source cutoff is below the target.  Letting the
cutoff tend to zero proves strict comparison.
\end{proof}

\begin{cor}\label{finite-coefficient-radius-crossed-comparison}
Assume the hypotheses of Theorem~\ref{orbit-radius-crossed-product}
except for COS, and suppose that $\operatorname{rc}(A)<\infty$.
If the action has $(\lambda,m)$-COS for some
$\lambda\in(0,1]$, $m\in\mathbb N$, strong URP in the sense
of Definition~\ref{strong-urp-definition}, or WRP with Cuntz control
in the crossed product in the sense of
Definition~\ref{cuntz-controlled-wrp-definition}, then $B$ has
strict comparison of positive elements.
\end{cor}

\begin{proof}
Exactness passes from $B$ to $A$.
Proposition~\ref{orbit-finite-radius-zero} gives
$\operatorname{orc}(A,\alpha)=0$.
In the strong URP and Cuntz-controlled cases,
Proposition~\ref{cuntz-controlled-wrp-orbit-zero-cos} supplies
$(1,1)$-COS.
Theorem~\ref{orbit-radius-crossed-product}
therefore gives $\operatorname{rc}(B)=0$ in all three cases.
A strict tracial gap becomes uniformly positive
after fixing a source cutoff, by compactness and continuous cutoff
functions as in Lemma~\ref{orbit-trace-gap}.
Since $B$ has $q$-comparison for every $q>0$, each source cutoff
is below the target.  Letting the cutoff tend to zero proves strict
comparison.
\end{proof}

\begin{rem}\label{cuntz-controlled-wrp-quantifiers}
The proof only uses the additional tower condition for nonzero
reserves $z\in A_+$.  Thus this restricted condition suffices for
Proposition~\ref{cuntz-controlled-wrp-orbit-zero-cos}.
The ordinary WRP hypothesis supplies the coefficient reserves;
the particular towers used for comparison need only have a
Cuntz-small complement in $B_\omega$.
The reserve must be prescribed before the towers are chosen.
Comparison of all complements with a single fixed large element does
not suffice for this absorption argument.
\end{rem}

\begin{rem}[Exact coverage and the Fourier boundary]
\label{strong-urp-boundary-remark}
For full projection towers the untapered matrix tower algebra is
unital.  This is sufficient for the coefficient comparisons in
Proposition~\ref{cuntz-controlled-wrp-orbit-zero-cos}.  For general Fourier
inputs, however, the localization in
Lemma~\ref{orbit-tower-localization} still uses scalar collars
vanishing on the tower boundary.  If a nonzero projection level
is discarded, its norm is one even when its trace is small.
Thus exact coverage alone does not justify setting the localized
cutoff equal to $1$ in the proof of
Theorem~\ref{orbit-radius-crossed-product}.

For example, let $u$ implement a $\mathbb Z$-action with an
exact cyclic projection tower $p_i=u^ip_0u^{-i}\in A$,
$0\leq i<N$, whose sum is $1$.  Its unital matrix tower algebra is
\[
 D=\left\{\sum_{i,j=0}^{N-1}u^ia_{ij}u^{-j}:
                   a_{ij}\in p_0Ap_0\right\}.
\]
For $x\in D$, Fourier coefficients give
\[
 \mathbb E(p_0xp_{N-1}u^*)=0,
 \qquad
 \mathbb E(p_0up_{N-1}u^*)=p_0.
\]
Consequently $\operatorname{dist}(u,D)=1$.
The matrix entry joining the last level to the first requires the
return unitary $p_0u^Np_0$, which belongs to the crossed-product corner
rather than $p_0Ap_0$.
For this reason the general positive-orbit-radius bound retains
its COS hypothesis.  In the zero-orbit-radius case, including
finite coefficient radius when $\Gamma$ is infinite,
Proposition~\ref{cuntz-controlled-wrp-orbit-zero-cos} makes that hypothesis
automatic under strong URP.
\end{rem}

\section{Examples of orbit radius of comparison and COS}\label{cos-examples-section}

The following three families have zero orbit radius of comparison
and satisfy $(1,1)$-COS.  In the first two, the coefficient algebra
may have positive finite radius of comparison; we compute the
normalized comparison allowance in the orbit matrices and use full
towers to obtain COS.  In the third, the coefficient algebra has
radius of comparison zero, and COS follows from strict comparison
in the crossed product.  These examples include nonsimple coefficient
algebras, as allowed in Section~\ref{cos-section}.

\begin{example}[Cantor skew products with finite radius of comparison]
\label{cos-cantor-skew-example}
Let $X$ be the Cantor set, let $h:X\to X$ be a minimal homeomorphism,
and let $D$ be a simple, separable, unital, exact, stably finite
\Cstar{}-algebra with $\mathrm T(D)\ne\emptyset$ and
$\operatorname{rc}(D)<\infty$.  For a point-norm continuous map
$x\mapsto\beta_x$ from $X$ to $\operatorname{Aut}(D)$, put
\[
 A=\mathrm C(X,D),\qquad
 \alpha(f)(x)=\beta_x\bigl(f(h^{-1}(x))\bigr).
\]
Then the $\mathbb Z$-action generated by $\alpha$ is minimal, has
strong URP, and has $(1,1)$-COS.  Moreover,
\[
 \operatorname{rc}(A)=\operatorname{rc}(D),\qquad
 \operatorname{orc}(A,\alpha)=0.
\]
\end{example}

\begin{proof}
For every nonempty compact metrizable zero-dimensional space $Y$,
\begin{equation}\label{eq:zero-dimensional-fibre-radius}
 \operatorname{rc}(\mathrm C(Y,D))=\operatorname{rc}(D).
\end{equation}
Indeed, let $r>\operatorname{rc}(D)$, and let
$a,b\in \mathrm M_\ell(\mathrm C(Y,D))_+$ satisfy
$\mathrm d_\rho(a)+r<\mathrm d_\rho(b)$ for every tracial state $\rho$ of $\mathrm C(Y,D)$.
The evaluation traces $\sigma\circ\operatorname{ev}_y$, with
$y\in Y$ and $\sigma\in\mathrm T(D)$, give
\[
 \mathrm d_\sigma(a(y))+r<\mathrm d_\sigma(b(y)).
\]
Consequently $a(y)\precsim_D b(y)$ for each $y$.  Fix $\eps>0$.
For every $y$, choose $v_y\in \mathrm M_\ell(D)$ such that
\[
 \|a(y)-v_yb(y)v_y^*\|<\eps/2.
\]
A stabilized comparison witness gives such a $v_y$ after compression
to the $\ell\times\ell$ corner.
Continuity gives a neighborhood $U_y$ on which the same norm is
less than $\eps$.  Refine a finite subcover to a finite clopen
partition $V_1,\ldots,V_q$, with $V_j\subseteq U_{y_j}$, and put
$v=\sum_j1_{V_j}\otimes v_{y_j}$.  Then
$\|a-vbv^*\|<\eps$.  Letting $\eps$ tend to zero gives $a\precsim b$,
so $\operatorname{rc}(\mathrm C(Y,D))\leq\operatorname{rc}(D)$.
Conversely, $r$-comparison for $\mathrm C(Y,D)$ implies $r$-comparison for $D$:
apply it to constant functions, using restriction of traces and
evaluation at any $y\in Y$.  This proves the reverse inequality.

To construct full central towers, fix a finite
$K\subseteq\mathbb Z$ and $\eta>0$.  Choose $H$ so large that every
interval $\{0,\ldots,l-1\}$ with $l\geq H$ is $(K,\eta)$-invariant.
Since $h$ is aperiodic, there is a nonempty clopen set $Y\subseteq X$
with $Y\cap h^j(Y)=\emptyset$ for $1\leq j<H$.
The first-return function to $Y$ is locally constant and bounded,
by minimality and compactness.  Its finitely many values
$l_1,\ldots,l_S$ are at least $H$.  Let $Y_s$ be the clopen set on
which the return time is $l_s$.  The Kakutani--Rokhlin partition
\cite[proof of Lemma~3.1]{Put-PJM} gives
\[
 X=\bigsqcup_{s=1}^S\bigsqcup_{j=0}^{l_s-1}h^j(Y_s).
\]
For $p_s=1_{Y_s}\otimes1_D\in Z(A)$, the fibre automorphisms fix
scalar functions, so
\[
 \alpha^j(p_s)=1_{h^j(Y_s)}\otimes1_D,\qquad
 \sum_{s=1}^S\sum_{j=0}^{l_s-1}\alpha^j(p_s)=1_A.
\]
These central projection towers give strong URP
(Definition~\ref{strong-urp-definition}).

Simplicity of $D$ identifies the ideals of $\mathrm C(X,D)$ with open subsets
of $X$.  Minimality of $h$ consequently makes the action minimal.
We compute the orbit radius using
\eqref{eq:zero-dimensional-fibre-radius}.  Fix $t>0$ and
$R>\operatorname{rc}(D)$, and choose a finite
$L\subseteq\mathbb Z$ with $|L|>R/t$ and $0<\eta<1$.
Every nonempty $(L,\eta)$-invariant $F$ has
$N:=|F|\geq|L|$, since $\operatorname{Int}_L(F)\ne\varnothing$.
For the Fourier inputs and scalar collars in
Definition~\ref{orbit-radius-definition}, put
\[
 a_F=(R_{F,\Lambda}(x)-\eps)_+,\qquad
 b_F=(R_{F,\Lambda}(y)-\delta)_+.
\]
If the normalized rank inequalities
\eqref{orbit-local-rank-gap} hold, then the unnormalized matrix
extensions satisfy
\[
 \mathrm d_\tau(a_F)+R<\mathrm d_\tau(a_F)+Nt<\mathrm d_\tau(b_F),
 \qquad\tau\in\mathrm T(A).
\]
Since $R>\operatorname{rc}(A)$, this gives $a_F\precsim b_F$,
and hence \eqref{orbit-local-comparison}.  Thus the normalized
comparison allowance is $R/N$, which tends to zero with the orbit
height.  As $t>0$ was arbitrary,
$\operatorname{orc}(A,\alpha)=0$.  The strong URP case of
Proposition~\ref{cuntz-controlled-wrp-orbit-zero-cos} gives
$(1,1)$-COS.

By \cite[Theorem~5.11]{Toms-Comp-DS}, the fibre $D$ can be a simple
unital AH algebra with any prescribed finite positive radius of
comparison.  The same radius then occurs in $A$, while the system
has $(1,1)$-COS.
\end{proof}

\begin{example}[Profinite translations with fibres of finite radius of comparison]
\label{cos-profinite-example}
Let $\Gamma$ be a countably infinite, residually finite, discrete
amenable group, and choose a decreasing sequence of finite-index
normal subgroups $\Gamma_n$ with $\bigcap_n\Gamma_n=\{e\}$.
Let $D$ satisfy the hypotheses of
Example~\ref{cos-cantor-skew-example}, and set
\[
 K=\varprojlim\Gamma/\Gamma_n,\qquad A=\mathrm C(K,D),\qquad
 \alpha_g(f)(x)=f(g^{-1}x).
\]
Then the translation action is free on $K$ and minimal, has strong
URP, and has $(1,1)$-COS.  In addition,
$\operatorname{rc}(A)=\operatorname{rc}(D)$ and
$\operatorname{orc}(A,\alpha)=0$.
\end{example}

\begin{proof}
The canonical map $\Gamma\to K$ is injective with dense image.
Left translation is therefore free and minimal on $K$.  Since $D$
is simple, minimality also holds for the action on $A$.
The space $K$ is compact, metrizable, and zero-dimensional, so
\eqref{eq:zero-dimensional-fibre-radius} gives the radius equality.
For $R>\operatorname{rc}(D)$ and $t>0$, choose a finite
$L\subseteq\Gamma$ with $|L|>R/t$ and an invariance tolerance
less than one.  Every corresponding nonempty invariant set $F$
then has $|F|\geq|L|$, so $R/|F|<t$.  The orbit-matrix comparison
calculation in Example~\ref{cos-cantor-skew-example} applies
unchanged and proves $\operatorname{orc}(A,\alpha)=0$.

By \cite[Corollary~5.6]{Deninger-Schmidt-Folner}, there are finite
sets $F_n\subseteq\Gamma$, each a set of representatives for
$\Gamma/\Gamma_n$, forming a F\o lner sequence.  Fix a finite
$L\subseteq\Gamma$ and $\eta>0$, and choose $n$ so that $F_n$ is
$(L,\eta)$-invariant in the left F\o lner convention.  Let
$\pi_n:K\to\Gamma/\Gamma_n$ be the quotient map and put
\[
 U_n=\pi_n^{-1}(\Gamma_n),\qquad
 p_n=1_{U_n}\otimes1_D\in Z(A).
\]
The cosets $gU_n$, $g\in F_n$, partition $K$, and hence
\[
 \alpha_g(p_n)\alpha_{g'}(p_n)=0\quad(g\ne g'),\qquad
 \sum_{g\in F_n}\alpha_g(p_n)=1_A.
\]
For these towers the normalized comparison allowance is explicitly
\[
 \frac{R}{|F_n|}=\frac{R}{[\Gamma:\Gamma_n]}
 \longrightarrow0,
 \qquad R>\operatorname{rc}(D),
\]
because $\Gamma$ is infinite and
$\bigcap_n\Gamma_n=\{e\}$.
Their constant images give a full central projection tower in
$A_\omega$, establishing strong URP.  Together with the orbit-radius
computation, strong URP gives
$(1,1)$-COS by
Proposition~\ref{cuntz-controlled-wrp-orbit-zero-cos}.
Thus this construction also allows fibres of positive finite radius.
\end{proof}

\begin{example}[Minimal skew products with pure fibres]
\label{cos-pure-fibre-example}
Let $X$ be an infinite compact metrizable space, let $h:X\to X$ be
a minimal homeomorphism, and let $D$ be a simple, separable, unital,
exact, stably finite, pure \Cstar{}-algebra with $\mathrm T(D)\ne\emptyset$.
Here purity means that $\operatorname{Cu}(D)$ is almost unperforated
and almost divisible.  For any point-norm continuous map
$x\mapsto\beta_x\in\operatorname{Aut}(D)$, define
\[
 A=\mathrm C(X,D),\qquad
 \alpha(f)(x)=\beta_x\bigl(f(h^{-1}(x))\bigr).
\]
Then
\[
 \operatorname{rc}(A)=0,\qquad
 \operatorname{orc}(A,\alpha)=0,
\]
and $(A,\mathbb Z,\alpha)$ has $(1,1)$-COS.  No restriction on
the covering dimension of $X$ or the mean dimension of $h$ is needed.
\end{example}

\begin{proof}
By \cite[Corollary~5.9]{Seth-Vilalta-pure-functions}, the algebra
$A=\mathrm C(X,D)$ is pure.  We first show that its radius of comparison
is zero.  Suppose that $a,b\in \mathrm M_n(A)_+$ satisfy
$\mathrm d_\tau(a)<\mathrm d_\tau(b)$ for all $\tau\in\mathrm T(A)$.
Fix a tracial state on $D$, which is faithful by simplicity.
The associated evaluation traces show that $b(z)\ne0$ for every
$z\in X$.  Hence $b$ is full in $\mathrm M_n(A)$, and there is
$L\in\mathbb N$ with $[1_A]\leq L[b]$.

If $\lambda$ is a functional on $\mathrm{Cu}(A)$ with
$\lambda([b])=1$, then
$1/n\leq\lambda([1_A])\leq L$.
After normalization at $[1_A]$, it is the rank functional of a
bounded $2$-quasitracial state, which is a trace because $A$ is exact.
Thus $\lambda([a])<\lambda([b])$ for every such functional.
Almost unperforation of $\mathrm{Cu}(A)$ and
\cite[Proposition~5.2.20(3)]{APT-Cu-tensor} give $a\precsim_A b$.
This proves $\operatorname{rc}(A)=0$.
Proposition~\ref{orbit-finite-radius-zero} now gives
$\operatorname{orc}(A,\alpha)=0$.

The action is free because $h$ is minimal on the infinite space $X$.
The crossed product $B=A\rtimes_\alpha\mathbb Z$ is simple by
\cite[Proposition~1.6]{ABP-CXD} and pure by
\cite[Theorem~4.6]{ABMPS-pureness}.  Exactness of $D$ and amenability
of $\mathbb Z$ give exactness of $B$.  Averaging a trace on $A$ and
composing with the canonical expectation gives a trace on $B$,
faithful by simplicity.  Thus $B$ is stably finite, and purity and
exactness give strict comparison by traces.

Suppose $a,b\in A_+$ satisfy $\mathrm d_\rho(a)<\mathrm d_\rho(b)$ for every
$\rho\in\mathcal M_1(A,\mathbb Z)$.  Every $\tau\in\mathrm T(B)$
restricts to such an invariant trace, and functional calculus in
$A\subseteq B$ gives
\[
 \mathrm d_\tau(a)=\mathrm d_{\tau|_A}(a)<\mathrm d_{\tau|_A}(b)=\mathrm d_\tau(b).
\]
Strict comparison in $B$ yields $a\precsim_B b$, proving $(1,1)$-COS.

For a concrete noncommutative Furstenberg transformation, take
$X=S^1$, $h(z)=e^{2\pi i\theta}z$ with
$\theta\in\mathbb R\setminus\mathbb Q$, and $D=C_r^*(\mathbb F_2)$.
For its canonical generators $u_1,u_2$, let
\[
 \beta_z(u_i)=zu_i\quad(i=1,2),\qquad
 \alpha(f)(z)=\beta_z\bigl(f(e^{-2\pi i\theta}z)\bigr).
\]
The maps $\beta_z$ are the point-norm continuous gauge automorphisms
arising from the characters of $\mathbb F_2$ with these values on
the generators.  The algebra $D$ is pure by
\cite[Proposition~7.3]{ABMPS-pureness}, and this construction appears in
\cite[Example~7.8]{ABMPS-pureness}.  It gives a system with
$\operatorname{orc}(A,\alpha)=0$, $(1,1)$-COS, and nonnuclear fibre.
\end{proof}

\chapter{\texorpdfstring{$\mathcal Z$}{Z}-stability and weak tracial approximation}\label{section-Z-stability}

We study $\mathcal Z$-stability for \Cstar{}-algebras with weak tracial
approximation, where $\mathcal Z$ denotes the Jiang--Su algebra.

\section{Weak tracial \texorpdfstring{$\mathcal Z$}{Z}-absorption}
The results of this section give an abstract formulation of part of
\cite[Section~4]{Niu-MD-Z-absorbing}.

The following decomposition uses three nested spectral cutoffs.  The
middle cutoff decomposes the given elements, while the inner and outer
cutoffs give exact annihilation and support relations, respectively.

\begin{lem}\label{cut-down}
Let $\mathcal S$ be a class of \Cstar{}-algebras, and let
$A\in\mathrm{TA}_0(\mathcal S)$.  For every finite set
$\{f_1,f_2,\ldots,f_n\}\subseteq A$ and every $\eps>0$, there are a
\Cstar{}-subalgebra $C\subseteq A$ with $C\in\mathcal S$, positive contractions
$\tilde h,h\in C$ with $\norm{h}=1$, and elements
$f_i^{(0)},f_i^{(1)}\in A$, $1\leq i\leq n$, such that
\begin{enumerate}
\item $\norm{f_i-(f_i^{(0)}+f_i^{(1)})}<\eps$;
\item $f_i^{(1)}\in C$;
\item $h f_i^{(0)}=f_i^{(0)}h=0$;
\item
  $\max\{\norm{[f_i^{(1)},h]},
           \norm{[f_i^{(1)},h^{1/2}]}\}<\eps$;
\item $\tau(1-h^2)<\eps$ for every $\tau\in\mathrm T(A)$; and
\item $\tilde h h=h$ and
  $\tilde h f_i^{(1)}=f_i^{(1)}=f_i^{(1)}\tilde h$.
\end{enumerate}
\end{lem}

\begin{proof}
Put $M=\max\{1,\norm{f_1},\ldots,\norm{f_n}\}$ and
$x_i=f_i/M$.  Choose $0<\eta<1$ such that
\[
 3M\eta<\eps\quad\text{and}\quad M\eta<\eps.
\]
For a positive contraction $e$, set
\[
 q(e)=g_{1/4}(e),\qquad h(e)=g_{1/2}(e).
\]
By the standard functional-calculus commutator estimate, there is
$\delta>0$ such that, whenever $e$ is a positive contraction,
$\norm{x}\leq2$, and $\norm{[e,x]}<\delta$, one has
\begin{equation}\label{cut-down-fc-comm}
 \max\big\{
 \norm{[q(e)^{1/2},x]},
 \norm{[(1-q(e))^{1/2},x]},
 \norm{[h(e),x]},
 \norm{[h(e)^{1/2},x]}
 \big\}<\eta.
\end{equation}

Apply $A\in\mathrm{TA}_0(\mathcal S)$ to
$\{x_1,\ldots,x_n\}$ with tolerance smaller than
$\min\{\delta,\eta,\eps\}$.  We obtain $C\subseteq A$ with
$C\in\mathcal S$, a positive contraction $e\in C$ of norm one, and
$x'_1,\ldots,x'_n\in A$ such that
\begin{equation}\label{cut-down-ta-data}
 \norm{x_i-x'_i}<\eta,\qquad
 \norm{[e,x'_i]}<\delta,\qquad
 ex'_ie\in C,\qquad
 \mathrm d_\tau(1-e)<\eps
\end{equation}
for all $i$ and all $\tau\in\mathrm T(A)$.  Since $\norm{x_i}\leq1$ and
$\eta<1$, we also have $\norm{x'_i}<2$.

Set
\[
 q=g_{1/4}(e),\qquad h=g_{1/2}(e),\qquad
 \tilde h=g_{1/8}(e),
\]
and put
\begin{align*}
 f_i^{(0)}&=M(1-q)^{1/2}x'_i(1-q)^{1/2},\\
 f_i^{(1)}&=Mq^{1/2}x'_iq^{1/2}.
\end{align*}
The pointwise relations among the three taper functions give
\begin{equation}\label{nested-taper-relations}
 qh=h,\qquad \tilde hq=q,\qquad \tilde hh=h.
\end{equation}
The distinct cutoffs ensure these exact relations on the transition
regions.

The function
\[
 r(t)=\begin{cases}g_{1/4}(t)^{1/2}/t,&t>0,\\0,&t=0,
       \end{cases}
\]
is continuous, since it vanishes on a neighbourhood of zero.  Hence
\[
 f_i^{(1)}=M r(e)(ex'_ie)r(e)\in C.
\]
Moreover, by \eqref{cut-down-fc-comm},
\begin{align*}
 &\left\|x'_i-\left((1-q)^{1/2}x'_i(1-q)^{1/2}
                  +q^{1/2}x'_iq^{1/2}\right)\right\|\\
 &\hspace{25mm}\leq
 \norm{[(1-q)^{1/2},x'_i]}+\norm{[q^{1/2},x'_i]}
 <2\eta.
\end{align*}
Together with \eqref{cut-down-ta-data}, this yields
\[
 \norm{f_i-(f_i^{(0)}+f_i^{(1)})}<3M\eta<\eps.
\]

Since $q$, $h$, and $\tilde h$ commute, \eqref{nested-taper-relations}
also gives
\[
 h f_i^{(0)}=f_i^{(0)}h=0,\qquad
 \tilde h f_i^{(1)}=f_i^{(1)}=f_i^{(1)}\tilde h.
\]
For $s\in\{1/2,1\}$, again using \eqref{cut-down-fc-comm},
\[
 \norm{[f_i^{(1)},h^s]}
 =M\norm{q^{1/2}[x'_i,h^s]q^{1/2}}
 <M\eta<\eps.
\]
Finally, $\norm h=1$, and in the tracial von Neumann algebra associated to
$\tau$ one has
\[
 1-h^2\leq 1_{[0,1)}(e)=s(1-e).
\]
Consequently,
\[
 \tau(1-h^2)\leq \mathrm d_\tau(1-e)<\eps,
\]
which completes the proof.
\end{proof}

\begin{rem}[Nested taper functions]\label{nested-taper-remark}
Lemma \ref{cut-down} is an abstract version of
\cite[Lemma~4.1]{Niu-MD-Z-absorbing}.  The taper construction in the
proof of the cited lemma requires a distinction between the outer and
inner taper functions.  Using the same coefficients for $h_U$ and $h_V$
introduces their squares in the product and does not give
$h_Uh_V=h_V$.  Separate taper functions $(a_r)$ and $(b_r)$ satisfy the
required identity when
\[
 a_rb_r=b_r
\]
and both vary sufficiently slowly between adjacent tower layers.
For example, for $L\in\mathbb N$ and depths $0\leq r\leq2L+1$, take
\[
 a_r=\min\{(r-1)_+/L,1\},\qquad
 b_r=\min\{(r-L-1)_+/L,1\}.
\]
Both successive increments are at most $1/L$, and $a_r=1$ wherever
$b_r\ne0$.  At depth $2L+1$ both functions are one.  Choosing the
F\o lner shapes for this depth preserves the commutator and trace
estimates and gives the exact identity.  The nested cutoffs
$\tilde h$, $q$, and $h$ in the proof above are the functional-calculus
counterpart of this outer--inner separation.
\end{rem}

\begin{prop}[cf. Proposition 4.7 of \cite{Niu-MD-Z-absorbing}]\label{tracial-Z}
Let $\mathcal S$ be a class such that every $C\in\mathcal S$ can be
locally approximated by \Cstar{}-subalgebras of finite nuclear dimension.  Let
$A$ be a separable, simple, unital \Cstar{}-algebra such that
$A\in\mathrm{TA}_0(\mathcal S)$, and suppose that $A$ is tracially
$\bar m$-almost divisible for some $\bar m\in\{0,1,2,\ldots\}$ in the sense of
\cite[Definition~3.5(ii)]{Winter-Z-stable-02}.  It suffices to assume the
divisibility estimates in that definition only for tracial states.
Then, for every finite set
$\mathcal F\subseteq A$, every $\eps>0$, and every $k\in\mathbb N$, there
is a c.p.c.~order zero map $\phi:\mathrm M_k(\Comp)\to A$ such that
\begin{enumerate}
\item $\norm{[\phi(a),f]}<\eps$ for $f\in\mathcal F$ and
      $a\in \mathrm M_k(\Comp)$ with $\norm a\leq1$;
\item $\mathrm d_\tau(1_A-\phi(1_k))<\eps$ for every $\tau\in \mathrm T(A)$.
\end{enumerate}
\end{prop}

\begin{proof}
We first prove the corresponding statement with the second condition
replaced by
\begin{equation}\label{tracial-Z-trace-target}
 \tau(1_A-\phi(1_k))<\beta,\qquad \tau\in \mathrm T(A),
\end{equation}
for arbitrary $\gamma,\beta>0$, with commutators bounded by $\gamma$.
By decreasing $\gamma$ and $\beta$, we may assume that
$\gamma,\beta<1$.

Write $\mathcal F=\{f_1,\ldots,f_n\}$ and put
$M=\max\{1,\norm{f_1},\ldots,\norm{f_n}\}$.  Choose $\kappa>0$ such that
\begin{equation}\label{tracial-Z-kappa}
 2M\kappa<\frac{\gamma}{8}
 \quad\text{and}\quad
 \kappa<\frac{\beta}{8}.
\end{equation}
By weak stability of order zero maps from $\mathrm M_k(\Comp)$
\cite[Proposition~2.5]{KW-DR-QD}, there is $\delta>0$ such that every
c.p.c.~map $\rho:\mathrm M_k(\Comp)\to E$ into a \Cstar{}-algebra $E$ satisfying
\[
 \norm{\rho(a)\rho(b)}<\delta
\]
for orthogonal positive contractions $a,b\in \mathrm M_k(\Comp)$ is within
$\kappa$ of a c.p.c.~order zero map.

Apply Lemma \ref{cut-down} with a sufficiently small tolerance.  It gives
$C\in\mathcal S$, $h\in C_+$, and decompositions
$f_i\approx f_i^{(0)}+f_i^{(1)}$ with $f_i^{(1)}\in C$ and
$hf_i^{(0)}=f_i^{(0)}h=0$.  Since $C$ is locally approximated by
subalgebras of finite nuclear dimension, there are a subalgebra
$D\subseteq C$ of finite nuclear dimension, a positive contraction
$e\in D$, and $b_i\in D$ such that
\begin{align}
 \norm{f_i-f_i^{(0)}-b_i}&<\frac{\gamma}{16},
       \label{tracial-Z-local-1}\\
 \norm{e^{1/2}f_i^{(0)}}+\norm{f_i^{(0)}e^{1/2}}
      &<\frac{\gamma}{16},\label{tracial-Z-local-2}\\
 \norm{[e^{1/2},b_i]}&<\frac{\gamma}{32},
       \label{tracial-Z-local-3}\\
 \tau(1-e)&<\frac{\beta}{8},\qquad \tau\in \mathrm T(A).
       \label{tracial-Z-local-4}
\end{align}
These inequalities follow by choosing the tolerance in Lemma
\ref{cut-down} sufficiently small and approximating $h$ and the
$f_i^{(1)}$ in a subalgebra of $C$ of finite nuclear dimension.
Continuity of $t\mapsto t^{1/2}$ gives
\eqref{tracial-Z-local-2}--\eqref{tracial-Z-local-3}.

The unital algebra $D+\Comp1_A$ has finite nuclear dimension.  Applying
\cite[Lemma~5.11]{Winter-Z-stable-02} to it, with the finite set
consisting of $e$ and the $b_i$ (rescaled to contractions), gives a
c.p.c.~order zero map $\varphi:\mathrm M_k(\Comp)\to A$ such that, for every
contraction $a\in \mathrm M_k(\Comp)$,
\begin{align}
 \norm{[\varphi(a),e]}&<\delta,\label{tracial-Z-winter-1}\\
 \norm{[\varphi(a),b_i]}&<\frac{\gamma}{8},
       \label{tracial-Z-winter-2}\\
 \tau(1-\varphi(1_k))&<\frac{\beta}{8},\qquad \tau\in \mathrm T(A).
       \label{tracial-Z-winter-3}
\end{align}
In the proof of the cited lemma, divisibility gives order zero maps in
the hereditary subalgebras generated by coordinate positive
contractions, with lower bounds on $\mathrm T(A)$.  The remainder of that
proof uses these bounds only for the induced limit traces
$\mathrm T_\infty(A)$, together with order zero lifting, norm estimates, and
finite nuclear dimension of the approximating subalgebra.  Thus the
proof also applies when the divisibility estimates are assumed only
for tracial states, without identifying quasitraces with traces.
Define
\[
 \rho(a)=e^{1/2}\varphi(a)e^{1/2},\qquad a\in \mathrm M_k(\Comp).
\]
If $a,b$ are orthogonal positive contractions, then
\[
 \rho(a)\rho(b)
 =e^{1/2}\varphi(a)[e,\varphi(b)]e^{1/2},
\]
and hence \eqref{tracial-Z-winter-1} and the choice of $\delta$ give a
c.p.c.~order zero map $\psi:\mathrm M_k(\Comp)\to A$ with
\begin{equation}\label{tracial-Z-order-zero-perturbation}
 \norm{\psi(a)-\rho(a)}<\kappa,\qquad \norm a\leq1.
\end{equation}
For $\norm a\leq1$, equations \eqref{tracial-Z-local-1}--
\eqref{tracial-Z-winter-2} give
\begin{align*}
 \norm{[\rho(a),f_i]}
 &\leq 2\norm{f_i-f_i^{(0)}-b_i}
      +\norm{e^{1/2}f_i^{(0)}}+\norm{f_i^{(0)}e^{1/2}}\\
 &\qquad +2\norm{[e^{1/2},b_i]}
      +\norm{[\varphi(a),b_i]}.
\end{align*}
Together with \eqref{tracial-Z-kappa} and
\eqref{tracial-Z-order-zero-perturbation}, this is strictly less than
$\gamma$.  Moreover, traciality and positivity give
\begin{align*}
 1-\tau(\rho(1_k))
 &=1-\tau(e\varphi(1_k))\\
 &\leq\tau(1-e)+\tau(1-\varphi(1_k))<\frac{\beta}{4}.
\end{align*}
Thus \eqref{tracial-Z-order-zero-perturbation} implies
\eqref{tracial-Z-trace-target}.  This proves the trace version.

For the dimension-function estimate, choose $\eta>0$ such that
\[
 \sqrt\eta<\eps
 \quad\text{and}\quad
 2M\sqrt\eta<\frac{\eps}{2}.
\]
Also require $\eta<1$.  Apply the trace version with centrality
tolerance $\eps/2$ and trace tolerance $\eta$, obtaining a map $\psi$.
Following \cite[Lemma~4.6]{Niu-MD-Z-absorbing}, set
$r=\sqrt\eta$, $v=\psi(1_k)$, and
$f(t)=\min\{t/(1-r),1\}$.  If $\pi_\psi$ is the supporting
*-homomorphism, define
\[
 \phi(x)=\pi_\psi(x)f(v),\qquad x\in \mathrm M_k(\Comp).
\]
This is a c.p.c.~order zero map into $A$, since $f(t)/t$ extends
continuously to $t=0$.  Moreover, $\norm{\phi-\psi}\leq r$.
Writing $\mu_\tau$ for the spectral probability measure of $v$, we have
\[
 \mathrm d_\tau(1-f(v))=\mu_\tau([0,1-r))
 \leq\frac{\tau(1-v)}r<\frac\eta r=r.
\]
Consequently,
\[
 \mathrm d_\tau(1-\phi(1_k))<\sqrt\eta<\eps,\qquad \tau\in \mathrm T(A).
\]
The choice of $\eta$ preserves the required commutator estimate.
\end{proof}

Weak tracial $\mathcal Z$-absorption requires approximately central
order zero maps whose images of the unit have uniformly large trace.

\begin{defn}\label{WTZ-1}
Let $A$ be a unital \Cstar{}-algebra with $\mathrm T(A)\neq\varnothing$.  We say that
$A$ is \emph{weakly tracially $\mathcal Z$-absorbing} if, for every
$k\in\mathbb N$, every finite set $\mathcal F\subseteq A$, and every
$\eps>0$, there is a c.p.c.~order zero map
$\phi:\mathrm M_k(\Comp)\to A$ such that
\begin{enumerate}
\item $\norm{[\phi(a),f]}<\eps$ for $f\in\mathcal F$ and
      $a\in \mathrm M_k(\Comp)$ with $\norm a\leq1$;
\item $\tau(\phi(1_k))>1-\eps$ for every $\tau\in \mathrm T(A)$.
\end{enumerate}
\end{defn}

For hereditary approximation algebras, we use a relative nonunital
version in which traces are normalized at a local unit $c$.

\begin{defn}\label{WTZ-0}
A \Cstar{}-algebra $C$ is \emph{weakly tracially
$\mathcal Z$-absorbing relative to a local unit} if, for every
$k\in\mathbb N$, every $\eps>0$, every finite set $\mathcal F\subseteq C$,
and every norm-one positive contraction $c\in C$ satisfying
\[
 cf=fc=f,\qquad f\in\mathcal F,
\]
there is a c.p.c.~order zero map $\phi:\mathrm M_k(\Comp)\to C$ such that
\begin{enumerate}
\item $\norm{[\phi(a),f]}<\eps$ for $f\in\mathcal F$ and
      $a\in \mathrm M_k(\Comp)$ with $\norm a\leq1$;
\item $\abs{\tau(\phi(1_k))-1}<\eps$ for every bounded positive trace
      $\tau$ on $C$ satisfying $\tau(c)=1$.
\end{enumerate}
\end{defn}

\begin{lem}\label{WTZ-hereditary}
Let $A$ be a simple unital \Cstar{}-algebra that is weakly tracially
$\mathcal Z$-absorbing in the sense of Definition \ref{WTZ-1}.  Then every
nonzero hereditary \Cstar{}-subalgebra of $A$ satisfies Definition
\ref{WTZ-0}.  Consequently, so does every finite direct sum of matrix
algebras over hereditary \Cstar{}-subalgebras of $A$.

In particular, for a simple unital \Cstar{}-algebra with $\mathrm T(A)\ne\varnothing$,
Definitions \ref{WTZ-1}
and \ref{WTZ-0} are equivalent.
\end{lem}

\begin{proof}
Let $H\subseteq A$ be a nonzero hereditary \Cstar{}-subalgebra, let
$\mathcal F\subseteq H$ be finite, and let $c\in H$ be a norm-one positive
contraction satisfying $cf=fc=f$ for $f\in\mathcal F$.  Since $c$ is full
in $A$, there are $z_1,\ldots,z_l\in A$ such that
\[
 b:=\sum_{j=1}^l z_j^*cz_j\geq\frac12 1_A.
\]
Every bounded positive trace $\tau$ on $H$ has the standard lower
semicontinuous extension $\tilde\tau$ to the ideal generated by $H$, which
is all of $A$.  If $\tau(c)=1$, then this extension is bounded and
\begin{equation}\label{hereditary-trace-bound}
 \norm{\tilde\tau}=\tilde\tau(1_A)
 \leq 2\tilde\tau(b)
 \leq 2\sum_{j=1}^l\norm{z_j}^2=:L.
\end{equation}
This bound is uniform over traces normalized at $c$.

Put $M=\max\{1,\norm f:f\in\mathcal F\}$.  Choose $\kappa,\eta>0$
such that
\[
 2M\kappa<\frac\eps3
 \quad\text{and}\quad
 L(\kappa+\eta)<\frac\eps3.
\]
Choose an order zero stability tolerance $\delta>0$ corresponding to
$\kappa$.  Definition \ref{WTZ-1}, applied with a sufficiently small
centrality tolerance to $\mathcal F\cup\{c,c^{1/2}\}$ and with trace
tolerance $\eta$, gives a c.p.c.~order zero map
$\varphi':\mathrm M_k(\Comp)\to A$ such that
\begin{align*}
 \norm{[\varphi'(a),c]}&<\delta,\\
 \norm{[\varphi'(a),f]}+M\norm{[\varphi'(a),c^{1/2}]}
     &<\frac\eps3,
\end{align*}
for $f\in\mathcal F$ and $\norm a\leq1$, and
\[
 \sigma(1-\varphi'(1_k))<\eta,\qquad \sigma\in \mathrm T(A).
\]

The c.p.c.~map
\[
 \rho(a)=c^{1/2}\varphi'(a)c^{1/2}\in\overline{cAc}\subseteq H
\]
is $\delta$-order zero: if $a,b$ are orthogonal positive contractions,
then
\[
 \rho(a)\rho(b)=c^{1/2}\varphi'(a)[c,\varphi'(b)]c^{1/2}.
\]
Weak stability therefore gives a c.p.c.~order zero map
$\phi:\mathrm M_k(\Comp)\to\overline{cAc}$ with
$\norm{\phi-\rho}<\kappa$.  Since $c^{1/2}f=f=fc^{1/2}$,
\[
 \norm{[\rho(a),f]}
 \leq\norm{[\varphi'(a),f]}
      +\norm f\,\norm{[\varphi'(a),c^{1/2}]},
\]
so $\norm{[\phi(a),f]}<\eps$.

Let $\tau$ be a bounded positive trace on $H$ with $\tau(c)=1$.
If $s=\norm{\tilde\tau}$, then $s^{-1}\tilde\tau\in \mathrm T(A)$, and
\eqref{hereditary-trace-bound} gives
\[
 0\leq 1-\tau(\rho(1_k))
 =\tilde\tau\big(c(1-\varphi'(1_k))\big)
 \leq\tilde\tau(1-\varphi'(1_k))<L\eta.
\]
The perturbation from $\rho$ to $\phi$ changes this value by at most
$L\kappa$, proving Definition \ref{WTZ-0} for $H$.

Applying Definition \ref{WTZ-1} simultaneously to all
matrix entries and using $a\mapsto\varphi'(a)\otimes1_r$ shows that
$\mathrm M_r(A)$ satisfies Definition \ref{WTZ-1}.  The first part therefore
applies to every hereditary subalgebra of $\mathrm M_r(A)$, in particular to
$\mathrm M_r(H)$.  For finite direct sums, let $c=(c_j)$ be the prescribed
local unit.  Apply the relative property to
$c_j/\norm{c_j}$ when $c_j\neq0$ and multiply the resulting map by
$\norm{c_j}$.  Homogeneity of the normalized trace estimate and
$\sum_j\tau_j(c_j)=1$ give the required estimate.  If
$\norm{c_j}<1$, the relation $c_jf_j=f_j$ implies $f_j=0$.

Finally, Definition \ref{WTZ-0} with $c=1_A$ implies Definition
\ref{WTZ-1}; the converse is the case $H=A$ just proved.
\end{proof}

\begin{prop}\label{tracial-Z-for-Z}
Let $\mathcal S$ be a class of \Cstar{}-algebras, not necessarily unital, each
of which satisfies Definition \ref{WTZ-0}.  If the unital \Cstar{}-algebra
$A$ belongs to $\mathrm{TA}_0(\mathcal S)$, then, for every finite set
$\mathcal F\subseteq A$, every $\eps>0$, and every $k\in\mathbb N$, there
is a c.p.c.~order zero map $\phi:\mathrm M_k(\Comp)\to A$ such that
\begin{enumerate}
\item $\norm{[\phi(a),f]}<\eps$ for $f\in\mathcal F$ and
      $a\in \mathrm M_k(\Comp)$ with $\norm a\leq1$;
\item $\mathrm d_\tau(1_A-\phi(1_k))<\eps$ for every $\tau\in \mathrm T(A)$.
\end{enumerate}
In particular, $A$ is weakly tracially $\mathcal Z$-absorbing in the
sense of Definition \ref{WTZ-1}.
\end{prop}

\begin{proof}
As in Proposition \ref{tracial-Z}, we first obtain commutators bounded
by $\gamma$ and tracial defects bounded by $\beta$, for arbitrary
$\gamma,\beta>0$.  We may assume that $\gamma,\beta<1$.
Write $\mathcal F=\{f_1,\ldots,f_n\}$,
$M=\max\{1,\norm{f_1},\ldots,\norm{f_n}\}$, and choose $\kappa>0$ with
\[
 2M\kappa<\frac\gamma8
 \quad\text{and}\quad
 \kappa<\frac\beta8.
\]
Let $\delta>0$ be an order zero stability tolerance corresponding to
$\kappa$.

Apply Lemma \ref{cut-down} with a sufficiently small tolerance.  We obtain
$C\in\mathcal S$, positive contractions $h,\tilde h\in C$, and
$f_i^{(0)},f_i^{(1)}\in A$ such that
\begin{align}
 \norm{f_i-f_i^{(0)}-f_i^{(1)}}&<\frac\gamma{16},
       \label{WTZ-local-decomp}\\
 h^{1/2}f_i^{(0)}=f_i^{(0)}h^{1/2}&=0,
       \label{WTZ-local-orth}\\
 \norm{[h^{1/2},f_i^{(1)}]}&<\frac\gamma{32},
       \label{WTZ-local-comm}\\
 \tau(1-h)&<\frac\beta8,\qquad \tau\in \mathrm T(A),
       \label{WTZ-local-trace}
\end{align}
and
\[
 \tilde hh=h=h\tilde h,\qquad
 \tilde h f_i^{(1)}=f_i^{(1)}=f_i^{(1)}\tilde h.
\]
Apply Definition \ref{WTZ-0} in $C$, with local unit $\tilde h$ and with
the finite set $\{h,f_i^{(1)}:1\leq i\leq n\}$.  Decreasing its tolerance
if necessary, we obtain a c.p.c.~order zero map
$\varphi:\mathrm M_k(\Comp)\to C$ such that, for $\norm a\leq1$,
\begin{align}
 \norm{[\varphi(a),h]}&<\delta,\label{WTZ-local-map-1}\\
 \norm{[\varphi(a),f_i^{(1)}]}&<\frac\gamma8,
       \label{WTZ-local-map-2}
\end{align}
and
\begin{equation}\label{WTZ-local-map-3}
 \abs{\sigma(\varphi(1_k))-1}<\frac\beta8
\end{equation}
for every bounded positive trace $\sigma$ on $C$ with
$\sigma(\tilde h)=1$.

Fix $\tau\in \mathrm T(A)$ and set $s=\tau(\tilde h)$.  By
\eqref{WTZ-local-trace}, $s\geq\tau(h)>1-\beta/8$, so $s>0$; hence
$s^{-1}\tau|_C$ is one of the traces in \eqref{WTZ-local-map-3}.  It
follows that
\begin{equation}\label{WTZ-local-map-trace}
 \tau(1-\varphi(1_k))
 \leq\tau(1-\tilde h)+\abs{\tau(\tilde h)-\tau(\varphi(1_k))}
 <\frac\beta4.
\end{equation}

Set $\rho(a)=h^{1/2}\varphi(a)h^{1/2}$.  By
\eqref{WTZ-local-map-1}, this is a $\delta$-order zero c.p.c.~map, and
there is a c.p.c.~order zero map $\psi:\mathrm M_k(\Comp)\to A$ with
$\norm{\psi-\rho}<\kappa$.  Equations \eqref{WTZ-local-decomp}--
\eqref{WTZ-local-map-2} give, for $\norm a\leq1$,
\begin{align*}
 \norm{[\rho(a),f_i]}
 &\leq 2\norm{f_i-f_i^{(0)}-f_i^{(1)}}
      +2\norm{[h^{1/2},f_i^{(1)}]}
      +\norm{[\varphi(a),f_i^{(1)}]}.
\end{align*}
The contribution from $f_i^{(0)}$ vanishes by
\eqref{WTZ-local-orth}.  The displayed bounds and
$\norm{\psi-\rho}<\kappa$ give
$\norm{[\psi(a),f_i]}<\gamma$.  Furthermore,
\[
 1-\tau(\rho(1_k))
 \leq\tau(1-h)+\tau(1-\varphi(1_k))<\frac{3\beta}{8},
\]
so $\tau(1-\psi(1_k))<\beta$.

Finally, apply this trace version with centrality tolerance $\eps/2$ and
with a trace tolerance $\eta>0$ satisfying
$\sqrt\eta<\eps$ and
$2M\sqrt\eta<\eps/2$ and $\eta<1$.  The cutoff in the last paragraph
of the proof of Proposition \ref{tracial-Z} gives a c.p.c.~order zero
map $\phi$ with $\norm{\phi-\psi}\leq\sqrt\eta$ and
$\mathrm d_\tau(1-\phi(1_k))<\sqrt\eta$ for all $\tau\in \mathrm T(A)$.
The choice of $\eta$ gives both conclusions.
\end{proof}

\begin{cor}\label{WRP-WTZ}
Let $A$ be a separable, simple, unital, infinite-dimensional \Cstar{}-algebra,
let $\Gamma$ be a countable discrete amenable group acting on $A$, and put
$B=A\rtimes_{\alpha,r}\Gamma$.  Suppose that the action has WRP and that
$B$ is simple.  If $A$ is weakly tracially $\mathcal Z$-absorbing, then
$B$ is weakly tracially $\mathcal Z$-absorbing.  If, in addition, $B$ has
strict comparison, then $B$ is tracially $\mathcal Z$-absorbing; if $B$
is also nuclear, then $B\cong B\otimes\mathcal Z$.
\end{cor}

\begin{proof}
By Theorem \ref{WTA-C-dyn-1}, $B\in\mathrm{TA}_0(\mathcal S)$, where
$\mathcal S$ consists of finite direct sums of matrix algebras over
hereditary \Cstar{}-subalgebras of $A$.  Lemma \ref{WTZ-hereditary} and
Proposition \ref{tracial-Z-for-Z} show that $B$ is weakly tracially
$\mathcal Z$-absorbing.

Now assume that $B$ has strict comparison, and let $0\neq b\in B_+$.  The
function $\tau\mapsto \mathrm d_\tau(b)$ is strictly positive and lower
semicontinuous on the compact space $\mathrm T(B)$; hence
\[
 \mu:=\min_{\tau\in \mathrm T(B)}\mathrm d_\tau(b)>0.
\]
Apply Proposition \ref{tracial-Z-for-Z} with defect tolerance less than
$\mu$.  Strict comparison gives
$1-\phi(1_k)\precsim b$, so $B$ is tracially
$\mathcal Z$-absorbing.  In the nuclear case,
\cite{HO-Z} then gives $B\cong B\otimes\mathcal Z$.
\end{proof}

\begin{cor}\label{WRP-TAF-AH}
Let $A$ be a separable, simple, unital, infinite-dimensional \Cstar{}-algebra
which has real rank zero and is weakly tracially $\mathcal Z$-absorbing.
Let $\Gamma$ be a countably infinite discrete amenable group acting on
$A$, and put $B=A\rtimes_{\alpha,r}\Gamma$.  Suppose that the action has
WRP and that $B$ is simple, has strict comparison, is nuclear, and
satisfies the UCT.  Then $B$ is tracially AF, equivalently,
$\operatorname{TR}(B)=0$.  In fact, $B$ is isomorphic to a simple unital
AH algebra with slow dimension growth and real rank zero.
\end{cor}

\begin{proof}
Corollary \ref{WRP-WTZ} gives
$B\cong B\otimes\mathcal Z$.  Since WRP requires
$\mathrm T(A)\neq\emptyset$ and $\Gamma$ is amenable, the fixed-point theorem gives
$\mathcal M_1(A,\Gamma)\neq\emptyset$.  Proposition
\ref{trace-factor} identifies this simplex with $\mathrm T(B)$.  Since $B$ is
simple, every trace on $B$ is faithful; in particular, $B$ is stably
finite.

Since $\Gamma$ is infinite and $A$ has real rank zero,
Corollary \ref{wrp-rr0-affine} gives
\[
 \overline{\rho_B(\mathrm K_0(B))}=\operatorname{Aff}(\mathrm T(B)).
\]
The algebra $B$ is exact because it is nuclear.  Therefore
\cite[Theorem~7.2]{Ror-Z-stable}, applied to the simple, unital, finite,
$\mathcal Z$-stable algebra $B$, shows that
$\operatorname{RR}(B)=0$.

Finiteness and $\mathcal Z$-stability give stable rank one and weak
unperforation, while real rank zero gives the refinement property for
projections.  Thus $(\mathrm K_0(B),\mathrm K_0(B)_+,[1_B])$ is a countable simple
weakly unperforated ordered Riesz group with order unit.  The range
theorem allows torsion in $\mathrm K_0(B)$ and does not require this group to
be a dimension group.
Moreover, the canonical map
\[
 \mathrm T(B)\longrightarrow \mathrm S(\mathrm K_0(B),[1_B])
\]
is an affine homeomorphism.  Indeed, injectivity follows from the density
of $\rho_B(\mathrm K_0(B))$.  Conversely, a state on the ordered $\mathrm K_0$-group
extends to a quasitrace by the state-extension theorem of
\cite{BlaTrace}, and this quasitrace is a trace because $B$ is exact
\cite{Haagtrace}.

This homeomorphism preserves extreme points.  Moreover,
density of $\rho_B(\mathrm K_0(B))$ excludes
$\mathrm K_0(B)/\operatorname{Inf}(\mathrm K_0(B))\cong\mathbb Z$, where
$\operatorname{Inf}(\mathrm K_0(B))=\ker\rho_B$: a cyclic pairing range
containing the constant function one is discrete.  Thus the hypotheses
of the range theorem for simple AH algebras with slow dimension growth
\cite{Vill-EHS} are satisfied; see also \cite{L-N-Range}.
There is therefore a simple unital AH algebra $C$ with slow
dimension growth and an isomorphism
\[
 \operatorname{Ell}(C)\cong\operatorname{Ell}(B).
\]
The density of the $\mathrm K_0$-range implies that $C$ is infinite-dimensional.
The algebra $C$ is therefore $\mathcal Z$-stable by \cite{Toms-SDG};
it is also nuclear and satisfies the UCT.  By \cite{CETWW-dim-n},
both $B$ and $C$ have finite nuclear dimension.  By
\cite[Theorem~4.9]{EGLN-DR} and quasidiagonality of traces
\cite{TWW-QD}, both stably finite algebras belong to the class
$\mathcal N_1$.  The classification theorem for this class
\cite{GLN-TAS-1,GLN-TAS-2} therefore yields $B\cong C$.
Consequently, $C$ has real rank zero.
Finally, for a simple unital AH algebra whose $\mathrm K_0$-range is dense,
$\mathcal Z$-stability is equivalent to being tracially AF by
\cite[Theorem~3.11]{TW1}.  Thus $B$ is tracially AF and is isomorphic to an AH algebra
with slow dimension growth and real rank zero.
\end{proof}

\begin{thm}\label{TW-conj}
Let $\mathcal S$ be a class such that every $C\in\mathcal S$ can be
locally approximated by \Cstar{}-subalgebras of finite nuclear dimension.  Let
$A$ be a separable, simple, unital, infinite-dimensional, nuclear
\Cstar{}-algebra with $A\in\mathrm{TA}_0(\mathcal S)$.  Then
$A\cong A\otimes\mathcal Z$ if and only if
\begin{enumerate}
\item $A$ has strict comparison of positive elements; and
\item for every $f\in\mathrm{LAff}_{++}(\mathrm T(A))$ there is
      $x\in(A\otimes\mathcal K)_+$ such that
      $\mathrm d_\tau(x)=f(\tau)$ for every $\tau\in \mathrm T(A)$.
\end{enumerate}
Consequently, within this class, $\mathcal Z$-stability is also equivalent
to the canonical ordered-semigroup description
\[
 \mathrm{Cu}(A)\cong
 \mathrm V(A)\sqcup\mathrm{LAff}_{++}(\mathrm T(A)).
\]
\end{thm}

\begin{proof}
If $A$ is $\mathcal Z$-stable, then strict comparison and the stated
description of the Cuntz semigroup are standard; see
\cite{Ror-Z-stable,BPT-Cuntz,CEI-CuntzSG}.  We prove the converse.

We show that $A$ is tracially $0$-almost divisible.  Let
$a\in \mathrm M_N(A)_+$ be a positive contraction, let $k\in\mathbb N$, and let
$\eps>0$.  The case $a=0$ is immediate; assume that $a\neq0$.
Choose
\[
 0<\theta<\min\left\{\frac12,\frac{\eps}{2N}\right\}.
\]
The function
\[
 f(\tau)=\frac{1-\theta}{k}\mathrm d_\tau(a),\qquad \tau\in \mathrm T(A),
\]
belongs to $\mathrm{LAff}_{++}(\mathrm T(A))$.  By (2), there is a positive
contraction $x\in A\otimes\mathcal K$ such that
\begin{equation}\label{TW-rank-choice}
 \mathrm d_\tau(x)=\frac{1-\theta}{k}\mathrm d_\tau(a),\qquad \tau\in \mathrm T(A).
\end{equation}
Since $A$ is nuclear, every normalized quasitrace is a trace.  Thus
$k \mathrm d_\tau(x)<\mathrm d_\tau(a)$ for every normalized quasitrace, and strict
comparison gives
\begin{equation}\label{TW-cuntz-comparison}
 k[x]\leq[a]\quad\text{in }\mathrm{Cu}(A).
\end{equation}

For $0<\delta<1/2$, let $g_{\delta,2\delta}$ be zero on
$[0,\delta]$, one on $[2\delta,1]$, and linear in between.  As
$\delta\searrow0$,
\[
 k\tau(g_{\delta,2\delta}(x))\nearrow k \mathrm d_\tau(x).
\]
Since $\mathrm d_\tau(a)\geq\tau(a)$ and
$\mathrm d_\tau(a)\leq N$, equation \eqref{TW-rank-choice} gives
\[
 k \mathrm d_\tau(x)=(1-\theta)\mathrm d_\tau(a)
 >\tau(a)-\frac\eps2.
\]
For each fixed $\delta$, the element $g_{\delta,2\delta}(x)$ belongs to
the Pedersen ideal of $A\otimes\mathcal K$, so the function
$\tau\mapsto\tau(g_{\delta,2\delta}(x))$ is continuous.  Compactness of
$\mathrm T(A)$ and monotonicity as $\delta\searrow0$ therefore give one
$\delta>0$ such that
\begin{equation}\label{TW-uniform-cutoff}
 k\tau(g_{\delta,2\delta}(x))>\tau(a)-\eps,
 \qquad \tau\in \mathrm T(A).
\end{equation}

Apply \cite[Proposition~2.12]{Winter-Z-stable-02} in the ambient algebra
$A\otimes\mathcal K$ to \eqref{TW-cuntz-comparison}, with the function
$g_{\delta,2\delta}$.  Its hereditary target is
$\operatorname{her}(a)\subseteq\mathrm M_N(A)$.  We obtain a c.p.c.~order zero map
\[
 \psi:\mathrm M_k(\Comp)\longrightarrow\operatorname{her}(a)
\]
such that $\psi(e_{11})$ is Murray--von Neumann equivalent to
$g_{\delta,2\delta}(x)$.  Since a trace composed with an order zero map is
a trace,
\[
 \tau(\psi(1_k))
 =k\tau(\psi(e_{11}))
 =k\tau(g_{\delta,2\delta}(x))
 >\tau(a)-\eps
\]
for every $\tau\in \mathrm T(A)$.  This proves that $A$ is tracially
$0$-almost divisible.

Proposition \ref{tracial-Z}, with $\bar m=0$, gives approximately
central order zero maps with uniformly small dimension functions on
their complements.  Given $0\neq b\in A_+$, simplicity and compactness
of $\mathrm T(A)$ give
$\min_{\tau\in \mathrm T(A)}\mathrm d_\tau(b)>0$.  Choosing the defect smaller than this
minimum and using strict comparison yields
$1-\phi(1_k)\precsim b$.  Thus $A$ is tracially
$\mathcal Z$-absorbing, and nuclearity gives
$A\cong A\otimes\mathcal Z$ by \cite{HO-Z}.

Finally, the canonical Cuntz-semigroup description implies (1) and (2),
while $\mathcal Z$-stability gives that description by the results cited
in the first paragraph.
\end{proof}

Since locally finite nuclear dimension passes to hereditary subalgebras,
matrix algebras, and finite direct sums, Theorems \ref{WTA-C-dyn-1} and
\ref{TW-conj} give the following consequence.

\begin{thm}\label{wrp-local-nucdim-Cu-Z}
Let $A$ be a separable, unital \Cstar{}-algebra with $\mathrm T(A)\neq\varnothing$
and locally finite nuclear dimension, and let $\Gamma$ be a countably
infinite discrete amenable group acting on $A$.  Suppose that the action
has WRP and that $B=A\rtimes_{\alpha,r}\Gamma$ is simple.  Then $B$ is
$\mathcal Z$-stable if and only if
\[
 \mathrm{Cu}(B)
 \cong\mathrm V(B)\sqcup\mathrm{LAff}_{++}(\mathrm T(B))
\]
as ordered semigroups, via the canonical identification.
\end{thm}

\begin{proof}
Locally finite nuclear dimension passes to hereditary subalgebras as
follows.  Let $H\subseteq A$ be hereditary and let
$\mathcal F\subseteq H$ be finite.
Choose a positive contraction $h\in H$ such that
$h^{1/2}xh^{1/2}$ is close to $x$ for $x\in\mathcal F$.
For $0<\eta<\delta<1$, local approximation gives a subalgebra
$D\subseteq A$ of finite nuclear dimension, a positive contraction
$d\in D$ with $\norm{d-h}<\eta$, and $x_D\in D$ close to each
$x\in\mathcal F$.  Put
\[
 q=(d-\delta)_+,\qquad
 f(t)=\begin{cases}\sqrt{(t-\delta)_+/t},&t>0,\\0,&t=0,\end{cases}
 \qquad z=h^{1/2}f(d),\qquad a_0=z^*z.
\]
Since $q=f(d)df(d)$, we have
\[
 (1-\eta/\delta)q\leq a_0\leq(1+\eta/\delta)q,
 \qquad \operatorname{her}_A(a_0)=\operatorname{her}_A(q).
\]
Let $z=va_0^{1/2}$ be its polar decomposition in $A^{**}$.
Writing $c=1-\eta/\delta>0$, the inequality $cq\leq a_0$ gives
\[
 \norm{\bigl(v-z(a_0+1/n)^{-1/2}\bigr)q^{1/2}}^2
 \leq \frac{1}{cn},\qquad n\in\mathbb N.
\]
Thus $vq^{1/2}\in A$ by norm approximation.  In particular,
polar conjugation embeds $E=\operatorname{her}_D(q)$ into
$\operatorname{her}_A(zz^*)\subseteq H$.  Moreover,
\[
 \norm{vq^{1/2}-q^{1/2}}
 \leq\norm{a_0^{1/2}-q^{1/2}}+\norm{z-q^{1/2}}
 \leq2\sqrt\eta.
\]
Consequently, $v(q^{1/2}x_Dq^{1/2})v^*$ approximates $x$ if
$\norm{h^{1/2}xh^{1/2}-x}$, $\delta$, and then $\eta$ are sufficiently
small, since $\norm{q-h}\leq\delta+\eta$.
The algebra $E$ has finite nuclear dimension by
\cite[Proposition~2.5]{WZ-ndim}, so its image in $H$ is the required
local approximating subalgebra.  Matrix amplifications and finite
direct sums preserve this local property by entrywise approximation.

Theorem \ref{WTA-C-dyn-1} gives $B\in\mathrm{TA}_0(\mathcal S)$, where
$\mathcal S$ consists of finite direct sums of matrix algebras over
hereditary \Cstar{}-subalgebras of $A$.  Every member of $\mathcal S$ locally
has finite nuclear dimension.  Moreover, $A$ is nuclear, and amenability
of $\Gamma$ implies that $B$ is nuclear.  Averaging a tracial state on
$A$ gives an invariant trace, which extends to $B$ and is faithful by
simplicity.  The canonical group unitaries make $B$ infinite dimensional.
Theorem~\ref{TW-conj} gives the conclusion.
\end{proof}

\section{Realization of strictly positive affine rank functions}

\begin{defn}\label{defn-e1-e2}
Let $A$ be a unital \Cstar{}-algebra with $\mathrm T(A)\ne\varnothing$.
We say that $A$ has \emph{property (E1)} if, for every $\eps>0$ and
every $\alpha\in\operatorname{Aff}(\mathrm T(A))$ with $\alpha>0$, there is
$a\in \mathrm M_\infty(A)_+$ such that
\[
 \abs{\mathrm d_\tau(a)-\alpha(\tau)}<\eps,\qquad \tau\in \mathrm T(A).
\]

We say that $A$ has \emph{property (E2)} if the following holds.  Let
$(\alpha_n)_{n=1}^\infty$ be a sequence in
$\operatorname{Aff}(\mathrm T(A))$ such that
\[
 0<\alpha_1<\alpha_2<\cdots
\]
pointwise on $\mathrm T(A)$.  Put
\[
 \delta_n=\min_{\tau\in \mathrm T(A)}
       \big(\alpha_{n+1}(\tau)-\alpha_n(\tau)\big)>0
\]
and
\[
 \eta_1=\frac{\delta_1}{4},\qquad
 \eta_n=\frac14\min\{\delta_{n-1},\delta_n\},\quad n\geq2.
\]
Then there are $a_n\in \mathrm M_\infty(A)_+$ such that
\[
 \abs{\mathrm d_\tau(a_n)-\alpha_n(\tau)}<\eta_n,\qquad
 \tau\in \mathrm T(A),\ n\in\mathbb N,
\]
and
\[
 a_1\precsim a_2\precsim\cdots.
\]
\end{defn}

\begin{rem}
Suppose that $A$ has property (E1) and strict comparison by traces.
Apply (E1) to $\alpha_n$ with tolerance $\eta_n$.  The resulting
elements satisfy, for every $\tau\in \mathrm T(A)$,
\[
 \mathrm d_\tau(a_n)
 <\alpha_n(\tau)+\frac{\delta_n}{4}
 <\alpha_{n+1}(\tau)-\frac{\delta_n}{4}
 <\mathrm d_\tau(a_{n+1}).
\]
Strict comparison gives $a_n\precsim a_{n+1}$, so $A$ has property
(E2).  If strict comparison is formulated using normalized
2-quasitraces, the same conclusion holds whenever every normalized
2-quasitrace on $A$ is a trace; this is automatic when $A$ is exact.
\end{rem}

\begin{lem}\label{K-surj}
Let $A$ be a separable unital \Cstar{}-algebra satisfying property (E2).
Then, for every $s\in\mathrm{LAff}_{++}(\mathrm T(A))$, there is
$x\in(A\otimes\mathcal K)_+$ such that
\[
 \mathrm d_\tau(x)=s(\tau),\qquad \tau\in \mathrm T(A).
\]
\end{lem}
\begin{proof}
Fix $s\in\mathrm{LAff}_{++}(\mathrm T(A))$.  By the standard description of
lower semicontinuous affine functions on the metrizable compact simplex
$\mathrm T(A)$, there are strictly positive continuous affine functions
\[
 \alpha_1<\alpha_2<\cdots,\qquad
 \sup_n\alpha_n(\tau)=s(\tau),\quad \tau\in \mathrm T(A).
\]
Choose $(a_n)$ as in property (E2).  Thus
$[a_1]\leq[a_2]\leq\cdots$ in $\mathrm{Cu}(A)$.  By
\cite[Theorem~1(i) and Appendix~6]{CEI-CuntzSG}, this sequence has a
supremum in $\mathrm{Cu}(A)$, represented by some
$x\in(A\otimes\mathcal K)_+$:
\[
 [x]=\sup_n[a_n].
\]

For a fixed $\tau$, if $s(\tau)<\infty$, then
\[
 \eta_n\leq\frac14
   \big(\alpha_n(\tau)-\alpha_{n-1}(\tau)\big)\longrightarrow0.
\]
If $s(\tau)=\infty$, then, for $n\geq2$,
$\eta_n<\alpha_n(\tau)/4$, and hence
$\mathrm d_\tau(a_n)>3\alpha_n(\tau)/4\longrightarrow\infty$.
Consequently,
\[
 \sup_n \mathrm d_\tau(a_n)=s(\tau)
\]
in either case.  Since the functional on $\mathrm{Cu}(A)$ induced by
$\tau$ preserves suprema of increasing sequences,
\[
 \mathrm d_\tau(x)=\sup_n \mathrm d_\tau(a_n)=s(\tau),\qquad \tau\in \mathrm T(A).
\]
\end{proof}

\begin{prop}\label{almost-surj}
For $\eps\geq0$, let $\mathcal S_\eps$ be the class of \Cstar{}-algebras $C$
with the following property.  For every $c\in C_+$ for which there is
$e\in C_+$ satisfying $ec=c$, there is $b\in \mathrm M_\infty(C)_+$ such that
\[
 \abs{\sigma(c)-\mathrm d_\sigma(b)}\leq\eps,\qquad \sigma\in \mathrm T(C).
\]
Such an element $c$ is called \emph{compactly supported in $C$}.
Thus $C\in\mathcal S_0$ means that every compactly supported positive
element of $C$ has its trace function realized exactly by a rank
function.  In particular,
\[
 \mathcal S_0\subseteq\mathcal S_\eps\quad(\eps>0),
 \qquad
 \mathcal S_\eps\subseteq\mathcal S_\delta
 \quad(0\leq\eps\leq\delta).
\]

Let $A$ be a simple unital \Cstar{}-algebra.  If
\[
 A\in\mathrm{TA}_0(\mathcal S_\eps)
 \quad\text{for every }\eps>0,
\]
then $A$ has property (E1).  Moreover, if $D\subseteq A$ is a
\Cstar{}-subalgebra and the rank-approximating elements $b$ from all local
models can be chosen in $\mathrm M_\infty(D)$, then the elements realizing
property (E1) can be chosen in $\mathrm M_\infty(D)$.
\end{prop}

\begin{proof}
Fix $\alpha\in\operatorname{Aff}(\mathrm T(A))$ with $\alpha>0$ and
$\eps>0$.  By \cite[Corollary~3.10]{BPT-Cuntz}, simplicity of $A$
gives $a'\in A_+$ such that
\[
 \tau(a')=\alpha(\tau),\qquad \tau\in \mathrm T(A).
\]
Put $M=\max\{1,\norm{a'}\}$ and choose $\gamma>0$ such that
$(2M+4)\gamma<\eps$.

Apply Lemma \ref{alt-defn-was} to $a'$ with tolerance $\gamma$.  We
obtain $C\in\mathcal S_\gamma$ and a positive contraction $h\in C$ of
norm one such that
\[
 \operatorname{dist}(ha'h,C)<\gamma,\qquad
 \mathrm d_\tau(1-h)<\gamma,\quad \tau\in \mathrm T(A).
\]
There is $c_0\in C_+$ with
$\norm{ha'h-c_0}<2\gamma$.  Set $c=(c_0-\gamma)_+$.  Then
$\norm{c-c_0}\leq\gamma$, and functional calculus gives
$e\in C_+$ with $ec=c$; hence $c$ is compactly supported in $C$.
For $\tau\in \mathrm T(A)$,
\begin{align*}
 0\leq\tau(a')-\tau(ha'h)
 &=\tau\big(a'^{1/2}(1-h^2)a'^{1/2}\big)\\
 &\leq M\tau(1-h^2)
 \leq2M\tau(1-h)<2M\gamma.
\end{align*}
It follows that
\begin{equation}\label{local-affine-trace-approx}
 \abs{\tau(c)-\alpha(\tau)}<(2M+3)\gamma,\qquad \tau\in \mathrm T(A).
\end{equation}

Since $C\in\mathcal S_\gamma$, there is $b\in \mathrm M_\infty(C)_+$ whose
rank function approximates the trace function of $c$.  For
$\tau\in \mathrm T(A)$, put $t=\norm{\tau|_C}\leq1$.  If $t>0$, applying the
defining estimate to $t^{-1}\tau|_C$ and multiplying by $t$ gives
\[
 \abs{\mathrm d_\tau(b)-\tau(c)}\leq\gamma,\qquad \tau\in \mathrm T(A).
\]
The same estimate holds if $t=0$, since both terms vanish.
Together with \eqref{local-affine-trace-approx}, this gives
\[
 \abs{\mathrm d_\tau(b)-\alpha(\tau)}<(2M+4)\gamma<\eps,
 \qquad \tau\in \mathrm T(A).
\]
This proves property (E1).  The final assertion follows by taking the
same element $b$ in $\mathrm M_\infty(D)$.
\end{proof}

\begin{cor}\label{diag-rank-wrp}
Let a countably infinite discrete amenable group $\Gamma$ act freely and
minimally on a compact metrizable space $X$, and suppose that the action
has WRP.  Put $B=\mathrm C(X)\rtimes_{\mathrm r}\Gamma$.  Then, for every
strictly positive continuous affine function $\alpha$ on $\mathrm T(B)$ and
every $\eps>0$,
there are $m\in\mathbb N$ and $f_1,\ldots,f_m\in \mathrm C(X)_+$ such that,
with $a=\operatorname{diag}(f_1,\ldots,f_m)$,
\[
 \abs{\mathrm d_\tau(a)-\alpha(\tau)}<\eps,\qquad \tau\in \mathrm T(B).
\]
\end{cor}
\begin{proof}
Freeness and minimality imply that $B$ is simple.  Fix $\gamma>0$.
In the proof of Theorem \ref{WTA-C-dyn-1}, the tiling shapes may be
chosen so that their cardinalities $n_s$ are all greater than any
prescribed integer.  The resulting local model has the form
\[
 C=\bigoplus_{s=1}^S \mathrm M_{n_s}(H_s),
\]
where every $H_s$ is a hereditary subalgebra of $\mathrm C(X)$ and hence is
isomorphic to $\mathrm C_0(U_s)$ for an open set $U_s\subseteq X$.  Choose the
shapes so large that $1/n_s<\gamma$ for every $s$.  Lemma
\ref{small-exi} and homogeneity of its estimate for bounded traces show
that $C\in\mathcal S_\gamma$; the same estimate passes to finite direct
sums by decomposing a trace over the summands.  Consequently,
\[
 B\in\mathrm{TA}_0(\mathcal S_\gamma)
 \quad\text{for every }\gamma>0.
\]

The construction in Lemmas \ref{rank-f-ext} and \ref{small-exi} uses
diagonal scalar-valued functions.  Under the tower embeddings in the
proof of Theorem \ref{WTA-C-dyn-1}, all diagonal corners of each local
model lie in $\mathrm C(X)$.  Thus the rank-approximating elements required in
Proposition \ref{almost-surj} can be chosen as finite diagonal elements
of $\mathrm M_\infty(\mathrm C(X))$.  That proposition, applied with $D=\mathrm C(X)$, proves
the assertion.
\end{proof}

\begin{lem}\label{rank-f-ext}
Let $X$ be a locally compact metrizable space, and let
$\alpha:X\to\mathbb N_0$ be a bounded lower semicontinuous function
whose support is compact.  Then there are $m\in\mathbb N$ and a positive
diagonal element $f\in \mathrm M_m(\mathrm C_0(X))$ such that
\begin{equation}\label{rank-equ}
 \operatorname{rank}(f(x))=\alpha(x),\qquad x\in X.
\end{equation}
\end{lem}
\begin{proof}
Let $r=\max_{x\in X}\alpha(x)$.  If $r=0$, take
$f=0\in \mathrm C_0(X)=\mathrm M_1(\mathrm C_0(X))$.
Otherwise, for $1\leq k\leq r$, put
\[
 U_k=\{x\in X:\alpha(x)\geq k\}.
\]
Since $\alpha$ is integer-valued and lower semicontinuous, each $U_k$
is open.  Moreover, $U_k$ has compact closure because it is contained
in the compact support of $\alpha$.  As $X$ is metrizable, there is
$\varphi_k\in \mathrm C_0(X)_+$ whose cozero set is exactly $U_k$.  Set
\[
 f=\operatorname{diag}(\varphi_1,\ldots,\varphi_r).
\]
At a point $x$, precisely $\alpha(x)$ of the functions $\varphi_k(x)$
are nonzero, which proves \eqref{rank-equ}.
\end{proof}

\begin{lem}\label{small-exi}
Let $C=\mathrm M_n(\mathrm C_0(X))$, where $X$ is a locally compact metrizable space.
Then $C\in\mathcal S_{1/n}$.  In fact, if $c\in C_+$ is compactly
supported and $\sigma$ is a nonzero bounded positive trace on $C$, the element
$b$ below can be chosen diagonal and satisfies
\[
 \abs{\mathrm d_\sigma(b)-\sigma(c)}<\frac{\norm\sigma}{n}.
\]
For the zero trace, both quantities are zero.
\end{lem}

\begin{proof}
Let $c\in C_+$ be compactly supported, and choose $e\in C_+$ with
$ec=c$.  Define
\[
 q(x)=\left\lceil\operatorname{Tr}_n(c(x))\right\rceil,
 \qquad x\in X.
\]
The function $q$ is bounded, integer-valued, and lower semicontinuous.
It has compact support: if $c(x)\ne0$, then $e(x)c(x)=c(x)$ forces
$\norm{e(x)}\geq1$, and the set
$\{x:\norm{e(x)}\geq1\}$ is compact.  Lemma \ref{rank-f-ext} therefore
gives a positive diagonal matrix over $\mathrm C_0(X)$ whose rank at $x$ is
$q(x)$.  After padding by zero entries to make its matrix size a
multiple of $n$, regard it as an element $b\in \mathrm M_\infty(C)_+$.

First suppose that $\sigma\in \mathrm T(C)$.  There is a Borel probability
measure $\mu$ on $X$ such that
\[
 \sigma(c)=\frac1n\int_X\operatorname{Tr}_n(c(x))\,d\mu(x),
 \qquad
 \mathrm d_\sigma(b)=\frac1n\int_X\operatorname{rank}(b(x))\,d\mu(x).
\]
Since $0\leq\lceil t\rceil-t<1$ for $t\geq0$, we obtain
\[
 0\leq \mathrm d_\sigma(b)-\sigma(c)
 =\frac1n\int_X
   \big(q(x)-\operatorname{Tr}_n(c(x))\big)\,d\mu(x)
 <\frac1n.
\]
The last inequality follows because the integrand is strictly less
than one at every point.  Normalization and homogeneity give the
estimate for every nonzero bounded positive trace.
\end{proof}

\begin{thm}\label{homogeneous-TA-Z}
Let
\[
 \mathcal H=\big\{\mathrm M_n(\mathrm C_0(X)):n\in\mathbb N,
       \ X\text{ locally compact and metrizable}\big\}
\]
be the class of basic homogeneous \Cstar{}-algebras.  Let $A$ be a separable,
simple, unital, nuclear, non-type-I \Cstar{}-algebra.  If
$A\in\mathrm{TA}_0(\mathcal H)$, then
\[
 A\cong A\otimes\mathcal Z
 \quad\Longleftrightarrow\quad
 A\text{ has strict comparison of positive elements}.
\]
\end{thm}
\begin{proof}
The forward implication is standard: $A$ is finite because it has a
faithful trace, and a simple, exact, finite, $\mathcal Z$-stable
\Cstar{}-algebra has strict comparison by \cite{Ror-Z-stable}.  We prove the
converse.

For $N\in\mathbb N$, put
\[
 \mathcal H_{\geq N}
 =\{\mathrm M_n(\mathrm C_0(X))\in\mathcal H:n\geq N\}.
\]
We claim that
\begin{equation}\label{homogeneous-large-multiplicity}
 A\in\mathrm{TA}_0(\mathcal H_{\geq N})
 \quad\text{for every }N\in\mathbb N.
\end{equation}
Suppose otherwise.  For some $N\geq2$, there are a finite set
$\mathcal F_0\subseteq A$ and $\eps_0>0$ for which the defining
$\mathrm{TA}_0$ approximation cannot be obtained with an algebra in
$\mathcal H_{\geq N}$.  Let $(a_j)_{j=1}^\infty$ be dense in the unit
ball of $A$.  Applying $A\in\mathrm{TA}_0(\mathcal H)$ to
\[
 \mathcal F_j=\mathcal F_0\cup\{a_1,\ldots,a_j\}
 \quad\text{with tolerance}\quad
 \delta_j<\min\{\eps_0,1/j\},
\]
gives
\[
 C_j=\mathrm M_{n_j}(\mathrm C_0(X_j))\subseteq A,\qquad h_j\in(C_j)_+,
 \qquad \norm{h_j}=1.
\]
The choice of $(\mathcal F_0,\eps_0)$ implies that $n_j<N$ for every
$j$.  For $j\geq k$, let $a'_{k,j}$ be the approximant to $a_k$ in
the defining $\mathrm{TA}_0$ approximation.  Then
\[
 \norm{[h_j,a_k]}
 \leq \norm{[h_j,a'_{k,j}]}+2\norm{a_k-a'_{k,j}}
 <3\delta_j
\]
and, since $h_ja'_{k,j}h_j\in C_j$,
\[
 \operatorname{dist}(h_ja_kh_j,C_j)
 \leq\norm{h_j(a_k-a'_{k,j})h_j}<\delta_j.
\]
Density of $(a_k)_k$ now gives
\begin{equation}\label{homogeneous-central-cut}
 \begin{aligned}
  \norm{[h_j,a]}&\longrightarrow0,\\
  \operatorname{dist}(h_jah_j,C_j)&\longrightarrow0
 \end{aligned}
 \qquad (a\in A).
\end{equation}

Fix a free ultrafilter $\omega$, let $q_\omega:\ell^\infty(A)\to
A_\omega$ be the quotient map, and set
\[
 h=q_\omega((h_j)_j),\qquad
 Q=q_\omega\left(\prod_{j=1}^\infty C_j\right)\subseteq A_\omega.
\]
Equation \eqref{homogeneous-central-cut} shows that $h\in A'\cap
A_\omega$ and $hAh\subseteq Q$.  Hence
\[
 \Phi:A\longrightarrow Q,\qquad \Phi(a)=hah,
\]
is a c.p.c.~order zero map.  It is nonzero because
$\Phi(1_A)=h^2$ and $\norm h=1$.  The kernel of an order zero map is an
ideal, so simplicity of $A$ implies that $\Phi$ is injective.

Every $C_j$ is $(N-1)$-subhomogeneous.  Equivalently, the standard
polynomial of degree $2(N-1)$ vanishes on $C_j$.  The same polynomial
therefore vanishes on $Q$ and on
$D=\mathrm C^*(\Phi(A))$.  By the structure theorem for order zero maps
\cite{WZ-order-zero}, there is a supporting *-homomorphism
\[
 \pi_\Phi:A\longrightarrow \operatorname{Mult}(D)\subseteq D^{**}
 \quad\text{such that}\quad
 \Phi(a)=\pi_\Phi(a)\Phi(1_A).
\]
It is faithful because $\Phi$ is faithful.  Since the standard
polynomial is multilinear, its identity passes to $D^{**}$ by separate
weak* continuity.  Thus $\pi_\Phi(A)$ satisfies the same standard identity.
Kaplansky's characterization \cite{Kap-SCalg} implies that $A$ is
$(N-1)$-subhomogeneous, hence type I.  This contradicts the hypothesis
and proves \eqref{homogeneous-large-multiplicity}.

Let $\eps>0$ and choose $N$ with $1/N<\eps$.  Lemma
\ref{small-exi} and \eqref{homogeneous-large-multiplicity} give
\[
 A\in\mathrm{TA}_0(\mathcal H_{\geq N})
 \subseteq\mathrm{TA}_0(\mathcal S_\eps).
\]
Since this holds for every $\eps>0$, Proposition \ref{almost-surj}
gives property (E1).  Nuclearity implies exactness, so strict comparison
and the remark following Definition \ref{defn-e1-e2} give property
(E2).  Lemma \ref{K-surj} then shows that every function in
$\mathrm{LAff}_{++}(\mathrm T(A))$ is realized by the rank function of a
positive element of $A\otimes\mathcal K$.

Finally, every basic homogeneous algebra locally has finite nuclear
dimension.  Indeed, a finite subset of $\mathrm M_n(\mathrm C_0(X))$ is contained in
$\mathrm M_n(D)$, where $D\subseteq \mathrm C_0(X)$ is generated by finitely many
matrix entries.  The spectrum of the finitely generated commutative
algebra $D$ embeds in a finite-dimensional Euclidean space, so $D$, and
hence $\mathrm M_n(D)$, has finite nuclear dimension.  Theorem \ref{TW-conj}
now gives $A\cong A\otimes\mathcal Z$.
\end{proof}

\begin{rem}\label{homogeneous-tracial-Z}
The proof of the reverse implication in Theorem
\ref{homogeneous-TA-Z} gives tracial $\mathcal Z$-absorption without
nuclearity if strict comparison is assumed with respect to tracial
states.  Properties (E1) and (E2) and the rank-realization argument
remain valid.  The proof of Theorem
\ref{TW-conj} then gives the tracial $0$-almost-divisibility estimates
for tracial states.  The version of Winter's centralization lemma for
tracial states used in Proposition \ref{tracial-Z} gives approximately
central order zero maps.  Nuclearity is needed only to pass from
tracial $\mathcal Z$-absorption to tensorial
$\mathcal Z$-stability.
\end{rem}

\begin{lem}\label{Z-S0}
Let $A$ be a separable, simple, unital, exact, $\mathcal Z$-stable
\Cstar{}-algebra with $\mathrm T(A)\ne\emptyset$.  Then
\[
 \overline{cAc}\in\mathcal S_0,\qquad c\in A_+.
\]
\end{lem}
\begin{proof}
Put $H=\overline{cAc}$.  The assertion is immediate if $c=0$.  Let
$0\neq a\in H_+$ be compactly supported in $H$; the case $a=0$ is
handled by taking $b=0$.  The function
\[
 f:\mathrm T(A)\longrightarrow(0,\infty),\qquad f(\tau)=\tau(a),
\]
is continuous and strictly positive.  The Cuntz-semigroup description
for simple, unital, exact, finite, $\mathcal Z$-stable algebras
\cite[Theorem~5.5]{BPT-Cuntz} gives
$x\in(A\otimes\mathcal K)_+$ such that
\begin{equation}\label{Z-S0-rank}
 \mathrm d_\tau(x)=\tau(a),\qquad \tau\in \mathrm T(A).
\end{equation}

Every trace on $A$ is faithful.  Thus
$\tau\mapsto \mathrm d_\tau(c)$ is a strictly positive lower semicontinuous
function on the compact space $\mathrm T(A)$, and consequently
\[
 \mu:=\min_{\tau\in \mathrm T(A)}\mathrm d_\tau(c)>0.
\]
Choose $m\in\mathbb N$ with $m\mu>\norm a$.  Equation
\eqref{Z-S0-rank} gives
\[
 \mathrm d_\tau(x)=\tau(a)\leq\norm a
 <m \mathrm d_\tau(c)=\mathrm d_\tau(c^{\oplus m}),\qquad \tau\in \mathrm T(A).
\]
Strict comparison yields $x\precsim c^{\oplus m}$.  Moreover, $A$ has
stable rank one by \cite{Ror-Z-stable}.  The Hilbert-module picture of
the Cuntz semigroup \cite{CEI-CuntzSG} therefore gives a positive
element
\[
 b\in\operatorname{her}(c^{\oplus m})\subseteq \mathrm M_m(H)
 \quad\text{with}\quad [b]=[x].
\]

We now consider tracial states on the possibly nonunital algebra
$H$.  Let $\sigma\in \mathrm T(H)$.  Since $H$ is full in the simple unital
algebra $A$, let $\widetilde\sigma$ be the standard lower
semicontinuous extension of $\sigma$ to $A$.  Fullness of
$0\ne q\in H_+$ gives finitely many $w_i\in A$ with
$\sum_iw_i^*qw_i=1_A$.  Traciality gives
\[
 \widetilde\sigma(1_A)
 =\sum_i\sigma(q^{1/2}w_iw_i^*q^{1/2})
 \leq\sum_i\norm{w_i}^2\sigma(q)<\infty.
\]
Thus $\widetilde\sigma$ is bounded.  Set
$s=\widetilde\sigma(1_A)>0$.  Then
$\tau=s^{-1}\widetilde\sigma\in \mathrm T(A)$.  Using
\eqref{Z-S0-rank} and $[b]=[x]$, we obtain
\[
 \mathrm d_\sigma(b)=s \mathrm d_\tau(b)=s \mathrm d_\tau(x)
 =s\tau(a)=\sigma(a).
\]
Thus the defining estimate for $\mathcal S_0$ holds with equality.
\end{proof}

\begin{rem}[Hereditary rank realization]\label{hereditary-rank-scope}
Hereditary subalgebras of a nonsimple \Cstar{}-algebra of finite nuclear
dimension need not belong to $\mathcal S_0$.  For example, let
$X=[0,1]$.  The pointwise rank of every element of
$\mathrm M_\infty(\mathrm C(X))$ is integer-valued, whereas evaluation traces take the
value $\frac12$ on the constant positive element $\frac12 1$.
Thus finite nuclear dimension alone does not imply the rank realization
in Lemma \ref{Z-S0} for nonsimple coefficient algebras.
\end{rem}

The preservation of $\mathcal Z$-stability below was obtained by
Orovitz, Phillips, and Wang under their weak tracial Rokhlin property
\cite[Corollary~8.13]{OPW-Strict-comparison-crossed-products}.
Their formulation controls the defect in the coefficient algebra
in the Cuntz order.  The following theorem proves the same conclusion
under WRP, which controls the defect uniformly over invariant traces.

\begin{thm}[WRP preserves $\mathcal Z$-stability]\label{WRP-finite-nucdim-Z}
Let $A$ be a separable, simple, unital, nuclear, $\mathcal Z$-stable
\Cstar{}-algebra, and let $\Gamma$ be a countable discrete amenable group
acting on $A$.  If $(A,\Gamma)$ has WRP, then
\[
 B=A\rtimes_{\alpha,r}\Gamma
\]
is simple and $\mathcal Z$-stable.
In particular, the conclusion holds for a separable, simple, unital,
infinite-dimensional coefficient algebra of finite nuclear dimension
with WRP.
\end{thm}
\begin{proof}
WRP requires $\mathrm T(A)\neq\emptyset$, so simplicity implies that
$A$ is stably finite.  It is infinite dimensional since it is
$\mathcal Z$-stable.  The action is minimal because $A$ is simple,
and Corollary~\ref{wrp-simple} gives simplicity of $B$.
Nuclearity of $A$ and amenability of $\Gamma$ give nuclearity of $B$.

Since $A$ is stably finite, nuclear, and $\mathcal Z$-stable, it has
strict comparison and almost divisible Cuntz semigroup, as in the
proof of Corollary~\ref{Z=>COS}.
Corollary~\ref{Z=>SCP} therefore gives strict comparison for $B$.
Furthermore, $\mathcal Z$-stability implies tracial
$\mathcal Z$-absorption by \cite[Proposition~2.2]{HO-Z}.  Choosing the
Cuntz defect below a nonzero positive element with uniformly small
dimension functions gives weak tracial $\mathcal Z$-absorption, as in
Remark~\ref{ta1-to-ta0}.  Corollary~\ref{WRP-WTZ} and the strict
comparison and nuclearity of $B$ give $B\cong B\otimes\mathcal Z$.
For the final assertion, finite nuclear dimension gives nuclearity
and \cite[Corollary~7.3]{Winter-Z-stable-02} gives
$\mathcal Z$-stability of $A$.

The conclusion also follows from
\cite[Corollary~5.6]{Gardella-Geffen-Naryshkin-Vaccaro-Z},
since WRP implies the weak Rokhlin condition in
\cite[Definition~5.2]{Gardella-Geffen-Naryshkin-Vaccaro-Z}.
Indeed, choose exact tiling shapes $D_s$ that are
$(K^{-1}\cup\{e\},\eta)$-invariant.
Their common $K$-translate intersection is
$\operatorname{Int}_{K^{-1}\cup\{e\}}(D_s)$ and has cardinality at
least $(1-\eta)|D_s|$.
The levels $\alpha_g(c_s)$ are central, orthogonal, and covariant;
their defect lies in $J_{\mathcal M_1(A,\Gamma)}$ and has value zero
under every invariant ultralimit trace.
\end{proof}

\begin{rem}[Infinite groups]\label{wrp-Z-finite-radius-extension}
For infinite groups, Corollary~\ref{wrp-finite-radius-Z} extends this
conclusion to stably finite coefficient \Cstar{}-algebras with locally
finite nuclear dimension and finite radius of comparison.  These
algebras need not be simple or $\mathcal Z$-stable and may have positive
radius of comparison.  The action is assumed to be minimal and to
satisfy at least one of the following conditions: WRP and COS, strong
URP, or WRP with Cuntz control in the crossed product.
Theorem~\ref{WRP-finite-nucdim-Z} also applies to finite groups.
\end{rem}

\chapter{Real rank and stable rank}\label{section-ranks}

\section{Permanence of real rank zero in weak tracial approximation}

Zhao and Fang \cite{Zhao-Fang-RR0} proved that real rank zero passes
from a centrally large subalgebra to an infinite-dimensional simple
unital \Cstar{}-algebra, without assuming stable rank one for the
subalgebra.  The same permanence theorem appears in the earlier joint
preprint \cite[Theorem~3.11]{Zhao-Fang-Fan-large-preprint}.  Its proof,
particularly Step~III, uses mutually equivalent projections
and a self-adjoint matrix congruence to absorb the positive and negative
parts of the defect.  The following proof adapts this method to weak
tracial approximation.

The approximating algebras may vary and need not be unital or simple,
so the cited permanence theorems do not apply directly.  Lifting an order
zero map transfers fixed matrix units to an approximating algebra while
preserving a uniform lower trace bound.  Strict comparison then
compares the positive and negative parts of the defect with the
resulting projections, and a self-adjoint matrix congruence gives an
invertible approximation.

\begin{thm}\label{rr0-twa}
Let $\mathcal S$ be a class of \Cstar{}-algebras of real rank zero, and let $A$
be a simple unital \Cstar{}-algebra in $\mathrm{TA}_0(\mathcal S)$.  If $A$ has
strict comparison of positive elements, then $A$
has real rank zero.
\end{thm}

\begin{proof}
Hereditary subalgebras of a \Cstar{}-algebra of real rank zero have approximate
identities of projections, and real rank zero passes to unitizations
and corners \cite{Brown-Pedersen-RR0}.  The approximation hypothesis
includes $\mathrm T(A)\ne\emptyset$.  Since $A$ is simple,
every tracial state on $A$ is faithful.

\medskip
\noindent\emph{Step 1: matrix units in an annihilating hereditary algebra.}
We first prove that $A$ has property (SP).  Let $0\ne a\in A_+$ be a
contraction and put
\[
 m=\min_{\tau\in\mathrm T(A)}\tau(a)>0.
\]
Choose an approximation satisfying the $\mathrm{TA}_0$ conditions with
tolerance $\eta<m/16$.  Taking a self-adjoint approximant to $a$ and
then the positive part of its compression gives $c\in C_+$ with
\[
 \|c-hah\|<2\eta.
\]
Since $0\leq1-h^2\leq2(1-h)$, we have
\[
 \tau(hah)\geq\tau(a)-\tau(1-h^2)>m-2\eta.
\]
Consequently $\|c\|>m-4\eta>m/2$, and the perturbation lemma for
Cuntz comparison gives
\[
 (c-m/4)_+\precsim_A hah\precsim_A a.
\]
The nonzero hereditary algebra generated by $(c-m/4)_+$ inside $C$
contains a nonzero projection.  A projection Cuntz subequivalent to $a$
is Murray--von Neumann equivalent to a projection in
$\operatorname{Her}_A(a)$, so this hereditary algebra also contains a
nonzero projection.  This proves (SP).

If $A$ is finite-dimensional, it has real rank zero.  Suppose
therefore that $A$ is infinite-dimensional.  Fix $x=x^*\in A$ with
$\|x\|\leq1$, and fix $\varepsilon>0$.  We may assume that $x$ is
noninvertible.  Continuous functional calculus then provides a
self-adjoint contraction $x_0$ and a nonzero positive contraction $d$
such that
\[
 \|x-x_0\|<\varepsilon/4,\qquad x_0d=dx_0=0.
\]
Indeed, replace $t\mapsto t$ by a function vanishing on a sufficiently
small neighborhood of zero, and choose $d$ supported in that neighborhood.

There are four mutually orthogonal equivalent nonzero projections in
$\operatorname{Her}_A(d)$.  To prove this, choose four orthogonal
nonzero positive elements in this simple non-elementary hereditary
algebra.  Property (SP) gives a nonzero projection in the hereditary
algebra generated by each element.  Simplicity gives a common nonzero
positive Cuntz lower
bound for these four projections: for two nonzero positive elements
$b,c$, choose $0\ne z\in b^{1/2}Ac^{1/2}$ and use
$zz^*\sim z^*z$, and iterate.  Applying (SP) to this lower bound and
using comparison of projections gives equivalent subprojections of the
four original projections.
Let
\[
 \iota\colon \mathrm M_4\longrightarrow A
\]
be the resulting homomorphism, and write $e_{ij}=\iota(E_{ij})$, where
$E_{ij}$ are the standard matrix units.  Thus $x_0e_{ij}=e_{ij}x_0=0$.
Set
\[
 \delta=\min_{\tau\in\mathrm T(A)}\tau(e_{11})>0.
\]
The matrix units and the constant $\delta$ are fixed before the
approximating algebras are chosen.

\medskip
\noindent\emph{Step 2: equivalent projections in the approximating algebras.}
Choose $C_n,h_n$ satisfying the $\mathrm{TA}_0$ conditions for the finite set
$\{x_0,e_{ij}:1\leq i,j\leq4\}$ with tolerances $\eta_n\searrow0$.
Choose self-adjoint approximants $x_n$ to $x_0$ and approximants
$e'_{ij,n}$ to $e_{ij}$, so that
\[
 \begin{gathered}
 h_nx_nh_n,\ h_ne'_{ij,n}h_n\in C_n,\qquad
 \|x_n-x_0\|,\ \|e'_{ij,n}-e_{ij}\|\longrightarrow0,\\
 \|[h_n,x_0]\|,\ \|[h_n,e_{ij}]\|\longrightarrow0,\qquad
 \sup_{\tau\in\mathrm T(A)}\mathrm d_\tau(1-h_n)<\eta_n.
 \end{gathered}
\]
Put
\[
 q_n=g_{1/20}(h_n),\qquad r_n=g_{3/4}(h_n),\qquad
 D_n=\operatorname{Her}_{C_n}(r_n).
\]
Consider the norm sequence algebra
$A_\infty=\ell^\infty(A)/c_0(A)$.  If $h=[(h_n)]$ and
$r=[(r_n)]$, then $h$ commutes with $x_0$ and all the matrix units.
The map
\[
 \mu\colon \mathrm M_4\longrightarrow
       \prod_n D_n\big/\bigoplus_n D_n\subseteq A_\infty,
 \qquad \mu(z)=r\iota(z),
\]
is c.p.c. and order zero.  Its range lies in the indicated quotient, since
\[
 r_n^{1/2}e'_{ij,n}r_n^{1/2}\in D_n:
\]
the function $g_{3/4}(t)^{1/2}/t$, extended by zero at zero, is
continuous, so the compression belongs to $C_n$ and is supported in
$r_n$.

The cone over $\mathrm M_4$ is projective
\cite[Theorem~10.2.1]{Lor-Mon}.  Equivalently, c.p.c. order zero maps
from $\mathrm M_4$ lift across quotient maps.  Hence $\mu$ lifts to c.p.c.
order zero maps $\phi_n\colon \mathrm M_4\to D_n$ satisfying
\begin{equation}\label{rr0-cone-transfer}
 \|\phi_n(E_{ij})-r_n^{1/2}e_{ij}r_n^{1/2}\|
       \longrightarrow0.
\end{equation}
Write $w_n=\phi_n(1_{\mathrm M_4})$ and $a_n=\phi_n(E_{11})$.  The scalar
inequality $1-g_{3/4}(t)\leq4(1-t)$ and
\eqref{rr0-cone-transfer} imply, uniformly over $\mathrm T(A)$,
\begin{equation}\label{rr0-reserve-trace}
 \tau(a_n)\geq\delta-4\eta_n-o(1).
\end{equation}

Define self-adjoint elements
\[
 \begin{split}
 a_{0,n}&=(1-q_n)^{1/2}x_n(1-q_n)^{1/2},\\
 a'_{1,n}&=q_n^{1/2}x_nq_n^{1/2}\in C_n.
 \end{split}
\]
Functional calculus and the commutator estimates give
\[
 \|x_0-a_{0,n}-a'_{1,n}\|\longrightarrow0.
\]
Moreover, $a_{0,n}D_n=D_na_{0,n}=0$, because $(1-q_n)r_n=0$.
In the sequence algebra, $[(a'_{1,n})]=g_{1/20}(h)x_0$ and
$[(w_n)]=r\iota(1)$, so
\[
 \|w_na'_{1,n}\|+\|a'_{1,n}w_n\|\longrightarrow0.
\]
Fix $0<\theta<\delta/32$ and put
\[
 a_{1,n}=(1-g_\theta(w_n))a'_{1,n}(1-g_\theta(w_n))\in C_n.
\]
Since $g_\theta(t)/t$ is bounded and continuous at zero, it follows that
\begin{equation}\label{rr0-decomposition}
 \|x_0-a_{0,n}-a_{1,n}\|\longrightarrow0.
\end{equation}
Both summands annihilate
$\operatorname{Her}_{C_n}((w_n-2\theta)_+)$.

For all sufficiently large $n$, \eqref{rr0-reserve-trace} gives
$\tau(a_n)>3\delta/4$ for every $\tau\in\mathrm T(A)$.  Real rank
zero gives a projection
\[
 p_{1,n}\in\operatorname{Her}_{C_n}((a_n-2\theta)_+)
\]
such that
\[
 \|(a_n-2\theta)_+
       -p_{1,n}(a_n-2\theta)_+p_{1,n}\|<\delta/8.
\]
As $a_n$ is a contraction, it follows that
\begin{equation}\label{rr0-projection-reserve}
 \tau(p_{1,n})\geq\tau((a_n-2\theta)_+)-\delta/8
       >\delta/2,\qquad \tau\in\mathrm T(A).
\end{equation}
The supporting homomorphism of $\phi_n$ gives four orthogonal equivalent
projections $p_{i,n}$.  They belong to
$\operatorname{Her}_{C_n}((w_n-2\theta)_+)$.  These projections belong
to $C_n$: choose a continuous
function $f$ with $f(0)=0$ and $f(t)=t^{-1}$ for $t\geq2\theta$,
and set
\[
 t_{i,n}=\phi_n(E_{i1})f(a_n)p_{1,n},\qquad
 p_{i,n}=t_{i,n}t_{i,n}^*.
\]
The order zero relations give $t_{i,n}^*t_{i,n}=p_{1,n}$ and the
asserted orthogonality and support.  Thus
$\tau(p_{i,n})=\tau(p_{1,n})>\delta/2$ for all four projections.

\medskip
\noindent\emph{Step 3: an invertible self-adjoint approximation.}
Choose $n$ sufficiently large that $\eta_n<\delta/2$ and the norm in
\eqref{rr0-decomposition} is less than $\varepsilon/4$.  Suppress $n$
from the notation and put
\[
 P=p_1+p_2+p_3+p_4,\qquad Q=1_A-P.
\]
The elements $a_0,a_1$ belong to $QAQ$, with $a_1\in C$.  Also,
the positive and negative parts of $a_0$ are supported in $1-q$, so
\[
 \mathrm d_\tau((a_0)_\pm)\leq \mathrm d_\tau(1-q)
       \leq \mathrm d_\tau(1-h)<\eta_n<\tau(p_i).
\]
Strict comparison therefore gives
\[
 (a_0)_+\precsim_A p_1,\qquad
 (a_0)_-\precsim_A p_3.
\]
For any $\kappa>0$, choose
$b_+\in QAp_1$ and $b_-\in QAp_3$ such that
\[
 \|(a_0)_+-b_+b_+^*\|<\kappa,\qquad
 \|(a_0)_--b_-b_-^*\|<\kappa.
\]
These elements are obtained from the Cuntz comparison witnesses by multiplying
on the left by $Q$ and on the right by $p_1$ or $p_3$.

The algebra $B=C+\mathbb C1_A$ has real rank zero.  If $C$ is
nonunital, $B$ is a quotient of its unitization; if $C$ has unit $e$,
then $B=C\oplus\mathbb C(1_A-e)$, with the second summand omitted
when $e=1_A$.  Hence $QBQ$ also has real rank zero, since $P\in C$.
If $Q=0$, \eqref{rr0-decomposition} shows that $x$ is within
$\varepsilon/2$ of zero, which gives the required invertible self-adjoint
approximation.
Otherwise choose an invertible self-adjoint $v\in QBQ$ with
\[
 \|v-a_1\|<\kappa.
\]
Choose partial isometries $u_1,u_2\in C$ such that
\[
 u_1^*u_1=p_1,\quad u_1u_1^*=p_2,\qquad
 u_2^*u_2=p_3,\quad u_2u_2^*=p_4.
\]
For $t>0$ define
\[
 \begin{split}
 R_t&=Q+b_++b_-+tP,\\
 W_t&=v+p_1-p_3+t(u_1+u_1^*+u_2+u_2^*).
 \end{split}
\]
The element $R_t$ is invertible, with inverse
\[
 R_t^{-1}=Q+t^{-1}P-t^{-1}(b_++b_-).
\]
The element $W_t$ is self-adjoint and invertible, since its $Q$-corner
is $v$ and its two remaining matrix corners are
\[
 \begin{pmatrix}1&t\\t&0\end{pmatrix}
 \quad\hbox{and}\quad
 \begin{pmatrix}-1&t\\t&0\end{pmatrix}.
\]
Consequently $R_tW_tR_t^*$ is invertible and self-adjoint.  Multiplying
the corners gives
\[
 R_tW_tR_t^*\longrightarrow
       v+b_+b_+^*-b_-b_-^*\quad\hbox{in norm as }t\searrow0.
\]
All remaining terms contain a factor $t$: the off-diagonal
terms are
$t(b_+-b_-)+t^2(b_+u_1^*+b_-u_2^*)$ and their adjoints, while the
$P$-corner is
$t^2(p_1-p_3)+t^3(u_1+u_1^*+u_2+u_2^*)$.

Choose $\kappa<\varepsilon/12$ and then $t$ sufficiently small that the
error in the preceding norm limit is less than $\varepsilon/4$.
The approximation estimates give
\[
 \|x-R_tW_tR_t^*\|
 <\varepsilon/4+\varepsilon/4+3\kappa+\varepsilon/4
 <\varepsilon.
\]
Invertible self-adjoint elements are therefore dense in $A_{\mathrm{sa}}$,
which proves real rank zero.
\end{proof}

\begin{rem}[Local weak comparison]\label{rr0-local-comparison}
The proof uses only the two comparisons in Step~3.  Thus
Theorem~\ref{rr0-twa} also holds under local weak comparison by tracial
states: it suffices that there be a fixed $\gamma\geq1$ such that
\[
 \gamma\sup_{\tau\in\mathrm T(A)}\mathrm d_\tau(a)
       <\inf_{\tau\in\mathrm T(A)}\mathrm d_\tau(b)
 \quad\Longrightarrow\quad a\precsim_A b.
\]
Indeed, choose $n$ with $\gamma\eta_n<\delta/2$.  No separability,
nuclearity, or exactness is required in the proof.  If strict
comparison is formulated using quasitraces rather than tracial states,
exactness identifies quasitraces with traces and therefore makes that
formulation compatible with Definition~\ref{defn_WAS}.
\end{rem}

\section{Permanence of stable rank one in weak tracial approximation}

The following proposition isolates the use of tracial approximation.
Its comparison hypothesis concerns only the positive elements in the
subalgebra generated by the approximant to $d$ and the compressing element.

\begin{prop}\label{sr1-twa-0}
Let $\mathcal S$ be a class of \Cstar{}-algebras of stable rank one.  Let $A$ be a
simple finite unital \Cstar{}-algebra with $\mathrm T(A)\neq\varnothing$, and let
$a\in A$ and $d\in A_+$ satisfy
\[
        ad=da=0,\qquad \|d\|=1.
\]
Suppose that, for every $\eta>0$, there are a \Cstar{}-subalgebra $C\subseteq A$
with $C\in\mathcal S$, a positive contraction $h\in C$ of norm one,
$a'\in A$, and a positive contraction $d'\in A$ such that
\begin{enumerate}
\item[(a)] $\|a-a'\|<\eta$ and $\|d-d'\|<\eta$;
\item[(b)] $ha'h,hd'h\in C$;
\item[(c)] $\|[h,a']\|<\eta$ and $\|[h,d']\|<\eta$;
\item[(d)] $\mathrm d_\tau(1-h)<\eta$ for every $\tau\in\mathrm T(A)$;
\item[(e)] if $x,y\in \mathrm C^*(d',h,1_A)_+$ satisfy
\[
        \mathrm d_\tau(x)<\mathrm d_\tau(y),\qquad \tau\in\mathrm T(A),
\]
then $x\precsim_A y$.
\end{enumerate}
Then $a\in\overline{\mathrm{GL}(A)}$.
\end{prop}

\begin{proof}
Since $A$ is simple, every tracial state on $A$ is faithful.  Compactness of
$\mathrm T(A)$ therefore gives
\[
 \delta:=\min_{\tau\in\mathrm T(A)}
             \tau\bigl(g_{1/2}(d)\bigr)>0.
\]
Fix $\varepsilon>0$, and choose
\[
        0<\theta<\frac{\delta}{64}.
\]
Apply the hypothesis with $\eta>0$ sufficiently small, as specified in
the estimates that follow.  Put
\[
 q=g_{1/20}(h),\qquad r=g_{3/4}(h),
\]
and define
\[
 \begin{split}
 a_0&=(1-q)^{1/2}a'(1-q)^{1/2},\\
 a_1'&=q^{1/2}a'q^{1/2},\\
 D&=r^{1/2}d'r^{1/2}.
 \end{split}
\]
The functions $g_{1/20}(t)^{1/2}/t$ and
$g_{3/4}(t)^{1/2}/t$, extended by zero at the origin, are continuous on
$[0,1]$.  Consequently, condition (b) gives
\[
        a_1',D\in C.
\]
Since $(1-q)r=0$, we have
\begin{equation}\label{sr1-a0-D}
        a_0D=Da_0=0.
\end{equation}

For the fixed value of $\theta$, conditions (a)--(c), functional calculus,
and $ad=da=0$ give, for sufficiently small $\eta$,
\begin{equation}\label{sr1-cut-estimates}
 \left\|a-(a_0+a_1')\right\|<\frac{\varepsilon}{2},
 \qquad
 \left\|g_\theta(D)a_1'\right\|
 +\left\|a_1'g_\theta(D)\right\|<\frac{\varepsilon}{4}.
\end{equation}
The first estimate follows from
$[q,a']\approx0$ and
\[
 (1-q)^{1/2}a'(1-q)^{1/2}+q^{1/2}a'q^{1/2}\approx a'.
\]
For the second estimate, $Da_1'\approx0\approx a_1'D$ follows from
$d'a'\approx0\approx a'd'$ and the commutator estimates.  Since
$g_\theta(t)=k_\theta(t)t$ for a continuous function $k_\theta$ with
$\|k_\theta\|\leq 2/\theta$, decreasing $\eta$ gives the required bound.

The scalar inequality
\[
        1-g_{3/4}(t)\leq4(1-t),\qquad 0\leq t\leq1,
\]
together with condition (d), gives
\[
 \tau(D)=\tau(d'r)
 \geq \tau(d')-\tau(1-r)
 >\tau(d)-5\eta.
\]
Since $d\geq\frac12g_{1/2}(d)$, we may decrease $\eta$ so that
\begin{equation}\label{sr1-D-trace}
        \tau(D)>\frac{\delta}{4},
        \qquad \tau\in\mathrm T(A).
\end{equation}

Set
\[
        a_1=(1-g_\theta(D))a_1'(1-g_\theta(D))
        \quad\text{and}\quad
        y=(D-2\theta)_+.
\]
Then $a_1\in C$.  The second estimate in
\eqref{sr1-cut-estimates} gives $\|a_1-a_1'\|<\varepsilon/2$, and hence
\begin{equation}\label{sr1-close-a}
        \|a-(a_0+a_1)\|<\varepsilon.
\end{equation}
Moreover,
\begin{equation}\label{sr1-y-orth}
        y(a_0+a_1)=(a_0+a_1)y=0.
\end{equation}
Indeed, this follows for $a_0$ from \eqref{sr1-a0-D}, and for $a_1$ from
$y(1-g_\theta(D))=0$.

Consider the positive contractions
\[
        e=1-g_{1/3}(h)
        \quad\text{and}\quad
        b=2(e-\tfrac12)_+.
\]
We have
\[
 \mathrm d_\tau(e)\leq \mathrm d_\tau(1-h)<\eta,
\]
whereas \eqref{sr1-D-trace} gives
\[
 \mathrm d_\tau(y)\geq\tau(y)
 \geq\tau(D)-2\theta>\frac{\delta}{8}.
\]
Taking $\eta<\delta/8$, condition (e) gives $e\precsim_A y$.
The standard cut-down form of Cuntz subequivalence provides $v_0\in A$
such that
\[
        v_0^*v_0=(e-\tfrac12)_+,\qquad
        v_0v_0^*\in\operatorname{Her}_A(y).
\]
With $v=\sqrt2v_0$, we therefore have
\begin{equation}\label{sr1-v}
        v^*v=b,\qquad p:=vv^*\in\operatorname{Her}_A(y).
\end{equation}
Because $b\,y=0$, the positive contractions $b$ and $p$ are orthogonal.

Define
\[
 \rho(t)=
 \begin{cases}
 t^{-1/2}\sin(\pi t/2),&t>0,\\
 0,&t=0,
 \end{cases}
\]
and
\begin{equation}\label{sr1-unitary}
 w=\cos\bigl(\tfrac{\pi}{2}(p+b)\bigr)+\rho(p)v-\rho(b)v^*.
\end{equation}
By \cite[Lemma~3.2]{LN-sr1}, $w$ is a unitary.  Set
\[
 \begin{split}
 z_1&=1-g_{1/10}(h),\\
 z_3&=g_{2\theta}(D),\\
 z_2&=1-z_1-z_3=g_{1/10}(h)-g_{2\theta}(D).
 \end{split}
\]
Since $D$ is supported in $r=g_{3/4}(h)$, the element $g_{1/10}(h)$
acts as the identity on $\operatorname{Her}(D)$.  Thus $z_2$ is positive,
and
\[
        z_1+z_2+z_3=1,\qquad z_1z_3=0.
\]
Also, $by=z_1y=0$, while $z_3y=y$.  It follows from
\eqref{sr1-v} that
\begin{equation}\label{sr1-zp}
        z_1p=0,\qquad z_3p=p,\qquad z_2p=0.
\end{equation}
Since $bz_1=z_1$, the last assertion of \cite[Lemma~3.2]{LN-sr1} gives
\begin{equation}\label{sr1-rotate}
        wz_1w^*\in\operatorname{Her}_A(p).
\end{equation}

Let
\[
        x=w^*(a_0+a_1).
\]
Equations \eqref{sr1-y-orth} and \eqref{sr1-rotate} show that $z_1x=0$.
The definitions of $a_1$ and $z_3$, together with
\eqref{sr1-a0-D}, give $xz_3=0$.  The middle compression belongs to
$C+\mathbb C1_A$, as follows.
Taking adjoints in \eqref{sr1-unitary}, using
$p(a_0+a_1)=0$, and using $z_2v=0$, which follows from
\eqref{sr1-zp}, we obtain
\[
 \begin{split}
 z_2xz_2
 &=z_2\cos\bigl(\tfrac{\pi}{2}(p+b)\bigr)
                    (a_0+a_1)z_2\\
 &=z_2\cos\bigl(\tfrac{\pi}{2}b\bigr)(a_0+a_1)z_2\\
 &=z_2\cos\bigl(\tfrac{\pi}{2}b\bigr)a_1z_2.
 \end{split}
\]
The last equality follows from $b(1-q)=1-q$ and hence
$\cos(\pi b/2)(1-q)=0$.

Put $B=C+\mathbb C1_A$.  If $C$ is nonunital, $B$ is a quotient of
its unitization.  If $C$ is unital, $B$ is either $C$ or the direct
sum of $C$ and a one-dimensional algebra.  Thus $B$ has stable rank
one.  Moreover, $z_1,z_2,z_3\in B$, and the last display shows that
$z_2xz_2\in B$.
Lemma~5.4 of \cite{A-NCP-LAlg} therefore gives
\[
        x\in\overline{\mathrm{GL}(A)}.
\]
Since $w$ is unitary, $a_0+a_1\in\overline{\mathrm{GL}(A)}$.
Equation~\eqref{sr1-close-a} and the arbitrariness of $\varepsilon$ prove
the proposition.
\end{proof}

\begin{thm}\label{sr1-twa}
Let $\mathcal S$ be a class of \Cstar{}-algebras of stable rank one.  Let $A$ be
a simple unital \Cstar{}-algebra in $\mathrm{TA}_0(\mathcal S)$.  If $A$ has
strict comparison of positive elements, then $A$ has stable rank one.
\end{thm}

\begin{proof}
The tracial states on $A$ are faithful, so $A$ is finite.  Let $x\in A$ and
let $\varepsilon>0$.  We may assume that $x$ is not invertible.  Since a
one-sided invertible element in a finite unital \Cstar{}-algebra is invertible,
$x$ is not one-sided invertible.

R{\o}rdam's zero-divisor reduction
\cite[Proposition~3.2]{RorUHF} gives a two-sided zero divisor $x_0\in A$
such that
\[
        \|x-x_0\|<\frac{\varepsilon}{2}.
\]
The left and right annihilators of $x_0$ may be different.
By \cite[Lemma~3.5]{RorUHF}, there are a unitary $u\in A$ and a nonzero
positive element $d\in A$ such that
\[
        d(ux_0)=(ux_0)d=0.
\]
After rescaling, we may suppose that $\|d\|=1$.

Apply the $\mathrm{TA}_0(\mathcal S)$ approximation to the finite set
$\{ux_0,d\}$.  The approximant to $d$ may be chosen to be a positive
contraction.  Indeed, if $\|d-d_0'\|<\kappa$, put
\[
 s=\frac{d_0'+(d_0')^*}{2}
 \quad\text{and}\quad
 d'=\frac{s+\kappa1_A}{1+2\kappa}.
\]
Then $0\leq d'\leq1_A$, $\|d-d'\|<4\kappa$, and $hd'h\in C$ because
$hd_0'h\in C$ and $h^2\in C$.  The commutator estimate is preserved as
well.  Since $\kappa$ is arbitrary, condition (a) holds with any
prescribed tolerance.  Strict comparison of positive elements gives
condition (e) of Proposition~\ref{sr1-twa-0}, and the remaining
conditions follow from $\mathrm{TA}_0(\mathcal S)$.  Hence
\[
        ux_0\in\overline{\mathrm{GL}(A)}.
\]
Multiplication by the unitary $u$ shows that
$x_0\in\overline{\mathrm{GL}(A)}$.  Consequently
$\operatorname{dist}(x,\mathrm{GL}(A))<\varepsilon$.  Since $x$ and
$\varepsilon$ were arbitrary, $A$ has stable rank one.
\end{proof}

\begin{rem}
The proof of Proposition~\ref{sr1-twa-0} uses only the comparison
\[
 1-g_{1/3}(h)\precsim (r^{1/2}d'r^{1/2}-2\theta)_+,
 \qquad r=g_{3/4}(h).
\]
Thus the strict comparison hypothesis in Theorem~\ref{sr1-twa} can be
replaced by this local comparison condition.  More generally, it
suffices that there be a fixed $K\geq1$ such that the relevant elements
satisfy
\[
 K\mathrm d_\tau(e)<\mathrm d_\tau(y)\quad(\tau\in\mathrm T(A))
 \quad\Longrightarrow\quad e\precsim_A y.
\]
In the proof, $\mathrm d_\tau(e)<\eta$ and $\mathrm d_\tau(y)>\delta/8$, so choosing
$\eta<\delta/(8K)$ suffices.  This is a comparison with the single
target $y$; comparison with $y^{\oplus m}$ alone does not imply the
unitary rotation in that proof.
\end{rem}

\section{Zero-divisor reduction to the coefficient algebra}

R{\o}rdam's argument gives left and right annihilators in the crossed
product.  The results in this section replace them by a common positive
annihilator in the coefficient algebra.  In the norm ultrapower, the
Rokhlin levels are exactly orthogonal, and Lemma~3.2 of \cite{LN-sr1}
gives the required unitary rotations.  Throughout the section, write
\(B=A\rtimes_{\alpha,r}\Gamma\) and denote the canonical conditional
expectation by \(\mathbb E:B\to A\).

\begin{lem}[One-sided coefficient reduction]\label{pre-alg-red}
Let $(A,\Gamma)$ be a minimal \Cstar{}-dynamical system with WRP, where
$A$ is unital and separable, $\mathrm T(A)\ne\emptyset$, and
$\Gamma$ is countably infinite, discrete, and amenable.  If
$a\in B$ and $0\ne b\in B_+$ satisfy $ba=0$, then, for every
$\eps>0$, there exist a unitary $u\in B$ and a norm-one positive
contraction $h\in A$ such that
\[
             \|hua\|<\eps.
\]
\end{lem}

\begin{proof}
Faithfulness of $\mathbb E$ gives $\mathbb E(b)\ne0$.  After rescaling
$b$, we may suppose that $\|\mathbb E(b)\|=1$.  Fix
$\kappa>0$.  By approximating $b^{1/2}$ by a Fourier polynomial, then
squaring and rescaling, choose
\[
 b_0\in \mathrm C_{\mathrm c}(\Gamma,A)_+,
 \qquad \|\mathbb E(b_0)\|=1,
 \qquad \|b-b_0\|<\kappa.
\]
Let $F\subseteq\Gamma$ contain the Fourier supports of $b_0$ and
$b_0^2$.  Since $\Gamma$ is infinite, choose
$\gamma_0\notin F$.

By Lemma~\ref{WRP-Kishimoto}, the action satisfies Elliott's condition.
Proposition~\ref{Elliott-central-sequence}, applied to
$\mathbb E(b_0)^2$ and $F\cup\{\gamma_0\}$, gives a norm-one positive
contraction $c\in A'\cap A_\omega$ such that its indicated translates
are mutually orthogonal and
\begin{equation}\label{one-sided-preserve-norm}
       \|c\mathbb E(b_0)^2c\|=1.
\end{equation}
Centrality and tower orthogonality give
\begin{equation}\label{one-sided-fourier-cut}
 c b_0^2c=c\mathbb E(b_0^2)c\in A_\omega,
 \qquad
 \alpha_{\gamma_0}(c)b_0c=0.
\end{equation}
Put $v=u_{\gamma_0}cb_0\in B_\omega$.  Then
\[
 vv^*=u_{\gamma_0}c\mathbb E(b_0^2)cu_{\gamma_0}^*\in A_\omega,
 \qquad v^*v=b_0c^2b_0,
\]
and the second relation in \eqref{one-sided-fourier-cut} shows that
$vv^*\perp v^*v$.  Kadison's inequality and
\eqref{one-sided-preserve-norm} give $\|v\|\geq1$.  Thus, for
$x=v/\|v\|$,
\begin{equation}\label{one-sided-initial-small}
 xx^*\in A_\omega,\qquad xx^*\perp x^*x,\qquad
 \|x^*xa\|\leq \|b_0\|\,\|(b_0-b)a\|.
\end{equation}

Functional calculus gives a perturbation whose range modulus is
identically one on a nonzero spectral region.  Write the polar
decomposition in $B_\omega^{**}$
as $x=(xx^*)^{1/2}s$.  Choose a continuous function $f:[0,1]\to[0,1]$ which vanishes
near zero, is identically one near one, and is uniformly as close
to the identity function as desired.  Then
\[
 t=f((xx^*)^{1/2})s\in B_\omega,\qquad
 tt^*\in A_\omega,\qquad tt^*\perp t^*t,\qquad \|t\|=1.
\]
Choose a continuous function $g:[0,1]\to[0,1]$ with $g(1)=1$
whose support is contained in the region on which $f=1$, and put
\[
                 h=g((xx^*)^{1/2})\in (A_\omega)_+.
\]
Then $\|h\|=1$ and $(tt^*)h=h(tt^*)=h$.  Lemma~3.2 of
\cite{LN-sr1} supplies a unitary $w\in B_\omega$ such that
\[
       (t^*t)(w^*hw)=w^*hw=(w^*hw)(t^*t).
\]
Consequently,
\[
 \|hwa\|=\|(w^*hw)a\|\leq\|t^*ta\|.
\]
By \eqref{one-sided-initial-small}, the right-hand side is smaller than
$\eps$ when first $\kappa$ and then the functional-calculus error are
chosen sufficiently small.

Choose positive-contraction representatives of $h$ in $A$ and
representatives of $w$ in $B$.  Polar decomposition gives unitary
representatives on an $\omega$-large set.  The required estimate holds
on an $\omega$-large set of coordinates; choosing one such coordinate
and normalizing the positive contraction proves the lemma.
\end{proof}

The following lemma aligns two positive elements of the coefficient
algebra by a unitary in the crossed product.

\begin{lem}[Alignment of coefficient supports]\label{small-unitary}
Let $(A,\Gamma)$ satisfy Elliott's condition, where $A$ is unital and
separable and $\Gamma$ is countable and discrete.  Assume that
$B=A\rtimes_{\alpha,r}\Gamma$ is simple and that every nonzero
hereditary \Cstar{}-subalgebra of $A$ is infinite-dimensional.  If
$h_0,h_1\in A_+$ are norm-one positive contractions, then, for every
$\eps>0$, there exist a norm-one positive contraction $k\in A$ and a
unitary $u\in B$ such that
\begin{equation}\label{aligned-supports}
 \|h_0k-k\|<\eps,
 \qquad
 \|h_1uku^*-uku^*\|<\eps.
\end{equation}
\end{lem}

\begin{proof}
We first construct orthogonal hereditary subalgebras on which
multiplication by $h_0$ and $h_1$, respectively, is close to the
identity.  Fix $0<\eta<\eps/2$ and put
\[
 e_i=(h_i-(1-\eta))_+,\qquad i=0,1.
\]
These elements are nonzero.  To obtain nonzero orthogonal positive
elements $d_i\in\overline{e_iAe_i}$, $i=0,1$, first observe that when
$e_0e_1=0$ the normalized elements $e_0,e_1$ suffice.  Otherwise set
\[
              T=e_0^{1/2}e_1e_0^{1/2}.
\]
If $\operatorname{sp}(T)\cap(0,\infty)$ contains at least two points,
choose nonzero functions $f_0,f_1\in \mathrm C_0((0,\infty))_+$ with disjoint
supports and with $T f_i(T)\ne0$, and set
\[
 \begin{split}
 d'_0&=e_0^{1/2}f_0(T)e_0^{1/2},\\
 d'_1&=e_1e_0^{1/2}f_1(T)e_0^{1/2}e_1.
 \end{split}
\]
Then $d'_0d'_1=0$, and both elements are nonzero.  If the positive
spectrum of $T$ is the singleton $\{\lambda\}$, its spectral projection
$p$ belongs to $A$.  The hereditary algebra $pAp$ is
infinite-dimensional, so it contains nonzero orthogonal positive
elements $x_0,x_1$.  The same formulas, with $f_i(T)$ replaced by
$x_i$, again give nonzero orthogonal elements.  After normalization we
obtain $d_0,d_1$.  Since the support projection of $e_i$ is dominated
by $1_{(1-\eta,1]}(h_i)$,
\begin{equation}\label{di-supported-hi}
       \|(1-h_i)d_i\|\leq\eta,\qquad i=0,1.
\end{equation}

Simplicity of $B$ gives $x\in\mathrm C_{\mathrm c}(\Gamma,A)$ such
that $v=d_0xd_1\ne0$.  Faithfulness of $\mathbb E$ gives
$\mathbb E(vv^*)\ne0$; rescale $v$ so that
$\|\mathbb E(vv^*)\|=1$.  Let $F$ contain the Fourier support of
$vv^*$.  Proposition~\ref{Elliott-central-sequence}, applied to
$\mathbb E(vv^*)^2$, gives a tower base $c\in A'\cap A_\omega$ with
\[
              \|c\mathbb E(vv^*)^2c\|=1.
\]
As in \eqref{one-sided-fourier-cut}, exact tower orthogonality gives
\[
        cvv^*c=c\mathbb E(vv^*)c\in A_\omega.
\]
Put $z=cv\in B_\omega$.  Since $c$ commutes with $d_0$ and $d_1$,
\[
 zz^*\in\overline{d_0A_\omega d_0},\qquad
 z^*z\in\overline{d_1B_\omega d_1}.
\]
In particular, $zz^*\perp z^*z$.  The preceding norm equality also
implies $\|zz^*\|=1$.

Applying the functional-calculus construction in the proof of
Lemma~\ref{pre-alg-red} gives $t\in B_\omega$ and a norm-one
positive contraction $k_\omega\in A_\omega$ such that
\[
 tt^*\perp t^*t,\qquad
 k_\omega\in\overline{d_0A_\omega d_0},\qquad
 (tt^*)k_\omega=k_\omega,\qquad
 t^*t\in\overline{d_1B_\omega d_1}.
\]
Lemma~3.2 of \cite{LN-sr1} gives a unitary $w\in B_\omega$ for which
\[
       w^*k_\omega w\in\overline{d_1B_\omega d_1}.
\]
Together with \eqref{di-supported-hi}, this yields
\[
 \|h_0k_\omega-k_\omega\|\leq\eta,\qquad
 \|h_1w^*k_\omega w-w^*k_\omega w\|\leq\eta.
\]
Choose representatives of $k_\omega$ and $w^*$ in $A$ and $B$,
respectively, and use polar decomposition to obtain unitaries on an
$\omega$-large set.  Choosing a coordinate at which the estimates
hold and normalizing the positive element gives \eqref{aligned-supports}.
\end{proof}

\begin{prop}[Coefficient zero-divisor reduction]\label{alg-red}
Let $(A,\Gamma)$ be a minimal \Cstar{}-dynamical system with WRP, where
$A$ is unital and separable, $\mathrm T(A)\ne\emptyset$, and
$\Gamma$ is countably infinite, discrete, and amenable.  Assume that
every nonzero hereditary \Cstar{}-subalgebra of $A$ is infinite-dimensional.
Put $B=A\rtimes_{\alpha,r}\Gamma$.  If $a\in B$ is noninvertible,
then, for every $\eps>0$, there are unitaries $u_0,u_1\in B$, an
element $b\in\mathrm C_{\mathrm c}(\Gamma,A)$, and a norm-one positive
contraction $d\in A$ such that
\[
             \|u_0au_1-b\|<\eps,\qquad db=bd=0.
\]
\end{prop}

\begin{proof}
Corollary~\ref{wrp-simple} shows that $B$ is simple.  It has a faithful
tracial state: average a tracial state of $A$ over $\Gamma$, compose the
result with $\mathbb E$, and use Corollary~\ref{min-faithful-trace}.
Thus $B$ is finite, and a noninvertible element of $B$ is neither left
nor right invertible.
By \cite[Proposition~3.2]{RorUHF}, after an arbitrarily small
perturbation of $a$ we obtain a two-sided zero divisor $a_0\in B$ and
nonzero elements $r_0,r_1\in B_+$ such that
\[
                   r_0a_0=0=a_0r_1.
\]

Apply Lemma~\ref{pre-alg-red} to $(a_0,r_0)$ and to $(a_0^*,r_1)$.
For any prescribed $\delta>0$, it gives unitaries $v_0,v_1\in B$ and
norm-one positive contractions $h_0,h_1\in A$ such that
\begin{equation}\label{two-one-sided-cuts}
 \|h_0v_0a_0\|<\delta,\qquad
 \|a_0v_1^*h_1\|<\delta.
\end{equation}
Lemma~\ref{small-unitary} gives a unitary $w\in B$ and a norm-one
positive contraction $k\in A$ such that
\[
 \|h_0k-k\|<\delta,\qquad
 \|h_1wkw^*-wkw^*\|<\delta.
\]
Using these relations and \eqref{two-one-sided-cuts}, we obtain
\begin{equation}\label{common-approx-annihilator}
 \|k v_0a_0v_1^*w\|<(1+\|a_0\|)\delta,\qquad
 \|v_0a_0v_1^*wk\|<(1+\|a_0\|)\delta.
\end{equation}

Fix $0<\sigma<1/4$ and put $q=g_\sigma(k)$.  Since
$g_\sigma(t)=t\ell_\sigma(t)$ for a continuous function
$\ell_\sigma$, \eqref{common-approx-annihilator} shows that both
$qv_0a_0v_1^*w$ and $v_0a_0v_1^*wq$ have arbitrarily small norm when
$\delta$ is sufficiently small.  Choose
$z\in\mathrm C_{\mathrm c}(\Gamma,A)$ close to
$v_0a_0v_1^*w$ and set
\[
        b=(1-q)z(1-q),\qquad
        d=\frac{(k-\sigma)_+}{1-\sigma}.
\]
Then $d$ is a norm-one positive contraction, $d(1-q)=0=(1-q)d$, and
hence $db=bd=0$.  Taking the perturbation of $a$, then $\delta$, and
then the Fourier approximation error sufficiently small gives
\(\|v_0av_1^*w-b\|<\eps\).  Taking
$u_0=v_0$ and $u_1=v_1^*w$ proves the proposition.
\end{proof}

The following lemma gives orthogonal coefficient elements that are
equivalent in the crossed product.  They will absorb the multiplicity
in COS.

\begin{lem}[Equivalent coefficient reserves]
\label{coefficient-multiplicity-reserves}
Let $\alpha$ be a minimal action on a unital \Cstar{}-algebra $A$, and
suppose that every nonzero hereditary \Cstar{}-subalgebra of $A$ is
infinite-dimensional.  Put $B=A\rtimes_{\alpha,r}\Gamma$.
Let $0\ne h\in A_+$, and let $c\in A_+$ be a positive contraction
with $ch=h$.  For every $m\in\mathbb N$, there are a norm-one
positive contraction $a\in A$, elements $g_1,\ldots,g_m\in\Gamma$,
and monomials $t_i\in A u_{g_i}$ such that
\[
 t_i^*t_j=\delta_{ij}a,\qquad ct_i=t_i,
 \qquad t_i t_i^*\in\operatorname{Her}_A(h).
\]
The range elements $t_i t_i^*$ are mutually orthogonal.
\end{lem}

\begin{proof}
Choose mutually orthogonal nonzero positive elements
$h_1,\ldots,h_m\in\operatorname{Her}_A(h)$.  Begin with
$t_1=h_1^{1/2}$ and $a=h_1$.
Suppose that $t_1,\ldots,t_r$ satisfy $t_i^*t_i=a\ne0$
and $t_it_i^*$ belongs to the respective hereditary algebra
$\operatorname{Her}_A(h_i)$.

Minimality gives $g\in\Gamma$ and $x\in A$ such that
\[
 y=a x\alpha_g(h_{r+1})\ne0.
\]
Indeed, otherwise $a$ annihilates the sum of the ideals generated by
all the translates of $h_{r+1}$; its closure is a nonzero invariant
ideal and hence is $A$, contradicting $a\ne0$.
Put $z=yy^*\in\operatorname{Her}_A(a)$.
Write $t_i=V_i a^{1/2}$ for the polar decompositions in $B^{**}$,
and replace each $t_i$ by $V_i z^{1/2}$.
These elements still belong to $A u_{g_i}$.  To see this, approximate
$z^{1/2}$ by $g_\eps(a)z^{1/2}$; the element
$V_i g_\eps(a)$ equals $t_i f_\eps(a)$, where
$f_\eps(t)=g_\eps(t)t^{-1/2}$ for $t>0$ and $f_\eps(0)=0$.
The approximants belong to the closed monomial subspace
$A u_{g_i}$.

The new elements have source $z$, and their range elements remain
in the respective hereditary algebras.  Add the monomial
$u_g^*y^*\in A u_{g^{-1}}$.  Its source is $z$ and its
range element is
$\alpha_{g^{-1}}(y^*y)\in\operatorname{Her}_A(h_{r+1})$.
Orthogonality of the hereditary algebras gives the off-diagonal
relations.  Induction produces all $m$ elements.
Divide the common source by its norm and all $t_i$ by the square root
of that norm.  Since $ch=h$, multiplication by $c$ fixes
$\operatorname{Her}_A(h)$; the inclusion of every range element
in this hereditary algebra then gives $ct_i=t_i$.
\end{proof}

Combining the coefficient reduction with the local comparison argument
of Proposition~\ref{sr1-twa-0} gives the following theorem.
The equivalent coefficient elements are fixed before the approximation
and absorb any prescribed comparison multiplicity.

\begin{thm}\label{sr1-wrp-cos-one}
Let $(A,\Gamma)$ be a minimal \Cstar{}-dynamical system with WRP, where
$A$ is unital and separable, has stable rank one, and
$\mathrm T(A)\ne\emptyset$, and where $\Gamma$ is countably infinite,
discrete, and amenable.  Assume that every nonzero hereditary
\Cstar{}-subalgebra of $A$ is infinite-dimensional.  If $(A,\Gamma)$ has
$(\lambda,m)$-COS for some $\lambda\in(0,1]$ and
$m\in\mathbb N$, then
\[
           \operatorname{tsr}(A\rtimes_{\alpha,r}\Gamma)=1.
\]
\end{thm}

\begin{proof}
Put $B=A\rtimes_{\alpha,r}\Gamma$.  As in the preceding proposition,
$B$ is simple and finite.  Let $x\in B$ be noninvertible.  By
Proposition~\ref{alg-red}, it suffices, up to multiplication by unitaries
and an arbitrarily small perturbation, to consider an element $b\in B$ for which
there is a norm-one positive contraction $d_0\in A$ satisfying
\[
                         d_0b=bd_0=0.
\]
Put $d=g_{1/2}(d_0)$ and $h_0=(d_0-1/2)_+$.  Then $d$ has
norm one, $db=bd=0$, and $dh_0=h_0$.  Apply
Lemma~\ref{coefficient-multiplicity-reserves} to obtain
a norm-one positive contraction $a\in A$ and fixed monomials
$t_1,\ldots,t_m\in B$ such that
\[
                 t_i^*t_j=\delta_{ij}a,\qquad dt_i=t_i.
\]
All these elements are chosen before the tracial approximation.

Let $\mathcal S$ be the class in Theorem~\ref{WTA-C-dyn-1}.  Stable
rank one passes to hereditary \Cstar{}-subalgebras, matrix algebras, and
finite direct sums, so every member of $\mathcal S$ has stable rank
one.  The same theorem gives $B\in\mathrm{TA}_0(\mathcal S)$, with the
approximating contraction in $A$ and with the approximants to elements
of $A$ still in $A$.

It remains to establish the comparison used in the proof of
Proposition~\ref{sr1-twa-0}.  Since every trace on $B$ is faithful, compactness
gives fixed constants
\[
 \delta_d=\min_{\tau\in\mathrm T(B)}\tau(g_{1/2}(d))>0,
 \qquad
 \delta_a=\min_{\tau\in\mathrm T(B)}\tau(a)>0.
\]
Fix $0<\theta<\delta_d/64$ and put $\kappa=\delta_a/4$.
Apply Theorem~\ref{WTA-C-dyn-1} with the finite set enlarged to
include $t_1,\ldots,t_m$.  Use Lemma~\ref{alt-defn-was} to obtain
the exact compressed approximants, with $d'\in A$ chosen to be a
positive contraction as in the proof of Theorem~\ref{sr1-twa}.
The functional-calculus cutoff in that lemma preserves arbitrarily
small commutators with the fixed columns.  Choose the tolerances
sufficiently small for the estimates in the proof of
Proposition~\ref{sr1-twa-0} and also require
\[
             \eta<\min\{\delta_a/4,\lambda\delta_a/2\}.
\]
Following the notation of Proposition~\ref{sr1-twa-0}, put
\[
 r=g_{3/4}(h),\qquad D=r^{1/2}d'r^{1/2},\qquad
 e=1-g_{1/3}(h),\qquad y=(D-2\theta)_+.
\]
Both $e$ and $y$ belong to $A_+$.  The enlarged finite set ensures
small commutators of $h$ with every $t_i$ and, by functional calculus,
of $r^{1/2}=g_{3/4}(h)^{1/2}$ with every $t_i$.
We construct a coefficient element whose $m$-fold direct sum is Cuntz
subequivalent to $y$.

Put
\[
 q_0=4(h-3/4)_+,\qquad
 W=(t_1q_0,\ldots,t_mq_0)\in \mathrm M_{1,m}(B).
\]
Then $r^{1/2}q_0=q_0$ and
\[
           W^*W=(q_0aq_0)^{\oplus m},\qquad \|W\|\leq1.
\]
For each column, the exact relation $dt_i=t_i$ gives
\[
 \|(1-D)t_iq_0\|
 \leq2\|[r^{1/2},t_i]\|+\|d'-d\|.
\]
With the columns $t_i$ fixed, $\|(1-D)W\|$ tends to zero with the
approximation tolerances.  Set
\[
                         V=g_{4\theta}(D)W.
\]
The cutoff vanishes on $[0,2\theta]$, so
$VV^*\in\operatorname{Her}_B(y)$.  Moreover,
\[
 \|V-W\|\leq\frac{1}{1-4\theta}\|(1-D)W\|.
\]
Indeed, $(1-g_{4\theta}(t))/(1-t)$ extends continuously to
$[0,1]$ and has norm at most $(1-4\theta)^{-1}$.
Since $V$ and $W$ are contractions, a further reduction of the
preliminary tolerances therefore ensures
\[
                         \|V^*V-W^*W\|<\kappa.
\]
The positive-element perturbation lemma gives, with
$z=(q_0aq_0-\kappa)_+\in A_+$,
\[
       z^{\oplus m}\precsim_B V^*V
            \sim VV^*\precsim_B y.
\]
All comparisons here are understood in matrix algebras over $B$.

The target $z$ retains a uniform positive trace bound.  Since
$q_0$ is a contraction and its defining function takes the value one
at one,
\[
 \tau(1-q_0^2)\leq \mathrm d_\tau(1-h)<\eta,
 \qquad \tau\in\mathrm T(B).
\]
Consequently,
\[
 \mathrm d_\tau(z)\geq\tau(z)
 \geq\tau(q_0aq_0)-\kappa
 >\delta_a-\eta-\kappa>\frac{\delta_a}{2}.
\]
By the trace identification in Proposition~\ref{trace-factor},
these bounds also hold for every invariant trace $\rho$ on $A$.
The defect estimate now gives
\[
 \mathrm d_\rho(e)<\eta<\lambda\frac{\delta_a}{2}
                 <\lambda \mathrm d_\rho(z),\qquad
 \rho\in\mathcal M_1(A,\Gamma).
\]
Property $(\lambda,m)$-COS and the preceding comparisons give
\[
                    e\precsim_B z^{\oplus m}\precsim_B y.
\]
Since $e$ and $y$ belong to $B$, this is also the
unstabilized comparison $e\precsim_B y$.
The unitary rotation in the proof of Proposition~\ref{sr1-twa-0}
therefore gives $b\in\overline{\mathrm{GL}(B)}$.  Multiplying by the
inverses of the two unitaries and letting the initial perturbation tend
to zero proves that every noninvertible
element of $B$ lies in $\overline{\mathrm{GL}(B)}$.
\end{proof}

\begin{rem}
The columns and the constants $\delta_a,\theta,\kappa$ are fixed
before the approximation.  Approximate commutation with the columns,
followed by a fixed spectral cut, gives the comparison of the
$m$-fold source with $y$.  Theorem~\ref{sr1-wrp-cos-D} uses inner
subtowers to obtain the same conclusion without assuming stable rank
one for the coefficient algebra.
In Lemma~\ref{pre-alg-red}, infinitude of $\Gamma$ is used to choose
a tower level outside the Fourier support of $b_0$.
\end{rem}

\section{Stable rank of \texorpdfstring{$A\rtimes\Gamma$}{A crossed with Gamma}}

For a unital \Cstar{}-algebra $D$, write $\mathrm T_\omega(D)$ for the
set of \emph{limit traces} on $D_\omega$, that is, the traces of the
form
\[
 \tau_\omega([(x_j)])=\lim_{j\to\omega}\tau_j(x_j),
 \qquad \tau_j\in\mathrm T(D).
\]
The WRP defect is uniformly small over traces on $D$, so every limit
trace annihilates it.  An arbitrary trace on the norm ultrapower need
not annihilate this defect; thus the estimates below concern
$\mathrm T_\omega(D)$ rather than the possibly larger space
$\mathrm T(D_\omega)$.

\begin{lem}\label{cut-pres-trace-lbd}
Let $A$ be a unital \Cstar{}-algebra, let $a\in A_+$, and suppose that
there are $x_1,\ldots,x_m\in A$ such that
\[
                 \sum_{i=1}^m x_i a x_i^*=1.
\]
For every $0<\eps<1$ there is $\eta>0$ with the following property.
If $c\in A_+$ and $\tau\in\mathrm T(A)$ satisfy
\[
 \tau\bigl(|x_i c-cx_i|^2\bigr)\leq\eta\tau(c^2),
 \qquad 1\leq i\leq m,
\]
then
\[
 \tau(cac)\geq
 \frac{1-\eps}{\sum_{i=1}^m\|x_i\|^2}\,\tau(c^2).
\]
\end{lem}

\begin{proof}
Put $S=\sum_i\|x_i\|^2$.  The fullness identity implies $S>0$.
Choose $\eta>0$ so small that
\begin{equation}\label{cut-pres-trace-choice}
 2\|a\|\eta^{1/2}\sum_{i=1}^m\|x_i\|\leq\eps.
\end{equation}
Fix $c$ and $\tau$ as in the statement, and write
$\|z\|_{2,\tau}=\tau(z^*z)^{1/2}$.  For
$r_i=cx_i$ and $s_i=x_ic$, traciality gives
\[
 \|r_i\|_{2,\tau},\ \|s_i\|_{2,\tau}
       \leq\|x_i\|\tau(c^2)^{1/2},
 \qquad
 \|r_i-s_i\|_{2,\tau}\leq\eta^{1/2}\tau(c^2)^{1/2}.
\]
Since
\[
 r_iar_i^*-s_ias_i^*=(r_i-s_i)ar_i^*+s_i a(r_i-s_i)^*,
\]
the tracial Cauchy--Schwarz inequality yields
\[
 \left|\tau(cx_i ax_i^*c)-\tau(x_i c a c x_i^*)\right|
 \leq 2\|a\|\|x_i\|\eta^{1/2}\tau(c^2).
\]
Summing and applying \eqref{cut-pres-trace-choice} gives
\begin{align*}
 S\tau(cac)
 &\geq\sum_{i=1}^m\tau(x_i c a c x_i^*)\\
 &\geq\sum_{i=1}^m\tau(cx_i a x_i^*c)-\eps\tau(c^2)\\
 &=(1-\eps)\tau(c^2),
\end{align*}
which is the desired estimate.
\end{proof}

\begin{lem}\label{trace-comm}
Let $A$ be a unital \Cstar{}-algebra, let $\tau\in\mathrm T(A)$, and let
$E\in A$ be a positive contraction satisfying
$\tau(E-E^2)=0$.  If $a\in A_+$ commutes with $E$ and
$f\colon\mathbb R\to\mathbb R$ is continuous with $f(0)=0$, then
\[
                        \tau(f(EaE))=\tau(Ef(a)E).
\]
\end{lem}

\begin{proof}
Let $\pi_\tau$ be the GNS representation associated to $\tau$.
The positive element $E-E^2$ belongs to the trace kernel, and hence
\[
                         p:=\pi_\tau(E)
\]
is a projection.  Since $E$ commutes with $a$, the projection $p$
commutes with $\pi_\tau(a)$.  Functional calculus and $f(0)=0$
give
\[
 \pi_\tau(f(EaE))
 =f(p\pi_\tau(a)p)
 =p f(\pi_\tau(a))p
 =\pi_\tau(Ef(a)E).
\]
Applying the GNS trace proves the assertion.
\end{proof}

The following lemma is a grouped form of Lemma~\ref{tower-order-zero-map}
and a noncommutative counterpart of \cite[Lemma~2.18]{LN-sr1}.

\begin{lem}\label{existence-0-map}
Let $F\subseteq\Gamma$ be finite and let $c\in A'\cap A_\omega$ be
a positive contraction such that the translates $g(c)$, $g\in F$,
are mutually orthogonal.  Let
$F_1,\ldots,F_r\subseteq F$ be pairwise disjoint nonempty sets of the
same cardinality.  There is a c.p.c. order zero map
\[
 \phi\colon \mathrm M_r(\Comp)\longrightarrow(A\rtimes\Gamma)_\omega
\]
such that
\[
                  \phi(e_{ii})=\sum_{g\in F_i}g(c),
                  \qquad 1\leq i\leq r.
\]
\end{lem}

\begin{proof}
Write $|F_i|=q$ and enumerate
$F_i=\{g_{i,1},\ldots,g_{i,q}\}$.  If $c\ne0$, apply Lemma~\ref{tower-order-zero-map} to the normalized base $c/\|c\|$ and
the shape $\bigsqcup_iF_i$, ordered lexicographically, and multiply
the resulting map by $\|c\|$.  This gives a c.p.c. order zero map
$\psi$ from $\mathrm M_{rq}(\Comp)$ with diagonal images $g_{i,j}(c)$.
If $c=0$, take $\psi=0$ instead.  Under the identification
$\mathrm M_{rq}\cong \mathrm M_r\otimes \mathrm M_q$, let
\[
 \iota\colon \mathrm M_r\longrightarrow \mathrm M_r\otimes \mathrm M_q,
 \qquad \iota(x)=x\otimes1_q.
\]
Then $\phi=\psi\circ\iota$ is c.p.c. and order zero, and
\[
 \phi(e_{ii})
 =\sum_{j=1}^q\psi(e_{(i,j),(i,j)})
 =\sum_{j=1}^q g_{i,j}(c).
\]
\end{proof}

\begin{lem}\label{pre-ord0-map}
Let $(A,\Gamma)$ be a unital \Cstar{}-dynamical system, where $\Gamma$ is
countably infinite, discrete, and amenable, and put
$B=A\rtimes_{\alpha,r}\Gamma$.  Let $a,b\in A_+$ and suppose that
\[
 \sum_{i=1}^p x_i a x_i^*=1,
 \qquad
 \sum_{j=1}^q y_j b y_j^*=1
\]
for some $x_1,\ldots,x_p,y_1,\ldots,y_q\in B$.  Set
\[
 \delta=\frac{1}{2\left(\sum_{i=1}^p\|x_i\|^2+
                              \sum_{j=1}^q\|y_j\|^2\right)}.
\]
Fix $\lambda\in(0,1]$ and a finite symmetric set
$K\subseteq\Gamma$ containing $e$, and fix a finite set
$G\subseteq\Gamma$ containing $e$.
Put $\beta=\min\{\delta/16,1/4\}$.
Then there are a finite set
$L\subseteq\Gamma$, a number $\eta>0$, and $n\in\mathbb N$ such that
whenever $F$ is $(L,\eta)$-invariant and $c\in A'\cap A_\omega$ is
a positive contraction whose $F$-translates are mutually orthogonal,
there is
a c.p.c. order zero map
\[
                \phi\colon \mathrm M_{n^2}(\Comp)\longrightarrow B_\omega.
\]
Writing
\[
 h_F=\sum_{g\in F}g(c),\qquad
 e_i=\phi(e_{ii}),\qquad
 E_k=\sum_{j=1}^n e_{n(k-1)+j},\qquad 1\leq k\leq n,
\]
one has $e_i\in A'\cap A_\omega$ and:
\begin{enumerate}
\item
\[
 \phi(1)=\sum_{i=1}^{n^2}e_i\leq h_F,\qquad
 \tau\left(h_F-\phi(1)\right)
 \leq\frac{\lambda\delta}{16n}\tau(h_F),\qquad
 \tau\in\mathrm T_\omega(B);
\]
\item for $1\leq k\leq n$ and $\tau\in\mathrm T_\omega(B)$,
\[
 \tau(E_k a E_k)\geq\delta\tau(E_k^2),\qquad
 \tau(E_k b E_k)\geq\delta\tau(E_k^2);
\]
\item $E_{k_1}\perp\alpha_g(E_{k_2})$ whenever
$g\in K$ and $k_1\ne k_2$;
\item $\alpha_g(E_k)\leq h_F$ for $g\in K$ and $1\leq k\leq n$.
\item Put $R_k=\sum_{j=4}^n e_{n(k-1)+j}$.  There are inner
subtower sums $P_k,T_k\in A'\cap A_\omega$ such that
\[
 P_k\leq R_k,\qquad T_k\leq E_k,
\]
and, for every $g\in G$, the levels of $\alpha_g(P_k)$ belong to
the subtower defining $R_k$, and those of $\alpha_g(T_k)$ belong
to the subtower defining $E_k$.  Moreover,
\[
 \tau(R_k-P_k)\leq\beta\tau(E_k),\qquad
 \tau(E_k-T_k)\leq\beta\tau(E_k),
 \qquad \tau\in\mathrm T_\omega(B).
\]
Consequently,
\[
 g_{1/2}(R_k)\alpha_g(g_1(P_k))=\alpha_g(g_1(P_k)),\qquad
 g_{1/2}(E_k)\alpha_g(g_1(T_k))=\alpha_g(g_1(T_k))
 \quad(g\in G).
\]
\end{enumerate}
The integer $n$ may be required to exceed any prescribed bound.
\end{lem}

\begin{proof}
We first derive a uniform commutator estimate from Følner invariance.
The estimate applies to tower sums whose bases are arbitrary positive
contractions, including zero.  Given a finite set
$\mathcal X\subseteq B$ and $\zeta>0$, there are a finite symmetric
set $L_0\subseteq\Gamma$ and $\eta_0>0$ such that, whenever $H$ is
$(L_0,\eta_0)$-invariant and $c\in A'\cap A_\omega$ is a positive
contraction with orthogonal $H$-translates, its support
$h_H=\sum_{g\in H}g(c)$ satisfies
\begin{equation}\label{tower-relative-commutator}
 \tau\bigl(|[x,h_H]|^2\bigr)\leq\zeta\tau(h_H^2),\qquad
 x\in\mathcal X,\quad\tau\in\mathrm T_\omega(B).
\end{equation}
Indeed, first consider a Fourier monomial $x=d u_s$.  Since $h_H$
commutes with $A$,
\[
 \|[d u_s,h_H]\|_{2,\tau}
 \leq\|d\|\,\|\alpha_s(h_H)-h_H\|_{2,\tau}.
\]
If $s,s^{-1}\in L_0$, tower orthogonality and traciality give
\begin{align*}
 \|\alpha_s(h_H)-h_H\|_{2,\tau}^2
 &=2\tau(h_H^2)-2\tau(h_H\alpha_s(h_H))\\
 &\leq2\bigl(|H|-|H\cap sH|\bigr)\tau(c^2)\\
 &\leq2\eta_0|H|\tau(c^2)
 =2\eta_0\tau(h_H^2).
\end{align*}
The middle inequality uses the common levels indexed by $H\cap sH$;
all other terms have nonnegative trace because they are products of
positive elements.  The estimate does not require $h_H$ and
$\alpha_s(h_H)$ to commute.
Approximate the finitely many elements of $\mathcal X$ by Fourier
polynomials.  The preceding estimate, the triangle inequality, and
\[
 \|[x-z,h_H]\|_{2,\tau}
 \leq2\|x-z\|\tau(h_H^2)^{1/2}
\]
then prove \eqref{tower-relative-commutator} after enlarging $L_0$ and
decreasing $\eta_0$.

Apply Lemma~\ref{cut-pres-trace-lbd}, with $\eps=1/2$, first to
$(a,x_1,\ldots,x_p)$ and then to $(b,y_1,\ldots,y_q)$.  For the
resulting tolerances, apply \eqref{tower-relative-commutator}.  Thus
$L_0$ and $\eta_0$ may be chosen so that every
$(L_0,\eta_0)$-invariant $H$ and central positive contraction $c$
with orthogonal $H$-translates satisfy
\begin{equation}\label{two-tower-trace-lower-bounds}
 \tau(h_H a h_H)\geq\delta\tau(h_H^2),\qquad
 \tau(h_H b h_H)\geq\delta\tau(h_H^2),\qquad
 \tau\in\mathrm T_\omega(B).
\end{equation}
Enlarge $L_0$ to contain $G$, and decrease $\eta_0$ so that
$\eta_0\leq\beta$.  These changes preserve
\eqref{two-tower-trace-lower-bounds}.

Choose $n>3$ such that
\begin{equation}\label{large-n-0}
 \frac{1}{n-3}<\frac{\delta}{24},\qquad
 \frac{\lambda\delta}{16n}<\frac12,\qquad
 \frac3n<\frac{\delta}{16},\qquad
 \frac{32}{n}<\frac{\lambda\delta}{4}.
\end{equation}
Put $\theta=\lambda\delta/(32n)$.  We use the left-Følner,
right-translate versions of Lemmas~7.1 and~7.2 of \cite{LN-sr1}.
They follow from the stated right-Følner, left-translate versions by
passing to inverses (equivalently, to the opposite group).  In
particular, the conclusion corresponding to Lemma~7.2 is
$K H_i\subseteq F_i$, consistent with the fact that
$\alpha_g$ moves a tower level indexed by $r$ to the level indexed by
$gr$.

By the exact tiling theorem \cite[Theorem~4.3]{DHZ-tiling}, in the
same inverse convention there is a right-translate tiling of $\Gamma$
with shapes $P_1,\ldots,P_t$ that are
$(L_0,\eta_0)$-invariant.  Apply the inverse version of
\cite[Lemma~7.2]{LN-sr1} to $\theta,n,K$ and these shapes, obtaining
a pair $(L_1,\eta_1)$.  Use the exact tiling theorem again to obtain a
right-translate tiling with $(L_1,\eta_1)$-invariant shapes
$Q_1,\ldots,Q_s$.  Finally, the inverse version of
\cite[Lemma~7.1]{LN-sr1}, applied to $\theta,n$ and the $Q_l$, gives
a pair $(L,\eta)$.

Fix $F$ and $c$ as in the statement.  Lemma~7.1
provides $R_1\subseteq F$ such that
\begin{equation}\label{small-R-1}
 |R_1|<\theta|F|,
\end{equation}
and $F\setminus R_1$ is tiled by right translates of the $Q_l$, with
the multiplicity of every shape divisible by $n$.  Distributing each
kind of tile evenly gives a decomposition
\[
                  F\setminus R_1=F_1\sqcup\cdots\sqcup F_n
\]
in which the $F_i$ are $(L_1,\eta_1)$-invariant and have equal
cardinality.  Lemma~7.2 then supplies sets $H_i\subseteq F_i$ such
that
\begin{equation}\label{trimmed-row-properties}
 K H_i\subseteq F_i,\qquad |H_1|=\cdots=|H_n|,\qquad
 |H_i|>(1-\theta)|F_i|,
\end{equation}
and every $H_i$ is tiled by right translates of the $P_l$, with all
multiplicities divisible by $n$.  If
\[
 R_2=(F\setminus R_1)\setminus\bigsqcup_{i=1}^nH_i,
\]
then
\begin{equation}\label{small-R-2}
                         |R_2|<\theta|F|.
\end{equation}

Divide, for each $i$, the tiles of every shape occurring in $H_i$
evenly into $n$ groups.  This gives
\[
                  H_i=H_{i,1}\sqcup\cdots\sqcup H_{i,n}.
\]
All $n^2$ sets $H_{i,j}$ have the same cardinality.  Moreover, each
is $(L_0,\eta_0)$-invariant, since left invariance is preserved by right
translation and disjoint unions.  Define
\[
 e_{n(i-1)+j}=\sum_{g\in H_{i,j}}g(c),\qquad 1\leq i,j\leq n.
\]
Lemma~\ref{existence-0-map} gives a c.p.c. order zero map
$\phi\colon \mathrm M_{n^2}(\Comp)\to B_\omega$ with these diagonal images.
Each $e_l$ belongs to $A'\cap A_\omega$.

For each cell define
\[
 J_{i,j}=\operatorname{Int}_G(H_{i,j}).
\]
Since $G\subseteq L_0$ and $H_{i,j}$ is
$(L_0,\eta_0)$-invariant, we have
\[
 |H_{i,j}\setminus J_{i,j}|<\beta|H_{i,j}|,
 \qquad GJ_{i,j}\subseteq H_{i,j}.
\]
Set
\[
 P_i=\sum_{j=4}^n\sum_{g\in J_{i,j}}g(c),\qquad
 T_i=\sum_{j=1}^n\sum_{g\in J_{i,j}}g(c).
\]
These are the required inner subtower sums.  Every trace on
$B_\omega$ has the same value on all translates of $c$, so
the cardinality estimates give the trace bounds in (5).
The absorption identities in (5) follow level by level from
$g_{1/2}(t)g_1(t)=g_1(t)$ and tower orthogonality.

The levels omitted from $\phi(1)$ are indexed by $R_1\sqcup R_2$.
Thus, for $\tau\in\mathrm T_\omega(B)$, inner invariance of $\tau$
and \eqref{small-R-1}--\eqref{small-R-2} give
\begin{align*}
 \tau(h_F-\phi(1))
 &=|R_1\sqcup R_2|\tau(c)\\
 &\leq2\theta|F|\tau(c)
 =\frac{\lambda\delta}{16n}\tau(h_F).
\end{align*}
For each $i$, one has
$E_i=\sum_{g\in H_i}g(c)=h_{H_i}$.  Equation
\eqref{two-tower-trace-lower-bounds} therefore gives
\[
 \tau(E_i a E_i)\geq\delta\tau(E_i^2),\qquad
 \tau(E_i b E_i)\geq\delta\tau(E_i^2),
\]
and hence assertion (2).

Finally, \eqref{trimmed-row-properties} and the disjointness of the
$F_i$ show that $gH_{k_2}$ is disjoint from $H_{k_1}$ whenever
$g\in K$ and $k_1\ne k_2$.  Tower orthogonality proves (3).  The
inclusion $gH_k\subseteq F_k\subseteq F$ gives (4).
\end{proof}

The hypothesis on hereditary subalgebras in Proposition~\ref{alg-red}
follows from minimality and WRP when $\Gamma$ is infinite.

\begin{lem}\label{no-atoms}
Let $(A,\Gamma)$ be a minimal \Cstar{}-dynamical system with WRP, where
$A$ is unital and separable, $\mathrm T(A)\ne\varnothing$, and
$\Gamma$ is countably infinite, discrete, and amenable.  Then every
nonzero hereditary \Cstar{}-subalgebra of $A$ is infinite-dimensional.
\end{lem}

\begin{proof}
Amenability gives an invariant tracial state $\rho$ on $A$.  The trace
$\rho\circ\mathbb E$ on $A\rtimes_{\alpha,r}\Gamma$ is faithful by
Corollary~\ref{min-faithful-trace}; in particular, $\rho$ is faithful.

Suppose that $p\in A$ is a nonzero minimal projection.  Given $N$, choose
distinct elements $g_1=e,g_2,\ldots,g_N\in\Gamma$.  Lemma~\ref{WRP-Kishimoto} and Proposition~\ref{Elliott-central-sequence},
applied to $p$ and $F=\{g_1,\ldots,g_N\}$, give a tower $(c,F)$ such
that
\[
                         \|cpc\|=1.
\]
Since $c\in A'\cap A_\omega$ and
$pA_\omega p=(pAp)_\omega=\Comp p$, we have $cpc=p$ and hence
$cp=pc=p$.  Applying $g_i$ gives
\[
             g_i(c)g_i(p)=g_i(p)g_i(c)=g_i(p).
\]
The tower levels $g_i(c)$ are mutually orthogonal, so the projections
$g_1(p),\ldots,g_N(p)$ are mutually orthogonal.  Invariance of $\rho$
therefore gives
\[
                         N\rho(p)\leq1.
\]
This is impossible for arbitrarily large $N$, since $\rho(p)>0$.
Thus $A$ has no minimal projections.

If a nonzero hereditary \Cstar{}-subalgebra $H\subseteq A$ were
finite-dimensional, it would contain a minimal projection $p$.  Since
$pAp\subseteq H$, this projection would also be minimal in $A$, a
contradiction.
\end{proof}

The equivalent coefficient elements also give the property~(D)
argument for arbitrary fixed comparison multiplicity.

\begin{thm}\label{sr1-wrp-cos-D}
Let $(A,\Gamma)$ be a minimal \Cstar{}-dynamical system with WRP, where
$A$ is unital and separable, $\mathrm T(A)\ne\varnothing$, and
$\Gamma$ is countably infinite, discrete, and amenable.  If
$(A,\Gamma)$ has $(\lambda,m)$-COS for some $\lambda\in(0,1]$
and $m\in\mathbb N$, then
\[
              \operatorname{tsr}(A\rtimes_{\alpha,r}\Gamma)=1.
\]
\end{thm}

\begin{proof}
Put $B=A\rtimes_{\alpha,r}\Gamma$.  By Corollary~\ref{wrp-simple}, $B$ is simple.  Let $z\in B$ be noninvertible and
fix $\eps>0$.  Lemma~\ref{no-atoms} permits us to apply Proposition~\ref{alg-red}.  We obtain unitaries $u_0,u_1\in B$, an element
$x\in\mathrm C_{\mathrm c}(\Gamma,A)$, and a norm-one positive
contraction $b_0\in A$ such that
\begin{equation}\label{D-zero-divisor-reduction}
       \|u_0zu_1-x\|<\eps,\qquad b_0x=xb_0=0.
\end{equation}

The hereditary algebra $\overline{b_0Ab_0}$ is infinite-dimensional,
so it contains orthogonal norm-one positive contractions $c_0,d_0$.
Put
\[
 c=g_{1/2}(c_0),\qquad d=g_{1/2}(d_0),\qquad
 h_c=(c_0-1/2)_+,\qquad h_d=(d_0-1/2)_+.
\]
These satisfy $ch_c=h_c$ and $dh_d=h_d$; all four elements
annihilate $x$ on both sides.
Apply Lemma~\ref{coefficient-multiplicity-reserves}
separately to these two corners.  Before choosing any tower parameters,
we obtain norm-one positive contractions $a,b\in A$ and fixed
monomials
\[
 t_i\in A u_{\gamma_i},\qquad \ell_i\in A u_{\nu_i},
 \qquad 1\leq i\leq m,
\]
such that
\[
 t_i^*t_j=\delta_{ij}a,\qquad ct_i=t_i,
 \qquad
 \ell_i^*\ell_j=\delta_{ij}b,\qquad d\ell_i=\ell_i.
\]
Their range elements are mutually orthogonal within each family and
belong to $\operatorname{Her}_A(h_c)$ or
$\operatorname{Her}_A(h_d)$, respectively.
The nonzero elements $g_{1/2}(a)$ and $g_{1/2}(b)$ are full in
the simple unital algebra $B$.  Hence there are
$v_1,\ldots,v_p,w_1,\ldots,w_q\in B$
such that
\begin{equation}\label{D-fullness-identities}
 \sum_{i=1}^p v_i g_{1/2}(a)v_i^*=1,
 \qquad
 \sum_{j=1}^q w_j g_{1/2}(b)w_j^*=1.
\end{equation}
Put
\[
 \delta=\frac{1}{2\left(\sum_i\|v_i\|^2+\sum_j\|w_j\|^2\right)}.
\]
Write
\[
                     x=\sum_{g\in K}x_gu_g,
\]
where $K$ is finite, symmetric, and contains $e$.

Put
\[
 G=\{e,\gamma_1,\ldots,\gamma_m,\nu_1,\ldots,\nu_m\}.
\]
Apply Lemma~\ref{pre-ord0-map} to $g_{1/2}(a)$,
$g_{1/2}(b)$, $\lambda$, $K$, and $G$, and denote the resulting data by
$L,\eta,n$.  By the choice made in \eqref{large-n-0},
\begin{equation}\label{D-large-n}
 \frac3n<\frac{\delta}{16},\qquad
 \frac{32}{n^2}<\frac{\lambda\delta}{4n}.
\end{equation}
Choose, by the exact tiling theorem
\cite[Theorem~4.3]{DHZ-tiling}, an exact tiling whose shapes
$F_1,\ldots,F_S$ are $(L,\eta)$-invariant.  WRP gives bases
$r_s\in A'\cap A_\omega$ such that all the levels $g(r_s)$,
$g\in F_s$, are mutually orthogonal and, with
\[
 Q_s=\sum_{g\in F_s}g(r_s),\qquad Q=\sum_{s=1}^S Q_s,
\]
one has $1-Q\in J_{\mathcal M_1(A,\Gamma)}$.

For each $s$, Lemma~\ref{pre-ord0-map} gives a c.p.c. order zero map
\[
                   \Phi_s:\mathrm M_{n^2}(\Comp)\longrightarrow B_\omega.
\]
The ranges of these maps are mutually orthogonal, so
$\Phi=\sum_s\Phi_s$ is again c.p.c. and order zero.  Put
\[
 e_i=\Phi(e_{ii}),\qquad H=\Phi(1),\qquad
 E_k=\sum_{j=1}^n e_{n(k-1)+j}.
\]
Then $e_i\in A'\cap A_\omega$, $H\leq Q$, and, for every
$\tau\in\mathrm T_\omega(B)$,
\begin{align}
 \tau(Q-H)&\leq\frac{\lambda\delta}{16n},
       \label{D-discarded-trace}\\
 \tau(E_k g_{1/2}(a)E_k)&\geq\delta\tau(E_k^2),
       \label{D-row-c-trace}\\
 \tau(E_k g_{1/2}(b)E_k)&\geq\delta\tau(E_k^2).
       \label{D-row-d-trace}
\end{align}
Proposition~\ref{trace-factor} identifies the coordinate traces with
invariant traces on $A$.

The image of $Q$ in the GNS representation of any such $\tau$ is the
unit.  Since the tower levels are orthogonal and $Qg(r_s)=g(r_s)^2$,
every tower level is a projection in that representation.  The same is
therefore true of $e_i$, $H$, and $E_k$.  A trace composed with an order
zero map has equal values on the diagonal matrix units, and hence
\begin{equation}\label{D-row-mass}
 \tau(E_k^2)=\frac1n\tau(H)
 \geq\frac1n\left(1-\frac{\lambda\delta}{16n}\right)
 >\frac{1}{2n}.
\end{equation}

For $0<t<1$, let $\Phi_t=g_t(\Phi)$ and write
\[
 e_i^{(t)}=\Phi_t(e_{ii}),\qquad
 E_k^{(t)}=\sum_{j=1}^n e_{n(k-1)+j}^{(t)},\qquad
 R_k^{(t)}=\sum_{j=4}^n e_{n(k-1)+j}^{(t)}.
\]
Set
\begin{align*}
 \phi&=\Phi_{1/4},& h&=\phi(1),\\
 s_k&=\sum_{j=1}^3e_{n(k-1)+j}^{(1/4)},
 &\widehat E_k&=E_k^{(1/4)},\\
 \widetilde c_k&=g_{1/2}(R_k^{(1/2)}cR_k^{(1/2)}),
 &c_k&=g_{1/4}(R_k^{(1/2)}cR_k^{(1/2)}),\\
 \widetilde d_k&=g_{1/2}(E_k^{(1/2)}dE_k^{(1/2)}),
 &d_k&=g_{1/4}(E_k^{(1/2)}dE_k^{(1/2)}).
\end{align*}
All these elements belong to $A_\omega$ except for the order zero map,
whose range lies in $B_\omega$.

For each WRP tower, take the inner subtower sums in
Lemma~\ref{pre-ord0-map}(5).  Summing over the towers gives
$P_k\leq R_k$ and $T_k\leq E_k$, where
$R_k=E_k-\sum_{j=1}^3e_{n(k-1)+j}$.
Put
\[
 q_k^c=g_1(P_k),\qquad q_k^d=g_1(T_k),\qquad
 A_k=q_k^c a q_k^c,\qquad B_k=q_k^d b q_k^d.
\]
Define the inner and outer cutoffs of the comparison targets by
\[
 \widetilde a_k=g_{1/2}(A_k),\qquad
 \widehat a_k=g_{1/4}(A_k),\qquad
 \widetilde b_k=g_{1/2}(B_k),\qquad
 \widehat b_k=g_{1/4}(B_k).
\]

We first compare the $m$-fold direct sums with the hereditary
subalgebras generated by $c_k$ and $d_k$.  Set
\[
 v_{i,k}=t_iq_k^c,\qquad w_{i,k}=\ell_iq_k^d,
 \qquad 1\leq i\leq m.
\]
The source relations are
\[
 v_{i,k}^*v_{j,k}=\delta_{ij}A_k,\qquad
 w_{i,k}^*w_{j,k}=\delta_{ij}B_k.
\]
Each range element belongs to $A_\omega$: if, for example,
$t_i=z_i u_{\gamma_i}$ with $z_i\in A$, then
\[
 v_{i,k}v_{i,k}^*=z_i\alpha_{\gamma_i}((q_k^c)^2)z_i^*.
\]
The range elements are orthogonal within each family.  The inner subtower
identities in Lemma~\ref{pre-ord0-map}(5), together with centrality
of the tower sums relative to $A$, give
\[
 R_k^{(1/2)}v_{i,k}=v_{i,k},\qquad cv_{i,k}=v_{i,k},\qquad
 E_k^{(1/2)}w_{i,k}=w_{i,k},\qquad dw_{i,k}=w_{i,k}.
\]
Thus $R_k^{(1/2)}cR_k^{(1/2)}$ fixes every $v_{i,k}$,
and $E_k^{(1/2)}dE_k^{(1/2)}$ fixes every $w_{i,k}$.
Their cutoff functions also fix these columns.  In particular, the
orthogonal equivalent range elements belong to the hereditary
algebras of $\widetilde c_k$ and $\widetilde d_k$, respectively.
Writing the columns as row matrices gives
\[
 \widehat a_k^{\oplus m}\precsim A_k^{\oplus m}
   \sim\sum_{i=1}^m v_{i,k}v_{i,k}^*\precsim c_k,
 \qquad
 \widehat b_k^{\oplus m}\precsim B_k^{\oplus m}
   \sim\sum_{i=1}^m w_{i,k}w_{i,k}^*\precsim d_k.
\]
These stabilized Cuntz comparisons have targets in the respective
hereditary subalgebras.

For the lower trace bounds on the targets, observe that in the GNS
representation of every limit trace, the tower sums are projections,
so their cutoffs
$q_k^c,q_k^d$ have the same images as $P_k,T_k$.
All these projections commute with $a$ and $b$.
By Lemma~\ref{trace-comm}, replacing a row by an inner subtower
decreases its trace by at most the trace of the omitted levels.  With
$\beta=\delta/16$, Lemma~\ref{pre-ord0-map}(5),
\eqref{D-row-c-trace}--\eqref{D-row-mass}, and
\eqref{D-large-n} therefore give
\begin{align}
 \tau(\widetilde a_k)
 &=\tau(P_k g_{1/2}(a)P_k)\nonumber\\
 &\geq\left(\delta-\frac3n-\beta\right)\tau(E_k^2)
 >\frac{7\delta}{16n}>\frac{\delta}{3n},
       \label{D-inner-c-lower}\\
 \tau(\widetilde b_k)
 &=\tau(T_k g_{1/2}(b)T_k)\nonumber\\
 &\geq(\delta-\beta)\tau(E_k^2)
 >\frac{15\delta}{32n}>\frac{\delta}{3n}.
       \label{D-inner-d-lower}
\end{align}
Since $a,b$ are contractions of norm one, the fullness identities
give $\delta\leq1/4$ and the value of $\beta$ in Lemma~\ref{pre-ord0-map} is $\delta/16$.

The corresponding dimension-function estimates at the coordinates
follow from the spectral cutoffs.  Choose positive-contraction representatives
$e_i=[(e_{i,j})]$ in $A$.  Since
\[
 \mathrm d_\rho(g_{1/4}(e_{i,j}))
 \leq 8\rho(e_{i,j}),
\]
subadditivity of dimension functions and equality of the diagonal
traces of $\Phi$ imply
\begin{equation}\label{D-small-sk-coordinate}
 \lim_{j\to\omega}\sup_{\rho\in\mathcal M_1(A,\Gamma)}
 \mathrm d_\rho(s_{k,j})\leq\frac{24}{n^2},
\end{equation}
where $(s_{k,j})$ is the corresponding representative of $s_k$.
Equations \eqref{D-inner-c-lower} and \eqref{D-inner-d-lower} imply,
for suitable positive representatives, that on an $\omega$-large set
of coordinates
\begin{equation}\label{D-target-coordinate-lower}
 \inf_{\rho\in\mathcal M_1(A,\Gamma)}
 \rho(\widetilde a_{k,j})>\frac{\delta}{4n},
 \qquad
 \inf_{\rho\in\mathcal M_1(A,\Gamma)}
 \rho(\widetilde b_{k,j})>\frac{\delta}{4n}.
\end{equation}

For the defect, choose positive-contraction representatives $H_j$ of
$H$.  The support of $1-g_{1/4}(H_j)$ is contained in the spectral
region $[0,1/4)$, and therefore
\[
 \mathrm d_\rho(1-g_{1/4}(H_j))
 \leq\frac{16}{9}\rho((1-H_j)^2).
\]
In every limit-trace GNS representation, $H$ is a projection and
\eqref{D-discarded-trace} gives
\[
 \tau((1-H)^2)=\tau(1-H)\leq\frac{\lambda\delta}{16n}.
\]
Consequently,
\begin{equation}\label{D-small-defect-coordinate}
 \lim_{j\to\omega}\sup_{\rho\in\mathcal M_1(A,\Gamma)}
 \mathrm d_\rho(1-g_{1/4}(H_j))
 \leq\frac{\lambda\delta}{9n}.
\end{equation}
The uniform estimates in
\eqref{D-small-sk-coordinate}--\eqref{D-small-defect-coordinate}
follow by contradiction.  Traces witnessing a failure at successive
coordinates define a limit trace in $\mathrm T_\omega(B)$ that
contradicts the ultrapower estimates.

By \eqref{D-large-n}, \eqref{D-small-sk-coordinate},
\eqref{D-target-coordinate-lower}, and
\eqref{D-small-defect-coordinate}, on an $\omega$-large set one has,
simultaneously for every invariant trace $\rho$,
\[
 \mathrm d_\rho(s_{k,j})<\lambda \mathrm d_\rho(\widetilde a_{k,j}),
 \qquad
 \mathrm d_\rho(1-g_{1/4}(H_j))
       <\lambda \mathrm d_\rho(\widetilde b_{k,j}).
\]
Property $(\lambda,m)$-COS gives, coordinatewise in $B$,
\[
 s_{k,j}\precsim\widetilde a_{k,j}^{\oplus m},\qquad
 1-g_{1/4}(H_j)\precsim\widetilde b_{k,j}^{\oplus m}.
\]
Choose positive-contraction representatives of $A_k,B_k$, and obtain
both their inner and outer cutoffs by coordinatewise functional
calculus.  The relations
\[
 \widehat a_{k,j}\widetilde a_{k,j}=\widetilde a_{k,j},\qquad
 \widehat b_{k,j}\widetilde b_{k,j}=\widetilde b_{k,j}
\]
then hold exactly at every coordinate.  Apply Lemma~\ref{comp-local-global} in $\mathrm M_m(B)$, identifying
$(\mathrm M_m(B))_\omega$ with $\mathrm M_m(B_\omega)$ and embedding scalar sources
in the first corner.  It gives
\[
 s_k\precsim\widehat a_k^{\oplus m}\precsim c_k,\qquad
 1-h\precsim\widehat b_k^{\oplus m}\precsim d_k.
\]
The last comparisons follow from the preceding $m$-fold comparisons.
Since the first and last elements belong to $B_\omega$, these
stabilized comparisons are also comparisons in $B_\omega$.
Thus
\begin{equation}\label{D-two-comparisons}
                         s_k\precsim c_k,
 \qquad                 1-h\precsim d_k.
\end{equation}

The remaining $\mathcal D_0$ relations of
\cite[Definition~5.1]{LN-sr1} follow from the cutoff identities and
tower orthogonality.  The cutoff functions satisfy
\[
                    g_{1/4}(t)g_{1/2}(t)=g_{1/2}(t),
\]
and hence
\[
             \widehat E_kc_k=c_k,\qquad
             \widehat E_kd_k=d_k.
\]
The subtower defining $c_k$ is orthogonal to the first three diagonal
entries, so $c_k\perp s_k$; also $c_k\perp d_k$ because $c\perp d$.
Both elements lie in $\overline{b_0B_\omega b_0}$.
The trace lower bounds show that $A_k$ and $B_k$ are nonzero;
hence there are nonzero columns $v_{i,k}$ and $w_{i,k}$ fixed by
$c_k$ and $d_k$.  Since $c_k,d_k$ are contractions, this proves
that they have norm one.  Summing the relations in
Lemma~\ref{pre-ord0-map}(3)--(4) over the WRP towers gives
\[
 \widehat E_{k_1}\alpha_g(\widehat E_{k_2})=0
 \quad(g\in K,\ k_1\ne k_2).
\]
Since the diagonal entries commute with $A$, the Fourier expansion of
$x$ gives
\[
                   \widehat E_{k_1}x\widehat E_{k_2}=0
                   \quad(k_1\ne k_2).
\]
Finally, with $p=q=n$, $r=3$, and $l=1$, the numerical condition in
the definition is $3n>2n$.  Together with
\eqref{D-zero-divisor-reduction} and \eqref{D-two-comparisons}, these
relations show that the constant copy of $x$ is a $\mathcal D_0$-element
of $B_\omega$.  By \cite[Proposition~5.4]{LN-sr1},
\[
                         x\in\overline{\mathrm{GL}(B_\omega)}.
\]

This ultrapower conclusion implies the same assertion in $B$.  If an
invertible $y\in B_\omega$ is close to $x$, lift $y$ and $y^{-1}$ to
bounded sequences $(y_j)$ and $(z_j)$ in $B$.  On an $\omega$-large
set, both $y_jz_j$ and $z_jy_j$ are within $1/2$ of the unit, so
$y_j$ is invertible; at one such coordinate it is as close to $x$ as
desired.  Thus $x\in\overline{\mathrm{GL}(B)}$.  Equation
\eqref{D-zero-divisor-reduction}, followed by multiplication by the
two unitaries, shows that $z$ is arbitrarily close to an invertible
element.  Hence $B$ has stable rank one.
\end{proof}

\begin{rem}
The equivalent coefficient elements and their monomial witnesses are
chosen before the fullness constant $\delta$ and the integer $n$.
Their group elements determine the inner subtowers.  The comparison
targets have trace at least $\delta/(3n)$, whereas $s_k$ has
dimension-function bound $24/n^2$.  The $m$ equivalent orthogonal
copies lie in a single hereditary subalgebra.  These comparisons
absorb the multiplicity in $(\lambda,m)$-COS without cancellation of
an amplified Cuntz comparison or an additional divisibility hypothesis.
\end{rem}

\section{$\mathcal Z$-stability of $A \rtimes \Gamma$ where $\mathrm{rc}(A) < \infty$}\label{finite-radius-regularity-section}

Corollary~\ref{finite-coefficient-radius-crossed-comparison} and
Theorem~\ref{sr1-wrp-cos-D} give Cuntz-semigroup regularity for crossed
products whose coefficient algebras have positive finite radius of
comparison.  Locally finite nuclear dimension then gives tensorial
absorption.

\begin{thm}[Finite coefficient radius of comparison]\label{finite-radius-crossed-Cu}
Let $A$ be a separable, unital, stably finite, nuclear \Cstar{}-algebra with
$\mathrm T(A)\neq\varnothing$, and let a countably infinite discrete amenable
group $\Gamma$ act minimally on $A$.  Suppose that at least one
of the following holds:
\begin{enumerate}
\item the action has WRP and $(\lambda,m)$-COS for some
      $\lambda\in(0,1]$ and $m\in\mathbb N$;
\item the action has strong URP;
\item the action has WRP with Cuntz control in the crossed product.
\end{enumerate}
If $\operatorname{rc}(A)<\infty$, then
\[
 B=A\rtimes_{\alpha,r}\Gamma
\]
is simple, has strict comparison of positive elements, and has stable
rank one.  Every $f\in\mathrm{LAff}_{++}(\mathrm T(B))$ is the rank function
of an element of $(B\otimes\mathcal K)_+$.  Moreover, the canonical
inclusion $\iota:B\to B\otimes\mathcal Z$, $b\mapsto b\otimes1$, induces
an isomorphism of ordered semigroups
\[
 \mathrm{Cu}(\iota):\mathrm{Cu}(B)
       \xrightarrow{\cong}\mathrm{Cu}(B\otimes\mathcal Z),
 \qquad
 \mathrm{Cu}(B)\cong
 \mathrm V(B)\sqcup\mathrm{LAff}_{++}(\mathrm T(B)).
\]
\end{thm}

\begin{proof}
Lemma~\ref{strong-urp-cuntz-controlled-wrp} and
Definition~\ref{cuntz-controlled-wrp-definition} show that each
alternative includes WRP.  In the second and third cases,
Propositions~\ref{orbit-finite-radius-zero} and
\ref{cuntz-controlled-wrp-orbit-zero-cos} give $(1,1)$-COS.
Thus all three alternatives imply WRP and COS.
Corollary~\ref{wrp-simple} gives simplicity of $B$.  Nuclearity of $A$
and amenability of $\Gamma$ give nuclearity, and hence exactness, of
$B$.  Amenability supplies an invariant trace on $A$, which extends to
a trace on $B$.  Simplicity makes this trace faithful, so $B$ is stably
finite.  The infinitely many canonical group unitaries are linearly
independent, so $B$ is non-elementary.

Corollary~\ref{finite-coefficient-radius-crossed-comparison} gives strict
comparison, and Theorem~\ref{sr1-wrp-cos-D} gives stable rank one.
Exactness identifies normalized
$2$-quasitraces with tracial states.  The rank-realization assertion is
therefore \cite[Theorem~8.11]{Thiel-sr1}; the canonical
Cuntz-semigroup isomorphism and its displayed description follow from
\cite[Corollary~8.12]{Thiel-sr1}.
\end{proof}

\begin{cor}[Finite coefficient radius of comparison and $\mathcal Z$-stability]
\label{wrp-finite-radius-Z}
Let $A$ be a separable, unital, stably finite \Cstar{}-algebra with
$\mathrm T(A)\neq\varnothing$, locally finite nuclear dimension, and
$\operatorname{rc}(A)<\infty$.  Let a countably infinite discrete
amenable group $\Gamma$ act minimally on $A$.  
Suppose that at least one
of the following holds:
\begin{enumerate}
\item the action has WRP and $(\lambda,m)$-COS for some
      $\lambda\in(0,1]$ and $m\in\mathbb N$;
\item the action has strong URP;
\item the action has WRP with Cuntz control in the crossed product.
\end{enumerate}
Then
\[
 B=A\rtimes_{\alpha,r}\Gamma
\]
is simple and $\mathcal Z$-stable.
\end{cor}

\begin{proof}
Locally finite nuclear dimension implies nuclearity of $A$.
Theorem~\ref{finite-radius-crossed-Cu} therefore applies and gives
simplicity of $B$ and the canonical ordered-semigroup description
\[
 \mathrm{Cu}(B)\cong
 \mathrm V(B)\sqcup\mathrm{LAff}_{++}(\mathrm T(B)).
\]
Theorem~\ref{wrp-local-nucdim-Cu-Z} gives
$B\cong B\otimes\mathcal Z$.
\end{proof}

\begin{rem}[Local dimension and the coefficient radius of comparison]
\label{finite-radius-centrality-gap}
The nuclear dimension hypothesis in Corollary~\ref{wrp-finite-radius-Z}
is local: each finite set is approximated inside a subalgebra of finite
nuclear dimension, with no common bound on these dimensions.
The coefficient algebra may have positive finite radius of comparison
and need not have strict comparison or be $\mathcal Z$-stable.
The local hypothesis supplies the central divisibility step through
Theorem~\ref{wrp-local-nucdim-Cu-Z}; without it,
Theorem~\ref{finite-radius-crossed-Cu} still gives Cuntz-semigroup
regularity, but this argument does not establish tensorial
absorption.

For an infinite group, this extends the conclusion of
Theorem~\ref{WRP-finite-nucdim-Z} beyond $\mathcal Z$-stable
coefficients.  Indeed, a simple, nuclear, stably finite,
$\mathcal Z$-stable coefficient algebra has finite nuclear dimension
by \cite{CETWW-dim-n} and radius of comparison zero.  Simplicity makes
every action minimal, and Corollary~\ref{Z=>COS} supplies COS.
Theorem~\ref{WRP-finite-nucdim-Z} also applies to finite groups.
\end{rem}

\begin{example}[A Villadsen algebra with positive radius and strong URP]
\label{villadsen-wrp-cos-example}
Motivated by Hirshberg's construction in
\cite{Hirshberg-exotic-circle}, we construct a simple unital Villadsen
algebra $A$ of the first type and an automorphism $\alpha$ such that
\[
 0<\operatorname{rc}(A)<1,
 \qquad (A,\mathbb Z,\alpha)\text{ has strong URP}.
\]
The algebra $A$ is not $\mathcal Z$-stable, whereas
$A\rtimes_\alpha\mathbb Z$ is simple and $\mathcal Z$-stable.
The construction determines the algebra and action simultaneously,
using one connected homogeneous building block at each stage.
\end{example}

\begin{proof}
\emph{The equivariant inductive system.}
Start with $X_0=[0,1]^2$, $r_0=q_0=1$, $\sigma_0=\mathrm{id}$,
and $u_0=1$.  Inductively, write
\[
 A_n=\mathrm M_{r_n}(\mathrm C(X_n)),\qquad
 \alpha_n(a)(x)=u_na(\sigma_n^{-1}x)u_n^*,
 \qquad \sigma_n^{q_n}=\mathrm{id},\quad u_n^{q_n}=1,
\]
where $\sigma_n$ permutes the copies of the seed cube and $u_n$ is a
constant unitary.  Give all products the maximum product metric.
Choose a finite $\sigma_n$-invariant $1/(n+1)$-net $E_n\subset X_n$
consisting of disjoint orbits of length $q_n$, and put $l_n=|E_n|$.
Such nets exist: each nonidentity power $\sigma_n^j$, $0<j<q_n$,
is a coordinate permutation whose fixed-point set is closed with empty
interior.  The complement of their finite union is the dense full-period
locus.  Starting with $\sigma_0$, whose period is one, the coordinate
cycle defining $\sigma_{n+1}$ ensures that its period is $k_n$.

After choosing $E_n$, choose $k_n$ divisible by $q_n$ so large that
\begin{equation}\label{villadsen-small-evaluations}
 k_n\geq n+2,\qquad \frac{l_n}{k_n}<2^{-n-2}.
\end{equation}
Set $X_{n+1}=X_n^{k_n}$ and $r_{n+1}=r_n(k_n+l_n)$.
List each evaluation orbit in the order
$t_j=\sigma_n^{-j}(t_0)$, $0\leq j<q_n$, and define
\begin{equation}\label{villadsen-equivariant-map}
 \phi_n(a)(x)=\operatorname{diag}\bigl(
 a(x_0),\ldots,a(x_{k_n-1}),\ (a(t))_{t\in E_n}\bigr).
\end{equation}
Thus every coordinate projection occurs once, followed by the point
evaluations.  Define
\[
 (\sigma_{n+1}(x))_j=\sigma_n(x_{j+1}),
 \qquad j\in\mathbb Z/k_n\mathbb Z.
\]
Let $P_{n+1}$ be the permutation matrix shifting the $k_n$ coordinate
block indices by $-1$ and, separately, shifting the indices in each
evaluation orbit by $-1$.  With each block of size $r_n$, put
\[
 u_{n+1}=\operatorname{diag}(u_n,\ldots,u_n)P_{n+1},
 \qquad q_{n+1}=k_n.
\]
The two factors commute.  Since $q_n$ divides $k_n$, both
$\sigma_{n+1}^{k_n}=\mathrm{id}$ and $u_{n+1}^{k_n}=1$.
Conjugation by $P_{n+1}$ sends block $j+1$ to block $j$.
For a coordinate block this gives
$u_na(\sigma_n^{-1}x_j)u_n^*$; for an evaluation block it gives
$u_na(t_{j+1})u_n^*=\alpha_n(a)(t_j)$.  Hence
\begin{equation}\label{villadsen-equivariance}
 \alpha_{n+1}\phi_n=\phi_n\alpha_n.
\end{equation}
The injective unital maps $\phi_n$ therefore define a Villadsen algebra
$A=\varinjlim(A_n,\phi_n)$ and an automorphism $\alpha$ of $A$.

\emph{Simplicity and the coefficient radius of comparison.}
For each fixed earlier stage, every composite coordinate projection
from $X_n$ to that stage is surjective and $1$-Lipschitz.  Its image of
$E_n$ is consequently a $1/(n+1)$-net.  The evaluation sets thus
satisfy the density condition for simplicity of the Villadsen
limit: every nonzero finite-stage positive element eventually acquires
a nonzero constant diagonal summand and becomes full.  In particular,
$A$ is simple.  As a unital AH algebra it has tracial states, which are
faithful by simplicity; thus $A$ is stably finite.  Moreover, $A$ has
locally finite nuclear dimension.

The cube is a solid seed space, so
\cite[Theorem~3.4]{ELN-Vill} gives
\begin{equation}\label{villadsen-positive-radius}
 \operatorname{rc}(A)
 =\frac12\lim_{n\to\infty}\frac{\dim X_n}{r_n}
 =\prod_{n=0}^{\infty}\frac{k_n}{k_n+l_n}
 \in(0,1).
\end{equation}
Positivity follows from \eqref{villadsen-small-evaluations}, and the
product is less than one because every $l_n$ is positive.  Since a
simple stably finite $\mathcal Z$-stable algebra has strict comparison,
$A$ is not $\mathcal Z$-stable.

\emph{Central projection towers.}
Consider the units of all diagonal blocks at stage $n+1$, including
the coordinate and point-evaluation blocks.  The automorphism
permutes them in cycles of lengths $k_n$ and $q_n$, respectively,
and both lengths tend to infinity.  Every block-unit projection
commutes with $\phi_n(A_n)$.

Fix $N\geq2$.  Every integer at least $N(N-1)$ is of the form
$aN+b(N+1)$ with $a,b\geq0$.  For all sufficiently large $n$,
partition each cycle into segments of lengths $N$ and $N+1$, in
the direction of the action.  Sum the initial projections of the
segments according to their length.  The two resulting bases have
shapes $\{0,\ldots,N-1\}$ and $\{0,\ldots,N\}$; all translated
levels are orthogonal and their sum is exactly $1$.
With $N$ fixed and $n\to\infty$, these bases commute with every
fixed earlier stage, so they define projection towers in
$A_\omega\cap A'$ with the same two shapes and exact coverage.
For any finite $K\subset\mathbb Z$ and $\eta>0$, choosing $N$
large enough makes both interval shapes $(K,\eta)$-invariant.
This proves strong URP in Definition~\ref{strong-urp-definition}.

Since $A$ is simple, the action is minimal.  Strong URP gives WRP
by Lemma~\ref{strong-urp-cuntz-controlled-wrp}, while
\eqref{villadsen-positive-radius} and
Propositions~\ref{orbit-finite-radius-zero} and
\ref{cuntz-controlled-wrp-orbit-zero-cos} give $(1,1)$-COS.
Finally, the locally finite nuclear dimension of $A$ and the
strong URP case of Corollary~\ref{wrp-finite-radius-Z} show that
$B=A\rtimes_\alpha\mathbb Z$ is simple and $\mathcal Z$-stable.
\end{proof}

\begin{rem}[Hirshberg's central projection towers]
\label{villadsen-strong-urp}
The underlying $\mathbb Z$-action in Hirshberg's construction
\cite[proof of Claim~(A) in Theorem~1]{Hirshberg-exotic-circle}
has full cyclic towers given by the central projections of the
finite-stage fibers, with lengths $2^n$.  Partitioning those cycles
into segments of the same two fixed lengths gives strong URP
by the same argument.  His coefficient algebra also has finite
radius of comparison: its homogeneous stages $C_n$ have base
dimension $2s(n)$ and matrix size $r(n)$, with $s(n)/r(n)\leq1$,
so \cite[Theorem~5.1]{Toms-Comp-DS} gives
\[
 \operatorname{rc}(C_n)\leq\frac{2s(n)-1}{2r(n)}\leq1.
\]
All quotients of the stages are stably finite, and the simple limit
is stably finite.  The inductive-limit inequality
\cite[Propositions~3.2.3 and 3.2.4(iii)]{RC-5-person}
therefore gives $\operatorname{rc}(C)\leq1$.
Thus Proposition~\ref{cuntz-controlled-wrp-orbit-zero-cos}
gives COS from finite coefficient radius of comparison, as in
Example~\ref{villadsen-wrp-cos-example}.
\end{rem}
\chapter{Examples}\label{section-examples}

This chapter gives examples of \Cstar{}-algebras with weak tracial approximation
and \Cstar{}-dynamical systems with the weak Rokhlin property. The main results are as follows.

First, any simple unital stably finite classified \Cstar{}-algebra can be weakly tracially approximated by (non-unital) interval algebras: 
\begin{overviewthm}[Theorem~\ref{classifiable-interval-approximation}]
Let $A$ be a unital, separable, simple, infinite-dimensional, nuclear,
stably finite, $\mathcal Z$-stable \Cstar{}-algebra satisfying the UCT.
Let $\mathcal I_{\mathrm h}$ and $\mathcal I_{\mathrm o}$ be the classes
of finite direct sums
\[
 \bigoplus_{j=1}^{r}\mathrm M_{n_j}(\mathrm C_0([0,1)))\quad\hbox{and}\quad
 \bigoplus_{j=1}^{r}\mathrm M_{n_j}(\mathrm C_0((0,1))),
 \qquad r,n_1,\ldots,n_r\in\mathbb N,
\]
respectively.  Then
\[
 A\in\mathrm{TA}_1(\mathcal I_{\mathrm h})
       \cap\mathrm{TA}_1(\mathcal I_{\mathrm o}).
\]
In particular, $A$ can be weakly tracially approximated by either class.
\end{overviewthm}

Then we have examples of \Cstar{}-dynamical systems which have the weak Rokhlin property. If the \Cstar{}-algebra is a bundle over a topological space $X$ and the action preserve $\mathrm{C}(X)$ and has the ordinary URP, the the action on the \Cstar{}-algebra also has the URP. 
\begin{overviewthm}[Theorem~\ref{central-topological-urp}]
Let $A$ be a unital separable \Cstar{}-algebra with
$\mathrm T(A)\ne\emptyset$, and let $\alpha$ be an action of a
countable discrete amenable group $\Gamma$ on $A$.
Suppose that $A$ contains an $\alpha$-invariant unital subalgebra
$\mathrm C(X)\subseteq Z(A)$, where $X$ is compact Hausdorff.
If the induced action on $X$ has topological URP, then
$(A,\Gamma)$ has URP as in Definition~\ref{URP-defn-C-dyn}
and hence WRP.

In particular, this applies to a unital separable \Cstar{}-algebra that is
a continuous field over $X$ and has a nonempty tracial state space, provided that
the action preserves its canonical central copy of $\mathrm C(X)$ and the
induced action on $X$ has topological URP.
\end{overviewthm}

If the action on a \Cstar{}-algebra $A$ is \Wstar{}-free and has only finitely many extreme invariant tracial state, using Rokhlin property in \Wstar{}-dynamical systems, we have that  the action has the WRP.
\begin{overviewthm}[Theorem~\ref{unique-invariant-trace-wrp}]
Let $A$ be a unital separable \Cstar{}-algebra with
$\mathrm T(A)\ne\varnothing$, and let a countable discrete amenable
group $\Gamma$ act on $A$.  If the action has a unique invariant
tracial state and is \Wstar{}-centrally free, then $(A,\Gamma)$ has WRP.
\end{overviewthm}

If the coefficient \Cstar{}-algebra $A$ is nuclear, the example above can be generalized to the case that the invariant states form a Bauer simplex with finite dimensional boundary. Note that in both cases, the \Cstar{}-algebra $A$ is allowed to be a commutative \Cstar{}-algebra.
\begin{overviewthm}[Theorem~\ref{finite-dimensional-trace-wrp}]
Let $A$ be a unital separable nuclear \Cstar{}-algebra with
$\mathrm T(A)\ne\emptyset$, and let
$\alpha\colon\Gamma\to\operatorname{Aut}(A)$ be an action
of a countable discrete amenable group.
Suppose that $\alpha$ is \Wstar{}-centrally free and that
$S=\mathrm T^\alpha(A)$ is a Bauer simplex with
$\dim(\partial_eS)<\infty$.
If $\Gamma$ is infinite, then $\alpha$ has WRP.
The same conclusion holds for finite $\Gamma$ when $A$ has
no finite-dimensional quotients.
\end{overviewthm}

\section{Classifiable \Cstar{}-algebras}\label{interval-models-section}

\begin{thm}\label{classifiable-interval-approximation}
Let $A$ be a unital, separable, simple, infinite-dimensional, nuclear,
stably finite, $\mathcal Z$-stable \Cstar{}-algebra satisfying the UCT.
Let $\mathcal I_{\mathrm h}$ and $\mathcal I_{\mathrm o}$ be the classes
of finite direct sums
\[
 \bigoplus_{j=1}^{r}\mathrm M_{n_j}(\mathrm C_0([0,1)))\quad\hbox{and}\quad
 \bigoplus_{j=1}^{r}\mathrm M_{n_j}(\mathrm C_0((0,1))),
 \qquad r,n_1,\ldots,n_r\in\mathbb N,
\]
respectively.  Then
\[
 A\in\mathrm{TA}_1(\mathcal I_{\mathrm h})
       \cap\mathrm{TA}_1(\mathcal I_{\mathrm o}).
\]
In particular, $A$ can be weakly tracially approximated by either class.
\end{thm}

\begin{proof}
The half-open interval assertion follows from the inductive-limit model of \cite{GLN-TAS-1}.  The open interval assertion
then follows from $\mathcal Z$-absorption.

\smallskip
\noindent\emph{(1) The inductive-limit model and its finite-dimensional quotient.}
By \cite[Theorem~A]{CETWW-dim-n}, $A$ has finite nuclear dimension.
Hence \cite[Theorem~4.10]{EGLN-DR} gives
$\operatorname{gTR}(A\otimes\mathcal Q)\leq1$, where $\mathcal Q$
is the universal UHF algebra.  The range theorem
\cite[Theorem~13.50]{GLN-TAS-1} applies because the Elliott invariant
of $A$ is a simple weakly unperforated invariant in the sense of
\cite[Subsection~13.7]{GLN-TAS-1}.  It provides a unital simple,
$\mathcal Z$-stable model $B$ with the same Elliott invariant as $A$
and with $\operatorname{gTR}(B\otimes\mathcal Q)\leq1$.
Thus \cite[Theorem~29.8]{GLN-TAS-2} gives $A\cong B$.
The model is a unital simple inductive limit
\[
 B=\varinjlim(B_n,\xi_{n,n+1})
\]
with unital injective connecting maps.  We use the building blocks in
\cite[Subsections~13.27--13.29]{GLN-TAS-1} and the connecting maps
constructed in \cite[Subsections~13.43--13.47]{GLN-TAS-1}.
Write
\[
 F_n=P_n\mathrm M_{N_n}(\mathrm C(X_n))P_n\oplus
          \bigoplus_{j=2}^{p_n}\mathrm M_{r_{n,j}},\qquad
 Q_n=\bigoplus_{j=1}^{p_n}\mathrm M_{r_{n,j}},
 \quad r_{n,1}=\operatorname{rank}(P_n).
\]
Here $X_n=[0,1]\vee X'_n$, with $1$ attached to $X'_n$ and with
distinguished point $\theta_n=0$.  The projection $P_n$ is constant on
the attached interval.  Let
\[
 \pi_n:F_n\longrightarrow Q_n,\qquad
 \pi_n(g,a_2,\ldots,a_{p_n})
       =(g(\theta_n),a_2,\ldots,a_{p_n}).
\]
For a finite-dimensional algebra $E_n$ and unital homomorphisms
$\beta_{n,0},\beta_{n,1}:Q_n\to E_n$, the building block is
\[
 B_n=\{(f,g)\in \mathrm C([0,1],E_n)\oplus F_n:
       f(e)=\beta_{n,e}(\pi_n(g)),\ e=0,1\}.
\]
Set $q_n(f,g)=\pi_n(g)$ and $J_n=\ker q_n$.  The endpoint
multiplicities in this construction can be chosen strictly positive,
as allowed in \cite[Subsections~13.12--13.13]{GLN-TAS-1},
so that $\beta_{n,0}$ is injective.
In the notation of \cite{GLN-TAS-1}, $F_n$, $Q_n$, and $q_n$
correspond to $A_{X,n}$, $F_n$, and $\pi^e_{A_n}$, respectively;
$J_n$ is the sum of the two coordinate ideals denoted there
by $I_n$ and $J_n$.
This description includes the homogeneous summand required for general
\Ksym{}-theory, so the blocks $B_n$ need not be one-dimensional NCCW
complexes.

\smallskip
\noindent\emph{(2) Uniform concentration at the point quotient.}
We claim that
\begin{equation}\label{GLN-point-ideal-small}
 \sup_{\tau\in \mathrm T(B)}\|\tau|_{J_n}\|\leq\delta_n,\qquad
 \delta_n=6\cdot2^{-2n}\longrightarrow0,
\end{equation}
where $B_n$ is identified with its image in $B$.  We derive this bound
from the estimates for the connecting maps.

Let $\varphi_{n,n+1}$ be the connecting map before the small
modification that makes the limit simple.  This map induces a homomorphism
$Q_n\to Q_{n+1}$ and hence maps $J_n$ into $J_{n+1}$.
For a homomorphism $\alpha:D\to D'$, write $\alpha^\sharp$ for the
induced map $\operatorname{Aff}(\mathrm T(D))\to
\operatorname{Aff}(\mathrm T(D'))$.
The claim in the proof of \cite[Theorem~13.42]{GLN-TAS-1}
gives a positive unital linear map
$\Gamma_{n+1}:\operatorname{Aff}(\mathrm T(Q_{n+1}))\to
\operatorname{Aff}(\mathrm T(B_{n+1}))$ such that
\[
 \left\|\Gamma_{n+1}q_{n+1}^\sharp
           \varphi_{n,n+1}^\sharp(v)
       -\varphi_{n,n+1}^\sharp(v)\right\|
       \leq4\cdot2^{-2n}\|v\|.
\]
For $x\in(J_n)_+$ with $\|x\|\leq1$, put
$\widehat x(\rho)=\rho(x)$.  Since
$q_{n+1}\varphi_{n,n+1}(x)=0$, the $\Gamma_{n+1}$ term in the
preceding estimate vanishes for $v=\widehat x$.  Thus
\[
\sigma(\varphi_{n,n+1}(x))\leq4\cdot2^{-2n},
       \qquad \sigma\in \mathrm T(B_{n+1}).
\]
The modified map $\xi_{n,n+1}$ satisfies a similar
uniform estimate on the positive unit ball of $J_n$.
Write $\varphi_{n,n+1}=(\varphi_E,\varphi_X)$ and
$\xi_{n,n+1}=(\xi_E,\xi_X)$ for their line and homogeneous
coordinates.  Equation~(e13.77) of \cite{GLN-TAS-1} gives
\[
 \bigl|\rho(\xi_X(x))-\rho(\varphi_X(x))\bigr|
       \leq2\cdot2^{-2n},\qquad \rho\in \mathrm T(F_{n+1}).
\]
Composing such a trace with the homogeneous quotient map gives a
trace on $B_{n+1}$, so $\rho(\xi_X(x))\leq6\cdot2^{-2n}$.
On the middle portion of each line, equation~(e13.82) expresses
$\xi_E(x)$ as a reparametrized evaluation of $\varphi_E(x)$;
its normalized matrix traces are therefore at most $4\cdot2^{-2n}$.
On the two endpoint collars, equations~(e13.86) and~(e13.87) bound
the difference from the corresponding endpoint trace of
$\varphi_E(x)$ by $2\cdot2^{-2n}$.  Thus these traces are at most
$6\cdot2^{-2n}$ as well.
The estimates for the modification hold on all of $B_n$, so the
bounds above are uniform on the positive unit ball of $J_n$.
Since every trace on the subhomogeneous block $B_{n+1}$ is
an integral of normalized traces of its irreducible representations,
we obtain
\[
 \sigma(\xi_{n,n+1}(x))\leq\delta_n,
       \qquad \sigma\in \mathrm T(B_{n+1}),\quad x\in(J_n)_+,\ \|x\|\leq1.
\]
Restricting a limit trace to $B_{n+1}$ and taking the supremum over
such $x$ proves \eqref{GLN-point-ideal-small}.  Consequently, for every
positive contraction $v\in J_n$,
\begin{equation}\label{GLN-point-ideal-rank-small}
 \mathrm d_\tau(v)=\lim_{m\to\infty}\tau(v^{1/m})\leq\delta_n,
       \qquad \tau\in \mathrm T(B),
\end{equation}
because all $v^{1/m}$ belong to $J_n$.  The bound on the norm of
the trace restriction controls all these roots simultaneously.

\smallskip
\noindent\emph{(3) The endpoint collar.}
Fix a finite subset $\mathcal G\subset B_n$ and a tolerance $\eta>0$.
Choose $0<a<1/4$ so small that the line coordinates of every
$x=(f,g)\in\mathcal G$ are within $\eta$ of their endpoint values on
$[0,2a]$ and $[1-2a,1]$, and that the first coordinate of $g$ on the
attached interval $[0,2a]\subset X_n$ is within $\eta$ of
$g(\theta_n)$.
For each $x$, a scalar cutoff gives $x'\in B_n$ satisfying
$\|x-x'\|<\eta$, $q_n(x')=q_n(x)$, and
\[
 f'(t)=\beta_{n,0}(q_n(x))\quad(0\leq t\leq a),\qquad
 f'(t)=\beta_{n,1}(q_n(x))\quad(1-a\leq t\leq1),
\]
while the first point coordinate is the constant matrix
$q_n(x)_1$ on $[0,a]\subset X_n$.
The modified coordinates agree with the originals outside the
larger collars and satisfy the same boundary conditions.

There is an injective homomorphism
\[
 \iota:\mathrm C_0([0,1),Q_n)\longrightarrow B_n
\]
defined as follows, extending every $b\in \mathrm C_0([0,1),Q_n)$ by $b(1)=0$.
Its line coordinate is
\[
 t\longmapsto
 \begin{cases}
   \beta_{n,0}(b(t/a)),&0\leq t\leq a,\\
   0,&a\leq t\leq1-a,\\
   \beta_{n,1}(b((1-t)/a)),&1-a\leq t\leq1.
 \end{cases}
\]
Its first point coordinate is $b_1(u/a)$ on the attached interval
$0\leq u\leq a$ and zero elsewhere on $X_n$; its remaining point
coordinates are $b_j(0)$.
These formulas are compatible at the endpoints and at the edges
of the collars, and injectivity of $\beta_{n,0}$ implies injectivity of $\iota$.
Moreover,
\begin{equation}\label{GLN-collar-quotient}
 q_n(\iota(b))=b(0).
\end{equation}

Choose a positive contraction $\chi\in \mathrm C_{\mathrm c}([0,1))$ with $\chi(0)=1$
and put
\[
 C=\iota(\mathrm C_0([0,1),Q_n)),\qquad
 h=\iota(\chi\,1_{Q_n}).
\]
Then $C\in\mathcal I_{\mathrm h}$, $\|h\|=1$, and $h$ is central in
$B_n$: each of its fiber values is a scalar multiple of the fiber
unit.  The constant coordinates on the collars give
\begin{equation}\label{GLN-collar-compression}
 [h,x']=0,\qquad
 hx'h=\iota(\chi^2q_n(x))\in C,
       \qquad x\in\mathcal G.
\end{equation}
Moreover, $q_n(h)=1_{Q_n}$, so that $1-h\in J_n$ and
\[
 \mathrm d_\tau(1-h)\leq\delta_n,\qquad \tau\in \mathrm T(B).
\]
Fix a finite set $\mathcal F=\{b_1,\ldots,b_m\}\subset B$
and $\varepsilon>0$.
Choose $n$ so that $\delta_n<\varepsilon$ and choose
$g_i\in B_n$ with $\|b_i-g_i\|<\varepsilon/4$.
Apply the preceding construction to $\mathcal G=\{g_1,\ldots,g_m\}$,
with tolerance $\eta<\varepsilon/4$, and denote the
resulting elements by $g_i'$.
Then $C$, $h$, and $g_i'$ satisfy
\[
 \|b_i-g_i'\|<\varepsilon/2,\qquad
 [h,g_i']=0,\qquad hg_i'h\in C,\qquad
 \mathrm d_\tau(1-h)<\varepsilon\quad(\tau\in \mathrm T(B)).
\]
These are the four conditions of Definition~\ref{defn_WAS}.
Thus $B\in\mathrm{TA}_0(\mathcal I_{\mathrm h})$.

\smallskip
\noindent\emph{(4) Open intervals.}
Choose a positive contraction
$z\in\mathcal Z$ with spectrum $[0,1]$ whose spectral measure
$\mu$ for the unique trace $\tau_{\mathcal Z}$ satisfies
$\mu(\{0,1\})=0$.
To construct such an element, take a nonscalar selfadjoint element
in the projectionless algebra $\mathcal Z$.  Its spectrum is an
interval.  After normalization, choose
two interior points $c<d$ that are not atoms of its spectral
measure, and apply the piecewise linear function with successive
values $1/2,0,1,1/2$ at $0,c,d,1$.  Its only minimum and maximum
occur at $c$ and $d$, respectively.

For $b\in \mathrm C_0((0,1),Q_n)$, extended by zero at both endpoints,
define
\[
 \widetilde\Theta(b)(s)=
       b\bigl(s+(1-s)z\bigr),\qquad 0\leq s\leq1.
\]
Continuous functional calculus is applied entrywise in the
$Q_n$ coordinate.  This is a homomorphism into
$\mathrm C_0([0,1),Q_n)\otimes\mathcal Z$, since its value at $s=1$ is
zero.  Put
\[
 \Theta=(\iota\otimes\operatorname{id}_{\mathcal Z})
                  \widetilde\Theta.
\]
It is injective, since
\[
 (q_n\otimes\operatorname{id}_{\mathcal Z})(\Theta(b))=b(z)
\]
and $z$ has full spectrum.  Given $\alpha>0$, choose a positive
contraction $\gamma\in \mathrm C_{\mathrm c}((0,1))$ of norm one which is equal to
one on $[\nu,1-\nu]$, with $\nu>0$ so small that
\[
 \mathrm d_{\tau_{\mathcal Z}}(1-\gamma(z))
        =\mu(\{t:\gamma(t)<1\})<\alpha.
\]
Set
\[
 D=\Theta(\mathrm C_0((0,1),Q_n)),\qquad
 k=\Theta(\gamma\,1_{Q_n}).
\]
The same constant coordinates on the collars give
\begin{equation}\label{GLN-open-compression}
 [k,x'\otimes1_{\mathcal Z}]=0,\qquad
 k(x'\otimes1_{\mathcal Z})k
       =\Theta(\gamma^2q_n(x))\in D,\quad x\in\mathcal G,
\end{equation}
and $\|k\|=1$.

Every trace on $B\otimes\mathcal Z$ is
$\tau\otimes\tau_{\mathcal Z}$ for some $\tau\in \mathrm T(B)$.
Its restriction to $B_n\otimes\mathcal Z$ splits into a positive
trace supported on $J_n\otimes\mathcal Z$, of norm at most
$\delta_n$, and a positive trace factoring through
$Q_n\otimes\mathcal Z$.
This decomposition is induced by the central
support of the ideal in the bidual.  On the quotient,
$k$ is $1_{Q_n}\otimes\gamma(z)$.
The dimension function is additive on the two trace components,
and the quotient component has norm at most one.  Hence
\begin{equation}\label{GLN-open-defect}
 \mathrm d_{\tau\otimes\tau_{\mathcal Z}}(1-k)
       \leq\delta_n+
             \mathrm d_{\tau_{\mathcal Z}}(1-\gamma(z))
       <\delta_n+\alpha.
\end{equation}

To verify the approximation property, fix
$\mathcal F=\{b_1,\ldots,b_m\}\subset B$ and $\varepsilon>0$.
Choose $n$ and $g_i\in B_n$ such that
$\delta_n<\varepsilon/2$ and
$\|b_i-g_i\|<\varepsilon/4$.
Apply the collar construction to $\mathcal G=\{g_1,\ldots,g_m\}$
with tolerance $\eta<\varepsilon/4$, obtaining $g_i'$, and choose
$0<\alpha<\varepsilon/2$ in the construction of $D$ and $k$.
Since $B$ is $\mathcal Z$-stable,
\cite[Theorem~2.2]{TW-D} gives an isomorphism
$\rho:B\otimes\mathcal Z\to B$ with
\[
 \|\rho(b_i\otimes1_{\mathcal Z})-b_i\|
       <\varepsilon/4\quad(1\leq i\leq m).
\]
Put $D'=\rho(D)$, $h_{\mathrm o}=\rho(k)$, and
$b_i'=\rho(g_i'\otimes1_{\mathcal Z})$.
Then $D'\in\mathcal I_{\mathrm o}$, $\|h_{\mathrm o}\|=1$, and
\[
 \|b_i-b_i'\|
 \leq \|b_i-\rho(b_i\otimes1_{\mathcal Z})\|
       +\|b_i-g_i\|+\|g_i-g_i'\|
 <3\varepsilon/4.
\]
Equation~\eqref{GLN-open-compression} gives
$[h_{\mathrm o},b_i']=0$ and
$h_{\mathrm o}b_i'h_{\mathrm o}\in D'$.
Every pullback trace $\tau\circ\rho$ is a product trace, so
\eqref{GLN-open-defect} gives
$\mathrm d_\tau(1-h_{\mathrm o})<\delta_n+\alpha<\varepsilon$
for every $\tau\in \mathrm T(B)$.
Thus $B\in\mathrm{TA}_0(\mathcal I_{\mathrm o})$.

\smallskip
\noindent\emph{(5) Cuntz comparison of the defect.}
The algebra $A\cong B$ has strict comparison by
\cite[Corollary~4.6]{Ror-Z-stable}.  For a given nonzero positive contraction
$a_0\in A$, simplicity and compactness of $\mathrm T(A)$ give
\[
 m_0=\min_{\tau\in \mathrm T(A)}\tau(a_0)>0.
\]
Given also a finite subset of $A$ and a tolerance $\varepsilon>0$,
apply the approximation property for either interval class with
tolerance $\min\{\varepsilon,m_0/2\}$.
Write $h$ for the resulting positive contraction, using
$h_{\mathrm o}$ in the open-interval case.
Then $\mathrm d_\tau(1-h)<m_0\leq \mathrm d_\tau(a_0)$ for all $\tau\in \mathrm T(A)$,
so strict comparison gives $1-h\precsim a_0$.
The other three defining conditions still hold with tolerance
$\varepsilon$.
Thus $A$ belongs to both $\mathrm{TA}_1$ classes.
\end{proof}

\section{\Cstar{}-bundles}\label{central-topological-towers-section}

\begin{thm}\label{central-topological-urp}
Let $A$ be a unital separable \Cstar{}-algebra with
$\mathrm T(A)\ne\emptyset$, and let $\alpha$ be an action of a
countable discrete amenable group $\Gamma$ on $A$.
Suppose that $A$ contains an $\alpha$-invariant unital subalgebra
$\mathrm C(X)\subseteq Z(A)$, where $X$ is compact Hausdorff.
If the induced action on $X$ has topological URP, then
$(A,\Gamma)$ has URP as in Definition~\ref{URP-defn-C-dyn}
and hence WRP.

In particular, this applies to a unital separable \Cstar{}-algebra that is
a continuous field over $X$ and has a nonempty tracial state space, provided that
the action preserves its canonical central copy of $\mathrm C(X)$ and the
induced action on $X$ has topological URP.
\end{thm}

\begin{proof}
We identify $\mathrm C(X)$ with the given central subalgebra of $A$.
Write the action on $X$ on the left, so that
$\alpha_g(f)(x)=f(g^{-1}x)$ and
$\operatorname{supp}(\alpha_g(f))=g\operatorname{supp}(f)$.
Tower shapes follow the corresponding left-multiplication convention
of Chapter~\ref{section-rokhlin}.

Fix a finite set $K\subseteq\Gamma$, a finite set
$\mathcal F\subseteq A$, and a tolerance $\varepsilon>0$.
It suffices to consider $0<\varepsilon<1$; enlarge $K$ to
contain the identity if necessary.
Topological URP gives nonempty open bases $U_1,\ldots,U_S$ and
$(K,\varepsilon/2)$-invariant shapes $F_1,\ldots,F_S$ such that
the sets $gU_s$, $g\in F_s$, are pairwise disjoint and
\[
 \mu\left(X\setminus
       \bigcup_{s=1}^S\bigcup_{g\in F_s}gU_s\right)
       <\varepsilon/2,\qquad
 \mu\in\mathcal M_1(X,\Gamma).
\]
Empty bases can be discarded.

We next choose compact subsets $B_s\subset U_s$ whose translated
levels have complement of measure less than $\varepsilon$ for every
$\mu$.  For each invariant probability measure $\mu$, regularity
and normality give open sets $V_s^\mu$ with
$\overline{V_s^\mu}\subset U_s$ and
\[
 \mu\left(\bigcup_{s=1}^S\bigcup_{g\in F_s}gV_s^\mu\right)
       >1-\varepsilon.
\]
The measure of this open union is lower semicontinuous in the
weak* topology.  Thus this inequality holds on an open
neighborhood of $\mu$ in $\mathcal M_1(X,\Gamma)$.
Compactness of that space gives finitely many such families.
Take $B_s$ to be the union of their closures, adding a point of
$U_s$ if needed to make $B_s$ nonempty.  Then
\[
 \mu\left(X\setminus
       \bigcup_{s=1}^S\bigcup_{g\in F_s}gB_s\right)
       <\varepsilon,\qquad
 \mu\in\mathcal M_1(X,\Gamma).
\]

Choose $f_s\in \mathrm C(X)_+$ with $0\leq f_s\leq1$, equal to one
on $B_s$, and supported in $U_s$.
Viewed in $A$, these are norm-one positive contractions.
Their translated levels are mutually orthogonal, and
\[
 q=\sum_{s=1}^S\sum_{g\in F_s}\alpha_g(f_s)
\]
is a positive contraction in $\mathrm C(X)$, equal to one on
$\bigcup_{s=1}^S\bigcup_{g\in F_s}gB_s$.
Every level is central in $A$, so the sum of commutator norms
in Definition~\ref{URP-defn-C-dyn} is zero
for every element of $\mathcal F$.

For $\tau\in\mathrm T(A\rtimes\Gamma)$, the restriction
$\tau|_A$ is $\alpha$-invariant.  Its restriction to $\mathrm C(X)$
is therefore an invariant probability measure $\mu_\tau$.
Functional calculus inside $\mathrm C(X)$ and monotone convergence give
\begin{equation}\label{central-topological-defect}
 \mathrm d_\tau(1-q)
   =\lim_{m\to\infty}\int_X(1-q)^{1/m}\,d\mu_\tau
   =\mu_\tau(\{x:q(x)<1\})
   <\varepsilon.
\end{equation}
Set $c_s=f_s$ and $\Gamma_s=F_s$.
The shapes are $(K,\varepsilon)$-invariant, the $c_s$ have norm
one, their translated levels are pairwise orthogonal, and the
sum of commutator norms is zero.
Equation~\eqref{central-topological-defect} gives
the uniform dimension-function bound on the complement.
Thus Definition~\ref{URP-defn-C-dyn} holds, and
Proposition~\ref{urp-implies-wrp} gives WRP.

A unital continuous-field algebra has its canonical central
$\mathrm C(X)$-action, so the stated special case follows.
\end{proof}

\section{Central sequences and freeness}\label{central-sequences-freeness}

We introduce the von Neumann algebra conventions for the Rokhlin
theorem in the next section.  Throughout this section, $M$ has separable
predual, and $\omega$ is the fixed free ultrafilter on $\mathbb N$.
Write $Z(M)$ for the center and $M_*$ for the predual.
All sequences in $M$ are assumed to be bounded in operator norm.
For $\psi\in M_*$ and $x,y\in M$, our module conventions are
\[
 (x\psi)(y)=\psi(yx),\qquad(\psi x)(y)=\psi(xy),
 \qquad [x,\psi]=x\psi-\psi x.
\]
For a normal positive functional $\phi$, put
\[
 \|x\|_\phi=\phi(x^*x)^{1/2},\qquad
 \|x\|_\phi^\sharp=
       \left(\frac{\phi(x^*x)+\phi(xx^*)}{2}\right)^{1/2}.
\]
If $\phi$ is faithful, $\|\cdot\|_\phi^\sharp$ induces the
$\sigma$-strong* topology on bounded subsets of $M$.

\subsection{Centralizing sequences and ultrapowers}

\begin{dfn}\label{centralizing-sequences}
A bounded sequence $(x_n)$ in $M$ is \emph{$\omega$-central} if
\[
 [x_n,a]\longrightarrow0\quad\hbox{$\sigma$-strong* along $\omega$}
       \qquad(a\in M),
\]
and is \emph{$\omega$-centralizing} if
\[
 \lim_{n\to\omega}\|[x_n,\psi]\|=0\qquad(\psi\in M_*).
\]
Replacing the ultralimit by the ordinary limit gives the usual
notions of a central sequence and a centralizing sequence.
Centralizing sequences are also called strongly central sequences.
\end{dfn}

\begin{lem}\label{central-versus-centralizing}
Every $\omega$-centralizing sequence is $\omega$-central.
If $\phi$ is a faithful normal state, an $\omega$-central sequence
is $\omega$-centralizing precisely when
$\lim_{n\to\omega}\|[x_n,\phi]\|=0$.
In particular, the two notions agree when $M$ is finite.
The same assertions hold for ordinary limits.
\end{lem}

\begin{proof}
Fix a faithful normal state $\phi$.
The linear space $M\phi=\{a\phi:a\in M\}$ is norm-dense in $M_*$
by faithfulness of $\phi$.  The identity
\[
 [x_n,a]\phi=[x_n,a\phi]-a[x_n,\phi]
\]
shows that its left side tends to zero in predual norm when
$(x_n)$ is centralizing.  For $b\in M$,
$\|b\|_\phi^2\leq\|b\|\,\|b\phi\|$; applying this to
$b=[x_n,a]$ proves strong convergence of the commutators to zero.
Since $[x_n^*,\psi]=-[x_n,\psi^*]^*$, the adjoint sequence is also
centralizing.  Applying the same argument to $x_n^*$ proves
strong* centrality.
Conversely,
\[
 \|[x_n,a\phi]\|
 \leq\|[x_n,a]\|_\phi+\|a\|\,\|[x_n,\phi]\|
\]
proves centralization on the dense space $M\phi$, and boundedness
extends it to all of $M_*$.  In the finite case choose a faithful
normal tracial state, for which $[x_n,\phi]=0$ identically.
\end{proof}

\begin{dfn}\label{ocneanu-ultrapowers}
Inside $\ell^\infty(\mathbb N,M)$, define
\begin{align*}
 I_\omega(M)&=\{(x_n):x_n\longrightarrow0
                          \ \sigma\hbox{-strong* along }\omega\},\\
 \mathcal C_\omega(M)&=\{(x_n):(x_n)\hbox{ is $\omega$-centralizing}\},\\
 \mathcal M^\omega(M)&=\{x:xI_\omega(M)\subseteq I_\omega(M),\
                          I_\omega(M)x\subseteq I_\omega(M)\}.
\end{align*}
The last algebra is the two-sided normalizer of $I_\omega(M)$.
One has
$I_\omega(M)\subseteq\mathcal C_\omega(M)\subseteq\mathcal M^\omega(M)$,
and $I_\omega(M)$ is a closed ideal of $\mathcal M^\omega(M)$.
The \emph{Ocneanu ultrapower} and the \emph{asymptotic centralizer}
are, respectively,
\[
 M^\omega=\mathcal M^\omega(M)/I_\omega(M),\qquad
 M_\omega=\mathcal C_\omega(M)/I_\omega(M).
\]
The algebra $M_\omega$ is also called the central sequence algebra
in Ocneanu's sense.  Constant sequences identify $M$ with a subalgebra of
$M^\omega$.
For $\psi\in M_*$, write
\[
 \psi^\omega([(x_n)])=\lim_{n\to\omega}\psi(x_n),
 \qquad \psi_\omega=\psi^\omega|_{M_\omega}.
\]
\end{dfn}

\begin{prp}[Ocneanu]\label{ocneanu-ultrapower-facts}
The algebra $M^\omega$ is a von Neumann algebra, and $M$ and
$M_\omega$ are von Neumann subalgebras.  Moreover,
\[
 M_\omega\subseteq M'\cap M^\omega,\qquad M\cap M_\omega=Z(M).
\]
For every faithful normal state $\phi$ on $M$, $\phi^\omega$ is
a faithful normal state on $M^\omega$, and $\phi_\omega$ is a
faithful normal tracial state on $M_\omega$.

The canonical normal conditional expectation is
\[
 \Psi\colon M^\omega\longrightarrow M,\qquad
 \Psi([(x_n)])=\mathop{\mathrm{weak}^*\!-\!\lim}_{n\to\omega}x_n.
\]
Its restriction to $M_\omega$ takes values in $Z(M)$, and
$\psi_\omega=\psi\circ\Psi|_{M_\omega}$.  Thus $\psi_\omega$
depends only on $\psi|_{Z(M)}$.

If $M$ is finite and $\tau$ is a faithful normal tracial state,
write $\|x\|_{2,\tau}=\tau(x^*x)^{1/2}$.  Then
\[
 M^\omega=
 \ell^\infty(\mathbb N,M)/
       \{(x_n):\lim_{n\to\omega}\|x_n\|_{2,\tau}=0\},
 \qquad M_\omega=M'\cap M^\omega.
\]
\end{prp}

\begin{proof}
The von Neumann algebra, normality, and faithfulness assertions
are the ultrapower facts of \cite[Section~5.1]{Ocneanu-thesis}.
For centralizing representatives $(x_n)$ and $(y_n)$,
\[
 |\psi(x_ny_n-y_nx_n)|
 \leq\|[x_n,\psi]\|\,\|y_n\|\longrightarrow_\omega0.
\]
Hence $\Psi(xy)=\Psi(yx)$ on $M_\omega$, proving the trace
assertion.  Since $\Psi$ is $M$-bimodular and
$M_\omega\subseteq M'\cap M^\omega$, its image there lies in
$Z(M)$.  The formulas for $\psi_\omega$ follow from the definition.
In the finite case $I_\omega(M)$ is the tracial $2$-norm ideal
of all of $\ell^\infty(\mathbb N,M)$, and
Lemma~\ref{central-versus-centralizing} identifies $M_\omega$ with
$M'\cap M^\omega$.
\end{proof}

\subsection{Lifted actions and freeness}

\begin{dfn}\label{liftable-vn-actions}
An automorphism $\theta$ of $M$ induces automorphisms
\[
 \theta^\omega([(x_n)])=[(\theta(x_n))],\qquad
 \theta_\omega=\theta^\omega|_{M_\omega}.
\]
These automorphisms are called \emph{liftable}.
More generally, suppose that $\theta_n\in\operatorname{Aut}(M)$
converges to $\theta\in\operatorname{Aut}(M)$ along $\omega$ in
the $u$-topology, that is,
\[
 \lim_{n\to\omega}\|\psi\circ\theta_n-\psi\circ\theta\|=0
       \qquad(\psi\in M_*).
\]
Then $[(x_n)]\mapsto[(\theta_n(x_n))]$ defines an automorphism
of $M^\omega$ leaving $M_\omega$ invariant
\cite[Section~5.2]{Ocneanu-thesis}.  Such automorphisms, and their
restrictions to $M_\omega$, are called \emph{semiliftable}.
\end{dfn}

For an action $\alpha\colon\Gamma\to\operatorname{Aut}(M)$,
write $\alpha^\omega$ and $\alpha_\omega$ for the induced actions.
The expectation is equivariant:
$\Psi\circ\alpha_g^\omega=\alpha_g\circ\Psi$.
In particular, if $\phi|_{Z(M)}$ is $\alpha$-invariant, then
$\phi_\omega$ is an $\alpha_\omega$-invariant trace.

\begin{dfn}\label{central-freeness-vn}
An automorphism of a von Neumann algebra is \emph{properly outer}
if its restriction to every nonzero invariant central corner
is outer.
An automorphism $\theta\in\operatorname{Aut}(M)$ is
\emph{centrally trivial} if $\theta_\omega=\operatorname{id}$.
For separable predual this is equivalent to
$\theta(x_n)-x_n\to0$ $\sigma$-strong* for every ordinary
centralizing sequence $(x_n)$.
It is \emph{properly centrally nontrivial} if
$\theta|_{pM}$ is not centrally trivial for every nonzero
$\theta$-invariant projection $p\in Z(M)$.

An action $\alpha$ on $M$ is \emph{free} if $\alpha_g$ is
properly outer for every $g\ne e$, and is \emph{centrally free}
if $\alpha_g$ is properly centrally nontrivial for every $g\ne e$.
\end{dfn}

\begin{dfn}[Ocneanu's relative outerness]\label{relative-strong-freeness}
An automorphism $\vartheta\in\operatorname{Aut}(M_\omega)$ is
\emph{strongly outer} if its restriction to
$\mathcal N'\cap M_\omega$ is properly outer for every countable
$\vartheta$-invariant self-adjoint subset
$\mathcal N\subseteq M_\omega$.
Equivalently, one may use every invariant von Neumann subalgebra
generated by a countable set.
Self-adjointness ensures that the relative commutant
is a von Neumann algebra.

An action $\beta$ on $M_\omega$ is \emph{strongly free} if each
$\beta_g$, $g\ne e$, is strongly outer in this sense.
An action $\alpha$ on $M$ is strongly free if $\alpha_\omega$
is strongly free.
\end{dfn}

\begin{lem}[Ocneanu]\label{cf-2-sf}
Let $\vartheta\in\operatorname{Aut}(M_\omega)$ be semiliftable,
represented by automorphisms $\theta_n$ converging to $\theta$
along $\omega$ in the $u$-topology.
If $\theta$ is properly centrally nontrivial, then $\vartheta$
is strongly outer.
In particular, a centrally free action $\alpha$ on $M$ induces
a strongly free action $\alpha_\omega$ on $M_\omega$.
\end{lem}

\begin{proof}
The first assertion is the lemma in
\cite[Section~5.6]{Ocneanu-thesis}.
Apply it to the constant sequence $\theta_n=\alpha_g$ for
each $g\ne e$ to obtain the second assertion.
\end{proof}

\subsection{Invariant traces and GNS actions}

Let $A$ be a unital separable \Cstar{}-algebra, and let a countable
discrete group $\Gamma$ act on $A$.  Write
\[
 \mathrm T^\alpha(A)=
       \{\tau\in\mathrm T(A):\tau\circ\alpha_g=\tau
                                 \text{ for all }g\in\Gamma\}.
\]
For $\tau\in\mathrm T^\alpha(A)$, put $M_\tau=\pi_\tau(A)''$.
Trace invariance gives a normal action
$\alpha_\tau\colon\Gamma\to\operatorname{Aut}(M_\tau)$ satisfying
\[
 (\alpha_\tau)_g(\pi_\tau(a))=\pi_\tau(\alpha_g(a)).
\]
It is implemented on the GNS space by
$[a]\mapsto[\alpha_g(a)]$ and preserves the faithful normal
tracial state $\bar\tau$ on $M_\tau$.
The GNS homomorphism $\pi_\tau$ need not be faithful on $A$.

\begin{dfn}\label{wstar-freeness-cstar}
The action $\alpha$ on $A$ is \emph{\Wstar{}-centrally free} if
$\alpha_\tau$ is centrally free for every
$\tau\in\mathrm T^\alpha(A)$.
It is \emph{\Wstar{}-strongly free} if $\alpha_\tau$ is strongly free,
in the sense of Definition~\ref{relative-strong-freeness},
for every $\tau\in\mathrm T^\alpha(A)$.
\end{dfn}

Lemma~\ref{cf-2-sf} shows that \Wstar{}-central freeness implies
\Wstar{}-strong freeness.
Consequently, the induced action on $(M_\tau)_\omega$ is properly outer
on the relative commutant of every countable invariant self-adjoint
set, as required by the small-tower lemma in the next section.
When $\Gamma$ is amenable and $\mathrm T(A)\ne\varnothing$,
the invariant trace space is nonempty.

\section{Exact Rokhlin towers in von Neumann algebras}

We prove an exact tower theorem for prescribed tiling shapes and then
lift these towers to \Cstar{}-algebras with finitely many extreme invariant
traces.  The von Neumann algebra argument rests on Ocneanu's small-tower lemma.

\begin{dfn}\label{exact-group-tiling}
An \emph{exact tiling} of a discrete group $\Gamma$ consists of finite
nonempty shapes $D_1,\ldots,D_S\subseteq\Gamma$, with $e\in D_s$,
and center sets $C_s\subseteq\Gamma$, such that
\[
 \Gamma=\bigsqcup_{s=1}^S\bigsqcup_{c\in C_s}D_sc.
\]
Thus the translated shapes $D_sc$ partition $\Gamma$.
For a finite set $P\subseteq\Gamma$ containing $e$, put
\[
 \operatorname{Int}_P F=\{g\in F:Pg\subseteq F\}.
\]
A finite nonempty set $F$ is $(P,\delta)$-invariant if
$|F\setminus\operatorname{Int}_P F|<\delta|F|$.
\end{dfn}

\begin{thm}[Downarowicz--Huczek--Zhang]
Let $\Gamma$ be a countable discrete amenable group.  For every finite
$P\subseteq\Gamma$ containing $e$ and every $\delta>0$, there is an
exact tiling all of whose shapes are $(P,\delta)$-invariant.
\end{thm}

\begin{proof}
For infinite $\Gamma$, apply \cite[Theorem~4.3]{DHZ-tiling} with
$K=P$ and a tolerance smaller than $\min\{1,\delta/|P|\}$.
Then \cite[Lemma~2.6]{DHZ-tiling} converts its symmetric-difference
invariance into the stated interior-boundary estimate.
For finite $\Gamma$, use the single shape $\Gamma$ and center set
$\{e\}$.
\end{proof}

\begin{dfn}\label{exact-wstar-rokhlin}
Let $\alpha\colon\Gamma\to\operatorname{Aut}(M)$ be an action on a
von Neumann algebra with separable predual, and let $\omega$ be a
free ultrafilter on $\mathbb N$.  Write $\beta=\alpha_\omega$ for
the induced action on $M_\omega$.
The action has the \emph{\Wstar{}-Rokhlin property} if, for every exact
tiling with shapes $D_1,\ldots,D_S$, there are projections
$p_1,\ldots,p_S\in M_\omega$ such that the projections
$\beta_k(p_s)$, $k\in D_s$, $1\leq s\leq S$, are pairwise
orthogonal and
\[
 \sum_{s=1}^S\sum_{k\in D_s}\beta_k(p_s)=1.
\]
\end{dfn}

For the next theorem and its proof, write $R=M_\omega$,
$t=\phi_\omega$, and $\beta=\alpha_\omega$.  Proposition~\ref{ocneanu-ultrapower-facts}
and equivariance of $\Psi$ give a faithful normal
$\beta$-invariant tracial state $t$ on $R$.

\begin{thm}\label{C:W*-Rokhlin-property}
Let $\Gamma$ be a countable discrete amenable group, let $M$ have
separable predual, and let $\alpha\colon\Gamma\to\operatorname{Aut}(M)$
be centrally free.  Suppose that $\phi$ is a faithful normal state
on $M$ and $\phi|_{Z(M)}$ is invariant under $\alpha$.
Then $\alpha$ has the \Wstar{}-Rokhlin property.

More precisely, for every exact tiling with shapes $D_1,\ldots,D_S$
and every countable subset $\mathcal N\subseteq M_\omega$, the
projections in Definition~\ref{exact-wstar-rokhlin} can be chosen
in $\mathcal N'\cap M_\omega$.
\end{thm}

The following form of Ocneanu's small-tower lemma applies.  By
Lemma~\ref{cf-2-sf}, the induced action $\beta$ in
Theorem~\ref{C:W*-Rokhlin-property} is strongly free.

\begin{lem}[Ocneanu's small-tower lemma]\label{tiny-tower}
Under the hypotheses of Theorem~\ref{C:W*-Rokhlin-property}, let
$K\subseteq\Gamma\setminus\{e\}$ be finite, let $\eta>0$, and let
$\mathcal N\subseteq R$ be a countable $\beta$-invariant subset.
There is a finite partition of unity by projections
$e_0,e_1,\ldots,e_q\in\mathcal N'\cap R$ such that
\[
 t(e_0)<\eta,\qquad e_j\beta_k(e_j)=0
       \quad(1\leq j\leq q,\ k\in K).
\]
\end{lem}

\begin{proof}
Apply the small-tower lemma in
\cite[Section~6.3]{Ocneanu-thesis} in the relative commutant of the
countably generated invariant von Neumann algebra generated by
$\mathcal N\cup\mathcal N^*$, choosing the remainder $2$-norm
tolerance smaller than $\eta^{1/2}$.
Strong freeness supplies proper outerness on this relative
commutant, and $t$ is the invariant trace used in that lemma.
When $K=\varnothing$, take $e_0=0$ and $e_1=1$.
\end{proof}

\begin{lem}\label{complete-tile-folner}
Fix an exact tiling with shapes $D_1,\ldots,D_S$.
For every finite $P\subseteq\Gamma$ containing $e$ and every
$\delta>0$, there is a finite nonempty $(P,\delta)$-invariant set
$H$ which is a union of complete tiles of this fixed tiling.
\end{lem}

\begin{proof}
Put $L=\bigcup_sD_sD_s^{-1}$.  Choose $0<\eta<1$ with
$\eta/(1-\eta)<\delta$, and use amenability to choose a finite
nonempty $(LP,\eta)$-invariant set $F$.
Let $H$ be the union of all tiles of the fixed tiling contained in
$F$.  If $g\in\operatorname{Int}_{LP}F$ and $p\in P$, then the
unique tile containing $pg$ is contained in $Lpg\subseteq F$.
Thus $Pg\subseteq H$, and
\[
 \operatorname{Int}_{LP}F\subseteq\operatorname{Int}_P H,
 \qquad |H|\geq(1-\eta)|F|>0.
\]
Since $H\subseteq F$, it follows that
\[
 \frac{|H\setminus\operatorname{Int}_P H|}{|H|}
 \leq\frac{|F\setminus\operatorname{Int}_{LP}F|}
              {(1-\eta)|F|}
 <\frac{\eta}{1-\eta}<\delta.
\]
\end{proof}

\begin{lem}[Increasing the tower mass]\label{exact-tower-improvement}
Assume the hypotheses of Theorem~\ref{C:W*-Rokhlin-property},
and fix an exact tiling with shapes $D_1,\ldots,D_S$ and a countable
$\beta$-invariant subset $\mathcal N\subseteq R$.
Suppose that $E_1,\ldots,E_S\in\mathcal N'\cap R$ are projections
whose levels $\beta_k(E_s)$, $k\in D_s$, are pairwise orthogonal.
If
\[
 b:=t\left(\sum_{s=1}^S\sum_{k\in D_s}\beta_k(E_s)\right)<1,
\]
then there is another such family $E_1',\ldots,E_S'$ with strictly
larger total trace.
\end{lem}

\begin{proof}
Write $q=\sum_{s,k\in D_s}\beta_k(E_s)$, and choose
$0<\eps<\min\{1/16,1-b\}$.  Put
\[
 U=\bigcup_sD_s,\quad L=\bigcup_sD_sD_s^{-1},\quad
 P=U^{-1}U,\quad m=\sum_s|D_s|,\quad
 \delta=\frac{\eps}{16m}.
\]
By Lemma~\ref{complete-tile-folner}, choose a
$(PL\cup U,\delta)$-invariant union of complete tiles
\[
 H=\bigsqcup_s\bigsqcup_{c\in L_s}D_sc,\qquad L_s\subseteq C_s.
\]
Retain the centers
\[
 L_s'=\{c\in L_s:D_sc\subseteq\operatorname{Int}_P H\},
 \qquad J=\bigsqcup_sD_sL_s'.
\]
If $g\in\operatorname{Int}_{PL}H$, every point of its tile belongs
to $Lg$ and hence to $\operatorname{Int}_P H$.  Therefore this
tile is retained, and
\begin{equation}\label{retained-complete-tiles}
 |H\setminus J|<\delta|H|.
\end{equation}

Choose $\eta>0$ such that $b<(1-\eps)(1-\eta)$.
Apply Lemma~\ref{tiny-tower} with $K=H^{-1}H\setminus\{e\}$,
enlarging $\mathcal N$ by all translates of the $E_s$.
The resulting partition $e_0,\ldots,e_r$ commutes with every
$\beta_g(E_s)$, and the projections $\beta_h(e_j)$, $h\in H$,
are pairwise orthogonal for $j\geq1$.
Set
\[
 x=\frac1{|H|}\sum_{h\in H}\beta_{h^{-1}}(q).
\]
Then $0\leq x\leq1$, $t(x)=b$, and $x$ commutes with each $e_j$.
There is a nonzero $f=e_j$, $j\geq1$, such that
$t(fx)\leq(1-\eps)t(f)$.  Otherwise,
\[
 b\geq t((1-e_0)x)>(1-\eps)(1-t(e_0))
                  >(1-\eps)(1-\eta),
\]
contrary to the choice of $\eta$.  Define
\[
 f'=\sum_{h\in H}\beta_h(f),\qquad
 \rho=t(f')=|H|t(f)>0.
\]
This is a projection commuting with every translate of every $E_s$,
and invariance of $t$ gives
\begin{equation}\label{packing-small-intersection}
 t(f'q)=|H|t(fx)\leq(1-\eps)\rho.
\end{equation}

For each $s$, define a base projection from the retained
centers:
\[
 S_s=\sum_{c\in L_s'}\beta_c(f).
\]
The levels $\beta_k(S_s)$, $k\in D_s$, are pairwise orthogonal,
since they are distinct levels indexed by the retained tiles in
$H$.  For every $k\in D_s$ and $l\in U$, we also have
\begin{equation}\label{new-levels-domination}
 \beta_k(S_s)\leq\beta_l(f').
\end{equation}
Indeed, $kc\in\operatorname{Int}_P H$ and $l^{-1}\in P$ imply
$l^{-1}kc\in H$, for each retained center $c$.
In particular, $S_s\leq f'$.

Set
\[
 E_s'=(1-f')E_s+S_s.
\]
Both summands are orthogonal projections in $\mathcal N'\cap R$.
Their translated levels are pairwise orthogonal.  The existing
levels remain orthogonal after multiplication by the commuting
projections, the new levels are orthogonal as above, and the
mixed products vanish by
\eqref{new-levels-domination}.

To compare the total traces, equation
\eqref{retained-complete-tiles} gives
\[
 A_0:=\sum_{s,k\in D_s}t(\beta_k(S_s))
       =|J|t(f)>(1-\delta)\rho.
\]
For $k\in U$, invariance of $H$ gives $|kH\triangle H|<2\delta|H|$.
Cancellation of the common summands, followed by the triangle
inequality in $\|\cdot\|_{1,t}$, yields
\[
 \|\beta_k(f')-f'\|_{1,t}
       \leq |kH\triangle H|t(f)<2\delta\rho.
\]
Consequently the removed mass satisfies
\begin{align*}
 B_0&:=\sum_{s,k\in D_s}
             t(\beta_k(f')\beta_k(E_s))\\
    &\leq t(f'q)+2m\delta\rho
     \leq(1-\eps+2m\delta)\rho.
\end{align*}
The mass of the new family is therefore
\[
 b'=b-B_0+A_0
    >b+(\eps-\delta-2m\delta)\rho
    >b+\frac{\eps}{2}\rho>b,
\]
which proves the assertion.
\end{proof}

\begin{proof}[Proof of Theorem~\ref{C:W*-Rokhlin-property}]
Adjoin adjoints to $\mathcal N$ and then take its countable
$\beta$-orbit.
Fix the prescribed tiling, and let $b_*$ be the supremum of the
total trace of orthogonal tower families with these shapes and
bases in $\mathcal N'\cap R$.  The zero family is admissible, and
every admissible family has total trace at most one.

The supremum is attained.  Indeed, choose admissible
families $(E_s^{(r)})_{s=1}^S$ with masses tending to $b_*$,
and regard $z_s=(E_s^{(r)})_r$ as elements of $\ell^\infty(R)$.
Let $C$ be the separable unital $\beta$-invariant \Cstar{}-algebra generated by
the $z_s$, their translates, and the constant sequences from
$\mathcal N$.  The index selection lemma
\cite[Section~5.5]{Ocneanu-thesis} gives an equivariant homomorphism
$\Theta\colon C\to R$ which fixes these constant sequences and
satisfies
\[
 t(\Theta((x_r)_r))=\lim_{r\to\omega}t(x_r)
       \quad((x_r)_r\in C).
\]
The lemma applies because $M_*$ is separable, $\Gamma$ is
countable, and $\beta$ is liftable.
The projection, orthogonality, and relative-commutation relations
hold coordinatewise in $C$ and are preserved by $\Theta$.
Thus the projections $\Theta(z_s)$ form an admissible family of
mass $b_*$.

If $b_*<1$, Lemma~\ref{exact-tower-improvement} increases the mass
of this maximizing family, a contradiction.  Hence $b_*=1$.
Its sum of levels is a projection $q$ with $t(q)=1$.
Faithfulness of $t$ implies $q=1$, proving the theorem.
\end{proof}

\begin{thm}\label{finite-invariant-traces-wrp}
Let $A$ be a unital separable \Cstar{}-algebra with
$\mathrm T(A)\ne\varnothing$, and let a countable discrete amenable
group $\Gamma$ act on $A$.  Suppose that the action is \Wstar{}-centrally
free and that $\mathcal S=\mathrm T^\alpha(A)$ has finitely many extreme
points.  Then $(A,\Gamma)$ has WRP.
\end{thm}

\begin{proof}
Amenability ensures that $\mathcal S$ is nonempty.  Write
$\partial_{\mathrm e}\mathcal S=\{\tau_1,\ldots,\tau_n\}$ and put
$\tau=n^{-1}\sum_i\tau_i$.  Every $\sigma\in\mathcal S$ is a convex
combination of the $\tau_i$, so, for $a\in A$,
\begin{equation}\label{finite-trace-domination}
 \|a\|_{2,\tau}^2
 \leq\sup_{\sigma\in\mathcal S}\|a\|_{2,\sigma}^2
 \leq n\|a\|_{2,\tau}^2.
\end{equation}
In the norm ultrapower
$A_\omega=\ell^\infty(A)/c_\omega(A)$, this proves
$J_\tau=J_{\mathcal S}$.

Let $N=\pi_\tau(A)''$, with its faithful normal tracial state
$\bar\tau$, and let $N^\omega$ be its tracial ultrapower.
The GNS homomorphism induces an equivariant homomorphism
\[
 \Phi\colon A_\omega\longrightarrow N^\omega.
\]
Kaplansky density, applied coordinatewise, and bounded positive
lifting through $A\to\pi_\tau(A)$ show that $\Phi$ is surjective.
Its kernel is $J_\tau$.
This construction applies even when $\tau$ is not faithful on $A$,
since it uses the GNS homomorphism and requires no inclusion $A\subset N$.

The induced action on $N$ is centrally free by hypothesis and
preserves $\bar\tau$.  Given an exact tiling with shapes
$D_1,\ldots,D_S$, Theorem~\ref{C:W*-Rokhlin-property} gives
projections $p_s\in N'\cap N^\omega$ with pairwise orthogonal
translated levels and
\[
 \sum_{s=1}^S\sum_{k\in D_s}\alpha_k(p_s)=1.
\]
Here $N_\omega=N'\cap N^\omega$, since $N$ is finite.

To lift the towers centrally, choose sequences of positive contractions
$b^{(s)}\in\ell^\infty(A)$ mapping to $p_s$, and set
\[
 C=\mathrm C^*(A,\alpha_g(b^{(s)}):g\in\Gamma,\ 1\leq s\leq S),
 \qquad
 I=C\cap\{(x_j):\lim_{j\to\omega}\|x_j\|_{2,\tau}=0\}.
\]
The constants identify $A$ with a subalgebra of $\ell^\infty(A)$;
$C$ is separable, and $I$ is an invariant ideal.
Centrality and orthogonality in $N^\omega$ imply
\[
 [b^{(s)},a]\in I,\qquad
 \alpha_k(b^{(s)})\alpha_l(b^{(r)})\in I
\]
for $a\in A$ and distinct levels $(s,k)\ne(r,l)$.
Choose a quasicentral approximate identity $(e^{(m)})_m$ for $I$
in $C$.  Every translate of this approximate identity is again
quasicentral for $I$ in $C$.  Thus the cutoffs
\[
 (1-e^{(m)})b^{(s)}(1-e^{(m)})
\]
have centrality defects tending to zero in norm, and the products
of distinct translated cutoffs also tend to zero in norm.
The images of the cutoffs in $N^\omega$ remain $p_s$.

Apply Kirchberg's epsilon test \cite[Lemma~3.1]{KR-CenSeq}, as in
the proof of Lemma~\ref{lift-towers-0}, including the coordinate
$2$-norm of the cutoff variable among the conditions.  It gives
a sequence $d=(d_j)$ of positive contractions with
$\lim_{j\to\omega}\|d_j\|_{2,\tau}=0$, such that
\[
 c_s=[(1-d)b^{(s)}(1-d)]\in A'\cap A_\omega
\]
are positive contractions, $\Phi(c_s)=p_s$, and their translated
levels are pairwise orthogonal.  The argument uses the trace-kernel
ideal pulled back to $\ell^\infty(A)$ and therefore applies
even when $\tau$ is nonfaithful on $A$.
Finally,
\[
 \Phi\left(1-\sum_{s=1}^S\sum_{k\in D_s}\alpha_k(c_s)\right)=0.
\]
The complement therefore belongs to $J_\tau=J_{\mathcal S}$,
proving WRP.
\end{proof}

\begin{thm}\label{unique-invariant-trace-wrp}
Let $A$ be a unital separable \Cstar{}-algebra with
$\mathrm T(A)\ne\varnothing$, and let a countable discrete amenable
group $\Gamma$ act on $A$.  If the action has a unique invariant
tracial state and is \Wstar{}-centrally free, then $(A,\Gamma)$ has WRP.
\end{thm}

\begin{proof}
Apply Theorem~\ref{finite-invariant-traces-wrp} to the singleton
invariant trace space.
\end{proof}

\begin{cor}\label{finite-invariant-traces-z}
Let $A$ be a separable simple unital infinite-dimensional
$\mathcal Z$-absorbing \Cstar{}-algebra with
$\mathrm T(A)\ne\varnothing$, and let a countable discrete amenable
group $\Gamma$ act on $A$.  Suppose the action is \Wstar{}-centrally free
and has finitely many extreme invariant traces.
Then $B=A\rtimes_{\alpha,\mathrm r}\Gamma$ is simple and weakly
tracially $\mathcal Z$-absorbing.
If $B$ also has strict comparison, then it is tracially
$\mathcal Z$-absorbing; if, in addition, $A$ is nuclear, then
$B\cong B\otimes\mathcal Z$.
\end{cor}

\begin{proof}
Theorem~\ref{finite-invariant-traces-wrp} gives WRP, and
Corollary~\ref{wrp-simple} gives simplicity of $B$, since $A$ is
simple.  $\mathcal Z$-absorption supplies asymptotically central
order zero maps from matrix algebras whose defects have uniformly
small trace.  Hence $A$ is weakly tracially $\mathcal Z$-absorbing.
Corollary~\ref{WRP-WTZ} gives the conclusion.
If $A$ is nuclear, amenability of $\Gamma$ makes $B$ nuclear,
giving the final assertion.
\end{proof}

\begin{cor}\label{finite-invariant-traces-sr1}
Let $(A,\Gamma)$ be a minimal \Cstar{}-dynamical system, where $A$ is
unital and separable, $\mathrm T(A)\ne\varnothing$, and $\Gamma$
is countably infinite, discrete, and amenable.
If the action is \Wstar{}-centrally free, has finitely many extreme
invariant traces, and has $(\lambda,m)$-COS for some
$\lambda\in(0,1]$ and $m\in\mathbb N$, then
\[
 \operatorname{tsr}(A\rtimes_{\alpha,\mathrm r}\Gamma)=1.
\]
\end{cor}

\begin{proof}
Combine Theorem~\ref{finite-invariant-traces-wrp} with
Theorem~\ref{sr1-wrp-cos-D}.
\end{proof}

\section{Trace-kernel ideals and crossed-product traces}
\label{trace-kernels-crossed-products}

In the section, let us prepare notations and lemmas for the next section. 

\subsection{Trace-kernel ideals and Rokhlin towers}

\begin{dfn}\label{D:trace-kernel}
Let $A$ be a unital \Cstar{}-algebra with $\mathrm T(A)\ne\emptyset$,
and let $\omega$ be a free ultrafilter on $\mathbb N$.  Put
\[
 c_\omega(A)=
 \{(a_n)\in\ell^\infty(\mathbb N,A):
                    \lim_{n\to\omega}\|a_n\|=0\},
 \qquad
 A_\omega=\ell^\infty(\mathbb N,A)/c_\omega(A),
\]
and identify $A$ with the constant sequences in $A_\omega$.
For a nonempty subset $S\subseteq\mathrm T(A)$, define
\[
 \|a\|_{2,S}=\sup_{\tau\in S}\tau(a^*a)^{1/2}
 \quad(a\in A),\qquad
 \|[(a_n)]\|_{2,S}=\lim_{n\to\omega}\|a_n\|_{2,S}.
\]
The latter seminorm is well defined on $A_\omega$.  The
\emph{trace-kernel ideal} associated to $S$ is
\[
 J_S=\{x\in A_\omega:\|x\|_{2,S}=0\}.
\]
Write $q_S\colon A_\omega\to A_\omega/J_S$ for the quotient map.
\end{dfn}

The trace property gives
$\|xy\|_{2,S}\leq\|x\|\|y\|_{2,S}$ and
$\|yx\|_{2,S}\leq\|x\|\|y\|_{2,S}$.
Consequently $J_S$ is a norm-closed two-sided ideal.
The supremum over $S$ is taken before the ultralimit; vanishing
against each fixed trace alone does not, in general, imply membership
in $J_S$.

\begin{dfn}[$\sigma$-ideals]\label{sigma-ideal-definition}
An ideal $I\triangleleft B$ is a \emph{$\sigma$-ideal}
if, for every separable \Cstar{}-subalgebra $C\subseteq B$, there is a
positive contraction $e\in I\cap C'$ such that $eb=b$ for all
$b\in I\cap C$.
\end{dfn}

\begin{prp}[Trace-kernel lifting]\label{trace-kernel-sigma}
The ideal $J_S$ is a $\sigma$-ideal of $A_\omega$.
If $A$ is separable, then $J_S\cap(A_\omega\cap A')$ is a
$\sigma$-ideal of $A_\omega\cap A'$, and
\[
 q_S(A_\omega\cap A')
       =q_S(A)'\cap(A_\omega/J_S).
\]
\end{prp}

\begin{proof}
The first assertion is \cite[Remark~4.7]{KR-CenSeq}, with the
constant family $S_n=S$; the relative assertion follows from
\cite[Proposition~4.5(i)]{KR-CenSeq}.
For the displayed equality, lift $y\in q_S(A)'\cap(A_\omega/J_S)$
to $b\in A_\omega$ and apply the $\sigma$-ideal property to
$C=\mathrm C^*(A,b)$.  Choose $e\in J_S\cap C'$ absorbing $J_S\cap C$.
Since $[b,a]\in J_S\cap C$ for $a\in A$, the element
$(1-e)b(1-e)$ commutes with $A$ and still maps to $y$.
\end{proof}

\begin{lem}[Orthogonal lifting of translated levels]
\label{trace-quotient-tower-lifting}
Suppose that $A$ is separable, a countable discrete group $\Gamma$
acts on $A$, and $S\subseteq\mathrm T(A)$ is nonempty and
$\alpha$-invariant.  Let $p_1,\ldots,p_R$ be positive contractions
in $q_S(A)'\cap(A_\omega/J_S)$, and let
$K_1,\ldots,K_R\subseteq\Gamma$ be finite.
If the translated levels
\[
 \{\alpha_g(p_s):1\leq s\leq R,\ g\in K_s\}
\]
are pairwise orthogonal, then there are positive contractions
$c_s\in A_\omega\cap A'$ with $q_S(c_s)=p_s$ whose corresponding
translated levels are pairwise orthogonal.
\end{lem}

\begin{proof}
Invariance of $S$ makes $J_S$ invariant, so the action descends
to the quotient.  Choose positive contraction lifts
$b_s\in A_\omega$ and put
\[
 C=\mathrm C^*(A,\alpha_g(b_s):g\in\Gamma,\ 1\leq s\leq R),
 \qquad I=J_S\cap C.
\]
The algebra $C$ is separable and invariant.  Choose a positive
contraction $e\in J_S\cap C'$ absorbing $I$, and set
$c_s=(1-e)b_s(1-e)$.  These elements have the prescribed images.
Since $[b_s,a]\in I$ for $a\in A$, they commute with $A$.

Each $\alpha_g(e)$ commutes with $C$ and absorbs $I$.
For two distinct levels, put
$b_{s,g}=\alpha_g(b_s)$ and $f_g=(1-\alpha_g(e))^2$.
Their product in the quotient is zero, so
$b_{s,g}b_{t,h}\in I$.  Therefore
\[
 \alpha_g(c_s)\alpha_h(c_t)
 =f_gb_{s,g}f_hb_{t,h}
 =f_gb_{s,g}b_{t,h}f_h=0.
\]
Thus the translated levels are orthogonal, even when $e$
is not invariant.
\end{proof}

\begin{cor}[WRP in the trace quotient]\label{wrp-trace-quotient}
Let $A$ be unital and separable with $\mathrm T(A)\ne\emptyset$,
let $\Gamma$ be countable, discrete, and amenable, and put
$S=\mathrm T^\alpha(A)=\mathcal M_1(A,\Gamma)$.
The action has WRP in the sense of
Definition~\ref{D:weak-Rokhlin-property} if and only if, for every
exact tiling with shapes $K_1,\ldots,K_R$, there are projections
$p_s\in q_S(A)'\cap(A_\omega/J_S)$ such that the levels
$\alpha_g(p_s)$, $g\in K_s$, are pairwise orthogonal and
\[
 \sum_{s=1}^R\sum_{g\in K_s}\alpha_g(p_s)=1.
\]
\end{cor}

\begin{proof}
The images of WRP contractions are orthogonal positive elements
summing to $1$, hence projections.  Conversely, lift the projections
by Lemma~\ref{trace-quotient-tower-lifting}.  The complement of the
sum of the lifted levels belongs to $J_S$, as in
Definition~\ref{D:weak-Rokhlin-property}.
\end{proof}

\subsection{GNS crossed products and conditional expectations}

For the remainder of this subsection, let $\Gamma$ be countable and
discrete, let $A$ be unital, and assume
$S=\mathrm T^\alpha(A)\ne\emptyset$.
Write
\[
 B=A\rtimes_{\alpha,\mathrm r}\Gamma,\qquad E\colon B\to A
\]
for the reduced crossed product and its canonical expectation.
When $\Gamma$ is amenable, the full and reduced crossed products
coincide.

\begin{lem}\label{L:W*-expectation}
For $\tau\in S$, let
\[
 M_\tau=\pi_\tau(A)'',\qquad
 N_\tau=\pi_{\tau\circ E}(B)''.
\]
Let $\bar\tau$ be the faithful normal tracial state on $M_\tau$,
and let $(\alpha_\tau)_g$ be the induced normal automorphisms.
The following assertions hold.
\begin{enumerate}
\item The map
\[
 [a u_g]\longmapsto[a]\otimes\delta_g
\]
extends to a unitary
$H_{\tau\circ E}\to H_\tau\otimes\ell^2(\Gamma)$.
In this representation,
\[
 \pi_{\tau\circ E}(a)=\pi_\tau(a)\otimes1,\qquad
 \pi_{\tau\circ E}(u_g)=U_{\tau,g}\otimes\lambda_g,
\]
where $U_{\tau,g}[a]=[\alpha_g(a)]$ and $\lambda$ is the left
regular representation.  In particular,
$\pi_{\tau\circ E}(A)''$ is canonically normally isomorphic
to $M_\tau$.

\item Under this identification,
\[
 N_\tau=M_\tau\rtimes_{\alpha_\tau}\Gamma.
\]
There is a faithful normal conditional expectation
$E_\tau\colon N_\tau\to M_\tau$ extending $E$, in the sense that
\[
 E_\tau(\pi_{\tau\circ E}(b))=\pi_\tau(E(b))
 \quad(b\in B).
\]
The trace $\bar\tau\circ E_\tau$ is faithful and normal.
Writing $v_{\tau,g}=\pi_{\tau\circ E}(u_g)$, the Fourier
coefficients
\[
 x_g=E_\tau(xv_{\tau,g}^*)\quad(x\in N_\tau)
\]
uniquely determine $x$.  The finite Fourier sums converge to $x$
in the $2$-norm of $\bar\tau\circ E_\tau$ as the finite subsets
of $\Gamma$ increase to $\Gamma$.

\item Put
\[
 \pi=\bigoplus_{\tau\in S}\pi_\tau,\qquad
 \pi_o=\bigoplus_{\tau\in S}\pi_{\tau\circ E},\qquad
 M=\pi(A)'',\quad M_o=\pi_o(A)'',\quad N=\pi_o(B)''.
\]
There is a canonical normal isomorphism $M\cong M_o$.
After identifying these algebras, the action extends normally
to an action $\beta$ on $M$, and
\[
 N=M\rtimes_\beta\Gamma.
\]
The canonical faithful normal expectation
$E_o\colon N\to M$ extends $E$.
Restriction to each direct summand gives surjective normal
homomorphisms
\[
 \rho_\tau^M\colon M\to M_\tau,\qquad
 \rho_\tau^N\colon N\to N_\tau.
\]
Each family is jointly faithful, and
\[
 E_\tau\circ\rho_\tau^N=\rho_\tau^M\circ E_o.
\]
\end{enumerate}
\end{lem}

\begin{proof}
Trace invariance gives
\[
 \langle[a u_g],[b u_h]\rangle
       =\delta_{g,h}\tau(b^*a).
\]
Thus the map in (1) is an isometry with dense range.
Left multiplication by $a$ and $u_g$ gives the two displayed
representation formulas.  Conjugating further by the unitary
\[
 \xi\otimes\delta_h\longmapsto
             U_{\tau,h}^*\xi\otimes\delta_h
\]
gives the usual regular representation:
the coefficient algebra acts on the $h$th summand by
$\pi_\tau(\alpha_{h^{-1}}(a))$, and the group acts by
$1\otimes\lambda$.  Taking weak closures proves the crossed-product
identification.

Compression to the identity summand gives the normal expectation
$E_\tau$; its faithfulness follows from positivity and the
corresponding diagonal compressions on all summands.
The invariant faithful trace on $M_\tau$ extends to
$\bar\tau\circ E_\tau$.
In its GNS Hilbert space the subspaces
$L^2(M_\tau,\bar\tau)v_{\tau,g}$ are mutually orthogonal and
span a dense subspace.  Orthogonal projection onto these subspaces
gives the stated Fourier coefficients and their $2$-norm
convergence.

For (3), take the direct sum of the unitaries in (1) and use
\[
 \bigoplus_{\tau\in S}
       (H_\tau\otimes\ell^2(\Gamma))
 \cong
 \left(\bigoplus_{\tau\in S}H_\tau\right)\otimes\ell^2(\Gamma).
\]
Here the coefficient algebra is $M\otimes1$, and the group
unitaries are $(\bigoplus_\tau U_{\tau,g})\otimes\lambda_g$.
This proves the normal isomorphism and the regular
crossed-product description, including the expectation.
The projections onto the direct summands commute with the
represented algebras.  Their restrictions therefore define the
normal coordinate homomorphisms, which are onto by taking weak
closures and jointly faithful by the direct sum representation.
The expectation identity holds on crossed-product polynomials
and hence, by normality, on all of $N$.
\end{proof}

\subsection{Proper outerness, trace extensions, and centers}

\begin{dfn}[\Wstar{}-proper outerness]\label{wstar-proper-outerness}
We call $\alpha$ \emph{\Wstar{}-properly outer} if
$(\alpha_\tau)_g$ is properly outer for every $\tau\in S$
and every $g\ne e$, using Definition~\ref{central-freeness-vn}.
\end{dfn}
\Wstar{}-central freeness implies \Wstar{}-proper outerness: an inner
automorphism is centrally trivial, also on any invariant
central corner.

The following intertwiner criterion will be used in the trace and
center calculations.

\begin{lem}\label{proper-outer-intertwiners}
Let $\theta$ be a properly outer automorphism of a von Neumann
algebra $P$.  If $b\in P$ satisfies
\[
 b m=\theta(m)b\quad(m\in P),
\]
then $b=0$.
\end{lem}

\begin{proof}
If $b\ne0$, write $b=v|b|$.
The intertwining relation shows that $b^*b$ and $bb^*$ are
central, so the initial and final projections
$p=v^*v$ and $q=vv^*$ are nonzero central projections.
Equivalent central projections are equal, hence $p=q$.
For $z\in Z(P)$, the intertwining relation gives
$zp=\theta(z)p$.  Taking $z=p$ yields $p\leq\theta(p)$,
whereas $z=\theta^{-1}(p)$ yields $p\leq\theta^{-1}(p)$.
Thus $\theta(p)=p$.
Cancelling $|b|$ on its support gives
$\theta(m)p=vmpv^*$ for $m\in P$.
Thus $\theta|_{pP}$ is inner, contradicting proper outerness.
\end{proof}

\begin{prp}\label{P:crossed_product_trace_properly_outer}
If $\alpha$ is \Wstar{}-properly outer, then
\[
 S\longrightarrow\mathrm T(B),\qquad \tau\longmapsto\tau\circ E
\]
is an affine homeomorphism.  Its inverse is restriction to $A$.
Amenability of $\Gamma$ is not required for this assertion
for the reduced crossed product.
\end{prp}

\begin{proof}
For $\tau\in S$, invariance and the trace identity on
crossed-product polynomials show that $\tau\circ E$ is a trace.
Conversely, let $\sigma\in\mathrm T(B)$ and put
$\tau=\sigma|_A\in S$.
In $Q=\pi_\sigma(B)''$, the induced trace $\bar\sigma$ is faithful
and normal.  The coefficient algebra
$P=\pi_\sigma(A)''$ is canonically normally isomorphic to $M_\tau$:
compression to $\overline{\pi_\sigma(A)\xi_\sigma}$ is its GNS
representation, and faithfulness of $\bar\sigma$ makes this
compression faithful.
There is a $\bar\sigma$-preserving normal conditional expectation
$F\colon Q\to P$, obtained by orthogonal projection in
$L^2(Q,\bar\sigma)$.

Put $V_g=\pi_\sigma(u_g)$.  Bimodularity gives
\[
 F(V_g)m=(\alpha_\tau)_g(m)F(V_g)\quad(m\in P).
\]
For $g\ne e$, proper outerness and
Lemma~\ref{proper-outer-intertwiners} imply $F(V_g)=0$.
Consequently $\sigma(a u_g)=0$ for $a\in A$ and $g\ne e$.
It follows on polynomials, and then by continuity, that
$\sigma=\tau\circ E$.
The displayed map and restriction are affine and weak*
continuous, proving the claim.
\end{proof}

\begin{lem}\label{L:center-inclusion}
Suppose that $\alpha$ is \Wstar{}-properly outer.
In the notation of Lemma~\ref{L:W*-expectation},
\[
 M'\cap N=Z(M),\qquad Z(N)=Z(M)^\beta.
\]
In particular $Z(N)\subseteq Z(M)$, and these conclusions
apply to \Wstar{}-centrally free actions.
\end{lem}

\begin{proof}
Let $x\in M'\cap N$.
For each $\tau\in S$, surjectivity of $\rho_\tau^M$ implies
$y=\rho_\tau^N(x)\in M_\tau'\cap N_\tau$.
Write $y_g=E_\tau(yv_{\tau,g}^*)$.  For $m\in M_\tau$,
bimodularity and covariance give
\[
 m y_g
 =E_\tau(y m v_{\tau,g}^*)
 =y_g(\alpha_\tau)_g(m).
\]
Taking adjoints and applying
Lemma~\ref{proper-outer-intertwiners} shows that $y_g=0$ for
$g\ne e$.
Uniqueness of Fourier coefficients therefore gives $y=E_\tau(y)$.
The expectation identity in Lemma~\ref{L:W*-expectation} and
joint faithfulness of the coordinate maps imply $x=E_o(x)\in M$.
Since $x$ commutes with $M$, it lies in $Z(M)$.
The reverse inclusion in $M'\cap N=Z(M)$ is immediate.
Finally, an element of $Z(M)$ is central in $N$ precisely when
it commutes with every implementing group unitary, equivalently
when it is fixed by $\beta$.
\end{proof}

\section{Weak Rokhlin property from invariant traces}
\label{invariant-traces-wrp-section}

Let us prove Theorem~\ref{finite-dimensional-trace-wrp} in this section. We first  construct a \Wstar{}-bundle over the extreme invariant traces and identify its ultrapower
with the uniform trace quotient. Then, we use invariant slicings to glue
fiber Rokhlin towers  into bundle towers.

\subsection{The invariant-trace bundle}\label{wstar-bundles}

We construct a continuous \Wstar{}-bundle over the extreme boundary
of the invariant trace space.  Since this space need not be a face of
$\mathrm T(A)$, we apply the tracial completion theorem to the
crossed product.  Its conditional expectation then determines
the coefficient bundle.

\subsubsection{Uniform tracial completions and fibers}

\begin{dfn}\label{D:wstar-bundle}
Let $K$ be a compact Hausdorff space.
A \emph{continuous tracial \Wstar{}-bundle} over $K$ consists of a unital
\Cstar{}-algebra $\mathcal M$, a specified unital inclusion
$\mathrm C(K)\subseteq Z(\mathcal M)$, and a faithful tracial conditional
expectation
\[
 E_{\mathcal M}\colon\mathcal M\longrightarrow \mathrm C(K)
\]
such that the closed unit ball $\mathcal M_1$ is complete for
\[
 \|x\|_{2,u}=\|E_{\mathcal M}(x^*x)\|^{1/2}.
\]
Here traciality means $E_{\mathcal M}(xy)=E_{\mathcal M}(yx)$,
and faithfulness means that $E_{\mathcal M}(x^*x)=0$ implies $x=0$.
The bundle is \emph{strictly separable} if $\mathcal M_1$ is
separable for this metric.
\end{dfn}

This is the convention of \cite[Definition~3.1]{BBSTWW}.
For $\lambda\in K$, put
\[
 \tau_\lambda=\operatorname{ev}_\lambda\circ E_{\mathcal M},
 \qquad
 \mathcal M_\lambda=\pi_{\tau_\lambda}(\mathcal M),\qquad
 \|x\|_{2,\lambda}=\tau_\lambda(x^*x)^{1/2}.
\]
The fiber $\mathcal M_\lambda$ is already a finite von Neumann
algebra, with its faithful normal GNS trace
\cite[Theorem~11]{Ozawa-Dixmier}.
The proof of that theorem also applies to nonmetrizable $K$;
see \cite[Section~3.1]{BBSTWW}.
In particular, factoriality is not part of the bundle definition.
On norm-bounded subsets, the strict topology from the
Hilbert $\mathrm C(K)$-module $L^2(\mathcal M,E_{\mathcal M})$ agrees
with the uniform $2$-norm topology.

For a unital \Cstar{}-algebra $D$ and a nonempty closed face
$F\subseteq\mathrm T(D)$, write
\[
 \pi_F=\bigoplus_{\tau\in F}\pi_\tau,\qquad
 N_F=\pi_F(D)''.
\]
Each $\tau\in F$ has a normal extension $\bar\tau$ to $N_F$.
The \emph{uniform tracial completion} $\overline D^{\,u}_F$
is realized inside $N_F$ by adjoining limits of norm-bounded
sequences from $\pi_F(D)$ which are Cauchy for
\[
 \|x\|_{2,F}=\sup_{\tau\in F}\bar\tau(x^*x)^{1/2}.
\]
Thus completion is taken on bounded sets.

\begin{thm}[Ozawa]\label{T:ozawa-bauer-completion}
Let $D$ be unital and separable, and let $F$ be a nonempty
closed face of $\mathrm T(D)$ with compact extreme boundary
$L=\partial_eF$.
There is a canonical unital isomorphism
\[
 \theta_F\colon \mathrm C(L)\longrightarrow
             \overline D^{\,u}_F\cap Z(N_F)
\]
such that
\[
 E_F=\theta_F^{-1}\circ
             \operatorname{ctr}_{N_F}|_{\overline D^{\,u}_F}
\]
makes $\overline D^{\,u}_F$ a continuous \Wstar{}-bundle over $L$.
Its expectation and uniform $2$-norm satisfy
\[
 E_F(\pi_F(d))(\tau)=\tau(d),\qquad
 \|\pi_F(d)\|_{2,u}=\|d\|_{2,F}
 \quad(d\in D,\ \tau\in L).
\]
Moreover, $\pi_F(D_1)$ is uniformly $2$-norm dense in
$(\overline D^{\,u}_F)_1$, and the fiber at $\tau\in L$
is canonically $\pi_\tau(D)''$.
\end{thm}

\begin{proof}
These are the Bauer-face conclusions of
\cite[Theorem~3]{Ozawa-Dixmier}.
The unit-ball density follows from the bounded convex
approximation in its proof.
\end{proof}

\subsubsection{Bundle actions and ultrapowers}

\begin{dfn}\label{D:bundle-action}
For a discrete group $\Gamma$, an action
$\beta\colon\Gamma\to\operatorname{Aut}(\mathcal M)$
on a \Wstar{}-bundle is a \emph{bundle action} if it preserves the
specified copy of $\mathrm C(K)$ and
\[
 E_{\mathcal M}\circ\beta_g
       =(\beta_g|_{\mathrm C(K)})\circ E_{\mathcal M}
 \quad(g\in\Gamma).
\]
It is \emph{fiberwise} if $\beta_g|_{\mathrm C(K)}=\operatorname{id}$
for every $g$.
In that case each $\tau_\lambda$ is invariant, and the
induced normal fiber action is determined by
\[
 \beta_{\lambda,g}(\pi_{\tau_\lambda}(x))
       =\pi_{\tau_\lambda}(\beta_g(x))
 \quad(x\in\mathcal M).
\]
\end{dfn}

\begin{dfn}\label{D:bundle-ultrapower}
Let $\omega$ be a free ultrafilter on $\mathbb N$.
The \Wstar{}-bundle ultrapower is
\[
 \mathcal M^\omega=
       \ell^\infty(\mathbb N,\mathcal M)/c_\omega(\mathcal M),
 \qquad
 c_\omega(\mathcal M)=
 \{(x_n)\in\ell^\infty(\mathbb N,\mathcal M):
                    \lim_{n\to\omega}\|x_n\|_{2,u}=0\}.
\]
Here $c_\omega(\mathcal M)$ is the uniform $2$-norm ideal.
The base $K^\omega$ is the ultracopower specified by
\[
 \mathrm C(K^\omega)=\mathrm C(K)_\omega,
\]
where the right hand side is the norm ultrapower.
The central inclusion of $\mathrm C(K^\omega)$ and the expectation
are induced coordinatewise:
\[
 E_{\mathcal M}^\omega([(x_n)])
                  =[(E_{\mathcal M}(x_n))],\qquad
 \|[(x_n)]\|_{2,u}=\lim_{n\to\omega}\|x_n\|_{2,u}.
\]
\end{dfn}

The resulting triple
$(\mathcal M^\omega,K^\omega,E_{\mathcal M}^\omega)$ is again
a continuous \Wstar{}-bundle \cite[Proposition~3.9]{BBSTWW}.
The base $K^\omega$ need not be metrizable.
The set-theoretic ultrapower $K_\omega$, formed by identifying
$(\lambda_n)$ and $(\mu_n)$ when
$\{n:\lambda_n=\mu_n\}\in\omega$, embeds densely into $K^\omega$
by the characters
\[
 [(f_n)]\longmapsto\lim_{n\to\omega}f_n(\lambda_n).
\]
A bundle action induces the coordinatewise action
$\beta_g^\omega([(x_n)])=[(\beta_g(x_n))]$;
a fiberwise action remains fiberwise.

\subsubsection{The coefficient bundle}

\begin{thm}\label{T:W*-bundle-const}
Let $A$ be unital and separable, and let
$\alpha\colon\Gamma\to\operatorname{Aut}(A)$ be an action of a
countable discrete group with
$S=\mathrm T^\alpha(A)\ne\emptyset$.
Suppose that $\alpha$ is \Wstar{}-properly outer and
$K=\partial_eS$ is compact.
Then there are a strictly separable continuous \Wstar{}-bundle
$(\mathcal M,K,E_{\mathcal M})$, a fiberwise action $\beta$,
and a unital equivariant homomorphism $\iota\colon A\to\mathcal M$
such that
\begin{enumerate}
\item $\iota(A_1)$ is uniformly $2$-norm dense in $\mathcal M_1$.
\item For $a\in A$ and $\tau\in K$,
\[
 E_{\mathcal M}(\iota(a))(\tau)=\tau(a),\qquad
 \|\iota(a)\|_{2,u}=\|a\|_{2,S}.
\]
\item There is a canonical normal isomorphism
\[
 \mathcal M_\tau\cong M_\tau=\pi_\tau(A)''
 \quad(\tau\in K)
\]
intertwining $\beta_\tau$ with $\alpha_\tau$.
\item The kernel of $\iota$ is
\[
 \ker\iota=\{a\in A:\|a\|_{2,S}=0\}.
\]
\end{enumerate}
If $\Gamma$ is amenable, the nonemptiness assumption on $S$
follows from $\mathrm T(A)\ne\emptyset$.
\end{thm}

\begin{proof}
Put $B=A\rtimes_{\alpha,\mathrm r}\Gamma$ and let
$E_\alpha\colon B\to A$ be the canonical expectation.
Proposition~\ref{P:crossed_product_trace_properly_outer} identifies
$F=\mathrm T(B)$ affinely and homeomorphically with $S$ by
$\tau\mapsto\tau\circ E_\alpha$.
Thus $S$ is a Choquet simplex, and this map identifies its
compact extreme boundary $K$ with $\partial_eF$.
In the representation of Lemma~\ref{L:W*-expectation}, write
\[
 \pi_o=\bigoplus_{\tau\in S}\pi_{\tau\circ E_\alpha},
 \qquad M=\pi_o(A)'',\qquad N=\pi_o(B)''.
\]
Let $E_o\colon N\to M$ be the normal extension of $E_\alpha$,
and let $\operatorname{ctr}_N$ be the center-valued trace of $N$.

By Theorem~\ref{T:ozawa-bauer-completion}, the completion
$\mathcal B=\overline B^{\,u}_F$ is a \Wstar{}-bundle over $K$.
We write
\[
 \theta\colon \mathrm C(K)\longrightarrow\mathcal B\cap Z(N),
 \qquad
 E_{\mathcal B}=\theta^{-1}\circ\operatorname{ctr}_N|_{\mathcal B}.
\]
On $N$, put
$\|x\|_{2,F}=\|\operatorname{ctr}_N(x^*x)\|^{1/2}$.

The expectations satisfy
\begin{equation}\label{Eq:two-expectations}
 \operatorname{ctr}_N\circ E_o=\operatorname{ctr}_N,
 \qquad E_o\circ\operatorname{ctr}_N=\operatorname{ctr}_N.
\end{equation}
For a normal state $f$ on $Z(N)$, the trace
$f\circ\operatorname{ctr}_N$ restricts to a trace on $B$.
Proposition~\ref{P:crossed_product_trace_properly_outer} implies
that this restriction equals its composition with $E_\alpha$.
Normality extends the equality to $N$.
Since normal states separate $Z(N)$, this proves the first
identity.
The second follows from
$Z(N)\subseteq Z(M)\subseteq M$, by
Lemma~\ref{L:center-inclusion}.

Kadison's inequality and \eqref{Eq:two-expectations} give
\[
 \operatorname{ctr}_N(E_o(x)^*E_o(x))
       \leq \operatorname{ctr}_N(E_o(x^*x))
       =\operatorname{ctr}_N(x^*x).
\]
Hence $E_o$ is contractive for $\|\cdot\|_{2,F}$.
For $x\in\mathcal B_1$, choose $b_n\in B_1$ such that
$\|\pi_o(b_n)-x\|_{2,F}\to0$.
Then
\[
 \|\pi_o(E_\alpha(b_n))-E_o(x)\|_{2,F}\longrightarrow0,
 \qquad \|E_\alpha(b_n)\|\leq1.
\]
Thus $E_o(\mathcal B)\subseteq\mathcal B\cap M$.
Set
\[
 \mathcal M=\mathcal B\cap M=E_o(\mathcal B),\qquad
 E_{\mathcal M}=E_{\mathcal B}|_{\mathcal M},\qquad
 \iota(a)=\pi_o(a).
\]
The unit ball of $\mathcal M$ is the fixed-point set of the
uniformly $2$-norm continuous contraction $E_o$ on
$\mathcal B_1$, and is therefore complete.
Moreover,
\[
 \theta(\mathrm C(K))=\mathcal B\cap Z(N)
                    \subseteq\mathcal M\cap Z(\mathcal M).
\]
Consequently $E_{\mathcal M}$ is a faithful tracial conditional
expectation onto this specified central copy of $\mathrm C(K)$.
Thus the construction defines a continuous \Wstar{}-bundle.

Applying the preceding approximation to $x\in\mathcal M_1$
gives contractions $E_\alpha(b_n)\in A_1$ approximating $x$.
This proves (1); separability of $A$ also proves strict
separability of the bundle.
The expectation formula in (2) follows from the corresponding
formula for $\mathcal B$.
Since $S$ is a Bauer simplex, an affine continuous function on $S$
has the same supremum on $S$ and on $K$.
Applying this to $\tau\mapsto\tau(a^*a)$ gives the norm identity.
Faithfulness of $E_{\mathcal M}$ then gives (4).

To identify the fibers, let $\rho_\tau^M\colon M\to M_\tau$
be the normal coordinate map of Lemma~\ref{L:W*-expectation}.
For $\tau\in K$ and $x\in\mathcal M$,
\[
 (\operatorname{ev}_\tau\circ E_{\mathcal M})(x)
             =\bar\tau(\rho_\tau^M(x)).
\]
Thus the bundle GNS representation is unitarily equivalent to
$\rho_\tau^M|_{\mathcal M}$ on $H_\tau$:
its cyclic subspace contains $\pi_\tau(A)\xi_\tau$, which is
dense in $H_\tau$.
Ozawa's fiber theorem shows that this image is a von Neumann
algebra.  It is contained in $M_\tau$ and contains
$\pi_\tau(A)$, so it equals $M_\tau$.

Finally, the normal action on $M$ is implemented in $N$ by the
canonical group unitaries.
Their inner automorphisms preserve $\mathcal B$, its uniform
$2$-norm, and $M$, so they restrict to an action $\beta$ on
$\mathcal M$.
They fix $Z(N)$ pointwise, hence fix $\theta(\mathrm C(K))$ pointwise.
The action is therefore fiberwise and $\iota$ is equivariant.
The coordinate maps intertwine this action with $\alpha_\tau$,
giving (3).
\end{proof}

The fibers in this theorem need not be factors:
an extreme invariant trace need not be extreme in $\mathrm T(A)$.
The map $\iota$ is injective whenever invariant traces separate
$A_+$; in particular, this holds for a minimal action, since
$\ker\iota$ is a proper invariant ideal.

\subsection{From bundle towers to \Cstar{}-algebra towers}
\label{rokhlin-bundle-to-cstar}

The bundle ultrapower is canonically isomorphic to the uniform trace
quotient.  This identifies Rokhlin towers in the bundle with WRP
towers in the \Cstar{}-algebra.  Throughout this subsection, $\omega$ is
the fixed free ultrafilter on $\mathbb N$, $A_\omega$ is the norm
ultrapower, and $\mathcal M^\omega$ is the \Wstar{}-bundle ultrapower.

\begin{dfn}\label{D:wstar-bundle-rokhlin}
Let $\beta\colon\Gamma\to\operatorname{Aut}(\mathcal M)$ be
a bundle action of a countable discrete group.
The action has the \emph{\Wstar{}-bundle Rokhlin property} if, for
every exact tiling of $\Gamma$ with shapes
$K_1,\ldots,K_R$, there are projections
$p_1,\ldots,p_R\in\mathcal M^\omega\cap\mathcal M'$ such that
\[
 \sum_{s=1}^R\sum_{g\in K_s}\beta_g^\omega(p_s)=1.
\]
The translated projections are automatically pairwise
orthogonal, since their sum is $1$.
\end{dfn}

The next lemma is the invariant-trace version of
\cite[Lemma~3.10]{BBSTWW}.
Its proof uses the bounded density established in
Theorem~\ref{T:W*-bundle-const} and the trace-kernel lifting
property from Section~\ref{trace-kernels-crossed-products}.

\begin{lem}[The uniform tracial quotient]\label{L:Sur}
Let $(\mathcal M,K,E_{\mathcal M})$, $\beta$, and
$\iota\colon A\to\mathcal M$ be as in
Theorem~\ref{T:W*-bundle-const}, and put
$S=\mathrm T^\alpha(A)$.
Then the coordinatewise map
\[
 \Phi\colon A_\omega\longrightarrow\mathcal M^\omega,
 \qquad
 \Phi([(a_n)])=[(\iota(a_n))]
\]
is an equivariant surjective homomorphism satisfying
\[
 \|\Phi(x)\|_{2,u}=\|x\|_{2,S},\qquad
 \ker\Phi=J_S.
\]
Moreover, its restriction gives an equivariant surjection
\[
 \Phi\colon A_\omega\cap A'
            \longrightarrow\mathcal M^\omega\cap\mathcal M'.
\]
In particular, there are canonical equivariant isomorphisms
\[
 A_\omega/J_S\cong\mathcal M^\omega,\qquad
 \frac{A_\omega\cap A'}{J_S\cap(A_\omega\cap A')}
                  \cong\mathcal M^\omega\cap\mathcal M'.
\]
\end{lem}

\begin{proof}
The inequality
$\|\iota(a)\|_{2,u}=\|a\|_{2,S}\leq\|a\|$
shows that the coordinatewise homomorphism annihilates
norm-null sequences, so $\Phi$ is well defined.
Equivariance and the displayed $2$-norm identity follow
from Theorem~\ref{T:W*-bundle-const}.
Since the uniform $2$-norm on $\mathcal M^\omega$ is
faithful, this also identifies its kernel with $J_S$.

To prove surjectivity, represent $y\in\mathcal M^\omega$
by a bounded sequence $(y_n)$ in $\mathcal M$, and choose
$C\geq1$ with $\|y_n\|\leq C$ for all $n$.
The uniform $2$-norm density of $\iota(A_1)$ in
$\mathcal M_1$ gives $a_n\in A$ such that
\[
 \|a_n\|\leq C,\qquad
 \|\iota(a_n)-y_n\|_{2,u}<1/n.
\]
Thus $\Phi([(a_n)])=y$.
Let
$\overline\Phi\colon A_\omega/J_S\to\mathcal M^\omega$
be the resulting equivariant isomorphism.

For the assertion concerning central sequences, observe that
if $y\in\mathcal M^\omega$ commutes with $\iota(A)$,
then, for $m\in\mathcal M_1$ and $a\in A_1$,
\[
 \|[y,m]\|_{2,u}
       \leq2\|y\|\|m-\iota(a)\|_{2,u}.
\]
Bounded density makes the right hand side arbitrarily
small.  Hence $[y,m]=0$, and
\[
 \mathcal M^\omega\cap\iota(A)'
       =\mathcal M^\omega\cap\mathcal M'.
\]
Proposition~\ref{trace-kernel-sigma} gives
\[
 q_S(A_\omega\cap A')
       =q_S(A)'\cap(A_\omega/J_S).
\]
Applying $\overline\Phi$ to this equality proves
surjectivity on the central sequence algebras.
The kernel of the restriction is
$J_S\cap(A_\omega\cap A')$, as asserted.
\end{proof}

\begin{lem}[Transfer of Rokhlin towers]\label{bundle-to-wrp}
In the setting of Theorem~\ref{T:W*-bundle-const}, suppose
also that $\Gamma$ is amenable.
Then $\alpha$ has WRP in the sense of
Definition~\ref{D:weak-Rokhlin-property} if and only if
$\beta$ has the \Wstar{}-bundle Rokhlin property of
Definition~\ref{D:wstar-bundle-rokhlin}.
\end{lem}

\begin{proof}
Fix an exact tiling with shapes $K_1,\ldots,K_R$.
Suppose first that $\beta$ has the \Wstar{}-bundle Rokhlin
property, and choose its projections $p_s$.
Lemma~\ref{L:Sur} identifies these with projections in
$q_S(A)'\cap(A_\omega/J_S)$.
Their translated levels are pairwise orthogonal.
Lemma~\ref{trace-quotient-tower-lifting} therefore gives
positive contractions $c_s\in A_\omega\cap A'$ satisfying
$\Phi(c_s)=p_s$ whose translated levels are pairwise
orthogonal in the norm ultrapower.
Equivariance gives
\[
 \Phi\left(1-\sum_{s=1}^R\sum_{g\in K_s}
                       \alpha_g^\omega(c_s)\right)=0.
\]
Since $\ker\Phi=J_S$, the lifted towers satisfy WRP.

Conversely, let $c_s$ be WRP contractions for this tiling.
By Lemma~\ref{L:Sur}, the positive contractions
$p_s=\Phi(c_s)$ belong to
$\mathcal M^\omega\cap\mathcal M'$.
Their translated levels are pairwise orthogonal and sum
to $1$, because the WRP defect belongs to $J_S$.
Orthogonal positive elements summing to $1$ are
projections: multiplying the sum by any one level shows
that this level equals its square.
Thus the translated levels, and hence the $p_s$, are
projections.
\end{proof}

\subsection{From fiber towers to bundle towers}
\label{rokhlin-fibers-to-bundle}

The local-to-global argument combines exact Rokhlin towers in the
fibers with invariant partitions which reproduce prescribed central
trace weights.

Throughout the bundle arguments, $(\mathcal M,K,E)$ is a strictly
separable continuous \Wstar{}-bundle and
$\beta\colon\Gamma\to\operatorname{Aut}(\mathcal M)$ is a fiberwise
action of a countable discrete group.
Write $M_\lambda=\pi_\lambda(\mathcal M)$ for the fiber, and
$\bar\tau_\lambda$ for its faithful normal trace.
Strict separability gives each fiber separable predual.
The ultrapowers $(M_\lambda)^\omega$ are tracial, whereas
$\mathcal M^\omega$ is the uniform bundle ultrapower.
Put $E^\omega=E_{\mathcal M}^\omega$.
Elements of $\mathrm C(K)$ and of $\mathcal M$ are identified with their
constant images in the corresponding ultrapowers.

\subsubsection{Fiber quotients and invariant slices}

\begin{lem}[Lifting central sequences from a fiber]
\label{L:Sur_Fiber}
For every $\lambda\in K$, the coordinatewise fiber map gives
equivariant surjections
\[
 \Pi_\lambda\colon\mathcal M^\omega\longrightarrow(M_\lambda)^\omega,
 \qquad
 \Pi_\lambda\colon\mathcal M^\omega\cap\mathcal M'
          \longrightarrow(M_\lambda)^\omega\cap M_\lambda'.
\]
\end{lem}

\begin{proof}
The fiber map $\pi_\lambda\colon\mathcal M\to M_\lambda$ is
surjective, and
$\|\pi_\lambda(x)\|_{2,\bar\tau_\lambda}\leq\|x\|_{2,u}$.
Bounded lifting through this \Cstar{}-quotient proves surjectivity of
$\Pi_\lambda$; the inequality makes the map well defined.
Fiberwise equivariance makes both maps equivariant.

Let $x=[(x_n)]$ belong to the fiber central sequence algebra.
Choose lifts $a_n\in\mathcal M$ with uniformly bounded norms,
and fix a uniformly $2$-norm dense sequence
$(s_j)$ in $\mathcal M_1$.
Centrality of $x$ permits a choice of integers $k(n)$ tending
to infinity along $\omega$ such that
\[
 \|[x_n,\pi_\lambda(s_j)]\|_{2,\bar\tau_\lambda}<1/k(n)
       \quad(1\leq j\leq k(n))
\]
on an $\omega$-large set.
For these $n$, continuity of the fiber $2$-norms gives an open
neighborhood $U_n$ of $\lambda$ on which
\[
 \|[a_n,s_j]\|_{2,\mu}<2/k(n)
       \quad(\mu\in U_n,\ 1\leq j\leq k(n)).
\]
Choose $f_n\in \mathrm C(K)$ with $0\leq f_n\leq1$,
$f_n(\lambda)=1$, and $\operatorname{supp}(f_n)\subseteq U_n$.
On the remaining indices take $f_n=1$.
Centrality of $f_n$ gives
\[
 \|[f_na_n,s_j]\|_{2,u}\leq2/k(n)
       \quad(1\leq j\leq k(n))
\]
on the same $\omega$-large set.
Thus $[(f_na_n)]$ commutes with every $s_j$.
The uniform norm bound and $2$-norm density extend commutation
to all of $\mathcal M$, while
$\Pi_\lambda([(f_na_n)])=x$.
\end{proof}

\begin{dfn}\label{invariant-bundle-slicing}
An \emph{invariant $n$-slicing} consists of pairwise orthogonal
projections
\[
 e_1,\ldots,e_n\in
             (\mathcal M^\omega\cap\mathcal M')^{\beta^\omega}
\]
such that
\[
 \sum_{i=1}^n e_i=1,\qquad
 E^\omega(e_i x)=\frac1n E(x)
       \quad(x\in\mathcal M,\ 1\leq i\leq n).
\]
The expectation on the right is viewed as a constant element of
$\mathrm C(K^\omega)=\mathrm C(K)_\omega$.
Equivalently, the trace identity holds at every point of
$K^\omega$.
\end{dfn}

A unital homomorphism
$\mathrm M_n\to(\mathcal M^\omega\cap\mathcal M')^{\beta^\omega}$
gives an invariant $n$-slicing by taking its diagonal matrix units.
Indeed, traciality of $E^\omega$ and commutation with $x$
make $E^\omega(e_{ii}x)$ independent of $i$, and their sum is
$E(x)$.

\begin{lem}[Relative slices]\label{relative-bundle-slices}
If an invariant $n$-slicing exists, then, for every countable
subset $\mathcal S\subseteq\mathcal M^\omega$, there are
pairwise orthogonal projections
\[
 f_1,\ldots,f_n\in
  (\mathcal M^\omega\cap\mathcal M'\cap\mathcal S')^{\beta^\omega}
\]
with
\[
 \sum_i f_i=1,\qquad
 E^\omega(f_i x)=\frac1n E^\omega(x)
       \quad(x\in\mathcal S).
\]
\end{lem}

\begin{proof}
Adjoin to $\mathcal S$ the constant images of a countable
uniformly $2$-norm dense subset of $\mathcal M_1$, together
with $1$, adjoints, and all $\Gamma$-translates.
Enumerate the resulting countable set as $(x_j)$, and choose
bounded representatives $(x_{j,m})_m$ in $\mathcal M$.
Represent the given slicing projections by positive contraction
sequences $(e_{i,k})_k$.

For each $m$, the finitely many $x_{j,m}$, $j\leq m$, are
elements of $\mathcal M$.
The slicing identities and centrality give an index $k(m)$ such that
\[
 \|[e_{i,k(m)},x_{j,m}]\|_{2,u}<1/m,\qquad
 \left\|E(e_{i,k(m)}x_{j,m})-\frac1nE(x_{j,m})\right\|<1/m
       \quad(i\leq n,\ j\leq m).
\]
Choose this index so that the errors in the projection,
orthogonality, and partition relations are also below $1/m$,
as are the invariance errors for the first $m$ elements of a fixed
enumeration of $\Gamma$.  Each condition holds on an
$\omega$-large set, so their finite intersection is nonempty.
Then
$f_i=[(e_{i,k(m)})_m]$ satisfies all the asserted relations.
Commutation with the adjoined dense subset extends to
$\mathcal M$ by boundedness.
\end{proof}

\subsubsection{Partitions with prescribed central trace weights}

\begin{lem}\label{bundle-central-partitions}
Suppose that invariant $n$-slicings exist for every $n\geq1$.
Let $h_1,\ldots,h_l\in \mathrm C(K^\omega)_+$ satisfy $\sum_jh_j=1$,
and let $\mathcal S\subseteq\mathcal M^\omega$ be countable.
Then there are pairwise orthogonal projections
\[
 p_1,\ldots,p_l\in
  (\mathcal M^\omega\cap\mathcal M'\cap\mathcal S')^{\beta^\omega}
\]
such that
\[
 \sum_jp_j=1,\qquad
 E^\omega(p_jx)=h_jE^\omega(x)
       \quad(x\in\mathcal S,\ 1\leq j\leq l).
\]
These identities extend to the uniform $2$-norm closure of
$\mathcal S$.
\end{lem}

\begin{proof}
Adjoin $1$ and a countable uniformly $2$-norm dense subset of
$\mathcal M_1$ to the given subset.
For $N\geq1$, Lemma~\ref{relative-bundle-slices} gives an
invariant relative slicing $e_{1,N},\ldots,e_{N,N}$ with
\[
 E^\omega(e_{r,N}x)=\frac1N E^\omega(x)
       \quad(x\in\mathcal S).
\]
Put $s_0=0$ and $s_j=\sum_{i=1}^j h_i$, so $s_l=1$.
For $0\leq t\leq1$, define the continuous scalar functions
\[
 F_{r,N}(t)=\min\{1,\max\{0,Nt-r+1\}\}
       \quad(1\leq r\leq N).
\]
They satisfy $\sum_{r=1}^N F_{r,N}(t)=Nt$.
Use functional calculus in the specified central algebra to set
\[
 d_{j,r}=F_{r,N}(s_j)-F_{r,N}(s_{j-1}),\qquad
 a_{j,N}=\sum_{r=1}^N d_{j,r}e_{r,N}.
\]
The coefficients $d_{j,r}$ are central positive contractions.
Consequently the $a_{j,N}$ are positive contractions, commute
with the specified elements, are invariant, and satisfy
\[
 \sum_j a_{j,N}=1,\qquad
 E^\omega(a_{j,N}x)=h_jE^\omega(x)
       \quad(x\in\mathcal S).
\]

We estimate their orthogonality in the uniform $2$-norm.
At a point $\lambda\in K^\omega$, the scalar $d_{j,r}(\lambda)$
is $N$ times the length of the intersection of
$[s_{j-1}(\lambda),s_j(\lambda)]$ with the $r$th interval
of the equal $N$-partition of $[0,1]$.
For $i\ne j$, both coefficients can be nonzero in at most one
such interval.
Since $E^\omega(e_{r,N})=1/N$, it follows that
\begin{equation}\label{central-partition-orthogonality-estimate}
 \|a_{i,N}a_{j,N}\|_{2,u}^2
 =\left\|\frac1N\sum_{r=1}^N(d_{i,r}d_{j,r})^2\right\|
 \leq\frac1N.
\end{equation}

A diagonal argument preserves the prescribed
$h_j$ and the specified elements, including their
nonconstant ultrapower representatives.
Fix representatives for them, and enumerate the countably many
commutation, invariance, partition, and expectation relations.
Choose $N_m>m^2$ and positive contraction representatives for
$a_{j,N_m}$.
There are nested sets $V_m\in\omega$, with
$V_m\subseteq\{k:k\geq m\}$, on which the errors in the first $m$
relations and in all orthogonality relations are below $2/m$.
For an index $k$, let
\[
 m(k)=\max\{m\leq k:k\in V_m\},
\]
using arbitrary contractions when this set is empty.
Then $m(k)\to_\omega\infty$.
At coordinate $k$, choose the $k$th coordinate of
$a_{j,N_{m(k)}}$, retaining the fixed representatives of the $h_j$
and the specified elements at that coordinate.
The resulting $p_j$ satisfy all the asserted relations.

They are positive contractions with pairwise zero products and
sum $1$, so $p_j^2=p_j$.
Finally, commutators with a fixed bounded element and the maps
$x\mapsto E^\omega(p_jx)$ and $x\mapsto h_jE^\omega(x)$
are uniformly $2$-norm continuous in $x$.
This proves the extension assertion.
\end{proof}

\begin{thm}[From fiber towers to bundle towers]
\label{fiber-to-bundle-rokhlin}
Suppose that every fiber action $\beta_\lambda$ has the
\Wstar{}-Rokhlin property of Definition~\ref{exact-wstar-rokhlin},
and that invariant $n$-slicings exist for all $n$.
Then $\beta$ has the \Wstar{}-bundle Rokhlin property.
\end{thm}

\begin{proof}
Fix an exact tiling with shapes $K_1,\ldots,K_R$.
Let $\mathcal F\subseteq\mathcal M_1$ be finite and let
$\varepsilon>0$.
At $\lambda\in K$, take fiber Rokhlin projections in
$(M_\lambda)^\omega\cap M_\lambda'$.
Spectral cutoffs give projection representatives in
$M_\lambda$.
Choose one coordinate and lift its base projections to
positive contractions $f_{\lambda,s}\in\mathcal M$.
The coordinate can be chosen so that, at $\lambda$,
\[
 \left\|1-\sum_{s=1}^R\sum_{g\in K_s}
                  \beta_g(f_{\lambda,s})\right\|_{2,\lambda}
                  <\varepsilon/2,\qquad
 \|[f_{\lambda,s},a]\|_{2,\lambda}<\varepsilon/2
       \quad(a\in\mathcal F).
\]
Also $\|f_{\lambda,s}^2-f_{\lambda,s}\|_{2,\lambda}=0$.
Continuity gives an open neighborhood $U_\lambda$ on which
all three kinds of errors are below $\varepsilon$.
Take a finite subcover $U_1,\ldots,U_L$, corresponding lifts
$f_{j,s}$, and a partition of unity
$h_1,\ldots,h_L\in \mathrm C(K)_+$ subordinate to this cover.

Apply Lemma~\ref{bundle-central-partitions} with a countable
uniformly $2$-norm dense subset of $\mathcal M_1$.
Its extension assertion gives invariant orthogonal projections
$p_j$ summing to $1$ and satisfying
\[
 p_j\in\mathcal M^\omega\cap\mathcal M',\qquad
 E^\omega(p_jx)=h_jE(x)\quad(x\in\mathcal M).
\]
Define positive contractions
\[
 q_s=\sum_{j=1}^L p_jf_{j,s}\in\mathcal M^\omega.
\]
For arbitrary $d_1,\ldots,d_L\in\mathcal M$, centrality and
orthogonality of the $p_j$ give
\begin{equation}\label{weighted-bundle-estimate}
 \left\|\sum_jp_jd_j\right\|_{2,u}^2
       =\left\|\sum_jh_jE(d_j^*d_j)\right\|.
\end{equation}
In each application below,
$\|d_j\|_{2,\mu}<\varepsilon$ on $U_j$.
Since $\operatorname{supp}(h_j)\subseteq U_j$,
the right hand side is at most $\varepsilon^2$.
Apply \eqref{weighted-bundle-estimate} to
\[
 f_{j,s}^2-f_{j,s},\qquad
 1-\sum_{s,g\in K_s}\beta_g(f_{j,s}),\qquad
 [f_{j,s},a]\ (a\in\mathcal F).
\]
Invariance of the $p_j$ is used in the middle expression.
We obtain
\[
 \|q_s^2-q_s\|_{2,u}\leq\varepsilon,\qquad
 \left\|1-\sum_{s,g\in K_s}\beta_g^\omega(q_s)\right\|_{2,u}
             \leq\varepsilon,\qquad
 \|[q_s,a]\|_{2,u}\leq\varepsilon\quad(a\in\mathcal F).
\]

Let $\mathcal F_m$ increase through a dense sequence in
$\mathcal M_1$ and let $\varepsilon=1/m$.
For each resulting family, choose a coordinate of positive contraction
representatives at which the finitely many displayed errors
are below $2/m$.
The sequence of these coordinates gives projections
$e_s\in\mathcal M^\omega\cap\mathcal M'$ satisfying
\[
 \sum_{s=1}^R\sum_{g\in K_s}\beta_g^\omega(e_s)=1.
\]
All translated levels are projections, and their sum is $1$,
so they are pairwise orthogonal.
Thus the bundle action has the \Wstar{}-bundle Rokhlin property.
\end{proof}

\subsubsection{Finite-dimensional gluing and the main conclusion}

\begin{lem}[Invariant diffuse gluing]
\label{invariant-diffuse-gluing}\label{equivariant-mcduff-gluing}
Suppose that the strictly separable bundle $(\mathcal M,K,E)$
has compact metrizable base of finite covering dimension, and
that $\beta$ is a fiberwise bundle action of a countable group,
in the sense of Definition~\ref{D:bundle-action}.
Assume that, for every $\lambda\in K$,
\[
 Z(M_\lambda)^{\beta_\lambda}=\mathbb C1
\]
and there is a positive contraction
\[
 z_\lambda\in
   ((M_\lambda)^\omega\cap M_\lambda')^{\beta_\lambda^\omega}
\]
with Lebesgue distribution on $[0,1]$ for the ultrapower trace.
Then invariant $n$-slicings exist for every $n\geq1$.
\end{lem}

\begin{proof}
Put $c=\dim K+1$.  We first construct a commuting tuple in
the bundle ultrapower from the fiber contractions.  Fix finite sets
$\mathcal F_m\subset\mathcal M_1$ increasing through a uniformly
$2$-norm dense sequence, finite sets $H_m\subset\Gamma$ increasing
to $\Gamma$, and numbers $\delta_m\downarrow0$.

\emph{Successive covering families.}
For each $m$, construct $c$ finite families successively.
To construct family $\ell$, at each $\lambda\in K$ choose a
positive contraction coordinate $t_\lambda\in M_\lambda$ of
$z_\lambda$ that is $\delta_m$-central for the
images of $\mathcal F_m$ and of every element in all earlier
families, $\delta_m$-invariant under $H_m$, and satisfies
\[
 \left|\bar\tau_\lambda(t_\lambda^p)-\frac1{p+1}\right|
      <\delta_m\quad(1\leq p\leq m).
\]
Lift $t_\lambda$ to a positive contraction $z\in\mathcal M$.
Choose the coordinate with a sufficiently small initial error.
Continuity then gives a neighborhood $U_\lambda$ on which
all these estimates hold with error $\delta_m$.
Choose a smaller neighborhood $V_\lambda$ with
$\overline{V_\lambda}\subset U_\lambda$ and a central function
$f_\lambda\in \mathrm C(K)_+$ satisfying
\[
 0\leq f_\lambda\leq1,\qquad
 f_\lambda|_{V_\lambda}=1,\qquad
 \operatorname{supp}(f_\lambda)\subset U_\lambda.
\]
Replace the lift by
\[
 b_\lambda=f_\lambda z+(1-f_\lambda)\tfrac12 1.
\]
Its commutators with the specified elements and its invariance
errors are below $\delta_m$ in the uniform $2$-norm:
outside $U_\lambda$ it is scalar, and the cutoff is central
and fixed by $\beta$. Its moments satisfy the stated estimates
on $V_\lambda$.  Choose a finite subcover consisting of sets
$V_\lambda$; the corresponding contractions $b_\lambda$ form
family $\ell$.
Every family covers all of $K$, and each later family uniformly
almost commutes with every earlier family.
This is the scalar version of the commuting covering-system
construction in \cite[Claim~8.4.1]{Gardella-Hirshberg-Vaccaro}.

After constructing all $c$ families, take a common
refinement of their neighborhood covers.  Finite covering
dimension gives a subordinate finite partition of unity $(h_j)$
with $c$ colors such that $h_jh_k=0$ for distinct indices of the
same color.  Such a partition is obtained from a refinement of
order at most $c$ by coloring the barycentric subdivision of its
nerve according to simplex dimension.
For each member of color $\ell$, choose $b_{\ell,j}$ from family
$\ell$ whose moment estimates hold throughout its support.  Put
\[
 a_{\ell,m}=\sum_{j\text{ of color }\ell}h_jb_{\ell,j},
 \qquad w_{\ell,m}=\sum_{j\text{ of color }\ell}h_j.
\]
Then $0\leq a_{\ell,m}\leq w_{\ell,m}$ and
$\sum_\ell w_{\ell,m}=1$.  At each fiber, at most one summand
of a given color is nonzero.  Consequently all centrality,
invariance, and cross-color commutator errors are bounded by
$\delta_m$, independently of the number of patches. Also, for
$1\leq p\leq m$,
\[
 \left\|E(a_{\ell,m}^p)
              -\frac{w_{\ell,m}^p}{p+1}\right\|<\delta_m.
\]
Thus, writing $Y=K^\omega$ and $\mathrm C(Y)=\mathrm C(K)_\omega$, the classes
$a_\ell=[(a_{\ell,m})]$ commute with each other, belong to
$(\mathcal M^\omega\cap\mathcal M')^{\beta^\omega}$, and satisfy
\begin{equation}\label{diffuse-gluing-marginals}
 0\leq a_\ell\leq w_\ell,\qquad
 \sum_\ell w_\ell=1,\qquad
 E^\omega(a_\ell^p)=\frac{w_\ell^p}{p+1}\quad(p\geq1),
 \quad w_\ell=[(w_{\ell,m})]\in \mathrm C(Y).
\end{equation}
These moment identities determine the distribution of
$(a_\ell)_y$ at every $y\in Y$: it is uniform on
$[0,w_\ell(y)]$, with the point mass at zero when $w_\ell(y)=0$.

\emph{A uniform diffuse variable.}
The maximum
\[
 a=\max_{1\leq\ell\leq c}a_\ell
\]
is defined by joint continuous functional calculus in the
commuting algebra generated by the $a_\ell$ and $\mathrm C(Y)$.
In each fiber the commuting tuple is a tuple of random
variables with the marginals in
\eqref{diffuse-gluing-marginals}.  The following estimate depends
only on these marginal distributions.
For an interval $[t,t+\delta]\subset[0,1]$, put
$r=\sqrt\delta$.  If $t\geq r$, a variable attaining the maximum
in this interval has weight at least $t$. The union bound gives
interval mass at most $c\delta/t\leq c\sqrt\delta$.
If $t<r$, choose a weight at least $1/c$; the event that the
maximum lies in the interval forces the corresponding variable,
after division by its weight, below $c(r+\delta)$.  This rescaled
variable is uniform on $[0,1]$.
In either case,
\begin{equation}\label{diffuse-gluing-concentration}
 \bar\tau_y(1_{[t,t+\delta]}(a_y))\leq2c\sqrt\delta.
\end{equation}
Intervals meeting the endpoints are handled by truncation.
In particular, all fiber distributions are atomless with a
common modulus of continuity.

Continuous ramps at $t$ are Cauchy in the uniform $2$-norm by
\eqref{diffuse-gluing-concentration}. Uniform $2$-norm
completeness of the ultrapower bundle therefore gives
$P_t=1_{[0,t]}(a)\in\mathcal M^\omega$.
Set $H_t=E^\omega(P_t)\in \mathrm C(Y)$. The same estimate gives
\[
 \|H_t-H_s\|\leq2c\sqrt{|t-s|},\qquad H_0=0,\quad H_1=1.
\]
Hence $H(y,t)=H_t(y)$ belongs to $\mathrm C(Y\times[0,1])$.
Joint continuous functional calculus with the central algebra
$\mathrm C(Y)$ defines $u=H(a)$. In the fiber at $y$, this is the
continuous distribution function of $a_y$ evaluated at $a_y$;
the probability integral transform gives
\[
 E^\omega(f(u))=\left(\int_0^1 f(t)\,dt\right)1
       \quad(f\in \mathrm C([0,1])).
\]
Thus $u\in(\mathcal M^\omega\cap\mathcal M')^{\beta^\omega}$
has uniform distribution in every fiber.

\emph{Mixed trace factorization.}
Strict separability gives countably many unitaries $v_j\in\mathcal M$
whose images generate every fiber as a von Neumann algebra:
use exponentials of a countable uniformly $2$-norm dense set of
selfadjoint contractions. Let $G$ be the countable group of
automorphisms generated by $\operatorname{Ad}(v_j)$ and
$\beta_g$, $g\in\Gamma$. Each automorphism fixes $\mathrm C(K)$ and
preserves $E$. Choose a symmetric full-support probability
$\mu$ on $G$ with $\mu(1)\geq1/2$, and set
$T=\sum_{\theta\in G}\mu(\theta)\theta$.
On each fiber $L^2$ space, $T$ is a selfadjoint contraction.
It is positive because the identity has weight at least $1/2$;
full support of $\mu$ identifies its fixed vectors as
\[
 L^2(Z(M_\lambda))^{\beta_\lambda}=\mathbb C1.
\]
For $x\in\mathcal M$, the continuous functions
\[
 \lambda\longmapsto
 E(|T^k x-E(x)|^2)(\lambda)
\]
decrease pointwise to zero. Dini's theorem gives
$\|T^k x-E(x)\|_{2,u}\to0$.
If $q\in(\mathcal M^\omega\cap\mathcal M')^{\beta^\omega}$,
traciality and invariance imply
$E^\omega(q\theta(x))=E^\omega(qx)$ for every $\theta\in G$.
Consequently
\[
 E^\omega(qx)=E^\omega(qT^kx)
       \longrightarrow E^\omega(q)E(x),
\]
so $E^\omega(qx)=E^\omega(q)E(x)$.

Finally, take the spectral projections of $u$ for the $n$ equal
subintervals of $[0,1]$. As above, uniform $2$-norm limits of
continuous ramps supply these projections inside $\mathcal M^\omega$.
They are invariant and commute with $\mathcal M$, are mutually
orthogonal, sum to $1$, and have expectation $1/n$.  The mixed
trace factorization therefore makes them an invariant $n$-slicing.
\end{proof}

\begin{cor}[Tensoring with a McDuff tracial algebra]
\label{mcduff-tensor-wrp}
Let $A$ be unital and separable with $\mathrm T(A)\ne\emptyset$,
and let a countable discrete amenable group $\Gamma$ act on $A$.
Suppose that $\alpha$ is \Wstar{}-centrally free and
$S=\mathrm T^\alpha(A)$ is Bauer.
Let $B$ be unital and separable with a unique tracial state
$\nu$, and suppose that $N=\pi_\nu(B)''$ is McDuff.
Then
\[
 \alpha\otimes\operatorname{id}_B
       \colon\Gamma\longrightarrow\operatorname{Aut}(A\otimes_{\min}B)
\]
has WRP.
\end{cor}

\begin{proof}
Every tracial state on $D=A\otimes_{\min}B$ is a product
$\rho\otimes\nu$.
Indeed, for $a\in A_+$ the positive trace
$b\mapsto\sigma(a\otimes b)$ is a scalar multiple of $\nu$;
linearity gives the product formula.
Hence the invariant trace space of
$\gamma=\alpha\otimes\operatorname{id}_B$ is canonically $S$,
and its GNS algebras are
\[
 \pi_{\rho\otimes\nu}(D)''
       =M_\rho\ \bar\otimes\ N.
\]
Uniqueness of $\nu$ makes $N$ a factor.

The action $\gamma$ is \Wstar{}-centrally free.
To see this, note that
$Z(M_\rho\bar\otimes N)=Z(M_\rho)\otimes1$.
A centrally trivial restriction of
$(\alpha_\rho)_g\otimes\operatorname{id}_N$ to a nonzero
invariant central corner would act trivially on the
embedded centralizing sequences from the corresponding
corner of $M_\rho$.
These sequences remain centralizing after tensoring with
$1_N$, as is seen from the product trace and density of
elementary tensors.
This would make the restriction of $(\alpha_\rho)_g$
centrally trivial, contradicting the hypothesis.

Construct the bundle $(\mathcal N,\partial_eS)$ for
$(D,\gamma)$ by Theorem~\ref{T:W*-bundle-const}.
For each $n$, McDuffness of $N$ gives a unital copy of
$\mathrm M_n$ in $N^\omega\cap N'$.
Choose bounded matrix-unit representatives and approximate
them coordinatewise by elements
$\pi_\nu(b_{ij,k})$, $b_{ij,k}\in B$, using Kaplansky density
and bounded lifting through the GNS homomorphism.
All matrix-unit and centrality errors tend to zero in
$\|\cdot\|_{2,\nu}$.
For $c\in B$,
\[
 \|1\otimes c\|_{2,S_D}=\|c\|_{2,\nu},
 \qquad S_D=\{\rho\otimes\nu:\rho\in S\}.
\]
Moreover, for $a\in A$ and $b\in B$,
\[
 \|[1\otimes b_{ij,k},a\otimes b]\|_{2,S_D}
       \leq\|a\|\|[b_{ij,k},b]\|_{2,\nu}.
\]
Thus their images give a unital homomorphism
\[
 \mathrm M_n\longrightarrow
        (\mathcal N^\omega\cap\mathcal N')^{\gamma^\omega}.
\]
Commutation with the whole bundle follows from its bounded
uniform $2$-norm density of the image of $D$.
The image is invariant because $\gamma$ fixes the second
tensor factor.

These matrix algebras supply invariant slicings.
The centrally free fiber actions have exact \Wstar{}-Rokhlin
towers by Theorem~\ref{C:W*-Rokhlin-property}.
Apply Theorem~\ref{fiber-to-bundle-rokhlin} and then
Lemma~\ref{bundle-to-wrp} to conclude.
\end{proof}

\begin{thm}\label{finite-dimensional-trace-wrp}\label{bauer-invariant-traces-wrp}
Let $A$ be a unital separable nuclear \Cstar{}-algebra with
$\mathrm T(A)\ne\emptyset$, and let
$\alpha\colon\Gamma\to\operatorname{Aut}(A)$ be an action
of a countable discrete amenable group.
Suppose that $\alpha$ is \Wstar{}-centrally free and that
$S=\mathrm T^\alpha(A)$ is a Bauer simplex with
$\dim(\partial_eS)<\infty$.
If $\Gamma$ is infinite, then $\alpha$ has WRP.
The same conclusion holds for finite $\Gamma$ when $A$ has
no finite-dimensional quotients.
\end{thm}

\begin{proof}
\Wstar{}-central freeness implies \Wstar{}-proper outerness.
Theorem~\ref{T:W*-bundle-const} therefore gives a strictly
separable bundle over $K=\partial_eS$, with a fiberwise
action $\beta$ and fibers
\[
 M_\tau=\pi_\tau(A)''\quad(\tau\in K).
\]
Write $\bar\tau$ for the faithful normal trace on $M_\tau$.
Extremality in the invariant trace simplex implies
\begin{equation}\label{invariant-fiber-center-ergodicity}
 Z(M_\tau)^{\alpha_\tau}=\mathbb C1.
\end{equation}
Indeed, a nontrivial invariant central projection would split
$\tau$ into two invariant tracial states. These are distinct:
Kaplansky density distinguishes their normal extensions on the
two central corners, contradicting extremality.
The central supports of the type decomposition are invariant.
Thus $M_\tau$ is either wholly of type $\mathrm{II}_1$ or
homogeneous of type $\mathrm I_r$ for one finite $r$.

In the type $\mathrm{II}_1$ case, nuclearity gives hyperfiniteness.
For every $N$, \cite[Corollary~5.5]{Gardella-Hirshberg-Vaccaro}
gives a unital homomorphism
\[
 \mathrm M_N\longrightarrow
       ((M_\tau)^\omega\cap M_\tau')^{\alpha_\tau^\omega}.
\]
The cited result applies also to nonfactorial hyperfinite
type $\mathrm{II}_1$ algebras.  The diagonal contraction with
$N$ equally weighted, equally spaced eigenvalues has moments
converging to those of Lebesgue measure.  Letting $N\to\infty$
and choosing coordinates diagonally for increasing finite
centrality, invariance, and moment conditions gives a positive
contraction in the displayed invariant central sequence algebra
with exact Lebesgue distribution.

In the type $\mathrm I_r$ case, write
\[
 M_\tau=L^\infty(Y_\tau,\nu_\tau)\ \bar\otimes\ \mathrm M_r.
\]
Equation~\eqref{invariant-fiber-center-ergodicity} says that the
induced probability-measure-preserving action on $Y_\tau$ is
ergodic.  It is also essentially free: if a nonidentity $g$ had
a positive-measure fixed-point set, its characteristic projection
would give an $\alpha_g$-invariant central corner on which
$\alpha_g$ fixes the center. An automorphism of a homogeneous
finite type I algebra fixing its center is inner, hence centrally
trivial, contrary to \Wstar{}-central freeness.

Suppose first that $\Gamma$ is infinite. Essential freeness
implies that $\nu_\tau$ is nonatomic, since an atom would have
infinitely many disjoint translates of equal positive measure.
By the hyperfiniteness theorem for amenable equivalence relations
\cite[Theorem~10]{CFW-Cartan}, the orbit relation on $Y_\tau$ is, modulo null
sets, the increasing union of finite probability-measure-preserving
subrelations $\mathcal R_k$.
Each finite Borel subrelation is smooth, so its quotient is a standard
Borel space.  The pushforward of $\nu_\tau$ is nonatomic: the
preimage of a quotient point is a finite class and has measure zero.
For any $N$, split this quotient probability space into $N$ sets
of measure $1/N$ and pull them back to central projections
$e_{1,k},\ldots,e_{N,k}\in M_\tau$.
They are an exact partition of $1$ and satisfy
\[
 \bar\tau(e_{i,k})=\frac1N,\qquad
 \|\alpha_g(e_{i,k})-e_{i,k}\|_{2,\bar\tau}^2
 \leq\nu_\tau\{y:(g^{-1}y,y)\notin\mathcal R_k\}\longrightarrow0.
\]
Their ultrapower classes are invariant central projections.
The corresponding equally spaced diagonal contractions, with
$N\to\infty$ and a diagonal choice of $k$, give a positive
contraction with Lebesgue distribution in the invariant central
sequence algebra of $M_\tau$.
If instead $\Gamma$ is finite, its ergodic action on $Y_\tau$
is a single finite orbit.  Then $M_\tau$ is finite-dimensional,
and the GNS homomorphism maps $A$ onto $M_\tau$.
This contradicts the assumed absence of finite-dimensional
quotients, excluding the type I case for a finite group.

Thus every fiber satisfies the hypotheses of
Lemma~\ref{invariant-diffuse-gluing}, which gives
invariant $n$-slicings for all $n$.
The centrally free fiber actions have exact \Wstar{}-Rokhlin towers
by Theorem~\ref{C:W*-Rokhlin-property}.
Theorem~\ref{fiber-to-bundle-rokhlin} gives bundle towers, and
Lemma~\ref{bundle-to-wrp} lifts them to WRP for $\alpha$.
\end{proof}

\begin{rem}[The finite-group restriction]\label{finite-group-diffuse-obstruction}
The quotient hypothesis cannot be dropped for arbitrary finite
groups.  Let $\mathbb Z/2\mathbb Z$ act antipodally on the circle,
and let $A=\mathrm C(S^1)$.  The action is free, hence \Wstar{}-centrally free,
and its invariant-measure simplex is Bauer with extreme boundary
$S^1/(\mathbb Z/2\mathbb Z)$ of dimension one. Nevertheless,
\[
 \|f\|_{2,S}\geq\frac1{\sqrt2}\|f\|\quad(f\in \mathrm C(S^1)),
\]
so $J_S=0$ in the norm ultrapower.
WRP for the single full-group tiling would give orthogonal
positive contractions $p,\alpha(p)$ with $p+\alpha(p)=1$;
hence $p$ is a projection. The norm ultrapower of $\mathrm C(S^1)$
has only the projections $0$ and $1$: uniform approximate
projections on the connected circle approach scalar projections,
and their scalar ultrapower is $\mathbb C$.
Neither projection can satisfy the displayed partition.
\end{rem}

\begin{cor}[Commutative coefficients]\label{commutative-bauer-urp}
Let a countably infinite discrete amenable group act freely on
a nonempty compact metrizable space $X$.
If its invariant-measure simplex is Bauer with finite-dimensional
extreme boundary, then the action has topological URP, and the
induced action on $\mathrm C(X)$ has WRP.
\end{cor}

\begin{proof}
$\mathrm C(X)$ is separable and nuclear, and its invariant traces are exactly
the invariant probability measures.  Freeness implies that every
invariant-measure GNS action on $L^\infty(X,\mu)$ is centrally
free: in an abelian von Neumann algebra a centrally trivial
automorphism is the identity, and a nonzero invariant central
corner on which a nonidentity $g$ acts identically would be a
positive-measure fixed-point set.
Apply Theorem~\ref{finite-dimensional-trace-wrp} and
Proposition~\ref{top-urp-wrp}.
\end{proof}


\bigskip
\begingroup\footnotesize
\noindent\mbox{Department of Mathematics, University of Toronto, Toronto, Ontario, Canada M5S~2E4}\par
\noindent\textit{Email address: }\texttt{elliott@math.toronto.edu}\par\smallskip
\noindent\mbox{School of Mathematics and Statistics, Northeast Normal University, Changchun, Jilin, China}\par
\noindent\textit{Email address: }\texttt{licg864@nenu.edu.cn}\par\smallskip
\noindent\mbox{Department of Mathematics and Statistics, University of Wyoming, Laramie, WY 82072, USA}\par
\noindent\textit{Email address: }\texttt{zniu@uwyo.edu}\par\smallskip
\noindent\mbox{Sanctuary Cognitive Systems Corp., 425 W 6th Ave, Vancouver, BC, Canada V5Y~1L3}\par
\noindent\textit{Email address: }\texttt{qingyun322@gmail.com}\par\smallskip
\endgroup

\end{document}